\documentstyle[12pt,amssymb,amscd]{amsart}
\input xy
\xyoption{all} \textwidth 7in \textheight
10in \topmargin -1cm \voffset -1cm \hoffset
-3cm
\def\fact{{\text{fact}}}
\def\bodya{{\cal B}_{\text{\bf presymm}}}
\def\bodyb{{\cal B}_{\text{\bf symm}}}
\def\follows{\Rightarrow}
\def\disjoint{{\sqcup}}
\def\bs{\backslash}
\def\A{{\frak A}}
\def\intr{\int}
\def\ins{{\frak ins}}
\def\inr{{\frak int}}
\def\la{\lambda}
\def\An{{\frak An}}
\def\P{{\cal P}}
\def\limdir{\text{limdir}}
\def\ev{\text{ev}}
\def\notsubset{\text{notsubset}}
\def\identical{=}
\def\ins{{\frak ins}}
\def\dirlim{\limdir}
\def\usual{{\text{usual}}}
\def\A{{\cal A}}
\def\B{{\cal B}}
\def\F{{\frak F}}
\def\a{{\frak a}}
\def\b{{\frak b}}
\def\e{{\frak e}}
\def\jk{{\frak jk}}
\def \q{{\frak q}}
\def \qq{\q}
\def\sym{{\text{sym}}}
\def\res{{\text{res}}}
\def\f{{\fun}^{\sim}}
\def \Im{\text{Im}}
\def\jj{{\frak jj}}
\def\fp{({\fun'})^{\sim}}
\def\r{{\frak r}}
\def \into{\to}
\def \comm{\text{\bf comm}}
\def\lis{{\bf lis}}
\def\cat{{\bf cat}}
\def\red{\text{red}}
\def\lism{{\bf lism}}
\def\complexes{\text{complexes}}
\def\op{\text{op}}
\def \ores{{\O^{\text{res}}}}
\def\asm{{\frak ope'}}
\def\ens{{\bf ens}}
\def\ske{{\bf ske}}
\def\ve{\varepsilon}
\def\vect{{\bf vect}}
\def \id{\text{Id}}
\def\Id{\id}
\def \as{{\frak as}}
\def\ope{{\frak ope}}
\def\fun{{\frak fun}}
\def\funp{({\frak fun}^{\psys})}
\def \Fun{{\frak Fun}}
\def\funct{{\frak Funct}}
\def\system{{{\cal SYS}}}
\def\psys{{\cal SYS}[p]}

\def\lie{{\frak lie}}
\def\li{\lie^1}
\def \ll{{{\li}_M}}
\def\psystem{{{\cal SYS}[p,\ll]}}
\def\psystemg{{{\cal SYS}[p,\lie_M]}}
\def\OPE{{\cal OPE}}

\def \OPER{{\cal OPE}^{\text{res}}}
\def\psystemb{\OPE[p,\lie_M]}
\def\psystemg{\psystemb}
\def\psis{\psystemb^{\text{res}}}

\def\univ{\text{univ}}
\def\psim{\psystemb^{\text{sym}}}
\def\psimr{\psystem^{\text{sym},\text{res}}}

\def\lisp{{[p,\lie_M,\S]}}
\def\lispq{{\lisp^{\text{res}}}}
\def\psml{\system[p,\lie_M,\S]}
\def\psmq{\psml^{\text{res}}}
\def\psmlh{\psml[[\h]]}
\def\psmqh{\psmq[[\h]]}
\def\coker{\text{Coker}}
\def\M{{\cal M}}
\def\gen{{\frak gen}}
\def\u{{\frak u}}
\def\Mg{{\frak M}}
\def \hh{{\frak h}}
\def \kk{{\frak kk}}
\def \II{{\frak II}}
\def \JJ{{\frak JJ}}
\def \GG{{\frak GG}}
\def \FF{{\frak FF}}
\def \m{{\frak m}}
\def \z{{\frak z}}
\def \N{{\cal N}}
\def \B{{\cal B}}
\def \G{{\cal G}}
\def \M{{\cal M}}
\def \O{{\cal O}}
\def \R{{\Bbb R}}
\def \Z{{\frak Z}}
\def \Co{{\Bbb C}}
\def \C{{\cal C}}
\def \i{{\frak i}}
\def \k{{\frak k}}
\def \ii{{\frak ii}}
\def \l{{\frak l}}

\def \L{{\frak L}}
\def \DDD{{\cal D}}
\def \D{{\cal D}}
\def \DD{\Delta}
\def \g{{\frak g}}
\def \ha{^\wedge}

\def \liminv{\text{liminv}}
\def \p{{\frak p}}
\def \pp{{\pi_p}}
\def \h{h}
\def \I{{\cal I}}
\def \s{{\frak s}}
\def \S{{\frak S}}
\def \U{{\cal U}}
\def \vs{\sigma}
\def\pt{\text{pt}}
\def\Hom{\text{Hom}}
\def \holi{{\frak holi}}
\def \La{{\holi_M}}
\def \Lab{{ {\frak holie}_M}}
\def \psl{\OPER[p,\Lab]}
\def \ps {\psl}
\def\ov{\overline}
\def\comp{{\bf comp}}
\def \symm{{\bf symm}}
\def\Rs{\R^\symm}
\def\h{\lambda}
\newtheorem{theorem}{Theorem}[section]
\newtheorem{Axiom}[theorem]{Axiom}
\newtheorem{Claim}[theorem]{Claim}

\newtheorem{Lemma}[theorem]{Lemma}
\newtheorem{Corollary}[theorem]{Corollary}
\newtheorem{Proposition}[theorem]{Proposition}
\newtheorem{Theorem}[theorem]{Theorem}

\newtheorem{Definition-Proposition}[theorem]
{Definition-Proposition}
\newtheorem{Condition}[theorem]{Condition}
\def\R{{\frak R}}
\def\Re{{\Bbb R}}
\def\Co{{\Bbb C}}
\def\Cm{\D\text{-sh}}
\title{A formalism for the  renormalization procedure\footnote{Partially supported by an NSF grant}}
\author{Dimitri Tamarkin}
\begin{document}

\maketitle

\centerline{INTRODUCTION} The purpose of
this paper is to analyze the procedure of
renormalization from the mathematical point
of view.  Our original motivation came from
trying to really understand the paper
\cite{CF}. This paper uses the so-called
Batalin-Vilkovitski formalism
\cite{BV},\cite{Sch}. Its main features
include:

1) given a QFT, one constructs a so-called
quantum Batalin-Vilkovitski bracket on the
space of observables. Using this bracket one
writes  a Master equation (a.k.a.
Maurer-Cartan equation);

2) every solution to this equation is
supposed to produce  a deformation of the
QFT.

It is the procedure of constructing such a
deformation that is called {\em
renormalizaton} in the current paper.

Unfortunately, the treatment in \cite{CF}
does not lead to a (mathematically)
non-contradictory definition of the
Batalin-Vilkovitski bracket or
renormalization (due to divergencies). The
goal of this paper is to begin filling this
gap up.

Before working with the QFT from \cite{CF}
(i.e. the  Poisson sigma model), it makes
sense to start with simpler theories and to
define the Batalin-Vilkovitski bracket and
the renormalization for them. In this paper
we do it for the theory of free boson in
$\Re^{2n}$, $n>1$. It turns out that the
construction generalizes more or less
straightforwardly to the situation in
\cite{CF}, which will be a subject of a
subsequent paper.

The author hopes that the constructions of
this paper will also work in a  more general
context.

We deal with QFT's via  a $\D_X$-module $M$
of observables of the theory ($X$ is the
space-time) and an OPE-product structure on
$M$. So, we start with a definition of an
OPE-product. To this end  one first has  to
prescribe possible singularities of these
OPE. We call such a prescription {\em a
system} (a precise definition is given
below). Given a system, we have a notion of
an OPE-algebra over this system.

  We then construct an
appropriate system   for the free scalar
boson Euclidean  theory in $\Re^{2n}$,
$n>1$, in which case the only possible
singularities are of the type: products  of
squares of Euclidean distances in the
denominator.  We denote this system by
$<\i>$. The cases $\Re^2,\Re^{2n+1}$ require
semi-integer powers  or logarithms, which
leads to  slightly more complicated
definitions. For simplicity we only work
with $\Re^{2n}$, $n>1$ throughout the paper.

We then show that the Batalin-Vilkovitski
bracket arises due to a certain additional
structure on the system. We call a system
with such a structure {\em pre- symmetric}.
The system $<\i>$ has no natural
pre-symmetric structure,  nevertheless we
construct a differential graded resolution
$<\R>\to <\i>$ which is pre-symmetric.
Furthermore, any OPE-algebra over $<\i>$ can
be lifted to an OPE-algebra over $<\R>$. The
building blocks for the  system $<\R>$ are
certain spaces of generalized functions. The
lifting procedure can be interpreted as a
regularization (i.e. passage from usual
functions with singularities to generalized
functions). It seems to be very similar to
the well-known Bogolyubov-Parasyuk-Hepp
procedure \cite{BP}. There also is some
affinity with the approach in \cite{Bor1}.

It is worth to mention that the homotopy
theory implies that, up-to homotopy, nothing
should depend on the choice of such a
lifting. What is not implied by the homotopy
theory is that we can always find an
"honest" lifting (as opposed to  a lifting
up-to higher homotopies). Furthermore, we
expect that the action of Hopf algebras
introduced in \cite{CKI} (see also
\cite{Bor1}, where a somewhat similar object
appears under the name of "the group of
renormalisations") should provide us with
("honest", not quasi-) isomorphisms between
the different liftings, which also looks a
bit different from what we are used to in
the homotopy theory.

 Next, we treat the
renormalization procedure. It turns out that
to accomplish such a procedure, one needs
certain additional properties of the system.
We call a system with these additional
properies {\em symmetric system}.
Unfortunately, the system $<\R>$ is  only
pre-symmetric, and not symmetric. The reason
is very simple: the renormalized OPE have
more sophisticated singularities. It turns
out, though, that there is a formal
"symmetrization" procedure, which produces a
symmetric system out of a pre-symmetric one.
So, starting from $<\R>$, we get a symmetric
system $<\Rs>$ and construct a renormalized
OPE in this system.

Morally, the system $<\Rs>$ is given in
terms of a $\D$-module whose solutions are
possible singularities of  the renormalized
OPE. Our last step is to interpret this
$\D$-module as a sub-module in the space of
real-analytic functions.

Our approach has to be compared with the
ones in \cite{Bor} and
\cite{CKI},\cite{CKII}. Our feeling is that
our approach is less general than the one in
\cite{Bor} (although, I believe, that they
become rather close, if one uses the
abstract definition of a system (see
\ref{systems}); the approach in
\cite{CKI},\cite{CKII} studies a concrete
renormalization procedure, nevertheless, it
seems that the Connes-Kreimer Hopf algebra
is a rather general phenomenon by means of
which one can identify different
regularizations (=liftings to $<\R>$) of an
OPE-algebra , as was mentioned above.

I hope that  the tools developed in this
paper can help complete  the project
described in \cite{Kor} in a mathematically
rigorous way . The major thing which is
predicted by physicists (i.e. in \cite{Kor})
and which is lacking in this paper is a
construction of a homotopy $d$-algebra
structure on the de Rham complex of the
$\D$-module of observables (we only
construct a Lie bracket).

The main technical tool that we use in this
paper is a $\D$-module structure on the
space of observables. The author started to
appreciate this structure in the process of
reading \cite{BD}.

 In
the case of the free boson the  module of
observables equals $\text{Sym}_{\O_X}
\D_X/\D_X\Delta$, where $X$ is the
space-time and $\Delta$ is the Laplacian.
This module is not free, which prompts using
resolutions and homological algebra.

 I would like to thank
D. Kazhdan for attracting interest to the
problem and valuable conversations without
which this paper would not be written.

 I also thank A.
Beilinson, P.Bressler, A. Cattaneo, V.
Drinfeld, G. Felder, M. Kontsevich,  Yu.
Manin, and B. Tsygan  for helping me and
sharing their ideas.
\section{Content of the paper}
 The paper consists of three parts. In the
 first part we introduce the notion of
 system and the structure of an OPE-algebra
 over a system. We then discuss a naive
 approach to renormalization, the naiveness
 being in ignoring all complications
 stemming from homological algebra. The rest
 of the paper is devoted to constructing a
 homotopically correct (=derived) version of
 this naive construction.  In the second
 part
 we explain the main steps in our
 construction with all technicalities
 omitted. The third part deals with these
 omitted technicalities.
 The first two parts occupy 35 pages, and the
 III-rd part occupies the remaining 75
 pages.

\bigskip

 \centerline{\Large PART I: Systems,
OPE, Naive renormalization}

\section{What is an  OPE?}
Before giving general definition of OPE, we
will introduce this notion  in the setting
of the theory of free boson. The general
notion of an OPE will be obtained via a
straightforward generalization.
\subsection{Notations}
We are going to consider the Euclidean
theory of free boson. Let $Y:=\Re^{2N}$ be
the space-time. We will prefer to work with
the complexification $X=\Co^{2N}$ viewed as
an affine algebraic variety  over $\Co$. Fix
a positively-definite quadratic form $q:Y\to
\Re$. Extend it to $X$ and denote the
extension by the same letter: $q:X\to \Co$.

For a finite set $S$ let $X^S$ be the
algebraic variety which is the product of
$\#S$ copies of $X$. Let $\DDD_{X^S}$ be the
sheaf of algebras of differential operators
on $X^S$. Let $\Cm_S$ be the category of
$\DDD_{X^S}$-sheaves, i.e.
non-quasi-coherent $\DDD_{X^S}$-modules. The
usage of non-quasi-coherent modules is
indispensable in the setting of this paper;
on the other hand, since we are not going to
use any of subtleties of the theory of
$\D$-modules, $\DDD_{X^S}$ sheaves won't
cause any discomfort.
\subsection{Extension of a $\D$-sheaf from a closed subvariety}
The material in this subsection is standard
and can be found for example in \cite{BD}.

Let $i:Y\to Z$ be a closed embedding of
algebraic variety and let $M$ be a
$\D_Y$-sheaf. Let $Y_n$ be the $n$-th
infinitesimal neighborhood of $Y$ in $Z$. It
is well known that $M$ is a crystal, i.e. it
naturally defines  an $\O_{Y_n}$-sheaf;
denote it by $N_n$. Set $i\ha
N:=\liminv_nN_n$; it is a topological
$\DDD_Z$-module, the topology is
$\I_Y$-adically complete, where $\I_Y$ is
the ideal of $Y$. There is a simple explicit
formula for $i\ha M$:
$$
i\ha Y\cong i_.\Hom_{\O_Y}(i^*\D_Z,M),
$$
where $i^*\D_Z$ is the quasi-coherent
inverse image of $\D_Z$ viewed as a
quasi-coherent $\O_Y$-module via the left
multiplication;  $i_.$ is the
sheaf-theoretic extension by zero; the
$\D_Z$-action on $i\ha Y$ is via the right
action on $i^*\D_Z$.

One can prove an analogue of Kashiwara's
theorem in this setting: the functor $i\ha$
is an equivalence of the following
categories: the first category is the
category of $\D_Y$-sheaves; the second
category is the category of $\D_Z$-sheaves
which are sheaf-theoretically supported on
$Y$ and are $\I_Y$-complete, the morphisms
are continuous morphisms. One of the
corollaries is the existence of natural maps
$i\ha(M)\otimes_{\O_Z} N\to
i\ha(M\otimes_{\O_Y} i^*N)$:  Kashiwara's
theorem implies that the right hand side is
the $\I_Y$-adic completion of the left hand
side.

 If
$i,k$ are consecutive embeddings, then $i\ha
k\ha\cong (ik)\ha$.
\subsubsection{} All our closed embeddings
are going to be the embeddings of a
generalized diagonal into some $X^S$. It is
convenient to describe them as surjections
$p:T\to S$. Each such a surjection produces
a closed embedding $i_p:X^S\to X^T$ in the
obvious way.
\subsubsection{}Another feature of the
 $\D$-modules theory  that will be
  used in this paper
is the existence of exterior product
functors
$$
\boxtimes_{a\in A}:\prod_{a\in
A}\Cm_{S_a}\to \Cm_{S},
$$
where $S_a, a\in A$ is a finite family of
finite sets and $S=\disjoint_{a\in A}S_a$.

The functor $i\ha$ is related with the
exterior product in the following way. Let
$p_a:T_a\to S_a$ be a family of projections.
Let $T=\disjoint_{a\in A} T_a$;
$S=\disjoint_{a\in A} S_a$; $p:T\to  S$;
$p=\disjoint_{a\in A} p_a$. We then have a
natural transformation $$\boxtimes_{a\in
A}\prod_{a\in A}i\ha_a\to
i\ha\boxtimes_{a\in A}
$$
both functors act from $\prod \Cm_{S_a}$ to
$\Cm_T$.
\subsection{Construction of functors which
are necessary to define an OPE } Let us now
take into account a specific feature of our
problem: the presence of the quadratic form
$q$ which describes the locus of
singularities of the corellators.   Let $S$
be a finite set and $s\neq t$ be elements in
$S$. Let $q_{st}:X^S\to \C$ be the function
$q(X_s-X_t)$, where $X_s$ are the
coordinates of a point on the $s$-th
component of $X^S$ ( $X_t$ are the
coordinates on the $t$-th copy of $X$.) Let
$\D_{st} $ be the divisor of zeros of
$q_{st}$. Denote by $Z_S:=X^S\backslash
(\cup_{st}\D_{ST})$. Let $j_S:Z_S\to X^S$ be
an open embedding. Set
$\B_S:=j_{S*}\O_{Z_S}$. $\B_S$ is a
$\D_{X^S}$-module.

For a projection $p:T\to S$ set
$\B_p:=\boxtimes_{s\in S}\B_{p^{-1}s}$;
$\B_p$ is a $\D_{X^T}$-module. Set
$\i_p:\C_S\to\C_T$,
\begin{equation}\label{functi}
\i_p(M)=i\ha(M)\otimes_{\O_{X^T}}\B_T.
\end{equation}

List the properties of these functors. First
of all they interact with the exterior
products in the same way as $i_p\ha$. The
behavior under compositions is different.
Let $$R\stackrel q\to T\stackrel p\to S$$ be
consecutive  surjections.  We  then have a
natural transformation $$\as_{pq}:\i_{pq}\to
\i_p\i_q,$$ which is not an isomorphism. Let
us construct $\as_{pq}$. We need an
auxiliary module
$\B_{p,q}=j_{p,q*}\O_{Z_{p,q}}$, where
$j_{p,q}:Z_{p,q}\to X^R$ is an open
subvariety defined by
$$Z_{s,t}=X^R\backslash(\cup_{q(s)\neq q(t)}
\D_{st}) .
$$
 It is clear that $\B_{pq}\cong
\B_{p,q}\otimes \B_q$ and that
$i_{p}^*\B_{p,q}\cong \B_p$. Here $i_q^*$ is
the inverse image for $\O_{X^{S}}$-coherent
sheaves.

Define $\as_{p,q}$ as  the composition
$$
i_{pq}\ha(M)\otimes \B_{pq}\cong i_q\ha
i_p\ha(M)\otimes \B_{p,q}\otimes \B_q\to
i_q\ha(i_p\ha(M)\otimes
i_q^{*}\B_{p,q})\otimes \B_q\cong
\i_q\i_p(M).
$$
\subsubsection{Co-associativity}\label{cosc}
The maps $\as_{p,q}$ have a co-associativity
property. Let
$$
U\stackrel r\to R\stackrel q\to T\stackrel p
\to S
$$
be a sequence of finite sets and their
surjections. We then have two
transformations from $\i_{pqr}\to
\i_r\i_q\i_p$:

the first one is given by
$$
\xymatrix{
\i_{pqr}\ar[rr]{\as_{r,pq}}&&\i_r\i_{pq}
\ar[rr]^{\id\times\as_{q,p}}&&\i_r\i_q\i_p}
$$
and the second one is given by
$$
\xymatrix{
\i_{pqr}\ar[rr]^{\as_{rq,p}}&&\i_{qr}\i_{p}
\ar[rr]^{\as_{r,q}\times\Id}&&\i_r\i_q\i_p}
$$
The co-associativity property  says that
these two transformations coincide.
\subsubsection{}\label{coexp} The maps $\as_{pq}$
interact with the exterior products in the
following way. Let
\begin{equation}\label{RTS} R_a\stackrel
{q_a}\to T_a\stackrel{p_a}\to S_a
\end{equation}
be a family of finite sets and their
surjections. Let $M_a\in \C_{S_a}$, $a\in A$
be arbitrary objects. Let
$$
R\stackrel q\to T\stackrel p\to S
$$
be the disjoint union of (\ref{RTS}) over
$A$.
 Let $$
 M=\boxtimes_{a\in A}M_a\in \C_S.
 $$
Then the following diagram is commutative.
\begin{equation}
\xymatrix{ \boxtimes_{a\in
A}\i_{p_aq_a}M_a\ar[d]\ar[r]&
\boxtimes_{a\in A}\i_{q_a}\i_{p_a}M_a\ar[r]&
\i_q\boxtimes_{a\in A}\i_{p_a}M_a\ar[r]&
\i_q\i_pM\\
 \i_{pq}M\ar[rrru]&&&}
 \end{equation}
\subsubsection{Abstract definition}\label{OPE}
\label{systems}\label{defsystem} We abstract
the properties of the functors $\i_p$.
Assume that for every surjection $p:S\to T$
of finite sets we are given functors
$\j_p:\Cm_T\to \Cm_S$ such that

1) If $p$ is a bijection, then $\j_p$ is the
equivalence of categories induced by $p$;

2) $\j_p$ interact with the exterior
products in the same way as $\i_p$; If all
$p_a$ are bijections, then the corresponding
transformation is the natural one.

3) Let
$$
R\stackrel q\to T\stackrel p\to S
$$
be a sequence of surjections.
 We then have transformations
$$
\as_{pq}:\j_{pq}\to\j_q\j_p.
$$
This transformation is an isomorphism if at
least one of $p,q$ is a bijection. If both
$p,q$ are bijections then the map $\as_{pq}$
is the natural isomorphism of the
corresponding equivalences.

4) The maps $\as_{pq}$ satisfy the
co-associativity property as in
(\ref{cosc}).

5) The maps $\as_{pq}$ interact with the
exterior product in the same way as in
(\ref{coexp}).

If all these properties are the case we say
that that we have {\em a system}.

The functors $\i_p$ and their
transformations form a system which we
denote by $<\i>$.
\subsubsection{Morphisms of
systems} Let $<\j>$, $<\k>$ be systems. {\em
A morphism  of systems} $F:<\j>\to <\k>$ is
a collection of transformations $F_p:\j_p\to
\k_p$ which commute with all elements of the
structure of system.
\subsection{Definition of OPE} With these
functors and their properties at hand we are
ready to define an OPE-algebra.

First of all, we need to fix a $\D_X$-module
$\M$ such that its sections are observables
of our theory. In the case of free boson we
set $\M=S_{\O_X}N$, where
$N=\D_X/\D_X\cdot\Delta$.

 As we know from physics, an OPE is a
prescription of maps
$$
\ope_S: \M^{\boxtimes_S}\to \i_{\pi_S}(\M),
$$
where $\pi_{S}:S\to \{1\}$ is the projection
onto a one-element set. These maps should be
equivariant with respect to bijections $S\to
S'$ of finite sets.

Let us formulate the conditions. It is
convenient to define  maps $\ope_p$ for an
arbitrary surjection $p:S\to T$,
$$
\ope_p: \M^{\boxtimes S}\to
\i_p(\M^{\boxtimes T})
$$ as the composition:
$$
\M^{\boxtimes S}\stackrel{\boxtimes_{t\in T}
\ope_{p^{-1}t}}\to \boxtimes_{t\in
T}\i_{\pi_{p^{-1}t}}(\M)\to
\i_p(\M^{\boxtimes T}).
$$

Let now $R\stackrel q\to T\stackrel p\to S$
be a sequence of surjections of finite sets.
We can define two maps
$$
\M^{\boxtimes R}\to \i_q\i_p \M^{\boxtimes
S}.
$$
The first one is induced by the map
$\as_{q,p}$:
$$
\M^{\boxtimes R}\stackrel{\ope_{pq}}\to
\i_{pq}\M^{\boxtimes
S}\stackrel{\as_{q,p}}\to
\i_q\i_p\M^{\boxtimes S};
$$

the second one is defined by:
$$
\M^{\boxtimes_R}\stackrel {\ope_{q}}\to
\i_q(\M^{T})\stackrel{\ope_p}\to
\i_q\i_p(\M^{\boxtimes S}).
$$
The axiom is that \begin{Axiom} These two
maps should coincide\end{Axiom}

\section{Additional features}
It turns out that the procedure of
renormalization depends on an additional
structure possessed by the  system $\i_{p}$,
which we are going to introduce. The
importance of this structure is not
restricted to the renormalization. The
author believes that this structure also
plays a key role in  formulation of the
quasi-classical correspondence principle
and in the
 connection between the Hamiltonian and Lagrangian formalism.
Thus, let us describe this structure.
\subsection{Preparation}
\subsubsection{The system $\l$} Let
 $\li$ be the operad which describes Lie algebras with
the bracket of degree $1$. Let
$\L(S):=\li(S)^*$ be the linear dual, here
$S$ is a finite set. Let $p:S\to T$ be a map
of finite sets. Set $$ \L(p):=\otimes_{t\in
T }\L(p^{-1}t).
$$

We then have maps $\L(p_1)\otimes \L(p_2)\to
\L(p_1\disjoint p_2)$ and $\L(rq)\to
\L(r)\otimes\L(q)$, where $p_i:S_i\to T_i$;
$r:T\to R$; $q:S\to T$ are maps of finite
sets.

Let now $p:S\to T$ be a surjection. Set
$\l_p:\Cm_{X^T}\to \Cm_{X^S}$;
$\l_p=(i_{p*})\otimes \L(p)$, where
$i_p:X^S\to X^T$ is an embedding determined
by $p$ and $i_{p*}$ is the correspondent
$\D$-module theoretic direct image.
 We then have natural maps
$$
\l_{p_1}(M_1)\boxtimes \l_{p_2}(M_2)\to
\l_{p_1\disjoint p_2} (M_1\boxtimes M_2)
$$
and $\l(rq)\to \l(q)\l(r)$, where
$p_i:S_i\to T_i$; $r:T\to R$; $q:S\to T$ are
maps of finite sets. and $M_i\in
\DDD_{X^{S_i}}$. These maps are induced by
the correspondent maps for $\L$.

Thus, the functors $\l$ possess the
structure which is similar to the one on
$\i$. One sees that all the properties for
$\i$ stated in \ref{defsystem} remain true
upon substituting $\l$ for $\i$. In other
words, $\l$ form a system which we denote by
$<\l>$.

An OPE-algebra  structure over the system
$<\l>$ on a $\D_X$-module $M$ is equivalent
to a *-Lie structure on $M[-1]$ as defined
in \cite{BD}. Let us recall the definition.
\subsubsection{Definition of *-Lie algebra
structure} A *-Lie stucture on a
$\D_X$-module $M$ is given by an
antisymmetric map $b:M\boxtimes M\to i_*M$,
where $i:X\to X\times X$ is the diagonal
embedding. The bracket $b$ is supposed to
satisfy an analogue of Jacobi identity.
\subsubsection{Quasi-isomorphisms of  systems}
Let $<\i>,<\j>$ be  systems and let
$F:<\i>\to<\j>$ be a morphism of systems.
$F$ is a  quasi-isomorphism if   for every
free  $\D_{X^S}$-module $M$ the induced map
$ \j_p(M)\to \i_p(M) $ is a
quasi-isomorphism for every surjection
$p:T\to S$.
\subsubsection{Definition of additional structure
I.}\label{stI} The most important part  of
our additional structure can be then
described as a choice of
 quasi-isomorphisms $<\R>\stackrel\sim\to <\i>$,
 $<\l>\stackrel\sim\to <\m>$  and a map
of systems $<\R>\to <\m>$:
$$\xymatrix{<\i> & <\l>\ar[d]\\
             <\R>\ar[u]\ar[r]&<\m>}
$$

There is even more structure on $<\R>$ which
we shall use.  This part is of some
importance, but not of principal importance,
and will be discussed later (see Sec.
\ref{pre-symm}).

In the rest of the Part I we ignore
homotopy-theoretical complications and
assume that we have a map systems $<\i>\to
<\l>$ (this helps to explain the ideas in a
simper way). A precise exposition will be
given in the subsequent parts of the paper.
Let us now discuss a motivation for the
introduced additional structure.
\subsection{Physical meaning}\label{physmot}
Physical meaning of the introduced
additional structure can be seen from
examining the case when $p:S\to \pt$, where
$S=\{1,2\}$ is a two element set. As a part
of our  structure, we have a  map
$$
\i_p\to i_{p*}\otimes \L(S).
$$
But $\L(S)=k[1]$, therefore, we simply get a
map \begin{equation}\label{mapl} \i_p\to
i_{p*}.
\end{equation}

Recall that
$$\i_p(M)=i\ha(M)\otimes_{\O_{X^S}} \B_S,$$
and one can show that
$$i_{p*}(M)=i\ha(M)\otimes_{\O_{X^S}}
i_{p*}\O_X.
$$

 Assume for simplicity that the map
 (\ref{mapl}) is induced in a natural way by
 a degree +1 map
 \begin{equation}\label{mapl1}
 \B_2\to i_{p*}\O_X
\end{equation}
(we keep in mind the above
identifications).

Such a map specifies an extension $\C_S$
fitting into exact sequence:
$$
0\to i_{p*}\O_X\to \C_S\to \B_S.$$

The meaning of $\C_S$ becomes clear, if we
come back to the real (versus complex)
picture. The global sections of $\B_S$
produce functions on the real part $Y^S$
with singularities on the diagonal. A global
section of $\C_S$ then has a meaning of
distribution on $Y^2$, whose restriction
onto the complement $Y^2\backslash Y$ is a
function from $\B_S$. If we take the space
$C'$ of  all such distributions, we shall
get a slightly larger extensions as the
kernel $C'\to \B_S$ consists of all
distributions supported on the diagonal,
which is larger than $i_{p*}\O_X$.
Nevertheless, it turns out that the space of
global sections of $\C_S$ can be defined as
a subspace of $C'$ (see \ref{cs}).

Set $\I_S:\D_{X}\to\D_{X^S}$ to be
$$
\I_S(M)=i\ha(M)\otimes_{\O_{X^S}}\C_S.
$$
For good $M$ (say flat as $\O_X$-modules),
we have an exact sequence
$$
0\to i_{p*}(M)\to \I_p(M)\to \i_p(M)\to 0.
$$

Let now $M$ be an $OPE$-algebra over $<i>$.
In particular, we have a map
$$
M^{\boxtimes S}\to \i_p(M).
$$
  We may now interpret the composition

$$
M^{\boxtimes S}\to \i_p(M)\to i_{p*}(M).
$$
  of
  this map with the map (\ref{mapl}).
  Assume that $M$ is a  complex of
   free $\D_X$ modules (bounded from above).
  Then we can lift the OPE-map to a map
  $$
  \mu:M^{\boxtimes S}\to \I_p(M)
  $$
with a non-zero differential and the desired
composition is equal to $d\mu$. The
procedure of lifting from $\i_p$ to $\I_p$
is nothing else but the {\em regularization
of divergences}. The map $\mu$ has the
meaning of the commutative product in the
Batalin-Vilkovitski formalism. Its
differential then has a meaning of the
Shouten bracket in the same formalism. This
simple physical argument suggest that the
map $d\mu$ should be a  *-Lie bracket of
degree +1.

\subsection{Geometrical meaning} We will hint
 at the  geometric  meaning of the additional
structure on $<\i>$.  Since our intention is
just to give a motivation,
 the  arguments will not be rigorous.

Recall that the functors $\i$ have been
constructed using the $\D_{X^S}$-modules
$\B_S$, which are defined as sheaves of
functions on certain affine varieties $Z_S$.
Therefore, the de Rham complex of $\B_S$
 computes the cohomology of $Z_S$ shifted by
$\dim Z_S=2ns$, where $2n=\dim_{\Bbb C} X$,
and $s=\# S$. Let
$\B_S':=\i_{\pp_S}(\O_X)=i\ha_{\pp_S}(\O_X)\otimes
\B_S$, Where $\pp_S:S\to \pt$ is the map
onto a point. The de Rham complex of $\B_S'$
computes the cohomology of the intersection
of $Z_S$ with a very small neighborhood of
the diagonal $X\subset X^S$.

On the other hand, $Z_S$ contains as its
real part the  space
 $Z_S^r:=Y^S\backslash
  (\cup_{s\neq t} \DD_{st})$, where $\DD_{st}$
  is the corresponding diagonal.
   Thus we have a map from the de Rham
cohomology of $\B_S'$ to the cohomology of
the intersection of $Z_S^r$ with a very
small neighborhood of the diagonal $Y\subset
Y^S$ in $Y^S$ which can be easily seen to be
the same as the cohomology of $Z_S^r$.  It
is well known that
$H^{(2n-1)(s-1)}(Z_S^r)\cong \L(S)[1-s]$,
where $2n=\dim Y$ and $s=\#S$. The shift on
the right hand side is made in such a way
that both sides have degree zero.

Let us slightly change our point of view.
Instead of  taking the full de Rham complex
let us pick a point $\vs\in S$ and let
$p_{\vs}:X^S\to X$ be the projection onto
 the correponding component. Let $p_{\vs*}(\B_S')$ be the fiber-wise de Rham
complex shifted by the dimension of the
fiber (in this case $H^0p_{\vs*}$ is the
usual $\D$-module theoretic direct image).

We see that the induced map $Z_S\to X$ is a
trivial fibration whose fiber $F_S$ is
homotopy equivalent to $Z_S$ and $\dim
F_S=\dim Z_S-2n$. Let $V$ be a small
neighborhood of $X\subset X^S$
 Then
$$
H^i(p_{\vs*}(\B_S'))\cong \O_X\otimes
H^{2n(s-1)+i}(Z_S\cap V)
$$
and we have a through  map
\begin{eqnarray*}
H^{1-s}(p_{\vs*}(\B_S))\cong \O_X\otimes
H^{(2n-1)(s-1)}(Z_S\cap V) \\
 \to H^{(2n-1)(s-1)}(Z_S^r\cap V)\otimes {\O_X}\to
 \L(s)[1-s]\otimes \O_X.
\end{eqnarray*}
Since $H^{>1-s}(p_{\vs*}(\B_S))=0$, we have
an induced map
$$
p_{\vs*}(\B_S')\to \O_X\otimes \L(S).
$$

It is well known that this map induces a map
$\B_S'\to i_*\O_X\otimes \L(S)$ in the
derived category of $\D_{X^S}$-sheaves.
Thus, the top  cohomology of the
configuration spaces can be interpreted as
maps $\B_S'\to i_*\O_X\otimes \L(S)$. These
maps can be extended to maps $\i_p(M)\to
\l_p(M)$ in the derived category of
$\DDD_{X^T}$-sheaves on $X^T$ for every free
$\DDD_{X^S}$-module $M$.

Of course, this argument is insufficient for
constructing a map of systems (as opposed to
 a collection of maps of  functors $\i_p\to
 \l_p$).

\section{Renormalization. 'Naive' version}
Here we will sketch a scheme for
renormalization ignoring
homotopy-theoretical problems. Although this
naive scheme is of purely heuristic value,
the correct renormalization scheme is in the
same relation to the naive one as derived
functors are to usual ones.

So, we shall simply assume that we are given
a map $<\i>\to <\l>$.

We start with defining the main ingredients.
\subsection{*-Lie structure on $M[-1]$}
Thus, we have a morphism of systems
$<\i>\to<\l>$. Assume that $M$ is an OPE
algebra over $\i$. Then it is also an OPE
algebra over $\l$. i.e. $M[-1]$ is a *-Lie
algebra. Let $\pi:X\to \pt$ be the
projection and denote by $\pi_*M$ the direct
image of $M$;
$$\pi_*M:=\omega_X\otimes_{\D_X}M.
$$
We know that $\g:=\pi_*M[-1]$ is then a DGLA
and this DGLA acts on $M$. Therefore, for
every surjection of finite non-empty sets
$p:S\to T$ we have a $\g$-action on
$$\hom(M^{\boxtimes S}, \i_p(M^{\boxtimes T})).
$$
A very important question for us is whether
the elements $\ope_p$ are $\g$-invariant. It
turns out that in general the answer is no.
We are going to impose an extra axiom which
would guarantee this property.
\def\j{{\i}}
\subsubsection{Extra axiom which ensures the
$\g$-invariance of $\ope_p$.}\label{invar}
Let $p:S\to T$ be a surjection of finite
sets as above. Pick an arbitrary element
$t\in T$; add one more element $\vs$ to $S$
and let $p_t:S\disjoint \{\vs\}\to T$ be a
map which extends $p$ in such a way that
$p_t(\vs)=t$. (this extra element is needed
to take into account the $\g$-action).  Let
$I:S\to S\disjoint \{\vs\}$ be the inclusion
and let $P:X^{S\disjoint \{\vs\}}\to X^S$ be
the natural projection corresponding  to
$I$. Let $P_*$ be the corresponding direct
image.
 We are going to
define several maps $P_*\j_{p_t}\to \j_p[1]$
as follows. Let $s\in S$ be such that
$p(s)=t$. Let $P_s:S\disjoint \vs \to S$ be
the map which is identity on $S$  and
$P_s\vs=s$. Then $p_t=pP_s$.
 We then have the following composition:
$$
P_*\j_{p_t}\to P_*\j_{P_s}\j_{p}\to
P_*\l_{P_s}\j_{p}.
$$
Note that
$$\l_{P_s}\cong
i_{P_s*}\otimes_{s'\in S}\L(P_s^{-1}s')\cong
i_{P_s*}[1].
$$
Thus, we can continue our composition:
$$
P_*\l_{P_s}\j_{p}\to P_*(i_{P_s*}[1])\j_p\to
\j_p[1],
$$
where we used the natural map
$$P_*i_{P_s*}\to \Id_{X^S}.$$
Let
$$
A_s:P_*\j_{p_t}\to \j_p[1]
$$ be the resulting
compositon.

There is one more way to decompose $p_t$.
Let $Q:=p\disjoint\Id:S\disjoint\{\vs\}\to
T\disjoint \{\vs\}$, Let $R:T\disjoint
\{\vs\}\to T$ be  the identity on $T$ and
let $R(\vs)=t$. Then again $p_t=RQ$.
Therefore, we have a composition:
$$
P_*\j_{p_t}\to P_*\j_{Q}\j_{R}.
$$
Let $P_T:X^{T\disjoint \{vs\}}\to X^T$ be
the natural projection. It is not hard to
see that we have an isomorphism
$$
P_*\j_Q\to \j_{P}P_{T*}.
$$
Thus, we continue as follows:
$$
P_*\j_{Q}\j_{R}\to
\j_pP_{T*}\j_R\to\j_pP_{T*}\l_R\cong
\j_pP_{T*}i_{R*}\to \j_p[1].
$$
Denote the composition of these maps by
$$
B_t:P_*\j_{p_t}\to \j_p[1].
$$
Let $C_t=B_t-\sum_{s\in S,p(s)=t}A_s$. Let
us show that the maps $C_t$ determine the
action of $\g$ on $\as_p$. Let $X\in \g$.
Let $L(X):=X.\ope_p$;
$L:\g\to\hom(M^{\boxtimes
S},\j_p(M^{\boxtimes T}))$.
\begin{Claim}\label{inv} $L$ is equal to the following
composition:
\begin{eqnarray*}
\g\stackrel{\oplus_{t\in T}\as_{p_t}}\to
\oplus_{t\in T} \g\otimes\hom(M^{ \boxtimes
S\disjoint\{\vs\}},\j_{p_t}(M^{\boxtimes
T}))\\
\to\oplus_{t\in T}\g\otimes \hom(P_*M^{
\mbox
S\disjoint\{\vs\}},P_*\j_{p_t}(M^{\boxtimes
T}))\cong \oplus_{t\in
T}\g[1]\otimes\hom(\g[1]\otimes M^{\boxtimes
S},P_*\j_{p_t}(M^{\boxtimes T}))\\
\to \oplus_{t\in T}\hom(M^{\boxtimes
S}[1],P_*\j_{p_t}(M^{\boxtimes
T}))\stackrel{\oplus C_t}\to
\hom(M^{\boxtimes S}[1],\j_{p}(M^{\boxtimes
T})[1]).
\end{eqnarray*}
 Therefore, if $C_t=0$ for all $t$, then
 $L=0$.
\end{Claim}
\pf Straightforward.
\endpf
\subsubsection{}  Call a system $<\j>$ endowed with a map
$<\j>\to <\l>$ {\em invariant} if all
$C_t=0$.
\subsubsection{Another axiom}
\label{pre-symm}\label{invsym} It turns out
that to construct a good theory one has to
introduce a one more  natural axiom on
$<\i>$. The importance  of this  axiom can
be fully appreciated only when one passes to
a more precise consideration.

Let us describe this axiom.
 Let $A=\{1,2\}$ be a two
element set. Let $q:A\to \pt$. Let $p:S\to
T$ be a surjection. Let $i:\pt \to A$ be the
inclusion in which a unique element $\pt$
goes to 1. For an  arbitrary injection
$j:U\to V$ let $p_j:X^V\to X^U$ be the
corresponding projection and $\p_i:\Cm_U\to
\Cm_V$ be the corresponding ${\cal
D}$-module-theoretic direct image.

We then construct two maps
$$
\p_{\Id_S\disjoint i}\i_{p\disjoint q}\to
\i_{p\disjoint \Id_{\pt}}.
$$
The first map is as follows:
\begin{eqnarray*}
 M_I:\p_{\Id_S\disjoint
i}\i_{p\disjoint q}\to \p_{\Id_S\disjoint
i}\i_{\Id_S\disjoint q}\i_{p\disjoint
\Id_{\pt}}\\
\to \p_{\Id_S\disjoint
i}\delta_{\Id_S\disjoint q}\i_{p\disjoint
\Id_{\pt}}\cong \i_{p\disjoint \Id_{\pt}}
\end{eqnarray*}
and the second one is:
\begin{eqnarray*}
M_{II}: \p_{\Id_S\disjoint i}\i_{p\disjoint
q}\to \p_{\Id_S\disjoint i}\i_{p\disjoint
\Id_{\pt}}\i_{\Id_S\disjoint q}\\
\to \p_{\Id_S\disjoint i}\i_{p\disjoint
\Id_{\pt}}\delta_{\Id_S\disjoint q}\cong
\i_{p\disjoint \Id_{\pt}} \\
\to
\i_{p\disjoint \Id_{\pt}}\p_{\Id_S\disjoint
i}\delta_{\Id_S\disjoint q},
\end{eqnarray*}
where we have used a natural isomorphism
$$
\p_{\Id_S\disjoint i}\i_{p\disjoint
\Id_{\pt}}\cong \i_{p\disjoint
\Id_{\pt}}\p_{\Id_S\disjoint i}.
$$

 Call a system $<\j>$ endowed with a map
$<\j>\to <\l>$ to be {\em pre-symmetric} if
 $M_I=M_{II}$ for all $p$.

 Finally, call a system {\em symmetric}, if
 it is  both pre-symmetric and invariant.
\subsubsection{} What is the situation with
the system $<\R>$ that we are going to
construct? It turns out, that up-to
homotopies, it is pre-symmetric, but not
symmetric. Pre-symmetricity is the
additional structure on $<\R>$ which was
mentioned in (\ref{stI}).

The above reasoning suggests that
renormalization is only possible in
symmetric (or, at least, invariant systems).
Therefore, a procedure of "fixing" $<\R>$
(which we call "symmetrization") is needed
to perform a renormalization. We shall
discuss a naive version of such a
symmetrization after a more detailed
explanation how renormalization goes on in a
symmetric system.
\subsection{Renormalization in a
symmetric system} As was mentioned, the
system $<\R>$ that we will construct in the
example of free boson is not symmetric.
Nevertheless, to appreciate the importance
of symmetricity, we will explain in the next
section that were $<\R>$ symmetric, the
renormalization of any OPE-algebra over
$<\R>$ could be defined in a very simple
fashion.

 Let $M$ be an
OPE-algebra in a {\em symmetric} system
$<\j>$. Then, in virtue of the map $\j\to
\l$, $M$ is also an OPE-algebra in $\l$,
i.e.  M[-1] is a *-Lie algebra.Let $\pp:X\to
\pt$ be the projection onto a point. Then
$\pp_*M[1]$ is a DGLA. Let $\h$ be a formal
variable (the 'interaction constant'). Pick
a Maurer-Cartan element
$$\S\in \h(\pp_*M[-1])[[\h]]^1=\h\pp_*M[[\h]];\
d\S+1/2[\S,\S]=0.
$$
This equation is called {\em quantum Master
equation.} Using $\S$ we can perturb the
differential on $M$; let
$M':=(M[[\h]],d+[\S,\cdot])$ be the
corresponding differential graded
$\D_{X[[\h]]}$-module.

The renormalization is the procedure of
constructing a $\Co[[\h]]$-linear OPE
structure over $<\j>$ on $M'$. In our
setting this procedure is trivial. Indeed,
since $M'=M[[\h]]$ as graded objects; the
OPE structure on $M$ gives rise to the maps
$$
\ope'_p:(M')^{\boxtimes_{\Co[[\h]]} S}\to
\j_p((M')^{\boxtimes_{\Co[[\h]]} T}).
$$
The $\l$-invariance of $\j$ and Claim
\ref{inv} imply that these maps are
compatible with the differential on $M'$.
Thus, $\ope'_p$  do define the renormalized
OPE on $M'$.
\def\ii{\i^{\text{symm}}}
\subsection{An idea how to fix
non-ivariance of $<\j>$: symmetrization. }
Let us try to define a system $\ii$ endowed
with a map $<\j>\to <\ii>$ such that in
$\ii$ all $C_t=0$. Then our OPE-algebra $M$
in $<\j>$ determines an OPE-algebra in
$<\ii>$ and the renormalization of this
algebra goes the way as was described above.

The obvious way to define $<\ii>$ is to
simply put
$$\ii_p:=\j_p/\text{Span}<\text{Im}
C_t>_{t\in T}.
$$
One checks that  the structure of system on
$\j$ is naturally transferred onto $\ii_p$.
\subsection{
  Summary}
Let us first summarize what we have done.

We start with a system $<\j>$ which is
quasi-isomorphic  to the original system
$<\i>$ and is endowed with a map $<\j>\to
<\l>$. We then construct a symmetric system
$<\ii>$ which fits into the diagram $<\j>\to
<\ii>\to <\l>$. Thereafter, having an OPE
algebra $M$ over $\j$, we observe that
$\pp_*M[-1]$ is a DGLA and we pick a
Maurer-Cartan element $\S\in
\h\pp_*M[[\h]][1]$. We then define the
$\D_{X[[\h]]}$-module $M'$ and define an OPE
structure on $M'$ over $\ii$.

What has to be done for this scheme to
really work?

Problem 1. We need to construct $\j$ with
the specified properties.

Problem 2. We have an OPE algebra $M$ over
$<\i>$ and a quasi-isomorphism $<\R>\to
<\i>$. We need to lift $M$ to an OPE algebra
over $<\R>$.

Problem 3. The passage from  $<\R>$ to
$<\Rs>$ is not stable under
quasi-isomorphism of systems. Thus we need
to develop a derived version of the map
$<\R>\mapsto <\Rs>$

Problem 4. After all, we get a renormalized
OPE-algebra in an abstract system $<\Rs>$.
To give a physical meaning to this system,
we have to find a construction which
transforms this OPE-algebra  into
OPE-products in terms of series of
real-analytic functions on $Y^S$.
\subsection{Plan for the future exposition}

The rest of the paper is devoted to  solving
these problems. As this involves a lot of
technicalities, we shall first retell the
content of the paper omitting them. Then the
detailed exposition, with proofs, will
follow.

First,  we shall formulate the list of
properties that the system $<\R>$, to be
constructed, should possess. These
properties form a homotopical variant of the
definition of the structure of pre-symmetric
system. Every system possessing  these
properties will be called pre-symmetric
(this should not lead to confusion with the
naive definition  of pre-symmetricity).

Secondly, we shall show how the
renormalization can be carried over for
OPE-algebras over a pre-symmetric system
$<\R>$ (including a construction for
symmetrization of $<\R>$ and a construction
of the renormalized  OPE-algebra over the
symmetrized system) . These steps constitute
a homotopically correct version of the above
outlined naive approach. Thereafter, we
construct a  pre-symmetric system $<\R>$
which is a resolution of the system $<\i>$.

To renormalize an OPE-algebra over $<\i>$
one has to be able to lift it to an
OPE-algebra over $<\R>$ so that the lifting
be compatible with the quasi-isomorphism  of
systems $<\R>\to <\i>$. This happens to be a
variant of the celebrated
Bogolyubov-Parasyuk theorem, saying that
such a lifting is always possible. An
analogous theorem can be shown by a
homotopy-theoretical non-sense, using the
quasi-isomorphicity of the map $<\R>\to
<\i>$; but for this to work one has to
replace the stucture of OPE-algebra up-to
higher homotopies.  Let us stress that
Bogolyubov-Parasyuk theorem produces a
lifting of  usual OPE-algebras, which is a
stronger statement. Homotopical approach, on
the other hand, provides for a homotopical
equivalence of two different liftings. These
homotopy-theoretic questions will be
discussed in a subsequent paper.

Finally, we solve the  Problem 4.

The exposition will be organized in such a
way that the most difficult technical
moments will be omitted at the "first
reading" , which is  the part II, and will
be discussed at the "second reading" (i.e.
the concluding part III).

\def\j{{\frak j}}

\bigskip

 \centerline{\bf \Large PART II: Exposition
without technicalities}

\def\Z{{ Y}}
\def \R{{\Bbb R}}
\def \aa{{\frak a}}
\def \a{{\frak a}}
\def \Re{\text{Re}}
\def \D{{\cal D}}
\def \u{{\frak u}}
\def\ii{{\frak ii}}
\def\ps{{\frak sys}^p}
\def \psm{{\frak psm}}
\def\system{{\frak sys}}
\def \Mg{{\frak M}}
\def\ope{{\frak ope}}
\def\h{\lambda}
\def \hh{{\frak h}}
\def \p{{\frak p}}
\def \g{{\frak g}}
\def \kk{{\frak kk}}
\def \f{{\frak f}}
\def \gen{{\frak gen}}
\def\univ{\text{univ}}
\def\overset{\supset}
\def\De{{\Delta}}
\def \S{{\frak D}}
\def \Co{{\Bbb C}}
\def \P{{\frak P}}
\def \Di{\text{Diag}}
\def\vect{\text{vect}}
\def \DD{{\frak D}}
\def \DDD{\Di}
\def \Gr{\text{Gr}}
\def \B{B}
\def\ha{{\wedge}}
\def \as{{\frak as}}
\def\Eq {\text{\bf Eq}}
\def \C{{\cal C}}
\def \U{{\cal U}}
\def \t{{\frak t}}
\def\disjoint {\sqcup}
\def \Ker{\text{Ker}}
\def\Prod{\prod}
\def\res{\text{res}}
\def\ve{\epsilon}
\def\circD{{\stackrel{\circ}{\Delta}}}
\def\Y{y}
\def\Diag{{\text{ Diag}}}
\def\Dg{{\frak Diag}}
\def\FM{{\frak FM}}
\def\zebra{{\bf Zebra}}
\def\B{{\cal B}}
\def\R{{\cal R}}
We shall pass from a naive approach to the
realistic one, in which the naive
definitions sketched above will be replaced
with appropriate homotopically correct
versions.

Our plan is as follows. In the following
section we give a homotopically correct
definition of pre-symmetric system.

Next we show show how, having an OPE-algebra
in a pre-symmetric system, one can
renormalize it.

 Next we have to show these definitions work
 in the example of free scalar boson. The
 major part of the required work is done in
 part III, in this part we only sketch the  main
 steps which are:

 1) we have to construct a pre-symmetric
 system $<\R>$ which maps
 quasi-isomorphically to the system $<\i>$;

 2) we have to show that every OPE-algebra
 over $<\i>$ lifts to an OPE-algebra over
 $<\R>$.

 Having done this we can apply the
 symetrization and renormalization
 procedures.

 3) And finally, we need to  be able to
 interpret the renormalized OPE in the
 symmetrized system in terms of expansions
 whose coefficients  are real-analytic
 functions on $Y^n$ without diagonals.

 So, let us follow our plan.

\section{Pre-symmetric systems}

In this section we shall give a homotopy
version of the notion of pre-symmetric
system (see \ref{invsym} for naive version).

The plan is as follows. We shall give two
slightly different (and slightly
non-equivalent) definitions of a homotopy
analog of a pre-symmetric system. Any
pre-symmetric system in the sense of the
first definition will naturally produce a
pre-symmetric system in the sense of the
second definition. The first definition is
given in terms of  functors $\R_p,\delta_p$,
in the second definition we  replace the
functors $\delta_p$ with functors of direct
image with respect to all projections
$X^S\to X^T$. We will see that the second
definition looks more natural. Moreover, the
second definition encloses all the structure
needed for symmetrization and
renormalization. So, we consider the second
definition as a more basic one. On the other
hand, to define a pre-symmetric system in
the example of free boson, we shall use the
first definition.

We start with formulation of the first
definition. First of all, we need to provide
for a homotopy-theoretical analog of a map
$<\i>\to <\l>$. This will be achieved via
replacement of $<\l>$ with a
quasi-isomorphic system
$<\l>\stackrel\sim\to<\m>$. We shall give
the definition of such an $<\m>$. A part of
a structure of pre-symmetric system on a
system $<\R>$  will then be a  map $<\R>\to
<\m>$. As was mentioned in \ref{invsym}, to
be pre-symmetric, the system $<\R>$ should
have additional properties. We will give
their  homotopical versions. This will
accomplish the first definition of a
pre-symmetric system. Finally, we formulate
the second definition  (which is essentially
a paraphrasing of the first definition in
terms of direct image functors with respect
to projections), it will then follow
automatically that every pre-symmetric
system in the sense of the first definition
gives rise to a pre-symmetric system in the
sense of the second definition.

\subsection{A homotopy version of the map
$<\i>\to <\l>$}\label{Rm} As was explained
above, the first step we need to do is to
endow the system $<\i>$ with a map of
systems $<\i>\to <\l>$. We shall do it in a
homotopical sense, i.e. we shall construct
systems $<\R>$ and $<\m>$ fitting into the
following commutative diagram:

$$
\xymatrix{ <\i> & <\l>\ar[d]^\sim\\
           <\R>\ar[u]^\sim\ar[r] &<\m>}
           $$
The vertical arrows should be
quasi-isomorphisms.

 Let us first define the
system $< \m>$.
\subsubsection{The system $<\m>$}

Let us define the complex $\m_p$ centered in
strictly  negative degrees by setting

\begin{equation}\label{forma1}
\m_p^{-n}=\oplus\delta_{p_1}\delta_{p_2}\cdots\delta_{p_n},
\end{equation}
where the direct sum is taken over all
diagrams
\begin{equation}\label{forma2}
S\stackrel{p_1}\to S/e_1\stackrel{p_2}\to
S/e_2\stackrel{p_3}\to\cdots\stackrel{p_{n-1}}\to
S/e_{n-1} \stackrel{p_n}\to T,
\end{equation}
where $\omega>e_1>e_2>\cdots>e_n>e$, where
$e$ is  the equivalence relation induced by
$p$ and $p_i$ are natural projections. The
differential is given by the alternated sum
$d=D_1-D_2+\cdots+(-1)^nD_{n-1}$ , where

$$D_i:
\delta_{p_1}\delta_{p_2}\cdots
\delta_{p_n}\to
\delta_{p_1}\delta_{p_2}\cdots
\delta_{p_{i-1}}\delta_{p_{i+1}p_i}\delta_{p_{i+2}}\cdots
\delta_{p_n}
$$
is induced by the isomorphism
$$\delta_{p_i}\delta_{p_{i+1}}\to
\delta_{p_{i+1}p_i},
$$

The maps $\as_{q,r}:\m_{rq}\to \m_q\m_r$ are
defined in the following natural way. Let
$p=rq$. Let $f$ be the equivalence relation
on $S$ determined by $q$ and $e$ be the
equivalence relation determined by $p$ so
that $f>e$. One can assume  that $S\stackrel
q\to S/f\stackrel r\to S/e$.

The map $\as_{q,r}$ restricted to
$$
\delta_{p_1}\delta_{p_2}\cdots \delta_{p_n}
$$
as in (\ref{forma1}), (\ref{forma2}), vanish
unless there exists a $k$ such that $e_k=f$,
in which case it  isomorphically maps this
term into
$$
(\delta_{p_1}\delta_{p_2}\cdots
\delta_{p_{k}})(\delta_{p_{k+1}}\cdots
\delta_{p_n})
$$

 The  factorization maps $$\fact:\boxtimes_{a\in
A}\m_{p_a}(M_a)\to \m_p(\boxtimes_a M_a)$$
are given by a "shuffle product". Here is
the construction.

Fix direct summands of $\m_{p_a}$:
$$\delta_{p_{1a}}\delta_{p_{2a}}\cdots\delta_{p_{n_aa}},
$$
where $p_{ia}:S_{ia}\to S_{i+1a}$,
 and define the
restriction of the factorization map onto
them.

Define {\em a shuffle} as a sequence
$$\alpha:=(a_1,a_2,\ldots,a_N)$$
where

- $a_k\in A$;

- $a_k$ enters  into the sequence
$a_1,a_2,\ldots a_N$ exactly $n_k$ times.

Given such a shuffle, let $\alpha_k(a)$ be
the number of times $a$ enters into the
subsequence $a_1,a_2,\cdots,a_k$.

Let $$S_k^\alpha:=\disjoint_{a\in A}
S_{\alpha_k(a)a}.$$

Define the map $$p_k^\alpha:S_k^\alpha\to
S_{k+1}^\alpha$$

as $$p_{\alpha_k(a)a}\disjoint_{a'\neq
a}\Id_{S_{\alpha_k(a')a'}}$$

We then have a natural map
\begin{eqnarray*}
\fact(\alpha):\boxtimes_{a\in
A}\delta_{p_{1a}}\delta_{p_{2a}}\cdots\delta_{p_{n_aa}}(M_a)
\stackrel\sim\to\\
\delta_{p_1^\alpha}\delta_{p_2^\alpha}\cdots
\delta_{p_N^\alpha}(\boxtimes_{a\in A}
M_a)\to \m_p(\boxtimes_a M_a).
\end{eqnarray*}

Set the  restriction of the map $\fact$ onto
$$
\boxtimes_{a\in
A}\delta_{p_{1a}}\delta_{p_{2a}}
\cdots\delta_{p_{n_aa}}(M_a)
$$
to be equal to
$$
\sum_\alpha
(-1)^{\text{sign}(\alpha)}\fact(\alpha),
$$
where $\text{sign}(\alpha)$ is the sign of
the shuffle.

 Denote by
$$
l^\m_p: \m_p\to\delta_p[1]
$$
the natural projection

Then a map  of systems $<\R>\to <\m>$ is
uniquely determined by  the knowledge of
compositions
\begin{equation}\label{lp1}
l_p:\R_p\to \m_p\to \l_p.\end{equation}

In the sequel we will work with these maps
rather than with the system $<\m>$.
\subsubsection{A quasi-isomorphism $<\l>\to
<\m>$} As a part of our program, we have to
define a quasi-isomorphism $<\l>\to <\m>$.
As it won't be used in the future,  we shall
give a very brief description.

It is not hard to see that the cohomology of
any complex $\m_p$ is concentrated in its
lowest degree (i.e. $\#T-\#S$, where $p:S\to
T$; and it is not hard to see that this
cohomology is isomorphic to $\l_p$, whence
the maps $\l_p\to\m_p$. The axioms for a map
of systems can  be easily  checked .
\subsubsection{First definition of
pre-symmetric system} As a part of the
structure of a pre-symmetric system (in the
sense of the first definition) we should
include maps (\ref{lp1}) which provide for a
homotopy-theoretical substitute for a map of
systems $<\i>\to <\l>$.
 To complete the definition we should
 add a structure which is a
homotopical analog  of properties
\ref{invsym}.  After we formulate this
structure, we will formulate the axioms
which should be satisfied by the elements of
the structure. This will complete the first
definition of a pre-symmetric structure.

 We shall start with the most natural
piece of structure. Let $\phi:S\to T$ and
$g:A\to B$ be surjections. Then we should
have a natural map
\begin{equation}\label{ptas}
\delta_{\Id_S\disjoint
\phi}\R_{\phi\disjoint \Id_B}\to
\R_{\phi\disjoint\Id_A}\delta_{\Id_T\disjoint
g}.
\end{equation}

Such a natural map also exists if one
replaces $<\R>$ with $<\i>$.

Indeed: $$\delta_{\Id_S\disjoint
\phi}\i_{\phi\disjoint \Id_B}(M)\cong
i_{\phi\disjoint
g}^\ha(M)\otimes_{\O_{X^{S\disjoint
A}}}({\cal  B }_{\phi}\boxtimes
i_{g*}\O_{X^B}),
$$
whereas

\begin{eqnarray*}
\i_{\phi\disjoint\Id_A}\delta_{\Id_T\disjoint
g}(M) \\ \cong \B_{\phi\disjoint
\Id_A}\otimes_{\O_{X^{S\disjoint A}}}
i_{\phi\disjoint
\Id_A}^\ha((\O_{X^T}\boxtimes
i_{g*}\O_{X_B})\otimes_{\O_{X^{T\disjoint
A}}} i^\ha_{\Id_T\disjoint g}(M)),
\end{eqnarray*}

and we see that  the right hand side in
(\ref{ptas}) is the completion of the left
hand side, whence the desired map.

The corresponding map for $\R$ is
constructed following the same principles.

The next piece of structure is more subtle
and is given by a family of maps
$$\R_{\phi\disjoint g}\to \R_{\phi\disjoint
\Id_A}\delta_{\Id_T\disjoint g},
$$
where $\phi:S\to T$ and $g:A\to B$ are
arbitrary surjections. The comparison of
this additional structure with the naive
structure  will be given after we list the
axioms satisfied by $\l_p$ and $s(\phi,g)$.
A presymmetric structure in the sense of the
first definition is then a collection of
maps $\l_p$ and $s(\phi,g)$ satisfying the
axioms formulated below.
\subsection{Axioms of the pre-symemtric
system (in the sense of the first
definition)}
\subsubsection{ Properties of the maps
$\l_p$}  The properties of $<\l_p>$ we are
going  to list  simply express the fact that
 the collection of maps $\l_p$ should
define a map of systems $<R>\to <\m>$.

{\bf Property 1} If $p$ is a bijection, then
$\l_p=0$.

{\bf Property 2} Let $f_i:S_i\to T_i$ be
nontrivial surjections. Then the composition
$$
\R_{f_1}(M_1)\boxtimes \R_{f_2}(M_2)\to
\R_{f_1\disjoint f_2}(M_1\boxtimes
M_2)\stackrel{\l_{f_1}\boxtimes
\l_{f_2}}\longrightarrow
\delta_{f_1\disjoint f_2}(M_1\boxtimes M_2)
$$
is zero.

If $f_1$ is a bijection, then the above
composition equals
$$
\R_{f_1}(M_1)\boxtimes
\R_{f_2}(M_2)\stackrel{a\boxtimes
\l_{f_2}}\longrightarrow
\delta_{f_1}(M_1)\boxtimes
\delta_{f_2}(M_2)\to \delta_{f_1\disjoint
f_2}(M_1\boxtimes M_2),
$$
where we used the isomorphism
$a:\R_{f_1}\to\delta_{f_1}$ for a bijective
$f_1$.

{\bf Property 3}  Define the differential
$d\l_p$.

Let $p:S\to T$ and let $e$ be the
equivalence relation on $S$ determined by
$e$. Let $e_1$ be a  strictly finer
nontrivial  equivalence relation. Set
$p_1:S\to S/e_1$, $p_2:S/e_1\to T$ to be the
natural projections so that $p=p_2p_1$. Set
$$
\l(e_1):\R_p\to
\R_{p_1}\R_{p_2}\stackrel{\l_{p_1},\l_{p_2}}\longrightarrow
\delta_{p_1}\delta_{p_2}\cong \delta_p.
$$

 We then have
 $$
 dl_p+\sum\limits_{e_1}\l(e_1)=0,
 $$
 where the sum  is taken over all nontrivial
 equivalence relations on $S$ which are
 strictly finer than $e$.
 \subsubsection{Properties of maps
 $s(\phi,g)$.}
 {\bf Property 1}

 The following diagram is
commutative:
$$
\xymatrix{ \R_{\phi\disjoint g\disjoint
h}\ar[rrd]^{s(\phi\disjoint g,h)}
\ar[rr]^{s(\phi,g\disjoint h)}&&
\R_{\phi\disjoint
g\disjoint\Id}\delta_{\Id\disjoint\Id\disjoint
h}
 \ar[d]^{s(\phi,g)}\\
&&
\R_{\phi\disjoint\Id\disjoint\Id}\delta_{\Id\disjoint
g\disjoint h}}
$$

{\bf Property 2}

 Assume that $\phi$ is not
bijective.  Then the  composition
$$
\R_{\phi\disjoint
g}\stackrel{s(\phi,g)}\longrightarrow
\R_{\phi\disjoint\Id}\delta_{\Id\disjoint
g}\stackrel{\l_{\phi\disjoint
\Id}}\longrightarrow \delta_{\phi\disjoint
\Id}\delta_{\Id\disjoint g}$$ equals
$$
\R_{\phi\disjoint
g}\stackrel{\l_{\phi\disjoint
g}}\longrightarrow \delta_{\phi\disjoint g}.
$$

If $\phi $ is bijective and $g$ is not, then
the above composition vanishes.

If both $\phi$ and $g$ are bijections, then
the above composition equals the natural
identification of the right and left hand
sides.

{\bf Property 3} Let $g=g_2g_1$, where
$g_1,g_2$ are surjections.
 Introduce a map
$$K(\phi_1,\phi_2,g_1,g_2):\R_{\phi_2\phi_1\disjoint
g_2g_1 }\to \R_{\phi_1\disjoint
g_1}\R_{\phi_2\disjoint g_2}
\stackrel{s(\phi_1,g_1),s(\phi_2,g_2)}{\longrightarrow}
\R_{\phi_1\disjoint\Id}\delta_{\Id\disjoint
g_1}\R_{\phi_2\disjoint\Id}\delta_{\Id\disjoint
g_2}\stackrel{(\ref{ptas})}\to
\R_{\phi_1\disjoint\Id}\R_{\phi_2\disjoint\Id}
\delta_{\Id\disjoint g_2g_1}.
$$

The property then says:
 The map
 $$\R_{\phi_2\phi_1\disjoint g}\stackrel{s(\phi_2\phi_1,g)}
 \longrightarrow
\R_{\phi_2\phi_1\disjoint\Id}\delta_{\Id\disjoint
g}\to
\R_{\phi_1\disjoint\Id}\R_{\phi_2\disjoint\Id}
\delta_{\Id\disjoint g}
$$
is equal to
$$
\sum_{g_2g_1=g} K(\phi_1,\phi_2,g_1,g_2),
$$
where the sum is taken over all diagrams
\begin{equation}\label{sett}
A\stackrel{g_1}\to A/{e_1}\stackrel{g_2}\to
B,
\end{equation}
where $e_1$ is an arbitrary equivalence
relation on $A$ such that $g$ passes through
$A/e$, and $g_1,g_2$ are the natural
surjections.

 {\bf Property 4}

 The following diagram is commutative:
$$
\xymatrix{ \R_{f_1\disjoint
g_1}(M_1)\boxtimes \R_{f_2\disjoint
g_2}(M_2)\ar[rr]^{s(\phi_1,g_1)\boxtimes
s(\phi_2,g_2)} \ar[d]&&\R_{f_1\disjoint
\Id}\delta_{\Id\disjoint g_1}(M_1)\boxtimes
\R_{f_2\disjoint \Id}\delta_{\Id\disjoint
g_2}(M_2)\ar[d]\\
\R_{f_1\disjoint f_2\disjoint g_1\disjoint
g_2}(M_1\boxtimes
M_2)\ar[rr]^{s(f_1\disjoint f_2,g_1\disjoint
g_2)} && \R_{f_1\disjoint
f_2\disjoint\Id}\delta_{\Id\disjoint
g_1\disjoint g_2}(M_1\boxtimes M_2)}
$$

{\bf Property 5} Denote
$$s(g_1,\phi,g_2):\R_{\phi\disjoint
g}\to \R_{\Id\disjoint g_1}\R_{\phi\disjoint
g_2}\stackrel{\l_{\Id\disjoint
g_1}}\longrightarrow \delta_{\Id\disjoint
g_1}\R_{\phi\disjoint\Id}\delta_{\Id\disjoint
g_2}\stackrel{(\ref{ptas})}\to
\R_{\phi\disjoint\Id}\delta_{\Id\disjoint
g};
$$
$$
s(\phi,g_1,g_2):\R_{\phi\disjoint g}\to
\R_{\phi\disjoint g_1}\R_{\Id\disjoint
g_2}\stackrel{\l_{\Id\disjoint g_2}}\to
\R_{\phi\disjoint\Id}\delta_{Id\disjoint
\g_1}\delta_{\Id\disjoint g_2}\to
\R_{\phi\disjoint\Id}\delta_{\Id\disjoint
g}.
$$

The property  asserts that

$$
ds(\phi,g)=\sum\limits_{g=g_2g_1}(s(g_1,\phi,g_2)-s(\phi,g_1,g_2))
$$
where the sum is taken over the same set as
in (\ref{sett}).

\subsubsection{ Comment on the meaning of
$s(\phi,g)$} To see this meaning consider a
special $g:A\to\pt$, where $A=\{1,2\}$, and
$\phi:S\to T$ is a surjection.   Calculate
the differential $ds(\phi,g)$.

It is equal to the difference $A-B$ of two
maps, where
$$
A:\R_{\phi\disjoint g}\to \R_{\phi\disjoint
\Id}\R_{\Id\disjoint g}\stackrel{
\l_{\Id\disjoint g}}\longrightarrow
\R_{\phi\disjoint \Id}\delta_{\Id\disjoint
g}
$$
and
$$
B:\R_{\phi\disjoint g}\to \R_{\Id\disjoint
g}\R_{\phi\disjoint \Id}\to
\delta_{\Id\disjoint g}\R_{\phi\disjoint
\Id}\to \R_{\phi\disjoint
\Id}\delta_{\Id\disjoint g}.
$$

Thus, the maps $s(\phi,g)$ provide for the
difference $A-B$ to be homotopy equivalent
to zero (up-to higher homotopies).

Let $j:S\disjoint \{1\}\to S\disjoint
\{1,2\}$ be the obvious inclusion.
 Composing $A-B$ with $\p_j$, we see that
$\p_jA=M_I$, $\p_jB=M_{II}$ as in Sec.
\ref{pre-symm}. Thus the maps $s(f,g)$ are
responsible for a homotopy analog of
pre-symmetricity of $<\R>$.

In the  next subsection the above described
structure will be reformulated in terms of
functors of direct image with respect to
projections. This will constitute a basis
for the further exposition.
\subsection{Reformulation in terms of
direct images with respect to projections:
second definition of a pre-symmetric system}
\label{defpresymm} Recall that the main
ingredient in the renormalization procedure
is an element of $p_*M$, where $p:X\to \pt$
is a projection. Thus we have to incorporate
into our picture direct images with respect
to projections. Let $i:S\to T$ be an
injection. It induces a projection
$p_i:X^T\to X^S$. Let
$\p_i:\Cm_{X^T}\to\Cm_{X^S}$ be the
corresponding $\D$-module theoretic direct
image. We want to incorporate it into our
picture and to describe the maps which can
be defined on superpositions of various
$\R_p$ and $\p_i$. These maps will be
derived from the maps $\l_p$ and $s(f,g)$.
Note that the direct images with respect to
injections are not applied, they are only
used to {\em produce} maps between different
iterations of $\R_p$ and $\p_i$.

Thus, we shall now describe these maps and
their properties.
\subsubsection{}
The map we shall describe here is somewhat
similar to (\ref{ptas}).

 Let $q:S\to T$ be an surjection
and $U$ be a finite set.
 Consider the following commutative diagram
\begin{equation}\label{pertas}
 \xymatrix{S\disjoint U\ar@{>>}[r]^p&
T\disjoint U\\
S\ar@{>>}[r]^q\ar@{^{(}->}[u]^i&
T\ar@{^{(}->}[u]^j}
\end{equation}
where $p=q\disjoint \Id$, and $i,j$ are the
natural injections. Then we have an
isomorphism
\begin{equation}\label{pertas1}
\p_i\R_p\to \R_q\p_j
 \end{equation}

  One can see that such an isomorphism
 is naturally defined, if we replace $\R$
 with $\i$.

\subsubsection{} Using the
maps $\l_p:\R_p\to\delta_p$, we can do the
following.

  Consider a commutative triangle
$$\xymatrix{S\ar@{>>}[r]^p& T\\
            R\ar@{^{(}->}[u]^i\ar@{^{(}->}[ur]^j&}
            $$
in which $i,j$ are injections and $p$ is a
proper surjection. We then have a degree +1
map
$$
L(i,p):\p_i\R_p\to  \p_j
$$
given by
$$
\p_i\R_p\to\p_i\delta_p\cong  \p_j
$$
\subsubsection{} Let us now "translate"
$s(f,g)$ into our new language.
  Consider a commutative
square
\begin{equation}
\xymatrix{ R\ar@{->>}[r]^p & T\\
           S\ar@{->>}[r]^q\ar@{^{(}->}[u]^i& P\ar@{^{(}->}[u]^j}
\end{equation}
in which $i,j$ are injections and $p,q$ are
proper surjections. Let $T_1=T\backslash
T_2$ be the subset of all $t\in T$ such that
$p^{-1}t\cap i(S)$ consists of $\geq 2$
elements.

 Call such a square {\em suitable}
if the following is satisfied:

$p^{-1}(T_1)\subset i(S)$, i.e.:
$$
\#(p^{-1}t\cap i(S))\geq 2\follows
p^{-1}(t)\subset i(S).
$$

 We then have a degree zero
map
$$
A(i,p,j,q):\p_i\R_p\to \R_q\p_j.
$$

Construction:  Decompose $T=T^1\disjoint
T^2$, where $T^1$ consists of all  $t\in T$
such that $p^{-1}t\subset i(S)$ (so that
$T_1\subset T^1$). Set
$$R^n=p^{-1}T^n,\quad S^n=i^{-1}R^n,
$$
etc., so that our suitable square splits
into a disjoint sum of two squares:

$$
\xymatrix{ R^n\ar@{->>}[r]^{p^n} & T^n\\
           S^n\ar@{->>}[r]^{q^n}\ar@{^{(}->}[u]^{i^n}&
           P^n\ar@{^{(}->}[u]^{j^n}}
$$
where $n=1,2$. It follows from the
definitions that $i^1,q^2$ are bijections so
that we may assume $S^1=R^1$, $S^2=P^2$,
$i^1=\Id,q^2=\Id$.

So, we have the following diagram:
\begin{equation}
\xymatrix{ S^1\disjoint R^2\ar@{->>}[r]^{p^1\disjoint p^2} & T^1\disjoint T^2\\
           S^1\disjoint S^2\ar@{->>}[r]^{q^1\disjoint \Id}\ar@{^{(}->}[u]^{\Id\disjoint i^2}
           & P^1\disjoint S^2\ar@{^{(}->}[u]^{j^1\disjoint j^2}}
\end{equation}

The desired map is then defined as follows:

$$
\p_{i^1\disjoint i^2}\R_{p^1\disjoint
p^2}\to \p_{\Id\disjoint
i^2}\R_{p^1\disjoint
\Id}\delta_{\Id\disjoint p^2}\to
\R_{p^1\disjoint \Id} \p_{\Id\disjoint
i^2}\delta_{\Id\disjoint p^2}\to
\R_{p^1\disjoint \Id}\p_{\Id\disjoint
p^2i^2}=\R_q\p_j.
$$

\subsubsection{Properties } The above
defined maps have the following properties,
easily derived from the ones of the maps
$\l_p,s(\phi,g)$. We shall now list them.

1. Let $$\xymatrix{R\ar@{>>}[r]^p&  T\\
        S\ar@{^{(}->}[u]^i\ar@{>>}[r]^q
                                 &P\ar@{^{(}->}[u]^j }
                                 $$
        be a suitable square and $q=q_2q_1$,
where $q_1,q_2$ are surjections.

Define the  set $X(q_1,q_2)$ of isomorphism
classes of  commutative diagrams
$$\xymatrix{R\ar@{>>}[r]^{p_1}&
U\ar@{>>}[r]^{p_2}& T\\
        S\ar@{^{(}->}[u]^i\ar@{>>}[r]^{q_1}
        \ar@{^{(}->}[u]^j&
        V\ar@{^{(}->}[u]^{j'}\ar@{>>}[r]^{q_2}&
         P\ar@{^{(}->}[u]^j }$$
         We will refer to such a diagram as
         $(p_1,p_2,j')$.
Both squares in every such  a diagram are
automatically suitable. Therefore, every
element $x:=(p_1,p_2,j')\in X(q_1,q_2)$
determines a map
$$
m_x:\p_i\R_p\to \p_i\R_{p_1}\R_{p_2}\to
\R_{q_1}\p_{j'}\R_{p_2}\to
\R_{q_1}\R_{q_2}\p_j.
$$

         Then  the
        composition
        $$
        \p_i\R_p\to \R_q\p_j\to
        \R_{q_1}\R_{q_2}\p_j$$
 equals
 $$
 \sum_{x\in X(q_1,q_2)} m_x
 $$

 2. Consider the following commutative
 diagram
 $$\xymatrix{ R\ar[r]^p              &T\\
S_1\ar[u]^{i_2}\ar[r]^q   &P_1\ar[u]^{j_2} \\
S\ar[u]^{i_1}\ar[r]^r     & P\ar[u]^{j_1}}
$$
in which both small squares are suitable.
Then the large square is also suitable and
the following maps coincide:
$$
\p_{i_2i_1}\R_p\to \R_r\p_{j_2j_1}
$$
and
$$
\p_{i_2i_1}\R_p\to \p_{i_1}\p_{i_2}\R_p\to
\p_{i_1}\R_q\p_{j_2}\to
\R_{r}\p_{j_1}\p_{j_2} \to \R_{r}\p_{j_2j_1}
$$

3. Consider the following commutative
diagram:
$$
\xymatrix{R\ar@{>>}[r]^p        &                 T\\
            S\ar@{>>}[r]^q\ar@{^{(}->}[u]^i &
                                        P\ar@{^{(}->}[u]^j\\
            Q\ar@{^{(}->}[u]^k\ar@{^{(}->}[ru]^l&}
            $$
where the upper square is suitable. Then the
following maps coincide:

$$
\p_{ik}\R_p\to \p_k\p_i\R_p\to
\p_k\R_q\p_j\to\p_{qk}\p_j=\p_l
$$
and $$\p_{ik}\R_p\to \p_{pik}=\p_l.$$

4. Let $$
\xymatrix{ R\ar@{->>}[r]^p & T\\
           S\ar@{->>}[r]^q\ar@{^{(}->}[u]^i& X\ar@{^{(}->}[u]^j}
$$
and

$$
\xymatrix{ R_1\ar@{->>}[r]^{p_1} & T_1\\
           S_1\ar@{->>}[r]^{q_1}\ar@{^{(}->}[u]^{i_1}& X_1\ar@{^{(}->}[u]^{j_1}}
$$

be suitable squares and let $s:S\to S_1$,
$r:R\to R_1$, $t:T\to T_1$, $x:X\to X_1$ be
bijections fitting the two squares into a
commutative cube. Then the map $A(i,p,j,q)$
can be expressed in terms of
$A(i_1,p_1,j_1,q_1)$ in the following
natural way:
$$
\p_i\R_p\cong
\p_s\p_{i_1}\p_{r^{-1}}\p_r\R_{p_1}\p_{t_1^{-1}}\cong
\p_s\p_{i_1}\R_{p_1}\p_{t_1^{-1}}
\stackrel{A(i_1,p_1,j_1,q_1)}\longrightarrow
\p_s\R_{q_1}\p_{j_1}\p_{t_1^{-1}}\cong
\p_s\R_{q_1}\p_{x^{-1}}\p_x\p_{j_1}\p_{t_1^{-1}}
\R_q\p_j.
$$

5. let $(i_k,p_k,j_k,q_k)$, $k\in K$ be a
collection of suitable squares. Let
$i_k:S_k\to R_k$; let $M_k$ be a collection
of $\D_{X^{S_k}}$-sheaves. Let
$i=\disjoint_{k\in K} i_k$,
$p=\disjoint_{k\in K} p_k$,
$j=\disjoint_{k\in K} j_k$,
$q=\disjoint_{k\in K} q_k$, and
$M=\boxtimes_{k\in K}M_k$. Then the square
$i,p,j,q$ is also suitable and the following
compositions coincide:
$$
\boxtimes_{k\in K}\p_{i_k}\R_{p_k}(M_k)\to
\boxtimes_{k\in K}R_{q_k}\p_{j_k}(M_k)\to
\R_q\p_j(M)
$$
and
$$
\boxtimes_{k\in K}\p_{i_k}\R_{p_k}(M_k)\to
\p_i\R_p(M)\to \R_q\p_j(M).
$$

6. Let $i_k:S_k\to R_k$, $k\in K$ be
injections and $p_k:R_k\to T_k$, $k\in K$ be
surjections such that $j_k:=p_ki_k$ are
injections. Let $M_k$ be
$\D_{X^{T_k}}$-modules. Let $i,j,p,M$ be
disjoint unions of the respective objects.

 Assume that at least
two of the maps $p_k$ are proper
surjections. Then the composition
$$
{\boxtimes_{k\in K}}\p_{i_k}\R_{p_k}(M_k)\to
\p_i\R_p(\boxtimes_k M_k) \to \p_j(M)
$$
vanishes.

If only one of the surjections $p_k$ is
proper, say $p_\kappa$, $\kappa\in K$, then
the above composition equals
\begin{eqnarray*}
\boxtimes_{k\in K}\p_{i_k}\R_{p_k}(M_k)=
\p_{i_\kappa}\R_{p_\kappa}(M_\kappa)
\boxtimes_{k\in K\bs \{\kappa\}}
\p_{i_k}\R_{p_k}(M_k)
\\
\stackrel{L(i_\kappa,p_\kappa)}\longrightarrow
\p_{j_\kappa}(M_\kappa) \boxtimes_{k\in K\bs
\{\kappa\}} \p_{i_k}\R_{p_k}(M_k)\to
\boxtimes_{k\in K}\p_{j_k}(M_k)\to \p_j(M)
\end{eqnarray*}

7. The diagram  (\ref{pertas}) is suitable,
and the corresponding map $A(i,p,j,q)$ is
the isomorphism (\ref{pertas1}).

 \subsubsection{Differentials} The
differential of the map $L(i,p)$ is computed
as follows. Consider the set of all
equivalence classes of decompositions
$p=p_2p_1$, where $p_1,p_2$ are surjections
and $p_1i$ is injection. We then have a map
$$l(p_1,p_2):\p_i\R_p\to
\p_i\R_{p_1}\R_{p_2}\to \p_{p_1i}R_{p_2}\to
\p_{p_2p_1i}=\p_{pi}.
$$
We then have
$$
dL(i,p)+\sum_{(p_1,p_2)}l(p_1,p_2)=0.
$$

2. Let
$$
\xymatrix{Q:R\ar@{>>}[r]^p&T\\
 S\ar@{>>}[r]^q\ar@{^{(}->}[u]^i&P\ar@{^{(}->}[u]^j}
 $$
 be a suitable square. Define two sets
 $L(Q)$ and $R(Q)$ as follows. The set
 $L(Q)$ is the set of all isomorphism
 classes of diagrams:
$$
\xymatrix{R\ar@{>>}[r]^{p_1}& R_1\ar@{>>}[r]^{p_2}&T\\
 S\ar@{>>}[rr]^q\ar@{^{(}->}[u]^i\ar@{^{(}->}[ur]^{i_1}&&
 P\ar@{^{(}->}[u]^j}
 $$
such that $p=p_1p_2$. It is clear that the
internal commutative square in this diagram
is also suitable.

Define the set $R(Q)$ as the set of
isomorphisms classes of diagrams

$$
\xymatrix{R\ar@{>>}[r]^{p_1}& R_1\ar@{>>}[r]^{p_2}&T\\
 S\ar@{>>}[rr]^q\ar@{^{(}->}[u]^i&&P\ar@{^{(}->}[u]^j
 \ar@{^{(}->}[ul]^{j_1}}
 $$
where $p=p_1p_2$. The internal square in
such a diagram is always suitable as well.

Every element $l:=(p_1,p_2,i_1)\in L(Q)$
determines a map
$$
f_l:\p_i\R_p\to \p_i\R_{p_1}\R_{p_2}\to
\p_{i_1}\R_{p_2}\to \R_{q}\p_j.
$$
Every element $r=(p_1,p_2,j_1)\in R(Q)$
determines a map
$$
g_r:\p_i\R_p\to \p_i\R_{p_1}\R_{p_2}\to
\R_q\p_{j_1}\R_{p_2}\to \R_q\p_j.
$$

We then have
$$
dA(i,p,j,q)=\sum_{l\in L(Q)}f_l-\sum_{r\in
R(Q)}g_r=0.
$$

This completes the list of properties.

\subsubsection{Second definition of a pre-symmetric system}
 Call a system $<\R>$
endowed with the above specified maps having
the above properties {\em a pre-symmetric
system} (in the sense of the second
definition). As we will  mainly use
pre-symmetric systems in the sense of the
second definition, we shall simply  refer to
them as pre-symmetric.

\section{Renormalization in pre-symmetric
systems} We are going to describe the
renormalization procedure for algebras over
 pre-symmetric systems.  The plan is as
follows.

First of all given an algebra $M$ over a
presymmetric system, we show that the direct
image $p_*M$ has an $L_\infty$-structure,
(here  $p:X\to \pt$). Next we have to show
how, given a solution to the Master
equation, one can deform the algebra $M$. As
in the naive approach, we see that to be
able to renormalize, one needs an extra
structure on our system, and we define this
structure (it is called symmetric). Next, we
show how the renormalization goes in
symmetric systems, and finally, we discuss a
procedure by means of which, given a
pre-symmetric system one can produce a
symmetric system (we call this procedure
{\em symmetrization}. So, the
renormalization of an algebra over  a
pre-symmetric system includes:

1) symmetrization of the system so  that we
get an OPE-algebra over a symmetric system;

2) renormalization in the symmetric system.

\subsubsection{ An $L_\infty$-structure on
$p_*M[1]$, where $M$ is an OPE-algebra over
$<\R>$}.

Let $M$ be an OPE-algebra over $<\R>$. We
are going to introduce an $L_\infty$
structure on $p_*M$, where $p:X\to \pt$ is
the projection Let $S$ be a finite set and
$i_S:\emptyset \to S$ be an embedding. Let
$\p_S:=\p_{i_S}$. It is clear that
$\p_{\pt}=p_{*}$ and that
$$\p_S(M^{\boxtimes
S})\cong (p_*M)^{\otimes S}.$$ Finally, set
$p_S:S\to \pt$.

 Define a  degree +1 map
$$C_S:(p_*M)^{\otimes_S}\to p_*M
$$ as the composition:
$$
(p_*M)^{\otimes_S}\cong \p_S(M^{\boxtimes
S})\stackrel{\ope_S}\to
\p_S\i_{p_S}(M)\stackrel{L(i_S,p_S)}\longrightarrow
\p_{i_{\pt}}M\cong p_*M.
$$

\begin{Claim} The maps $C_S$ endow $p_*M[1]$
with an $L_\infty$-structure.
\end{Claim}

\pf The key ingredient in the proof is

\begin{Lemma}
Let $q:S\to T$ be a surjection such that one
can decompose $S=S_1\disjoint S_2$,
$T=T_1\disjoint T_2$, $q=q_1\disjoint q_2$,
where $q_i:S_i\to T_i$, $i=1,2$ are both
non-bijective surjections. Then the
composition
$$
\p_S(M^{\boxtimes S})\to
\p_S\i_q(M^{\boxtimes T})\to
\p_T(M^{\boxtimes T})
$$
vanishes.
\end{Lemma}
\pf  Let $A=\{1,2\}$.  Then the above
composition equals:
\begin{eqnarray*}
\p_S(M^{\boxtimes S})\cong
\p_{i_{S_1}\disjoint i_{S_2}}(M^{\boxtimes
S_1}\boxtimes M^{\boxtimes S_2})\to
\p_{i_{S_1}\disjoint
i_{S_2}}(\i_{q_1}(M^{T_1})\boxtimes
\i_{q_2}(M^{T_2}))\\
\to \p_S\i_q(M^{\boxtimes
T_1}\boxtimes M^{\boxtimes T_2})\to
\p_T(M^{\boxtimes T}).
\end{eqnarray*}

Here $i_{S_1}:\emptyset\to S_1$,
$i_2:\emptyset\to S_2$.

  The composition of the two last arrows
   vanishes  by the Property 6 in the
previous subsection.
\endpf

The Claim now follows directly from the
formula of the differential of $L(i,p)$ .
\endpf
\subsubsection{Action of the DGLA $p_*M[1]$
on $M$}

Define the maps
$$
A_S:(p_*M)^{\otimes S}\otimes M\to M
$$

as follows. Let $S_0=S\disjoint \pt$. Let
$k:\pt\to S_0$ be the natural embedding. Let
$p_{S_0}:S_0\to \pt$.

We then set
$$
A_S: (p_*M)^{\otimes S}\otimes M\cong
\p_k(M^{\boxtimes
S_0})\stackrel{\ope_{S_0}}\to
\p_k\i_{p_{S_0}}(M)\stackrel{L(k,p_{S_0})}\to
M.
$$
It is not hard to see that the collection of
maps $A_S$ determines an $L_\infty$-action
of $p_*M[1]$ on $M$.

\subsection{Symmetric systems} Pre-symmetric systems do not
fit for renormalization. The reason is more
or less  the same as in the naive approach,
but let us reformulate it in terms of direct
images with respect to projections.

Let $p:S\to \pt$ and pick an element $s\in
S$. Let

 Let $S':=S\disjoint \{s\}$. Let $t\in S'$.
 Define $p_t:S\to S'$ as follows:
 $$
 p_t(r)=r$$
 if $r\neq s$;
 $$
p_t(s)=t.
$$

Let $p_s:S\to \{a,s\}$, where $a$ is an
abstract element, $a\neq s$, by setting
$$
p_s(t)=a$$ if $t\neq s$; $p_s(s)=s$.

Let $q:S'\to \pt$ and $r:\{a,s\}\to \pt$.
Let $i:S'\to S$, $j:\{s\}\to \{a,s\}$ be
natural embeddings.

We then have several maps
$$
\p_i\R_p\to \R_q.
$$

a) Let $t\in S'$. Set
$$
L_t:\p_i\R_p\to \p_i\R_{p_t}\R_q\to \R_q;
$$
Set
$$
R:\p_i\R_p\to \p_i\R_{p_s}\R_r\to
\R_{q}\p_j\R_r\to \R_q.
$$

Then luck of symmetricity manifests itself
in the fact that the difference
$$
R-\sum_{\t\in S'} L_t$$ is not homotopic to
0.

We thus need to add extra homotopies which
would take care about it. it turns out that
this can be accomplished in a very symple
way:

call a system $<\R>$ {\em symmetric} if the
maps $A(i,p,j,q)$ are defined for {\em all}
commutative squares
$$\xymatrix{
R \ar@{>>}[r]^p & T\\
S\ar@{^{(}->}[u]^i\ar@{>>}[r]^q &
P\ar@{^{(}->}[u]^j}
$$
where $p,q$ are both non-bijective
surjections (not necessarily suitable). The
properties remain the same as for
pre-symmetric system except that we drop the
suitability condition everywhere.

We shall demonstrate how the renormalization
goes in symmetric systems.

Let now $\S\in \h p_*M^0[[\h]]$ be a MC
element.   For a finite set $T$ set
$$\S_T:=\S^{\boxtimes T}\in
\p_T\h^{|T|}M^{\boxtimes T}[[\h]]$$
 Let $i:R\to S$ be an injection. Let
 $T=S\backslash i(R)$.
  We then have a map
$$M^{\boxtimes R}\to
\p_i\h^{|T|}M^{\boxtimes S}[[\h]]$$ defined
by:
$$
M^{\boxtimes R}\stackrel{\otimes \S_T}\to
M^{\boxtimes R}\otimes
\p_T\h^{|T|}M^{\boxtimes_ T}[[\h]] \cong
\p_i\h^{|T|}M^{\boxtimes S}[[\h]].$$

\subsubsection{} Let $i:S\to R$ be an
injection and $q:R\to T$ be a surjection
such that $p:=qi$ is a surjection.

We then have a map
$$
\ope(q,i):M^{\boxtimes S}\to
\p_iM^{\boxtimes R}\to \p_i\R_qM^{\boxtimes
T}\to \R_pM^{\boxtimes T}.
$$

set
\begin{equation}\label{local*}
\ope^r_p=\sum \ope(q,i),
\end{equation}

where the sum is taken over all isomorphism
classes of decompositions $p=qi$. Let
$M^r:=M[[\h]], d_\S$, where $d_\S$ is the
differential twisted by $\S$. Then
$(M',\ope^r)$ is the renormalized
OPE-algebra.

Note that the sum (\ref{local*}) is infinite
but it converges in the $\h$-adic topology.

\subsection{Symmetrization} Finally, we need
a method how, given a pre-symmetric system
one gets a symmetric system.

The idea is as follows. Let $f:S\to T$ be a
map of finite sets. Construct a category
$\bodya(f)$ whose objects are compositions
$\p_i\R_{p_1}\R_{p_2}\cdots\R_{p_n}$, where
$i$ is injective, $p_k$ are surjective, all
the maps are composable and
$$
p_np_{n-1}\cdots p_1i=f.
$$  The morphisms are all possible
morphisms one can get using the axioms of
pre-symmetric system. Given a pre-symmetric
system $<\R>$ and a $\DD_{X^T}$-sheaf $N$,
the application
\begin{eqnarray*}
\p_i\R_{p_1}\R_{p_2}\cdots\R_{p_n}\mapsto\\
\p_i\R_{p_1}\R_{p_2}\cdots\R_{p_n}N
\end{eqnarray*}
 produces a  functor
$$
\r^f(N):\bodya(f)\to \Cm_{X^T}.
$$

Let $\bodyb(f)$ be the same thing, but we
use axioms of a symmetric system. We then
have a tautological functor $R:\bodya(f)\to
\bodyb(f)$. One can construct a bifunctor
$$B:\bodya^{\op}(f)\times \bodyb(f)\to \complexes,$$
where $B(X,Y)=\hom_{\bodyb(f)}(X,Y)$.

Set $$\r_\symm^f(N):\bodyb(f)\to
\complexes$$ to be
$$
\r^f(N)\otimes_{\bodya(f)} B.
$$

{\bf Remark} Let $R^{-1}$ be a functor from:

the category of functors $\bodyb(f)\to
\complexes$

 to:

  the category of functors $\bodya(f)\to
 \complexes$

 which is the pre-composition with $R$.
One can show that $R^{-1}$ has a left
adjoint $R_!$ and that
$\r^f_\symm(N)=R_!\r^f(N)$.\endpf

We can now construct a system $<\Rs>$ which
is a  symmetrization of $\R$
 by setting $\Rs_p(N)=\r^p(N)(\R_p)$.
We have to say that the introduction of a
structure of system on the collection of
functors $<\Rs_p>$ is not at all a
consequence of a general non-sense. It turns
out that in order to define such a structure
one has to use certain specific features of
the categories $\bodya,\bodyb$.

We also have a natural map $<\R>\to <\Rs>$.
Therefore, given an OPE-algebra  over
$<\R>$, we can transform it into an
OPE-algebra over $<\Rs>$ and then
renormalize it.

We shall now give a more explicit
construction of  $<\Rs>$. In fact, the
resulting system $<\Rs>$ is isomorphic to
the above described one. This follows from a
more detailed study of the categories
$\bodya, \bodyb$ which id done in \ref{babb}
\def\Aut{{\text{Aut}}}
\def\N{\Rs}
\section{Explicit construction of
$<\N>$}\label{explicitrsymm}
\subsection{Main objects}
\subsubsection{Groupoid $C'_f$} Let $f:S\to T$
be a surjection. Define a groupoid $C'_f$
whose objects are diagrams
$$
\xymatrix{S\ar[r]^i &U\ar[r]^p & T,}
$$
where $i$ is injective, $p$ is surjective,
and $pi=f$.  Isomorphisms are morphisms of
these diagrams inducing identities on $S,T$.
\subsubsection{Groupoid $C_f$} Let $(i,p)\in
C_f$. Call $p$ {\em $i$-super-surjective} if
for every $t\in T$, the pre-image  $p^{-1}t$
either:

 contains at least two elements from
$i(S)$

or:

consists of one element from $i(S)$.

Let $C_f$ be the full sub-groupoid of $C'_f$
consisting of all pairs $(i,p)$, where $p$
is $i$-super-surjective.

\subsubsection{Functors $\M(i,p)$, $\M_f$}
For an object $(i,p)$ in $C_f$, set
$\M(i,p):= \p_i\R_p$. It is clear that $\M(
,)$ is a functor from $C_f$ to the category
of functors from the category of
$\D_{X^T}$-sheaves to  the category of
$\D_{X^S}$-sheaves. Set
$$\M_f:=\limdir_{C_f}\M(i,p).$$
Denote by  $I(i,p):\M(i,p)\to \M_f$ the
natural map. It is clear that $I(i,p)$
passes through $\M(i,p)_{\Aut_{C_f}(i,p)}$.
Furthermore, we have an isomorphism
\begin{equation}\label{bezdiff}
\oplus \M(i,p)_{\Aut_{C_f}(i,p)}\to \M_f,
\end{equation}
where the sum is taken over an arbitrary set
of representatives of isomorphism classes of
$C_f$.

\subsection{Differential} The symmetrized
resolution $\Rs_f$  is given by the functor
$\M_f$ as in (\ref{bezdiff}), on which a new
differential is introduced. This
differential is of the form $d+L+R$, where
$d$ is the differential on $\M_f$, and
 degree +1 endomorphisms $L,R:\M_f\to \M_f$
 shall be defined below.

\subsubsection{Map $L:\M_f\to \M_f$}
\subsubsection{Set $E_L(i,p)$}
Let
$$\xymatrix{S\ar[r]^i & U\ar[r]^p  &T}
$$
be an object in $C_f$. Define a  finite set
$E_L(i,p)$ whose elements are equivalence
relations $e$ on $U$ such that

1) $p$ passes through $U/e$;

2) the composition
$$\xymatrix{S\ar[r]^i & U\ar[r]  &U/e}
$$
is injective.

Let $\pi_e:U\to U/e$ be the natural
projection, let $p_e:U/e\to T$ be the map
induced by $p$, and $i_e=\pi_e i$.

It turns out that $(i_e,\pi_e)\in C_f$.
Indeed, $\pi_e^{-1}(t)$ is the quotient of
$p^{-1}t$ by $e$ and elements of $i(S)$ are
$e$-non-equivalent, which implies the
super-surjectivity.
\subsubsection{The map $L$}
 Define a map
$L_e:\M(i,p)\to \M(i_e,p_e)$ as follows:
$$
\p_i\R_p\to
\p_i\R_{\pi_e}\R_{p_e}\stackrel{L(i,\pi_e)}\longrightarrow
\p_{i_e}\R_{p_e}.
$$

Define a map $L(i,p):\M(i,p)\to \M_f$ by
setting
$$
L(i,p)=\sum\limits_{e\in
E_L(i,p)}I(i_e,p_e)L_e.
$$

It is easy to see that the collection of
maps $L(i,p)$ descends to a map $L:\M_f\to
\M_f$.
\subsection{Map $R:\M_f\to \M_f$}
\subsubsection{Set $E_R(i,p)$}
Let
$$\xymatrix{S\ar[r]^i & U\ar[r]^p  &T}
$$
be an object in $C_f$. Define a  finite set
$E_R(i,p)$ whose elements are equivalence
relations $e$ on $U$ such that

1) $p$ passes through $U/e$;

2) The restriction of $e$ on $S$ coincides
with the equivalence relation on $S$
determined by $f$.

Let $\pi_e:U\to U/e$. Let $T_e:=\Im(\pi_ei)$
and $V:=V_e:=\pi_e^{-1}T_e$ and
$W:=W_e:=U\backslash U_e$. Let
$e_V\text{(resp. }e_W)$ be the restriction
of $e$ on $V\text{(resp. }W)$.

It is clear that

1) $i(S)\subset V$;

2) The map $p^{V/e_V}:V/e_V\to T$  induced
by $p$ is bijective.

So, we have a diagram:

$$
\xymatrix{S\ar@{^{(}->}[r]^{i_e} & V
\ar@{>>}[r]^{\pi_e|_{V}}\ar@/^2pc/[rr]^{p_e}&
V/e_V\ar[r]^\sim_{p^{V/e_V}} &
T\\
&&&\\
&W\ar@{>>}[r]^{\pi_e|_{W}}&W/e_W\ar[uur]^{p^{W/e_W}}&}
$$

Elements of  $E_R(i,p)$ can be equivalently
defined as collections $(W,e_W)$, where
$W\subset U$, $W\cap i(S)=\emptyset$, and
$e_W$ is an equivalence relation on $W$ such
that $p|_W$ passes through $e_W$.  Indeed,
let $V:=U\backslash W$ and let $e_V$ be
induced on $V$ by $p|_V$. Set
$e:=e_V\disjoint e_W$. This establishes a
1-1 correspondence between different
descriptions of $E_R(i,p)$.

 Let us check that $(i_e,p_e)\in C_f$.
Indeed,  for every $t\in T$,
$p_e^{-1}t=p^{-1}t\cap V$. Since $V\supset
i(S)$, we have: if $p^{-1}t\cap i(S)$ has at
least two elements, then so does
$p_e^{-1}t$; otherwise $p^{-1}t$ consists of
exactly one element from $i(S)$ and
$p_e^{-1}t=p^{-1}t$.

We will now define a map $R_e:\M(i,p)\to
\M(i_e,p_e)$ To this end we shall consider a
diagram:
$$\xymatrix{
S\ar@{^{(}->}[r]^{i_e}\ar@{^{(}->}@/^2pc/[rr]^i
&
V\ar@{^{(}->}[r]^{I}\ar@{>>}[d]^{\pi_e|_V}\ar@{>>}@/_5pc/[ddr]_{p_e}
&
 U\ar@{>>}[d]^{\pi_e}\ar@{>>}@/^2pc/[dd]^p\\
        & V/e\ar@{^{(}->}[r]^J \ar[dr]^\sim_{p^{V/e_V}}&  U/e\ar@{>>}[d]^{r}\\
        &                       & T}
$$

We then observe that the square
$(I,\pi_e,J,\pi_e|_V)$ is  clearly suitable.
We therefore can define $R_e:M(i,p)\to
M(i_e,p_e)$  via the following chain of
maps:
\begin{eqnarray*}
R_e:\p_i\R_p\cong \p_{i_e}\p_I\R_{rp_e}\to
\p_{i_e}\p_I\R_{p_e}\R_r\to
\p_{i_e}\R_{\pi_e|_V}\p_J\R_r\\\to
\p_{i_e}\R_{\pi_e|_V}\p_{p^{V/e_V}}\cong
\p_{i_e}\R_{p_e}
\end{eqnarray*}

We then define $$R=\sum_{e\in R(i,p)}
I(i_e,p_e)R_e.$$
\subsubsection{Definition of the differential} We define
the differential on $<\Rs>$ as a sum
$d+L+R$.
\subsection{Asymptotic decomposition maps
$\as_{f_1,f_2}:\N_{f_2f_1}\to
\N_{f_1}\N_{f_2}$}

Suppose we have a  chain of surjections
 $$
\xymatrix{S\ar@{>>}[r]^{f_1}& R
\ar@{>>}[r]^{f_2} & T},
$$
so that  $f=f_2f_1$.

Let
$$ \xymatrix{ S\ar@{^{(}->}[r]^i&U\ar@{>>}[r]^p&T}
$$
be in $C_f$.
 The map
$\I(i,p):\M(i,p)\to \M_f$ determines a
similar map $\I(i,p)\to \N_f$. In order to
construct the map $\as_{f,g}$ we will first
define maps
$$
\as(i,p,f,g):\M(i,p)\to \M_{f_1}\M_{f_2}.
$$

Define the set $E(i,p,f_1,f_2)$ whose
elements are equivalence relations $e$ on
$U$ such that

1) $p$ passes through $U/e$

2) The restriction $e|_S$ coincides with the
equivalence relation on $S$ determined by
$f_1$.

Let $V_e\subset U$ be the set of all
elements which are equivalent (with respect
to $e$) to elements of $S$. Let
$W_e=U\backslash V_e$. Let $i_e:S\to V_e$;
$p_e:V_e\to V_e/e$, $j_e: V_e/e\to U/e$,
$q_e:U/e\to T$ be the map induced by $p$. We
then have the following commutative diagram:
$$
\xymatrix{
S\ar@{^{(}->}[r]^{i_e}\ar@/^2pc/@{^{(}->}[rr]^i
& V_e\ar@{^{(}->}[r]^I\ar@{>>}[d]^{p_e}
                  & U\ar@{>>}[d]^\pi\ar@/^2pc/@{>>}[dd]^p\\
&V_e/e\ar@{^{(}->}[r]^{j_e}&U/e\ar@{>>}[d]^{q_e}\\
&& T}
$$
It is easy to check that the square
$(I,\pi,j_e,p_e)$ is suitable. This allows
us to
 define a
map
$$
\as(i,p,e):\M(i,p)\to \M(i_e,p_e)\M(j_e,q_e)
$$
as follows:
\begin{eqnarray*}
\M(i,p)\cong \p_i\R_p\to
\p_{i_e}\p_I\R_{q_e\pi}\to
\p_{i_e}\p_I\R_\pi\R_{q_e}\\
\to
\p_{i_e}\R_{p_e}\p_{j_e}\R_{q_e}=\M(i_e,p_e)\M(j_e,q_e).
\end{eqnarray*}

Let
$$
\as(i,p,f_1,f_2):\M(i,p)\to
\M_{f_1}\M_{f_2}$$ be given by the formula:
$$
\as(i,p,f_1,f_2)=\sum\limits_{e\in
E(i,p,f_1,f_2)}\I_{(j_e,q_e)}\I_{(i_e,p_e)}\as(i,p,e).
$$

 This completes the definition of the map
 $\as_{f_1,f_2}$.
\subsection{Factorization maps}
Let $f_a:S_a\to T_a$, $a\in A$ be a family
of surjections

Let $(i_a,p_a)\in C_{f_a}$, $a\in A$, be a
family of objects. Let $i=\disjoint_{a\in A}
i_a$, $p=\disjoint_{a\in A} p_a$,
$f=\disjoint_{a\in A} f_a$.

Let $M_a\in \Cm_{X^{T^a}}$. Let
$M:=\boxtimes_{a\in A} M_a$.

 We then
have a natural map
$$
\boxtimes_a \p_{i_a}\R_{p_a}(M_a)\to
\p_i\R_p(M),$$ induced by the factorization
maps for $<\R>$. These maps give rise to the
factorization maps in $<\Rs>$.

\subsection{ Maps $L(i,f):\p_i\Rs_f\to
\p_j$} Let
$$
\xymatrix{R\ar@{>>}[r]^f& T\\
          S\ar@{^{(}->}[u]^i\ar@{^{(}->}[ur]^j&}
        $$
be a commutative diagram. The map
$L(i,p):\p_i\Rs_f\to \p_j$ is then defined
via maps

$$
\p_i\M(k,p)\cong
\p_{ki}\R_p\stackrel{L(ki,p)}\longrightarrow
\p_j,
$$
where $pk=f$.
\subsection{The maps
$A(i,p,j,q):\p_i\Rs_p\to \Rs_q\p_j$}
 Let
$$
\xymatrix{R\ar@{>>}[r]^f  & T\\
          S\ar@{^{(}->}[u]^i\ar@{>>}[r]^q  & P\ar@{^{(}->}[u]^j}
$$
be a commutative diagram. The maps
$$A(i,f,j,q):\p_i\Rs_f\to \Rs_q\p_j
$$

are defined as follows.

Let $(k,p)\in C_f$. Let $u=ki$. One can show
that there exists a unique, up-to an
isomorphism decomposition $u=u_2u_1$ into a
product of two injections such that in the
diagram
$$
\xymatrix{R \ar@{>>}[r]^p      &   T\\
          S_1\ar@{^{(}->}[u]^{u_2}\ar@{>>}[r]^{q_1}     &
          P\ar@{^{(}->}[u]^{j}\\
          S\ar@{^{(}->}[u]^{u_1}\ar@{>>}[ur]^{q} &}
$$
uniquely, up-to an isomorphism, constructed,
given a decomposition $u=u_2u_1$, the square
$$
(u_2,p,j,q_1)$$ is suitable, and in the pair
$$
(u_1,q_1),
$$
the map $q_1$ is super-surjective.

The map $A(i,f,j,q)$ goes as follows:
$$
\p_i\M(k,p)\cong \p_{ki}\R_p\cong
\p_{u_1}\p_{u_2}\R_p\stackrel{A(u_2,p,j,q_1)}\longrightarrow
\p_{u_1}\R_{q_1}\p_j\cong \M(u_1,q_1)\p_j.
$$

\section{Constructing the system $<\R>$ with
the above explained properties}

\subsection{Step 1: Spaces of generalized
functions $C_S$} Our motivation comes from
the construction in \ref{physmot}. In the
case  when $p:S\to \pt$, where $S$ has two
elements, this construction suggests that
one can replace $\i_p$ with a complex $0\to
i_{p*}\to \I_p\to 0$, where we put $\I_p$ in
degree 0. Denote this complex by $\R_p$.
 On the one
hand we have a map $\R_p\to\i_p$, so that
the induced map $\R_p(M)\to \i_p(M)$ is a
quasi-isomorphism for good $M$'s; on the
other hand we have a map $\R_p\to
i_{p*}\O_X$ of degree +1. Thus, $\R_p$ has
all the desired  properties.

Let us try to expand this construction to an
arbitrary case. It is natural to start
 with constructing
certain spaces of generalized functions
$\C_S$ on $X^S$  so that each $\C_S$  is a
sub- $\D_{X^S}$ submodule of the space of
complex-valued generalized functions on
$Y^S$ with compact support . In pursuit of
making $\C_S$ as small as possible we
construct $\C_S$ in such a way that they are
holonomic $\D_{X^S}$-modules; their
structure is as follows. Let  $D$ be a
generalized diagonal in $X^S$ and let
$\C_{S,{[D]}}$ be the maximal submodule
supported on $D$.  This defines a filtration
on $\C_S$ whose terms are labelled by the
ordered set of generalized diagonals in
$X^S$.
 The associated graded term
$$\C_{S,[D]}/
\text{span}_{E\subsetneq D}C_{S,[E]}\cong
i_{D*}\B_D,
$$
where $\B_D$ is the $\D_D$-module of all
meromorphic functions with singularities
along hyper-surfaces $q(X_i-X_j)=0$, where
$X_i|_D\neq X_j|_D$.

Construction of  such $\C_S$ is done by
means of certain analytical considerations.
Some of them  a very similar to standard
methods of regularization of divergent
integrals. The detailed exposition is in
\ref{genfun}-\ref{genfun1}.
\subsection{Step 2: Functors $\I_p$ and
their properties} Next we construct the
functors $\I_p$ out of $\C_S$ in the same
way as $\i_p$ was constructed out of $\B_S$:
let $p:S\to T$ be a surjection of finite
sets; set
$$
\C_p:=\boxtimes_{t\in T} \C_{p^{-1}t}.
$$
 Define $\I_p:\D_{X^T}\to \D_{X^S}$ by
$$
\I_p(M)=i_p^\ha(M)\otimes_{\O_{X^S}} \C_p.
$$

We then have natural maps $\I_p\to \i_p$. We
then ask ourselves whether $\I_p$ form a
system. The answer is no. It probably could
be yes if $\C_S$ would be a bit larger
subspace of generalized functions, because
we have a technique of asymptotic
decomposition of generalized functions due
to  J. Bernstein (unpublished). But there
are examples in which we see that already
for the set $S=\{1,2,3\}$ consisting of
three elements there are functions $f\in
C_S$, whose asymptotic decomposition near
the diagonal $X^1=X^2$ requires introduction
of such functions as $\log(X_1-X_3)$. For
example, let $Y=\Re^4$ and take
$$
f(X^1,X^2,X^3)=\frac
1{|X^1-X^3|^2|X^2-X^3|^2}.
$$

This is a locally  $L^1$-function,
therefore, it determines a generalized
function. Let us investigate its asymptotic
as $X^1$ approaches $X^2$. According to  J.
Bernstein, we should consider the following
expression:

$$
u(\la)=\int \frac
{g(X_1,X_1+(X_1-X_2)/\la,X_3)}{|X^1-X^3|^2|X^2-X^3|^2}d^4X^1d^4X^2d^4X^3,
$$
where $g$ is a compactly supported smooth
function and $\la$ is a small positive
parameter. Our goal is to  find an
asymptotic for $a(\la)$. Let $x=X_1$,
$a=X_2-X_1$, $b=X_3-X_1$. Let
$G(x,a,b):=g(x,x+a,x+b)$. We then have
$$
u(\la)=\int\frac
{G(x,a/\la,b)}{|b|^2|b+a|^2}d^4bd^4ad^4x
$$

One can show that
$$
u(\la)=C\int
G(x,a/\la,0)\ln(|a|^2)d^4ad^4x+v(\la),
$$
where $v(\la)$ is bounded as $\la\to +0$,
and $C$ is a constant.

This means that
$$
v(\la)=\int
g(X_1,X_1+(X_1-X_2)/\la,X_3)\big\{ \frac 1
{|X^1-X^3|^2|X^2-X^3|^2}-C\ln(|X_1-X_2|^2)\delta(X^1-X^3)
\big\} d^4X^1d^4X^2d^4X^3,
$$
is bounded as $\la\to +0$. This demonstrates
that, at least, we have to include
$\ln(|X_1-X_2|^2)$ into our picture to get
an asymptotic decomposition of $$ \frac 1
{|X^1-X^3|^2|X^2-X^3|^2}.
$$

The geometrical meaning of this phenomenon
is that the cohomology of the complex
variety which is the complement in
$\Co^4\times \Co^4$ to the set of complex
zeroes of $|Z_1-Z_2|^2=0$ differs from the
cohomology of the real part, which is ${\Bbb
R}^{4}\times {\Bbb R}^4$ minus the diagonal.
We need to add functions which would kill
the de-Rham cocycles which are non-trivial
on the complexification but become trivial
upon restriction to the real part.

Nevertheless, we have maps
\begin{equation}\label{I-dec}\label{I-asadd}
\I_{pq}\to \I_q\i_p
\end{equation}
 for all  surjections
$p,q$.

For certain $p,q$ we also have maps
\begin{equation}\label{I-decadd}
\I_{pq}\to \I_q\I_p.
\end{equation} Namely,
this happens if
$$q=q_1\disjoint \Id:S_1\disjoint S_2\to
R_1\disjoint S_2
$$ and
$$p=\Id\disjoint
p_1:R_1\disjoint S_2\to R_1\disjoint T_2,
$$
or if $p,q$ can be brought to this form via
conjugations by bijections. This
circumstance will play an important role in
the future steps, but now let us concentrate
only on the maps $\I_{pq}\to \I_q\i_p$. They
have  associativity properties similar to
those of $\i$ and they nicely behave with
respect to $\boxtimes$. They  are compatible
with the corresponding maps $\i_{pq}\to
\i_q\i_p$.

There is an additional feature stemming from
the fact that the submodule
$\C_{S,\Delta}\subset \C_S$, where
$\Delta\subset X^S$ is a generalized
diagonal, is isomorphic to
$i_{\Delta*}\C_\Delta$.

Let $p$ be a  surjection. Denote
$\delta_p:=i_{p*}$. We then have a natural
map
\begin{equation}\label{delI}
\delta_p\I_q\to \I_{qp},
\end{equation}
 whenever surjections $p,q$ are composable.
 These maps behave nicely with respect to
 the other parts of the structure.
\subsubsection{Iterations of functors $\I$
and $\i$} We will work with all possible
functors of the form
$$
\j^1_{p_1}\j^2_{p_2}\cdots \j^n_{p_n},
$$
where $p_i:S_i\to S_{i+1}$ are surjections
and $\j^s_{p_s}$ is either $\i_{p_s}$ or
$\I_{p_s}$.  Fix a surjection $p:S\to T$ and
consider the class  $\zebra_p$ of all such
compositions with $p_np_{n-1}\cdots p_1=p$
(in particular, $S_1=S$, $S_{n+1}=T$). The
asymptotic decomposition maps (\ref{I-dec})
and their compositions produce maps between
objects of $\zebra_p$ (warning: we exclude
the  maps (\ref{I-decadd})). For example, we
can construct a map $\I_{qrp}\to
\I_p\i_r\i_q$ as a composition:
$$
\I_{qrp}\to \I_p\i_{qr}\to \I_p\i_r\i_q.$$
We can also take another composition:
$$
\I_{qrp}\to \I_{rp}\i_q\to \I_p\i_r\i_q.
$$
The associativity property implies that
these compositions are equal.

On the other hand,  there is no way to
construct a map $\I_{qrp}\to \i_p\I_q\i_r$.

Thus, $\zebra_p$ is naturally a category.
Furthermore, it turns out that, because of
the associativity properties, there is at
most one arrow between different arrows,
i.e. $\zebra_p$ is equivalent to a poset
which will be denoted by $\zebra(p)$.
 Let
us describe it. First of all, each
isomorphism class in $\zebra_p$ does not
even form a set because of the indeterminacy
in the choice of intermediate sets $S_i$.
This can be easily resolved by demanding
each $S_i$ to be  $S/e_i$, where $e_i$ is an
equivalence relation on $S$. More precisely,
let $e$ be the equivalence relation on $S$
determined by $p:S\to T$, $T$ being
identified with $S/e$. Let $\Eq_e$ be the
poset of all equivalence relations on $S$
which are finer than $e$. Let us write
$e_1>e_2$ if $e_1$ is finer than $e_2$.
Denote by $\omega$ the trivial (the finest)
equivalence relation on $S$. An element of
$\zebra(p)$ is then a pair $F,\{\j^s\}$,
where $F=(\omega=e_1>\cdots e_{n+1}=e)$ is a
proper flag of equivalence relations and
$\{j^s\}_{s=1}^n$ is a  sequence of symbols
$\i$ or $\I$. It is convenient to visualize
an object of zebra as a subdivision of a
large segment into $n$ small subsegments;
the equivalence relations are associated
with the nodes ($e_s$ is associated with the
$s$-th node from the left) and $\j^s$
determines one of two colors of the small
segment between the $s$-th and the $s+1$-th
node.

To such data we assiociate the functor
$$[F,\{\j^s\}]= \j^1_{p_1}\j^2_{p_2}\cdots
\j^n_{p_n},
$$

where $p_i:S/e_i\to S/e_{i+1}$ is the
natural projection.  Let us describe the
order (we assume that an arrow $X\to Y$
exists iff $X\leq Y$).  We say that $X<Y$ if

1) the flag of $Y$ is a refinement of the
flag of $X$. Thus, each  small segment of
the flag of $X$ is then subdivided into
even smaller segments (call them
microscopic) of the flag of $Y$.

2) If a small segment of the flag of $X$ is
colored into the color $"\i"$, then all its
microscopic subsegments are also colored
into $"\i"$. If a small segment is colored
into $"\I"$, then the color of its leftmost
microscopic segment may be arbitrary, but
the colors of its remaining microscopic
segments must by $"\i"$. The detailed
exposition can be found in \ref{funI}.

\subsection{Step 3: OPE-algebras over the
collection of functors $\I_p$. The functors
$\P_p$}

Albeit the functors $\I_p$ do not form a
system it is still possible to make a
meaningful definition of an OPE-algebra over
a collection of functors $\I_p$, which we
will now do.

Let $M$ be a $\D_X$-module. An OPE-structure
over a collection of $\I_p$ is a collection
of maps
$$
\ope_{p_S}:M^{\boxtimes S}\to \I_{p_S}(M),
$$
where $p_S:S\to\pt$, with  certain
properties. To formulate them, we first form
maps
$$
\ope_p: M^{\boxtimes S}\to \I_p(M^{\boxtimes
T})
$$
for an arbitrary surjection $p:S\to T$, in
the same way as it was done in the
definition of an OPE-algebra over a system.

The natural maps $\I_p\to \i_p$ give rise to
maps $$\ope^{\i}_p:M^{\boxtimes S}\to
\i_p(M^{\boxtimes T}).
$$

 Let $p=p_np_{n-1}\cdots p_1$, where $p_i:S_i\to S_{i+1}$ and
 $\j^1,\j^2,\ldots,\j^n$ be as above.
 We can construct maps
 $$
 M^{\boxtimes S}\to \j^1_{p_1}\cdots
 \j^n_{p_n}M^{\boxtimes T}
 $$
 as follows:
 $$\xymatrix{
 M^{\boxtimes
 S}\ar[r]^{\ope^{\j^1}_{p_1}}&
 \j^1_{p_1}M^{\boxtimes S_2}\ar[r]^{\ope^{\j^2}_{p_2}}
 & \j^1_{p_1}\j^2_{p_2}M^{\boxtimes
 S_3}\cdots}
 $$

Thus for every object $X\in \zebra_p$, we
have a map
$$
\ope_X: M^{\boxtimes S}\to X(M^{\boxtimes
T}).
$$

Let $u:X\to Y$ be an arrow in $\zebra_p$. We
then have a composition
$$\ope_X\circ u(M^{\boxtimes T}):M^{\boxtimes S}\to
Y(M^{\boxtimes T}).
$$
 We demand that this composition be equal to
 $\ope_Y$. If this is the case, then we say
 that the maps $\ope_{p_S}$ define an
 OPE-algebra structure on $M$ over the
 collection $\I$.

 We can now do the following. Set
 $$\P_p(M)=\liminv_{X\in \zebra_p}
 X(M^{\boxtimes T}).$$
Then the above axiom implies that the maps
$\ope_X$ produce a map
$$
\ope^{\P}_p: M^{\boxtimes S}\to
\P_p(M^{\boxtimes T}).
$$

It is not hard to see that the functors
$\P_p$ form a system. Indeed: let $p=rs$.
Then $\P_s\P_r$ can be realized as an
inverse limit of $X(M^{\boxtimes T}$ over a
full subcategory (=subset with an induced
order) of $\zebra(p)$ formed by all $X$'s
whose flags contain the equivalence relation
on $S$ determined by $r$, whence a map
$$\P_p\to \P_s\P_r.$$

\subsubsection{ Example} Let $S=\{1,2,3\}$
and $p:S\to \pt$.  We have the following
equivalence relations on $S$:

a) the finest one $\omega$;

b) the relations $e_{ij}$, $i\neq j$,
$i,j\in \{1,2,3\}$, in which $i\sim j$, and
the remaining element is only equivalent to
itself;

c) the coarsest relation $\alpha$ in which
all elements are equivalent.

Let $S_{ij}:=S/e_{ij}$. Let $p_{ij}:S\to
S/e_{ij}$ and $q_{ij}:S/e_{ij}\to \pt$. Then
$\P_p$ is the inverse limit of the following
diagram:

$$
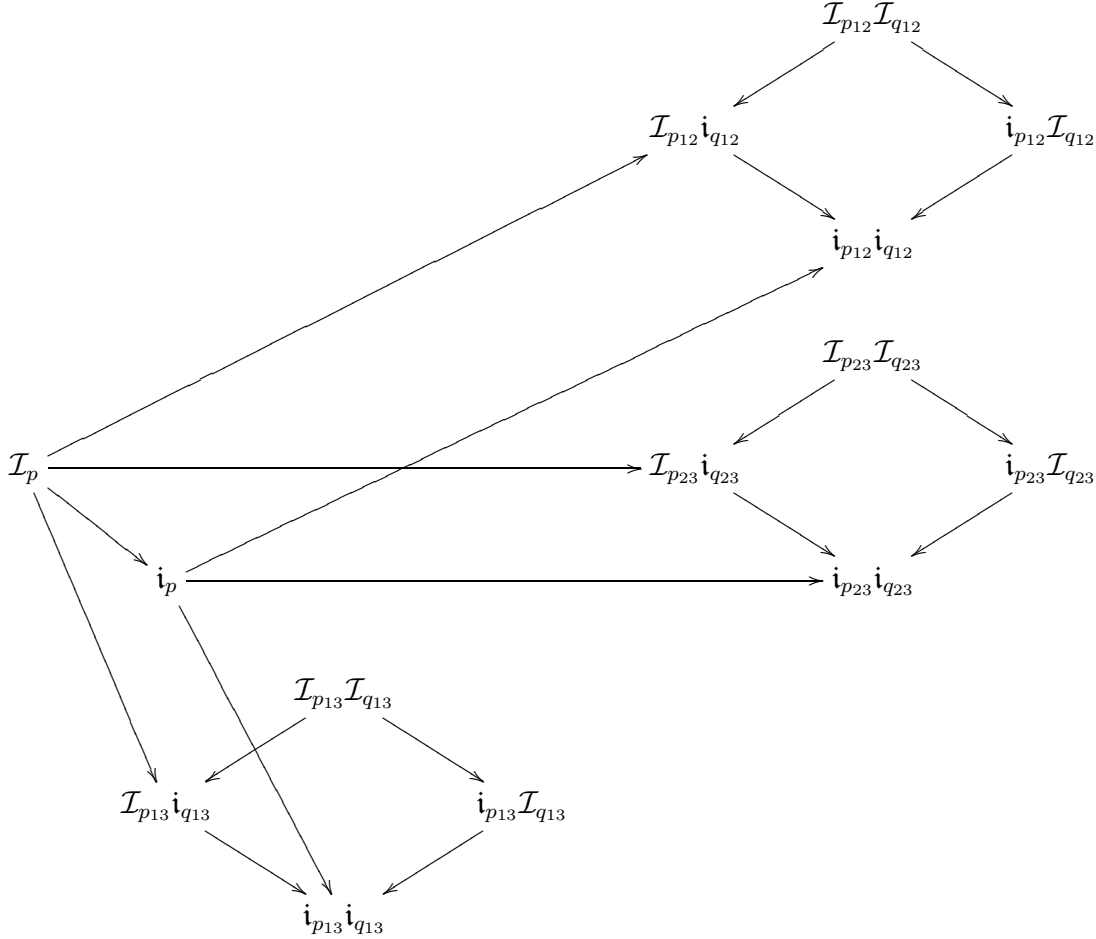
\xymatrix{
      & && & &     \I_{p_{12}}\I_{q_{12}}\ar[dl]\ar[dr] &&\\
&&&&\I_{p_{12}}\i_{q_{12}}\ar[dr]&&\i_{p_{12}}\I_{q_{12}}\ar[dl]&\\
      & && & &      \i_{p_{12}}\i_{q_{12}} &&\\
     & &&& &     \I_{p_{23}}\I_{q_{23}}\ar[dl]\ar[dr] &&\\
\I_p\ar[dr]\ar[rrrr]\ar[rrrruuu]\ar[rddd]&&&&\I_{p_{23}}\i_{q_{23}}\ar[dr]&&\i_{p_{23}}\I_{q_{23}}\ar[dl]&\\
       &\i_p\ar[rrrr]\ar[rrrruuu]\ar[rddd]& & & &     \i_{p_{23}}\i_{q_{23}} &&\\
      &&     \I_{p_{13}}\I_{q_{13}}\ar[dl]\ar[dr] &&&&&\\
&\I_{p_{13}}\i_{q_{13}}\ar[dr]&&\i_{p_{13}}\I_{q_{13}}\ar[dl]&&&&\\
  & &      \i_{p_{13}}\i_{q_{13}} &&&&&
        }
$$

This diagram is co-final to the sub-diagram:
\begin{equation}\label{colimit}\xymatrix{
&\I_{p_{12}}\i_{q_{12}}&\ar[l]^1\I_{p_{12}}\I_{q_{12}}\\
\I_p\ar[ru]\ar[r]\ar[rd] &
\I_{p_{23}}\i_{q_{23}}&\ar[l]^2\I_{p_{23}}\I_{q_{23}}\\
&\I_{p_{13}}\i_{q_{13}}&\ar[l]^3\I_{p_{13}}\i_{q_{13}}}
\end{equation}

We see that $\P_p$ is an extension of $\I_p$
by the kernels of the arrows 1,2,3, which
are $\I_{p_{ij}}\i_{q_{ij}*}$, where $i\neq
j$, $i,j=1,2,3$.
\subsubsection{ } The features of functoriality of the
collection of functors $\P_p$ are inherited
from those of the collection
$\delta_p,\i_p,\I_p$. The most important
ones  are the following ones:

1) the structure of system on the collection
of functors $\P_p$ ;

2) maps $\P_p\delta_q\P_r\to \P_{rqp}$,
where $p,q,r$ are surjections and $q$ is
{\em not} a bijection.

 Let us sketch the definition. First of all,
such a map is uniquely defined by
prescription all compositions

$$f_X:\P_p\delta_q\P_r\to \P_{rqp}\to X,
$$
 where $X$ runs through the set of all
  elements in $\zebra(p)$.

 Let $R:=rqp$; $R:S\to T$; let $Q=qp$.  Let
 $e$ (resp. $e_q$, resp. $e_p$) be the
 equivalence relation
 determined by $R$ (resp. $Q$, resp. $p$).
 It follows that$$
 \omega\geq e_p>e_q\geq e,$$
where $\omega$ is the trivial equivalence
relation on $S$.
 Without loss of generality, we may assume
that $p:S\to S/e_p$, $q:S/e_p\to S/e_{q}$,
$r:S/e_q\to S/e$ are  the natural
projections.

Let now $X$ be given by a flag
$$(\omega=f_1>f_2>\cdots f_{n+1}=e)$$
and a coloring $\j^1,\j^2,\ldots \j^n$.

The map $f_X$ is then specified by the
following conditions:

1) $f_X=0$ unless there exists a $k$ such
that $f_k=e_p>e_q\geq f_{k+1}$ and
$\j_k=\I_k$.

2) Assume that such a $k$ exists.  Let
$\rho:S/{f_{k+1}}\to T$ be the natural
projection. Let $\sigma:S/e_q\to
S/{f_{k+1}}$ so that $r=\rho\sigma$ and
$\sigma q:S/f_k\to  S/f_{k+1}$ is the
natural projection.
 Define
elements $X_r\in \zebra(r)$ and $X_\pi\in
\zebra(\pi)$  as follows:

$X_r$ is given by the flag
$\omega>f_1>\cdots>f_k=e_r$, and the
coloring  $(\j^1,\j^2,\ldots,\j^{k-1})$;

$X_p$ is given by the flag
$$\omega_{S/f_k}\geq
f_{k+1}/f_k>f_{k+2}/f_k>\ldots >e/f_k,
$$
of equivalnce relations  on $S/f_k$. It
follows that $X$ decomposes as
$X=X_p\I_{\sigma q}X_\rho$.

The map $f_X$ then goes as follows:
$$
\P_p\delta_q\P_r\to
\P_p\delta_q\P_\sigma\P_\rho\to X_p\delta_q
\I_\sigma X_\rho\to X_p\I_{\sigma
q}X_\rho=X.
$$
\subsubsection{Example} Let us come back to
our example $S=\{1,2,3\}$ and $p:S\to \pt$.
We know that $\P_p$ is the inverse limit of
the diagram (\ref{colimit}). Let us describe
the map
$$
\delta_{p_{12}}\P_{q_{12}}\to \P_q.
$$
First of all, $\P_{q_{12}}\to \I_{q_{12}}$
is an isomorphism.

We then have maps
$$
\delta_{p_{12}}\I_{q_{12}}\to
\I_{p_{12}}\I_{q_{12}}
$$
and
$$
\delta_{p_{12}}\I_{q_{12}}\to \I_q.
$$
The diagram
$$\xymatrix{
\delta_{p_{12}}\I_{q_{12}}\ar[r]\ar[d]&
\I_{p_{12}}\I_{q_{12}}\ar[d]\\
\I_p\ar[r]&\I_{p_{12}}\i_{q_{12}}}
$$
turns out to be commutative (this is hidden
behind the words "these maps behave well
with respect to the other elements of the
structure" after (\ref{delI})). Furthermore,
the compositions
$$
\delta_{p_{12}}\I_{q_{12}}\to \I_q\to
\I_{p_{23}}\i_{q_{23}},
\I_{p_{13}}\i_{q_{13}}
$$

as well as

$$
\delta_{p_{12}}\I_{q_{12}}\to
\I_{p_{12}}\I_{q_{12}}\to
\i_{p_{12}}\I_{q_{12}}
$$
all vanish, whence the desired map
$\delta_{p_{12}}\P_{q_{12}}\to \P_p$.

Consider now the map
$\P_{p_{12}}\delta_{q_{12}}\to  \P_p$.
Again, we have an isomorphism
$$
\P_{p_{12}}\to \I_{p_{12}}.
$$
 We also have a map
$$\I_{p_{12}}\delta_{q_{12}}\to
\I_{p_{12}}\I_{q_{12}},
$$
 the composition
$$
\I_{p_{12}}\delta_{q_{12}}\to
\I_{p_{12}}\I_{q_{12}}\to
\I_{p_{12}}\i_{q_{12}}
$$
being zero. Furthermore, the sequence
$$
0\to\I_{p_{12}}\delta_{q_{12}}\to
\I_{p_{12}}\I_{q_{12}}\to
\I_{p_{12}}\i_{q_{12}}\to 0
$$
is exact. Therefore, the map
$$
\I_{p_{12}}\delta_{q_{12}}\to \P_p$$
realizes an embedding of the kernel of the
arrow 1 in (\ref{colimit}) into   $\P_p$.

Describe the map $\delta_p\to \P_p$. It is
given by the inclusion $\delta_p\to \I_p$;
since the composition of this map with every
arrow coming out of $\I_p$ vanishes, this is
a well defined map. This map can be also
described as a composition:
$$
\delta_p\cong
\delta_{p_{12}}\delta_{q_{12}}\to
\delta_{p_{12}}\I_{q_{12}}\to \P_p.
$$

Finally,  the map
$$
\delta_{p_{12}}\delta_{q_{12}}\to \P_p$$ is
given by
$$
\delta_{p_{12}}\delta_{q_{12}}\to
\I_{p_{12}}\I_{q_{12}}$$ and is different
from the previous one!

The maps that we considered fit into a
commutative diagram
$$\xymatrix{
        &&&&   \P_p &&&&\\
&\delta_{p_{12}}\I_{q_{12}}\ar[urrr]
&\I_{p_{12}}\delta_{q_{12}}\ar[urr] &
\delta_{p_{23}}\I_{q_{23}}\ar[ur]
&&\I_{p_{23}}\delta_{q_{23}}\ar[ul]&
\delta_{p_{13}}\I_{q_{13}}\ar[ull]
&\I_{p_{13}}\delta_{q_{13}}\ar[ulll]&
      \\
&\ar[u]&
\delta_{p_{12}}\delta_{q_{12}}\ar[u]\ar[ul]&\ar[u]
&
\delta_{p_{23}}\delta_{q_{23}}\ar[ur]\ar[ul]&&\ar[u]&
\delta_{p_{13}}\delta_{q_{13}}\ar[u]\ar[ul]&\\
&&&&\delta_p\ar@{-}[ul]\ar@{-}[ulll]\ar@{-}[urr]&&&&}
$$

This diagram specifies a map from the direct
limit of its three lowest floors to $\P_p$.
It turns out that this map is an inclusion
whose cokernel is isomorphic to $\i_p$ via
the natural map $\P_p\to \I_p\to \i_p$.

This implies  that $\P_p$ has a three term
filtration  (the two lowest floors are
combined) whose successive quotients are

1)$$\delta_p
\oplus_{i<j}\delta_{p_{ij}}\delta_{q_{ij}};$$

2)$$\oplus_{i<j}
\delta_{p_{ij}}\i_{q_{ij}}\oplus
\i_{p_{ij}}\delta_{q_{ij}}
$$

3) $\i_p$.
 \subsubsection{Filtration on $\P$} The
 filtration on  functors $I_p$ define a
 filtration on $\P_p$. See \ref{Filt}-\ref{Filt1} for its
 description.  Its successive quotients are
 direct sums of the terms of the form
 $$
 \i_{p_1}\delta_{q_1}\i_{p_2}\delta_{q_2}\cdots
 \delta_{q_n}\i_{p_{n+1}}
 $$
 with fixed $n$.
 Here $p_{n+1}q_np_n\cdots
 q_1p_1=p$; all $p$'s and $q$'s are surjective
 and all $q$'s are not bijective.
\subsection{Resolution $\R$}
We are now ready to define the desired
resolution. The starting point is the maps
$\P_p\to \i_p$, which are surjections. Our
goal is to kill the kernel, which turns out
to be spanned by the images of all maps
$$
\P_a\delta_b\P_c\to \P_p,$$
 where $cba=p$.

Thus, it makes sense to assign
$$\R_p^{0}:=\P_p$$ and
$$
\R_p^{-1}:=\oplus \P_a\delta_b\P_c,
$$
where the direct sum is taken over all
sequences
$$
S\stackrel c\to S/e_1\stackrel b\to
S/e_2\stackrel a\to T,$$ where $ e_1>e_2>e$
are equivalence relations on $S$, $e$ is
determined by $p$, and $a,b,c$ are natural
projections. The differential is given by
the above described  maps
$\P_c\delta_b\P_a$.

The $n$-th term $\R^{-n}_p$ is given by the
direct sum of the terms
$$
\P_{p_1}\delta_{q_1}\P_{p_2}\delta_{q_2}\cdots
\delta_{q_n}\P_{p_{n+1}},
$$
where the sum is taken over all diagrams of
the form $$ S\stackrel{p_1}\to
S/e_1\stackrel{q_1}\to
S/f_1\stackrel{p_2}\to
S/e_2\stackrel{q_2}\to
S/f_2\stackrel{p_3}\to
\cdots\stackrel{q_n}\to
S/f_n\stackrel{p_{n+1}}\to T,$$ where
$$e_1>f_1\geq e_2>f_2\geq\cdots >f_n\geq e,$$
and $p_i,q_j$ are all natural projections.
The differential $d:\R^{-n}_p\to
\R^{-n+1}_p$ is given by the alternated sum
of maps induced by

a) $\P_{p_i}\delta_{q_i}\P_{p_{i+1}}\to
\P_{p_{i+1}q_ip_i}$ and

b)
$\delta_{q_i}\P_{p_{i+1}}\delta_{q_{i+1}}\to
\delta_{q_{i+1}p_{i+1}q_1}$, which are
non-zero iff $p_{i+1}=\Id$, in which case
they are natural isomorphisms.

One then has to check that $d^2=0$ and to
define on $\R_p$ a structure of system. For
all this we refer the reader to \ref{resdif}
\subsubsection{Example}
Let $S=\{1,2,3\}$. Then the complex $\R_p$
is depicted as follows:

$$
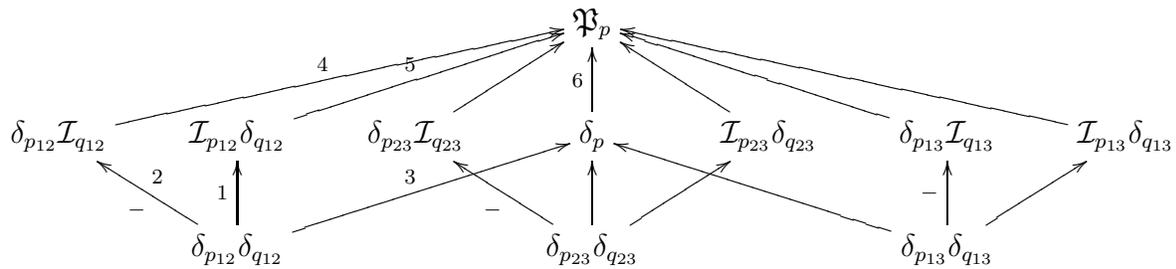
\xymatrix{
        &&&&   \P_p &&&&\\
&\delta_{p_{12}}\I_{q_{12}}\ar[urrr]^{4}
&\I_{p_{12}}\delta_{q_{12}}\ar[urr]^{5} &
\delta_{p_{23}}\I_{q_{23}}\ar[ur]
&\delta_p\ar[u]^{6}
&\I_{p_{23}}\delta_{q_{23}}\ar[ul]&
\delta_{p_{13}}\I_{q_{13}}\ar[ull]
&\I_{p_{13}}\delta_{q_{13}}\ar[ulll]&
      \\
&&
\delta_{p_{12}}\delta_{q_{12}}\ar[u]^{1}\ar[ul]^{-}_{2}\ar[urr]^{3}&
&
\delta_{p_{23}}\delta_{q_{23}}\ar[ur]\ar[ul]^{-}\ar[u]&&
\delta_{p_{13}}\delta_{q_{13}}\ar[u]^{-}\ar[ur]\ar[ull]&\\
}
$$
where all the arrows are the natural maps;
the arrows marked with $-$ are taken with
the negative sign. Let us check that
$d^2=0$. It suffices to check that
$$
d^2|_{\delta_{p_{12}}\delta_{q_{12}}}:
\delta_{p_{12}}\delta_{q_{12}}\to \P_p$$ is
zero. This reduces to checking that the
compositions

$$Ad^2|_{\delta_{p_{12}}\delta_{q_12}}:
\delta_{p_{12}}\delta_{q_{12}}\to
\P_p\stackrel A\to \I_p;
$$

$$
B_{ij}d^2|_{\delta_{p_{12}}\delta_{q_12}}:
\delta_{p_{12}}\delta_{q_{12}}\to
\P_p\stackrel{B_{ij}}\to
\I_{p_{ij}}\I_{q_{ij}}$$

do all vanish. Let us so do.

$Ad^2$. We have: $A42=A63$; $A51=0$. Hence
$Ad^2=A42-A63+A51=0$.

$B_{ij}d^2$. If $\{i,j\}\neq \{1,2\}$, then
all three maps $$
B_{ij}42=B_{ij}63=B_{ij}51=0.$$ Consider now
the remaining case $B_{12}d^2$. We then
have: $B_{12}42=B_{12}51$ and $B_{12}63=0$,
which implies that $B_{12}d^2=0$.

\def\m{{\frak m}}
\subsection{The system $<\m>$ and the map
$ <\R>\to <\m>$} Recall that the whole
purpose of constructing $<\R>$ was to
establish a link between the systems $<\i>$
and $<\l>$. Unfortunately, there is no
direct map $<\R>\to <\l>$. Instead, we shall
construct a map $m:<\R>\to <\m>$ satisfying
the properties described in \ref{Rm}.

 Define
a map $m_p:\R_p\to \m_p$ by the following
conditions.

1) $m_p$ vanishes on all terms
$$
\R_{p_1}\delta_{q_1}\R_{p_2}\delta_{q_2}\cdots
\delta_{q_n}\R_{p_{n+1}}
$$
where at least on $p_i\neq \Id$. Otherwise,
$m_p$ is the identical embedding onto the
term $$\delta_{q_1}\delta_{q_2}\cdots
\delta_{q_n}$$ of $\m_p$.

Denote by
$$
l_p:\R_p\to \m_p\to\delta_p[1]
$$
the natural composition.
\subsection{The  additional structure induced by the
maps (\ref{I-asadd})}\label{I.I} Recall that
the collection of maps $<\I_p>$ has a
functoriality $(\ref{I-asadd})$ which we
have never used. It turns out that this
additional functoriality yields an
additional structure on the system $<\R>$.

To obtain this additional structure one has
first to understand the additional structure
on the system $<\P>$ produced by these
functors. Consider some examples.

Let $A=\{1,2\}$ be a 2-element set and let
$g:A\to \pt$. Let $f:S\to T$ be a
surjection. Let $$f\disjoint g:S\disjoint A
\to T\disjoint \pt$$ be a disjoint union.

We may define two maps
$$n_1,n_2:\P_{f\disjoint g}\to \I_{f\disjoint
\Id_A}\I_{\Id_T\disjoint g}
$$
The map $n_1$ is just the natural projection
onto  a member of $\zebra_{f\disjoint g}$.
The map $n_2$ is the composition

$$
\P_{f\disjoint g}\to \I_{f\disjoint g}\to
\I_{f\disjoint \Id_A}\I_{\Id_T\disjoint g},
$$
where we first apply the natural projection
and then the map (\ref{I-asadd}).

It follows that the compositions of
$n_1,n_2$ with the map
$$\lambda:\I_{f\disjoint
\Id_A}\I_{\Id_T\disjoint g}\to
\I_{f\disjoint \Id_A}\i_{\Id_T\disjoint g}$$
do coincide, therefore the difference
$n_2-n_1$ determines a map to the kernel of
$\lambda$ , i.e. a map
$$
\xi(f,g):\P_{f\disjoint g}\to
\I_{f\disjoint\Id_A}\delta_{\Id_T\disjoint
g}.
$$
This is only true because of the special
form of $g$.

For a general $g:A\to B$ the kernel of
$\lambda$ is spanned by the images of all
maps
$$
\I_{f\disjoint\Id_A}\delta_{u_1}\I_{u_2}\to
\I_{f\disjoint \Id_A}\I_{\Id_T\disjoint g}
,$$ where $u_1,u_2$ are surjections, $u_1$
is not a bijection, and
$$u_2u_1=\Id_T\disjoint g.$$
So that the structure of $n_1-n_2$ becomes
more complicated.

  Nevertheless,
one can define maps $$
\xi(f,g):\P_{f\disjoint g}\to
\I_{f\disjoint\Id_A}\delta_{\Id_T\disjoint
g}
$$ for an arbitrary $g$ by means of
the following inductive process. Let
$|g|=|A|-|B|$. Since $g$ is a surjection,
$|g|\geq 0$. If $g$ is a bijection, then we
have a natural isomorphism
$$
\P_{f\disjoint g}\to
\P_{f\disjoint\Id_A}\delta_{\Id_T\disjoint
g},
$$
because $\Id_T\disjoint g$ is a bijection.

Set $\xi(f,g)$ to be the composition of this
isomorphism with the natural map
$$
\P_{f\disjoint\Id_A}\delta_{\Id_T\disjoint
g}\to
\I_{f\disjoint\Id_A}\delta_{\Id_T\disjoint
g}.
$$

Let us now assume that $\xi(f,g)$ is defined
for all $g$ with $|g|<N$. Define it for all
$g$ with $|g|=N$. Let $e$ be an equivalence
relation on $A$ induced by $g$. Let
$\omega_S>\ve\geq e$, let $h_\ve:A\to A/\ve$
and $k_\ve:S/\ve\to B$ so that $k_\ve
h_\ve=g$.

Define a map
$$
C(\ve):\P_{f\disjoint g}\to \P_{f\disjoint
k_\ve}\P_{\Id_T\disjoint \Id_{A/\ve}}\to
\I_{f\disjoint \Id_A}\delta_{\Id_T\disjoint
k_\ve}\I_{\Id_T\disjoint \Id_{A/\ve}}\to
\I_{f\disjoint \Id_A}\I_{\Id_T\disjoint g}.
$$

 set
 $$\xi'(f,g):=-\sum_\ve C(\ve).$$
If $\xi'(f,g)$ passes through
$\I_{f\disjoint \Id_A}\delta_{\Id_T\disjoint
g}$, it determines a map $$ \P_{f\disjoint
g}\to \I_{f\disjoint
\Id_A}\delta_{\Id_T\disjoint g},$$ which we
assign to be $\xi(f,g)$. It can be checked
that if this rule was obeyed when $\xi(f,g)$
was defined for all $g$ with $|g|<N$, then
$C'(f,g)$ passes through $\I_{f\disjoint
\Id_A}\delta_{\Id_T\disjoint g}$ and gives
rise to the map $\xi(f,g)$.

On the next step the maps $\xi(f,g)$ are
lifted to  maps $$c(f,g):\P_{f\disjoint
g}\to \P_{f\disjoint
\Id}\delta_{\Id\disjoint g},$$ which in turn
produce maps
$$
s(f,g):\R_{f\disjoint g}\to \R_{f\disjoint
\Id}\delta_{\Id\disjoint g}
$$
with a non-zero differential, which is
described in (\ref{diffs}) .
\subsubsection{} Thus, we described a  construction of a
pre-symmetric (up-to homotopies) system
$<\R>$.

\def\An{{\cal An}}
\def\Ga{\Gamma}
\def\Gl{{\Gamma^{\circ}}}
\def\Gll{{\Gamma^{\circ\circ}}}
\def\Al{{\A^{\circ}}}

\section{Realization of the system $<\Rs>$ in the spaces
of real-analytic functions}

Our answer to the renormalization problem l
is given in terms of a system $<\Rs>$. To be
able to get a physically meaningful answer
we need an OPE expansion in terms of series
of real-analytic functions on the $Y^S$
minus all generalized diagonals.

The nicest possible way to do it includes
constructing a system  which is explicitly
linked to the spaces of real-analytic
functions on $Y^S$ minus all generalized
diagonals and constructing a map from
$<\Rs>$ to this system. Unfortunately, we do
not know how to realize this project. The
problem is that  arbitrary real-analytic
functions do not have a good asymptotic
expansion in a neighborhood of generalized
diagonals, therefore, we cannot form a
system based on such spaces.

Let us describe a palliative measure we take
instead.

First of all, we shall work with spaces of
global sections rather than with sheaves.
So, whenever we use a notation for a sheaf,
it will actually mean the space of global
sections. If our sheaf is a
$\D_{X^S}$-module, then its space of global
sections is a module over the space of
global sections of $\D_{X^S}$. Whenever we
say " a $\D_{X^S}$-module", we will actually
mean  "a module over the space of global
sections of $\D_{X^S}$.

Let $Y\subset Y^S$ be the main diagonal. We
pick a  vector field  which contracts
everything to $Y$ and take analytic
functions on $Y^S$ minus the complement to
all generalized diagonals which are
generalized eigenvalues of this field.

Denote this space spanned by such functions
by $\Al_S$.  This space has a grading given
by the generalized eigenvalue. Let
$\Al_S^{\geq N}$ be the span of all elements
whose generalized eigenvalue is $\geq N$.

Then the spaces $\Al_S^{\leq N}:=
\Al_S/\Al_S^{\geq N}$ do not depend on a
choice of particular vector field.

Let $p:S\to T$ be a projection. We define a
functor $ \Al_p $ from the category of
$\D_{X^T}$-modules to the category of
$\D_{X^S}$-modules by the formula
$$
\Al_p(M)=\liminv_N
i_p^\ha(M)\otimes_{\O_{X^S}} \Al_S^{\leq N}.
$$

These functors do not form a system.
Nevertheless,  given  $a\in \Al_p(M)$, $b\in
\Al_{p_1}\Al_{p_2}(M)$, where $p=p_2p_1$,
one can  say whether $b$ is an asymptotic
decomposition of $a$ or not. The problem is
that not every $a$ has such a decomposition.

We define a functor
$$\Gll(p_1,p_2)\subset
\Al_p\oplus \Al_{p_1}\Al_{p_2}
$$
so that $\Gll(p_1,p_2)(M)$ consists of all
pairs $(a,b)$ such that $b$ is an asymptotic
decomposition of $a$. In other words,
instead of a  map $\Al_p\to
\Al_{p_1}\Al_{p_2}$ we have a
"correspondence" given by $\Gll(p_1,p_2)$.

Next, we construct maps $\int_p:\N_p\to
\Al_p$. We then show that these maps are
compatible with the correspondences
$\Gll(p_1,p_2)$ as follows:

Let $$ \int_{p_1,p_2}:\N_p\to
\N_{p_1}\N_{p_2}
\stackrel{\int_{p_1}\otimes\int_{p_2}}\longrightarrow
\Al_{p_1}\Al_{p_2}.
$$

We then show that
$$\int_p\oplus\int_{p_1,p_2}:\N_p\to
\Al_p\oplus \Al_{p_1}\Al_{p_2}
$$
passes through $\Gll(p_1,p_2)$.

This construction provides us with an OPE
product on $M$ in terms of series of
real-analytic functions on $Y^S$.

The construction of the maps $\int_p$
resembles the construction of the maps
\ref{I.I}, which is based on the maps
(\ref{I-asadd}).  The construction of
$\int_p$ is based on the existence of
asymptotic decompositions of generalized
functions from $\C_S$ near generalized
diagonals. Namely, let $p_S:S\to \pt$, and
let $p_S=p_2p_1$ be a decomposition. We
construct maps
$$
\C_S\to \A_{p_1}\I_{p_2},$$ where $\A_{p_1}$
is constructed in the same way as
$\Al_{p_1}$ but generalized functions which
are non-singular on the complement to
generalized diagonals  and are generalized
eigenvectors of the vector field which
shrinks everything to the main diagonal, are
used.

\bigskip

 \centerline{\bf \Large PART III. Technicalities}
In the concluding part of the paper we give
constructions and proof required for
everything in the previous part to work.
This includes

1) constructing the system $<\R>$ and
endowing it with a pre-symmetric structure;

2) Bogolyubov-Parasyuk lifting theorem

3) more details on the symmetrization
procedure and on the renormalization in
symmetric systems. To this end we need to
develop certain machinery ("pseudo-tensor
bodies").

4) real-analytic interpretation of the
symmetric system that we obtain from $<\R>$.

\section{Constructing the system $<\R>$}
\subsubsection{}
Let $\Z={\R}^N$, where $N$ is a fixed
natural even number. We fix the coordinates
$x^1,x^2,\ldots,x^N$ on $\Z$. For $x\in \Z$
we set $q(x)=\sum_{i=1}^N (x^i)^2$. Also we
take the standard orientation on $\Z$.
\subsubsection{}
Let $S$ be a finite set. Let $\Z^S$ be the
space of functions $S\to \Z$. Let
$[n]=\{1,2,\ldots,n\}$, then $\Z^{[n]}\cong
\Z^n$. Since $Z$ is even-dimensional, the
orientation on $\Z$ produces canonically an
orientation on $\Z^S$. Thus, $\Z^S$ will be
assumed to have an orientation.

Let $e$ be an equivalence relation on
$S$.Denote by $\Delta_e\subset \Z^S$ the
corresponding generalized diagonal
consisting of points $\Y:S\to \Z$ such that
$s\sim_e t\Rightarrow \Y(s)=\Y(t)$.

Let $f:S\to T$ be a map of finite sets. We
have an induced map $f^{\#}:\Z^T\to \Z^S$.
Let $p_e:S\to S/e$. Then $\Delta_e=\Im p_e$.
If $f$ is surjective, then $f$ identifies
$\Delta_e$ with $\Z^{S/e}$. We will use this
identification.

For two equivalence relations $e_1$ and
$e_2$ on a finite set $S$ we write $e_1\leq
e_2$ iff
 $s\sim_{e_2} t\Rightarrow s \sim_{e_1} t$.
We have $e_1\leq e_2$ iff
 $\Delta_{e_1}\subseteq \Delta_{e_2}$.
Denote by $\alpha$ the least equivalence
relation (i.e. in  every two points are
equivalent) and by $\omega$ the greatest
equivalence relation (i.e.  every two
distinct points are {\em not} equivalent).
Let $s,t\in S$ be distinct elements.

Let $T\subset S$. Denote by $e_T$ the
equivalence relation in which two distinct
elements are equivalent iff both of them are
in $T$. For example, $\omega=e_\emptyset$;
$\alpha=e_S$.
 Set $\De_T:=\De_{e_T}$; $\De_{st}:=\De_{\{s;t\}}$.
\subsubsection{} Denote
$$U_S=\Z^S-\cup_{e\neq \omega}\Delta_e.
$$
Obviously, a point $\Y:S\to \Z$ is in $U_S$
iff the map $\Y$ is injective.

\subsubsection{} Let $s,t\in S$ be distinct elements.
Denote by $q_{st}:\Z^S\to \R$ the function
defined according to the rule
\begin{equation}\label{qs1s2}
q_{st}(\Y)=q(\Y(s)-\Y(t)),
\end{equation}
where $\Y:S\to\Z$ is a point in $\Z^S$ and
 $q$ is the standard quadratic form on $\Z$.
Of course, the set of zeros of $q_{st}(\Y)$
is $\De_{st}$.
\subsubsection{} Denote by $B_S$ the space of functions
$U_S\to \C$ which can be expressed as a
ratio $P(Y)/Q(Y)$, where $P$ is an arbitrary
polynomial and $Q$ is a product of
nonnegative integer
 powers of $q_{st}$ for arbitrary $s,t$.
\subsubsection{} As usual, we denote by $\S_{\Z^S}$
the   space of compactly supported top forms
on $\Z$ and by $\S'_{\Z^S}$ the space of
distributions on $\S_{\Z^S}$. Any smooth
function on $\Z^S$ will be regarded as a
distribution in the usual way (recall that
the orientation on $\Z$ produces a canonical
orientation on $\Z^S$.)
\subsection{Quasi-polynomial
distributions}\label{t}\label{cs}
 We have a diagonal action of
the group $\R^N$ on $\Z^S$ by translations.
This induces an action of the abelian
$N$-dimensional Lie algebra
  $\t_N$ on $\S_{\Z^S}$; $\S'_{\Z^S}$.
Call a distribution $f$ quasi-polynomial if
there exists an $M$ such that $t^Mf=0$. Let
$\P_{\Z^S}$ be the subspace of all
quasi-polynomial distributions.

We have natural continuous  maps
$$
T_{S_1S_2}:\S'_{\Z^{S^1}}\otimes
\S'_{\Z^{S^2}}\to \S'_{\Z^{S^1\disjoint
S^2}}.
$$ which induce maps:
\begin{equation}\label{tensor}
T_{S_1S_2}:\P_{\Z^{S^1}}\otimes
\P_{\Z^{S^2}}\to \P_{\Z^{S^1\disjoint S^2}}.
\end{equation}
\subsection{Definition of subspaces $C_S\in
\P_{\Z^S}$}\label{genfun} We define these
subspaces recursively.

1) If $S$ is empty or has only one element,
we set $C_S:=\P_{\Z^S}$.

2) Suppose, we have already defined
$C_S\subset \P_{\Z^S}$ for all $S$ with at
most $m$ elements.  For an $S$ with $m+1$
elements, we say that a quasi-polynomial
distribution $f$ on $\Z^S$ is in $C_S$ iff
for any partition $S=S_1\disjoint S_2$,
there exists an integer $M$ such that
\begin{equation}\label{factor}
(\Prod_{s_1\in S_1;s_2\in S_2
}q_{s_1s_2}^M)f\in T_{S_1S_2}(C_{S_1}\otimes
C_{S_2}),
\end{equation}
where $T_{S_1S_2}$ is as in (\ref{tensor})
and $q_{s_1s_2}$ is as in (\ref{qs1s2}).
\subsection{Example} Let $S=\{1,2\}$. We will also use the
symbol $[2]$ for $\{1,2\}$. For $\Y:[2]\to
\Z$ we write $\Y_1:=\Y(1)$ and
$\Y_2:=\Y(2)$.
 Then
$f\in C_S$ iff $f$ is quasi-polynomial and
there exists an $M$ such that
$q(\Y_1-\Y_2)^Mf=P(\Y_1,\Y_2)$, where $P$ is
a polynomial.

Define a map $\pi:C_S\to B_S$ by $\pi
f=P/q^M$. It is clear that this map is well
defined and that
 $\Ker \pi$ consists
of all functions $f\in C_S$ supported on the
diagonal. Denote $C_{S,\Delta}:=\Ker \pi.$
We are going to describe this space.

  Let $s=(\Y,\Y)$ be a point on the diagonal.
Then on any relatively compact neighborhood
$U$ of $s$, any distribution supported on
the diagonal is of the form
\begin{equation}\label{di}
f=\sum_{\mu}f_{\mu}(\Y_1)\delta^{\mu}(\Y_2-\Y_1),
\end{equation}
where the sum is taken over a finite set of
multi-indices $\mu$ and $f_\mu$ are
distributions on $\Z$.

Suppose that $f\in C_{S,\Delta}$. Then
 $f$ is quasi-polynomial, $t^{M'}f=0$ for some $M'$.
Therefore, $t^{M'}f_\mu=0$ for all $\mu$
meaning that each $f_\mu$ is a polynomial of
degree less than $M'$. This immediately
implies that (\ref{di}) is true everywhere
for some polynomials $f_\mu$. Conversely, if
all $f_\mu$ are polynomials, then $f\in
C_{S,\Delta}$.

The map $\pi$ defines an injection
$C_S/C_{S,\Delta}\to B_S$. Let us show that
this is in fact a bijection. This means that
for any integer $M>0$ and any polynomial
$P(\Y_1,\Y_2)$ there exists a distribution
$F$ such that $Fq(\Y_1-\Y_2)^M=P$. It is
sufficient to do it for $P=1$.
 Let us construct such an $F$.
\subsubsection{}
To this end, take an $fd\Y_1d\Y_2\in
\S_{\Z^S}$, where $d\Y$ is the standard
volume form on $\Z=R^N$ and consider the
expression$$
Z(s,f)=\int_{\Z^2}f(\Y_1,\Y_2)q(\Y_1-\Y_2)^sd\Y_1d\Y_2.
$$
\begin{Claim}
This integral  uniformly converges on any
strip
 $\text{Re} s> K$, where $K>-N/2$.
\end{Claim}
\pf To show it, change the variables
$\Y=\Y_1$ $z=\Y_2-\Y_1$, and
$g(\Y,z)=f(\Y,\Y+z)$. Then
$$
Z(s,f)=\int_{\Z^2}g(\Y,z)q(z)^sd\Y dz.
$$

Let $S^{N-1}\subset \Z$ be the unit sphere
$q(\Y)=1$. Let $\alpha:R_+\times S^{N-1}\to
\Z$ be the map: $\alpha(r,n)=rn$. Let $dn$
be the measure on $S^{n-1}$
 determined by $q$. Then
$$
Z(s,f)=\int\limits_0^\infty
r^{2s+N-1}h_f(r)dr,
$$
where
$$
h_f(r)=\int\limits_{\Z\times
S^{N-1}}g(\Y,rn)d\Y dn.
$$

 Whence the statement.
\endpf
\subsubsection{} Since $Z(s,f)$ is (up to a shift)
the Mellin transform of $h_f$, we know that
$Z(s,f)$ has a meromorphic continuation to
the whole complex plane, the poles can only
occur at $s=-(N+k)/2$, $k=0,1,2, \ldots$ and
are of at most first order. Denote
$$U_M(f)=\res_{s=-M}\frac{Z(s,f)}{s+M}.
$$
\begin{Claim} $U_M$ is a distribution.
\end{Claim}
\pf Set $s'=2s+N-1$; $M'=-2M+N-1$ ($M'$
corresponds to $s=M$). Integration by parts
yields:
$$
Z(s,f)=\displaystyle\frac{
(-1)^P}{(s'+1)(s'+2)\cdots (s'+P)}
\int\limits_0^\infty(\frac {d^P}{dr^P}
h_f(r))r^{s'+P}dr,
$$
whenever $s'+P>0$ and $s'\neq
-1,-2,\ldots,-P+1$. Choose $P$ large enough
so that $M'+P>1$. Set
$$
l(s',r)=\displaystyle\frac{
(-1)^P}{(s'+1)(s'+2) \cdots (s'+P)}r^{s'+P}.
$$
and
$$
\lambda(r)=\res_{s'=M'}l(s',r)=r^(M'+P)C_{M',P},
$$
where $C_{M',P}$ is a constant. Thus,
$$
U_{M'}(f)=C_{M',P}\int\limits_0^\infty
r^{M'+P}\frac {d^P}{dr^P}h_f(r)dr.
$$
It is clear that the function $h_f$ is
smooth and rapidly decreasing as $r\to
\infty$. Furthermore, $f\mapsto h_f$ is a
continuous map from $\S_{\Z^S}$ to the space
of rapidly decreasing infinitely
differentiable functions on $[0,\infty]$, in
which the topology is given by the family of
seminorms
 $$
\|h\|_{K,L}=\max_{r}r^L|h^{(K)}(r)|.
$$
Since the map $f\mapsto h_f$  is continuous,
so is $U_M$, whence the statement.
\begin{Claim}

1.$q^M(Y_1-Y_2)U_M=1$

2. $\t.U_M=0$. (for $\t$ see \ref{t}).
\end{Claim}
\pf 1. \begin{eqnarray*}
U_M(q^M(Y_1-Y_2)f)=\res_{s=-M}\frac{Z(s,q^Mf)}{s+M}=
\res_{s=-M}\frac{Z(s+M,f)}{s+M}\\
=\res_{s=0}\frac{Z(s,f)}{s} =Z(0,f)=\int
fd\Y_1d\Y_2.
\end{eqnarray*}
2. Clear.
\endpf
\begin{Corollary} $U_M\in C_S$; $\pi(U_M)=1/q^M$.
Therefore, $\pi$ is surjective.
\end{Corollary}

Thus, we have an exact sequence:
\begin{equation}\label{ext}
0\to C_{S,\Delta}\to C_S\to B_S\to 0.
\end{equation}

 We see that this is an extension
of $\DD_{S}$-modules. From our description
of $C_{S,\De}$, it follows that
$C_{S,\De}=i_*\O_{\De}$, where $i:\D\to
\Z^2$ is the diagonal embedding. One can
show that this extension does not split. One
can construct a similar extension  when $N$
is odd, in which case it splits; the reason
is that the Green function for the Laplace
operator requires extraction a square root.

\section{Study of $C_S$}\label{genfun1}
\subsection{Action of differential operators}
Denote by $\DD_{\Z^S}$ the algebra of
polynomial differential operators on $\Z^S$.
it is clear that each $\P_{\Z^S}$ is a
$\DD_{\Z^S}$-module.
\begin{Claim} Each $C_S$ is a  $\DD_{\Z^S}$-submodule
of $\P_{\Z^S}$.
\end{Claim}
\pf This is obvious when $S$ has 0 or 1
element. For an arbitrary $S$ the proof can
be easily done by induction. Indeed, we only
need to check that for any $f\in C_S$ and
any polynomial differential operator $D$,
$Df$ satisfies (\ref{factor}). It suffices
to consider only operators of zeroth and
first order. If the order of $D$ is zero,
the statement is immediate. Assume that the
order of $D$ is 1 and $D1=0$. Let
 $$
Q_{S_1S_2}=\Prod_{s_1\in S_1;s_2\in S_2
}q_{s_1s_2}
$$
and $fQ^M\in T_{S_1S_2}(C_{S_1}\otimes
C_{S_2})$. It is immediate that the space on
the right hand side is a
$\DD_{\Z^S}$-submodule of $\P_{\Z^S}$. We
then have
$$
Q^{M+1}Df=D(Q^{M+1}f)-(M+1)Q^M(DQ)f\in
 T_{S_1S_2}(C_{S_1}\otimes C_{S_2}).
$$
\endpf

\subsection{Map $\pi:C_S\to B_S$ and its surjectivity}
\subsubsection{} Let $f\in C_S$.
\begin{Claim} There exists a natural number $M$
such that
\begin{equation}\label{reg}
(\prod_{\{s,t\}\subset S} q_{st})^Mf=P,
\end{equation}
where $P$ is a polynomial and the product is
taken over all 2-element subsets of $S$.
\end{Claim}
\pf  This is obvious wnen $S$ is empty or
has only one element. For general $S$ the
argument follows from (\ref{factor}) by
induction on the number of elements in $S$.
\endpf

Write
$$
\pi(f)=P/(\prod_{s\neq t} q_st)^M.
$$
It is clear that $\pi(f)$  depends only on
$f$ and that $\pi:C_S\to B_S$ is a
$\DD_{\Z^S}$-module map.
\begin{Proposition}\label{p} The map $\pi$ is surjective.
\end{Proposition}
\subsection{Proof of Proposition \ref{p}}
It is sufficient to construct for every
integer $M>0$ an $F\in C_S$ such that
$$
F(\prod_{\{s,t\}\subset S}q_{st})^M=1.
$$
This is what we are going to do.
\subsubsection{}For convenience, denote by
$P:=P_2(S) $ the set of all 2-element
subsets of $S$; for $T=\{s,t\}\in P$ write
$q_T=q_{st}.$ Denote $\U=\Co^{P}$; for $s\in
U$, write
$$q^s:=\prod_{T\in P}
q_T^{s_T},
$$
 It is clear that for every $s\in \U$, $q^s:\Z^S\to
\Co$ is an analytic function on $U_S$.
\subsubsection{}
Denote by $d\Y$ the standard volume form on
$\Z$; set
$$
\Omega:=\prod_{s\in S}d\Y_s.
$$
Note that the product does not depend on the
order of multiples. Let $f\Omega\in
\S_{\Z^S}$. Write
$$
Z(f,s)=\int_{\Z^S}(fq^s\Omega).
$$
This integral converges if  $\Re\; s_T>0$
for every $T\in P$.
\subsubsection{} \begin{Claim}\label{mero}
For any $f$,  $Z$ extends to a meromorphic
function on $\U$. It can only have poles of
the first order along the divisors of the
form
 $$
D(R,n):=\{(\sum_{T\subset R}2s_T)+n=0\},
$$
where $R\subset S$ is a subset with at least
2 elements; $T$ is an arbitrary 2-element
subset of $R$; $n>(\#R-1)(N-1)$ is a
positive integer.
\end{Claim}
\pf Let $\FM$ be the  real Fulton-MacPherson
compactification of $U_S$ so that we have a
surjection $P:\FM\to \Z^S$. Denote
$V=P^{-1}U_S$. We know that $P$ identifies
$V$ and $U_S$. The complement $
\FM\backslash V$ can be represented as
$\FM\backslash V=\cup_{R\subset S} \FM_R$,
where $\#R>1$ and each $\FM_R$ is a smooth
subvariety of codimension 1;
$P(\FM_R)=\D_R$, where $\D_R$ is the
diagonal given by the equivalence relation
$e_R$ on $S$ in which $x\sim_{e_R} y$ and
$x\neq y$ iff $x,y\in R$.

Let $P'(S)$ be the set of non-empty subsets
of $S$. Let $K\subset P'(S)$. Then
$$
\FM_K:=\cap_{R\in K}\FM_R\neq \emptyset
$$
 if and only for every $R_1,R_2$
from $K$, either one of them is inside the
other, or they do not intersect. In this
case we call $K$ {\em forest}. Let
$$
\FM_K^o:=\FM_K\backslash \cup_{L\overset
K,L\neq K}\FM_L.
$$
For every point $x\in \FM_K^o$, there exists
a neighborhood $W$ of $x$ and a
non-degenerate system of functions $t_R$,
$R\in K$  (i.e. all $dt_R$ are linearly
independent at every point $y\in W$) such
that $\FM_R$ is given by the equation
$t_R=0$.

\begin{Claim} 1) We have
$$
P^{-1}\Omega=\prod_{R\in
K}t_R^{(N-1)(\#R-1)}\omega,
$$
where $\omega$ is nondegenerate at $x$.

2)
$$
 P^{-1}q_{st}=(\prod_{\{s,t\}\subseteq R}t_{R}^2)u_{st},
$$
where $u_{st}(x)\neq 0$.
\end{Claim}

Without loss of generality we can assume
that

1)both $\omega$ and all $u_{st}$ do not
vanish on $W$;

2) $\phi:=P^{-1}f$ is supported on $W$.

\subsubsection{}
We have
\begin{equation}\label{Mellin}
Z(s,f)=\int_{\Z^S}(\prod_{R\in
K}t_R^{2s_R+(\#R-1)(N-1)})F(s,\Y)\phi\Omega,
\end{equation}
where $F(s,\Y)$ is an integer function in
$s$ and $s_R=\sum_{T\subset R}s_T$.
Therefore, $Z(s,f)$ can only have poles of
at most first order along the divisors
$D(R,n)$, where $n> (\#R-1)(N-1)$.
\endpf
\subsubsection{}
Let $M\in \U$  be such that all $M_T$ are
integer. Choose  an arbitrary total order
$<_P$ on $P_2(S)$ and a point $\ve\in \U$
such that

1) each $\lambda_T$ is positive real number;

\noindent 2)$$ \sum_{T\in P} \lambda_T<1;
$$
3)for all $T$,
$$\lambda_T>\sum_{T'<T}\lambda_{T'}.
$$
Let $C\subset \Co$ be the unit cirle. Then
for all $z\in C^{P_2(S)}$, $Z(s,f)$ is
regular at $M+\lambda z$. Set
$$
U(M,\lambda,f)=\frac 1{(2\pi i)^{\#
P}}\int_{C^P} Z(M+\lambda\ve,f)\prod_{T\in
P} \frac {d\ve_T}{\ve_T},
$$
Note that the sign of this integral is well
defined.

It is clear that $U(M,\lambda,f)$ is
independent of $\lambda$; we set
$U(M,<_P,f):=U(M,\lambda,f)$.
\begin{Claim}
$f\mapsto U(M,<_P)$ is a distribution.
\end{Claim}
\pf Let $P:\FM\to \Z^S$, $x\in \FM$ and a
neighborhood $W$ of $x$ be as in the proof
 of Claim (\ref{mero}).

It is sufficient to check that $U(M,<_P)$ is
continuous when restricted to  a subspace
$\S_W$ of densities $f$ such that $P^{-1}f$
is supported in $W$.

Let $s_R':=2s_R+(\#R-1)(N-1)$. Let $L_R$ be
arbitrary positive integers. Then we can
modify (\ref{Mellin}) a follows:
$$
Z(s,f)=\prod_{R\in K} \frac{
(-1)^L}{(s_R'+1)(s_R'+2)\cdots (s_R'+L_R)}
\int_W(\prod_Rt_R^{s_R'+L_R}
\partial_{t_R}^{L_R}(\phi(\Y)F(s,\Y))\omega,
$$
where we assume that we have extended the
set of functions $t_R$ to a coordinate
system on $W$ and that $\omega$ is the
standard density in this coordinate system.

Pick $L_R$ to be large enough. Then it is
immediate that

$$
U(M,<_P,f)=\int_W(\sum_{S\in K}
\partial_t^S A_S)f\omega,
$$
where $A_S$ are smooth functions on $W$.
Therefore, $U(M,>_P)$ is a distribution.
\endpf
Let $S=S_1\disjoint S_2$ so that
$\Z^S=\Z^{S_1}\times \Z^{S_2}$. Let $f_i\in
\S_{\Z^{S^i}}$; define $f_1\boxtimes
f_2:\Z^S\to \Co$ by $f_1\boxtimes
f_2(Y_1\times Y_2)=f_1(Y_1)f_2(Y_2)$, where
$Y_i\in \Z^{S_i}$. Assume that $M_T\geq 0$
whenever $T=\{s_1,s_2\}$ with $s_1\in S_1$,
$s_2\in S_2$.
\begin{Claim} \label{fc}
We have
$$
U(M,<_P,f_1\boxtimes
f_2)=U(M|_{S_1},<_{(P|_{P_2(S_1)})},f_1)
U(M|_{S_2},<_{(P|_{P_2(S_2)})},f_2).
$$
\end{Claim}
\pf Clear
\begin{Claim}
1) $q^LU(M,<_P)=U(M+L,<_P)$;

2) $q^M(U(M,<_P))=1$;

3) $U(M,<_P)\in C_S$.
\end{Claim}
\pf 1) Clear;

 2) follows from  1);

3) Note that $\t.U(M,<_P)=0$, therefore
$U(M,<_P)$ is quasi-polynomial. The property
(\ref{factor}) follows by induction from
Claim \ref{fc}.
\endpf
Thus, we have shown that $\pi:C_S\to B_S$ is
surjective.

\subsection{Filtration on $C_S$}
Let $\Diag^n\subset \Z^S$ be the union of
generalized diagonals of codimension $n$.
Let $F^nC_S:=C_{S,\Diag^n}\subset C_S$ be
the submodule consisting of distributions
supported on $\Diag^n$.  We will study this
filtration.

\subsubsection{}
 let $\De\subset \Z^S$ be a diagonal.
Let $i_\De:\De\to \Z^S$ be the corresponding
inclusion. Let $\DD_\De$ be the algebra of
polynomial differential operators on $\De$.
Let $\omega_\De$ (resp. $\omega_{\Z^S}$) be
the bundle of top forms on $\De$ (resp. on
$\Z^S$). It is well known that
$$
\D_{\De\to \Z}:=
\omega_{\D_\De}\otimes_{\O_{\De}}
\D_{\Z^S}\otimes_{\O_{\Z^S}}\omega_{\Z^S}^{-1}
$$
is a right $\D_\De$ and a left
$\D_{\Z^S}$-module. Let $M$ be a left
$\D_\De$-module. Set
$$
i_{\De*}M:=M\otimes_{\D_De}\D_{\De\to \Z}.
$$
For example, let $\S'_{\Z^S,\De}\subset
\S'_{\Z^S}$ be the submodule of
distributions $F$ such that

1)  $F$ is supported on $\De$

2) there exists an $M=M(F)$ such that $Fg=0$
for any smooth function vanishing on $\De$
at order $\geq M(F)$. (note that locally on
$\De$ the condition 2) is always true). We
have a natural isomorphism
$$
I_{\De}:=I_{\De,\Z^S}:i_{\De*}\S'_{\De}\stackrel\sim\to
\S'_{\Z^S,\De}.
$$
\begin{Claim}
1. Let $\De_1\subset \De_2\subset \De_3$.
Consider the composition
$$
i_{\De_1\De_3*}\C_{\De_1}\cong
i_{\De_2\De_3*}i_{\De_1\De_2*}\C_{\De_1}
\stackrel{I_{\De_1\De_2}}\to
i_{\De_2\De_3*}\C_{\De_2}\stackrel{I_{\De_2\De_3}}{\to}
\C_{\De_3}.
$$
It is equal to $I_{\D_1\D_3}$.
\end{Claim}
\subsubsection{}
Let $\De$ be given by an equivalence
relation $e$ on $S$. We then have an
isomorphism $\De\cong \Z^{S/e}$. Denote
$C_\D:=C_{\Z^{S/e}}$.
\begin{Proposition}\label{nos}
1) $$I_\De(i_{\De*}C_\De)\subset
C_{\Z^S;\De};
$$
2)$$I_\De|_{C_\De}:i_{\De*}C_\De\to
C_{\Z^S;\De}
$$ is an isomorphism.
\end{Proposition}
\pf 1) It suffices to show that
$I_\De(i_{\De*}C_\De)\subset C_S$. Let
$A\subset i_{\De*}C_\De$ be the subspace of
all elements annihilated by multiplication
by any function vanishing on $\De$. Let
$f\in \S_{\Z^S}$, $a\in C_\D$ and $u\in
\omega_\D
\otimes_{\O_{\Z^S}}\omega_{\Z^S}^{-1}$. Then
$uf|_\De\in \S_\De$  and
$I_\De(au)(f)=a(uf|_\De)$.

Using this formula and a simple induction,
we see that
 $I_\De(A)\subset
C_S$. It is also well known that
$i_{\De*}C_\De$ is generated by $A$. This
completes the proof of 1).

2) We need Lemma
\begin{Lemma} Let $U\subset \Z$ be a non-empty open set and assume
that $F\in C_S$ vanishes on $U$. Then $F=0$.
\end{Lemma}
{\em Proof of Lemma} The statement is
obvious when $S$ is empty or  has 1 element.
Let us now use induction. Let
$S=S_1\disjoint S_2$. We know that for some
$M$
$$
\prod_{s_1\in S_1,s_2\in
S_2}q_{s_1s_2}^MF\in
T_{S_1S_2}(C_{S_1}\otimes C_{S_2}).
$$
 There exist non-empty open sets $A_i\in \Z^{S_i}$
 such that $A_1\times A_2\subset U$.
Write:$$ \prod_{s_1\in S_1,s_2\in
S_2}q_{s_1s_2}^MF\in T_{S_1S_2}(\sum_i
a_i\otimes b_i),
$$
where $a_i\in C_{S_1}$, $b_i\in C_{S_2}$ and
$a_i$ are linearly independent. By induction
assumption, restrictions of $a_i$ onto $A_1$
are also linearly independent (because if
these restrictions are dependent, then the
same dependence holds for the whole
$Y^{S^1}$.) Therefore, there exist $p_i\in
\S_{A_1}$ such that $a_i(p_j)=\delta_{ij}$.
Let $q\in \S_{A_2}$. We know that
$F(p_i\boxtimes q)=0$. Therefore,
$b_i(q)=0$, since $q$ is arbitrary, $b_i$
vanishes on $A_2$, hence by induction
assumption, $b_i=0$. Therefore,
 $\prod_{s_1\in S_1,s_2\in S_2}q_{s_1s_2}^MF=0$.
Therefore, $F$ is supported on
$\D_{S_1S_2}:=\cap_{s_i\in S_i}\D_{s_1s_2}$,
hence on
$$
E=\cap_{S=S_1\disjoint S_2}\D_{S_1S_2}.
$$
Show that $E$ is the smallest diagonal
$\De\alpha$. Indeed it is clear that
$\De_\alpha\subset E$. If $\Y\notin
\De_\alpha$, then there exists a partition
$S=S_1\disjoint S_2$ such that $S_1,S_2$ are
non-empty and $Y_{s_1}\neq Y_{s_2}$ whenever
$s_i\in S_i$. Therefore $\Y\notin
\De_{S_1S_2}$, hence not in $E$.

Thus, $F$ is supported on $\De_\alpha$.
Since $F$ is quasi-polynomial and vanishes
on $U$, it also vanishes on $U+a$ for all
$a\in \De_\alpha$. Therefore, $F$ vanishes
on a neighborhood of $\De_\alpha$.
Therefore, $F=0$.

{\em Proof of Proposition \ref{nos} 2)}

1.Choose a relatively compact open set $U\in
\Z^S$. Then it is well known that there
exists $M$ such that $F$ is annihilated by
multilpication by any function vanishing on
$\De\cap U$ of order $\geq M$. In virtue of
the Lemma, this implies that $F$ is actually
annihilated by multiplication by
 any function
vanishing on $\De$ of order $\geq M$. 2.It
suffices to check that there exists $f\in
I_\De(i_{\De*}C_\De)$ such that $F-f$
vanishes on $U$. it is easy to see that the
latter is equivalent to the following: for
any polynomial $P$ vanishing on $\De$ of
order $M-1$, $Pf\in I_\De(A)$. This follows
easily by induction.
\endpf
\subsection{} Let $\Dg_n:=\Dg_n(S)$ be the set  (not the union!)
of all diagonals
 in $\Z^S$ of codimension $n$.
We have a map
$$
I_n:=\oplus_{\D\in\Dg_n} I_\De\oplus i:
\oplus_{\D\in\Dg_n} i_{\D*}C_\D\oplus
C_{S,\Diag_{n+1}}
 \to C_{S,\Diag_{n}}
$$
\begin{Claim}
1) $I_n$ is surjective;

2) if $\sum_{\De\in\Dg_n} f_\De+g\in \Ker
I_n$, where $f_\De\in i_{\De*}\De$, and
$g\in C_{S,\De_{n+1}}$
 then all $f_\De$ are supported
on $\De_{n+1}$.
\end{Claim}
\pf We need Lemma.
\begin{Lemma} Let $\De\in \Dg(S)_n$
There exists $M:P_2(S)\to {\Bbb Z}_{\geq 0}$
such that $q_\De:=q^M=0$ on any $\De'\in
\Dg_n$, $\De'\neq \De$ but $q^M\neq 0$ on
$\De$.
\end{Lemma}
\pf Set $M_{st}=1$ if $\De_{st}\supset \De$;
otherwise set $M_{st}=0$.
\endpf
{\em Proof of Claim}

1) Let $F\in \C_{S,\De_n}$. Then
$F(q_\De)^m\in C_{S,\De}$ for $m>>0$. We
have an isomorphism $C_{S,\De}\cong
i_{\De*}C_{\De}$. We also have the map
$\pi:C_\De\to B_\De$ which induces a map
$$
i_*(\pi):i_{\De*}C_{\De}\to i_{\De*}B_\De.
$$
 In particular
$[i_{\De*}(\pi)]F(q_\De)^m\in
i_{\De*}B_\De$. Since the multiplication
onto $q_\De$ is invertible on
$i_{\De*}B_\De$, there is an element $x\in
i_{\De*}B_\De$ such that
$$
x(\q_\De)^m=(i_{\De*}(\pi)(F)(q_\De)^m.
$$
Since $\pi$ is surjective, so is
$i_{\De*}\pi$. Pick a pre-image $x':=x'_\De$
of $x$ in $i_{\De*}C_{\De}$. Then
$$([i_{\De*}(\pi)](F(q_\De)^M-x'(q_\De)^M)=0.
$$
It then follows that $F-x'$ is supported on
the union of all $n$-dimensional diagonals
except $\D$. Since each $x'_\D$ is supported
on $\De$, we have:
$$
F-\sum_{\De}x'_\De
$$
is supported on $\De_{n+1}$.

Proof of 2). Let $\sum f_\De +g\in\Ker I_n$.
 It follows that $(q_\De)^mf_\De$ is supported
on $\De_{n+1}$ it is easy to check that if
$x\in \De$ and $(q_\De)^m(x)=0$, then $x\in
\De_{n+1}$. Therefore, $f_\De$ is supported
on $\De_{n+1}$.
\endpf
\subsubsection{}
\begin{Corollary}\label{filtr}
The map $\oplus_{\De\in \Eq(S)_n} I_\D$
induces an isomorphism:
$$
\oplus_{\De\in \Dg_n}i_*B_\D\to
C_{S,\Dg_n}/C_{S,\Eq(S)_{n+1}}.
$$
\end{Corollary}
\subsubsection{} Let $X:=\Co^N$ be the complexification
of $\Z$ viewed as an algebraic variety over
$\Co$. Let $\D_{X^S}$ be the sheaf of
differential operators on $X^S$. Then $C_S$
defines a $\D_{X^S}$-module $\C_S$ in a
usual way. The above claim implies that
$\C_S$ is a holonomic $\D_{X^S}$-module
(because each quotient
$C_{S,\Dg(S)_n}/C_{S,\Dg(S)_{n+1}}$
determines a
 holonomic $\D_{X^S}$-module).

\subsubsection{} Let $\Dg(S)$ be the set of diagonals
in $X^S$ ordered with respect to the
inclusion. We denote by the same symbol the
corresponding category. We have a functor
$\D\mapsto \C_{S,\D}$ from $\Eq(S)$ to the
category of $\DD_{X^S}$-modules.
\subsubsection{} Let $I$ be a small category and
$\C$ an abelian $k$-linear category. Let
$F:I\to\C$ be a functor. Let $I'$ be the
abelian category of functors
$I^{\text{op}}\to \vect$. For $A\in I'$ we
can form the Eilenberg-MacLane tensor
product $F\otimes_{I} A\in \C$. We call $F$
{\em perfect} if the functor $A\mapsto
F\otimes_{I} A$ is exact.
\begin{Claim} The functor $\De\mapsto \C_{S,\De}$
is perfect.
\end{Claim}
\pf Let $n>0$ be an integer and
$\De\in\Dg(S)$. Set
$F_n(e):=(\C_{S,\De})_{n}$. We see that
$F_n:\Dg(S)\to \D_{X^S}$  are subfunctors of
our functor $ F=F_0$.

It suffices to show that for every $n$,
$G_n:=F_n/F_{n+1}$ is perfect.

It follows that
$$
G_n(s)=\bigoplus\limits_{t\in \Dg_n;t\leq s}
G_t,
$$
where $G_t=i_{t*}\B_{\De_t}$. The structure
maps are the  obvious ones.

We have
$$
G_n\otimes_{\Dg(S)}
A=\bigoplus\limits_{t\in\Dg(S)_n}
A(t)\otimes G_t
$$
and we see that the functor
$$
A\mapsto \G_n \otimes_{\Dg(S)}$$ is exact.
Therefore, $G_n$ is perfect.
\endpf
\subsubsection{}
Let $S_a,a\in A$ be a finite family of
 finite sets. Then we have
a $\prod_{a\in A}\Dg(S_a)$-filtration on
$\prod_{a\in A} \C_{S_a}$ viewed as a
$\D_{X^{\disjoint_{a\in A }S_a}}$-module.
The same agument shows that the
corresponding functor from the category
$\prod_{a\in A}\Dg(S_a)$ to the category of
$\DD_{X^{\disjoint_{a\in A }S_a}}$-modules
is perfect.

\subsubsection{}
We are going to study how the map
$I_{\D_1\D_2}:i_{\D_1\D_2*}\C_{\D_1}\to
\C_{\D_2}$ is compatible with the
filtrations. Te answer is very simple: this
map induces an isomorphism
$$
i_{\D_1\D_2*}\C_{\D_1}\to \C_{\D_2,\D_1}.
$$

The filtration on the l.h.s induced by the
filtration on $\C_{\D_1}$ coincides with the
filtration induced by the one on
$\C_{\D_2,\D_1}$.
\section{Asymptotic maps}\label{asympt}
\subsection{Construction}
Let $\D_e\subset Y^S$ be a diagonal given by
an equivalence relation $e$ on $S$. Let
$p:S\to S/e$ be the canonical projection.
Let $S_i:=p^{-1}i$, $i\in S/e$. Denote
$$
\C_S^e:=i_{\D_e\Z^S}\ha(\B_\D)\otimes\boxtimes_{i\in
S/e} \C_{S_i}.
$$
The multiplication by $q_{st}$
 is invertible on $\C_S^e$ whenever $p(s)\neq p(t)$.
 Let $Q_e$ be the product of all such
$q_{st}$.
\subsubsection{}
We are going to construct  a map
$$
\as_{S,e}:\C_S\to \C_{S}^e
$$
as follows.

First of all it suffices to define a
corresponding map on the level of global
sections. Let $F\in C_S$.
 It follows from the defintion
that there exists an $M$ such that
$$
Q^M_eF\in T((\otimes_{i\in S/e }C_{S_i}),
$$
where the tensor product is taken over
$\Co$. where $T$ is the natural inclusion
$$
(\otimes_{i\in S/e}C_{S_i})\to C_S,
$$
induced by the superposition of maps from
(\ref{tensor}). On the other hand we have an
obvious map
$$
(\otimes_{i\in S/e}C_{S_i})\to C_S^e.
$$
Since the multiplication by $Q$ is
invertible on $C_S^e$, we have a well
defined map $\as'_{S,e}:C_S\to C_S^e $,
which determines the desired map $\as_S$.
\section{Properties of $\as_S$}
\subsection{Compatibility with the filtrations}.
\subsubsection{Filtration on $C_{S}^e$}
Let $f\geq e$ be an equivalence relation. It
can be equivalently described as a set of
equivalence relations $f_i$ on $S_i$. Set
$$(\C_{S}^e)_{f}:=i_{\D_e\Z^S}^\ha(\B_S)\otimes\boxtimes_{i}
\C_{S_i,\De_{f_i}} \subset \C_{S}^e.
$$
Thus we have a filtration of $\C_{S}^e$
indexed by the ordered set $\Dg(S)^{\geq e}$
of all equivalence relations which are
greater or equal than $e$. It is clear that
this filtration is perfect (i.e. the
corresponding functor
$$\Dg(S)^{\geq e}\to \D_{X^S}\text{-mod}
$$
is perfect.) We can also consider $\C_{S}^e$
as an object perfectly filtered by $\Dg(S)$
such that $(C_{S}^e)_f=0$ if $f$ is not
greater or equal to $e$.

We have an isomorphism
$$\Gr_f \C_{S}^e=i_{\D_f\Z^S*}\{{\i_{\D_e\D_f}}^\ha(\B_{S/e})
\otimes\boxtimes_i\B_{S_i/f_i}\}
$$
if $f\geq 0$ (otherwise the corresponding
element is zero).
\subsubsection{}
The map $\as_{S,e}$ is compatible with the
filtrations. Let $f\geq e$. The induced map
from $\Gr_{f}\C_S\cong i_{\D_f\Z^S}\B_{S/f}$
to $\Gr_f\C_{S}^e$ is induced by the
asymptotic map
$$
\B_S\to i_{\D_e\D_f}^\ha(\B_{S/e})
\otimes\boxtimes_i \B_{S_i/f_i}.
$$

\def\follows{\Rightarrow}
\def\overset{\supseteq}
\def\suset{\supset}
\def\liminv{\text{liminv}}
\def \varprojlim{\liminv}
\def\into{{\to}}

\def\D{\Delta}
\def\DD{{\cal D}}
\def\DDD{{{\cal D}\text{-mod}}}
\def \cD{\DD}
\def\as{{\frak As}}
\def\B{{\cal B}}
\def\C{{\cal C}}
\def\disjoint{\sqcup}
\def \E{{\cal E}}
\def \F{{\cal F}}
\def \G{{\cal G}}
\def \H{{\cal H}}
\def\FL{{\cal FL}}
\def\Funct{\text{\bf Funct}}
\def\Gr{\text{Gr}}
\def\i{{\frak  i}}
\def\I{{\cal I}}

\def\N{{\text{\bf Zebra}}}
\def\zebra{\N}
\def\segments{{\text{\bf Segments}}}
\def \Mr{{\cal M}^{\text{res}}}
\def\M{{\cal M}}
\def\oo{{\bf o}}
\def\P{{\cal P}}
\def\R{{\cal R}}
\def\ve{{\frak v}}
\def\vect{\text{\bf Vect}}
\def \d{{\frak d}}
\def \ha{\wedge}
\def \Eq{\text{Eq}}
\def \O{{\cal O}}
\def \NN{{\zebra(f,e)}}
\def\PP{{\cal R}}
\def\RR{{{\frak R}^{\text{\bf o}}}}
\def\Flags{{\text{\bf Flags}}}
\def\p{{\frac p}}

\section{Formalism
$\I,\i,\delta$}\label{funI}

In this section we will define functors
$\I,\i,\delta$. The functors $\i$ are the
same as the ones used to define an OPE (see
(\ref{functi})).  The functors $\delta$ are
the functors of direct image in the theory
of $\D$-modules.

The functors $\I$ are built from $\C_S$.

These functors will be used to construct a
required
 resolution of
the system $<\i>$.

\subsection{Main defintions}
\subsubsection{}
Let $\D_e\subset \D_f$ be two diagonals in
$X^S$ determined by the equivalence
relations $e\leq f$. Let $p:S/f\to S/e$ be
the canonical projection. Let
$(S/f)_i:=p^{-1}(i)$, $i\in S/e$. Set
$$
\I_{\D_1\D_2},\i_{\D_1\D_2},\delta_{\D_1\D_2}:\DDD_{X^{\D_1}}\to\DDD_{X^{\D_2}}.
$$
to be defined  by the formulas:

$$
\I_{\D_2\D_1}(M)=
i_{\D_1\D_2}(M)^{\ha}\otimes\boxtimes_{i\in
S/e}{\C_{(S/f)_i}};
$$

$$
\i_{\D_2\D_1}(M)=
i_{\D_1\D_2}^\ha(M)\otimes\boxtimes_{i\in
S/e}{\B_{(S/f)_i}};
$$

$$
\delta_{\D_2\D_1}(M) =i_{\D_1\D_2*}(M).
$$

Sometimes we will also use the notation
$\I_i,\i_i,\delta_i$, where $i:\D_2\to \D_1$
is the inclusion of the corresponding
diagonals.
\subsubsection{Exactness}
Let $T\in S$ be a subset and $p_T:X^S\to
X^T$ be the corresponding projection. Call
an $H\in \DDD_{X^S}$ {\em $T$-exact} if $H$
is locally free as a
$p_T^{-1}\O_{X^T}$-module. Let $i:\D_e\into
X^S$ be a diagonal and let $T\subset S$ be
such that the through map $T\into S\to S/e$
is a bijection.

Let $M\in \DDD_{\D_e}$. Write
$$
i_H(M)=i^\ha(M)\otimes H
$$
\begin{Claim}

1.Let  the functor $H$ be $T$-exact. Then
the functor
$$
i_H(\cdot):\DDD_{\D_e}\to\DDD_{X^S}
$$
is exact;

2. Let
 $$
0\to \H_1\to\H_2\to \H_3\to 0
$$
be an exact sequence of $T$-exact modules.
Then the sequence
$$
0\to i_{H_1}(M)\to i_{H_2}(M)\to
i{H_3}(M)\to 0
$$

is exact for all $M\in\DDD_{\D_e}$.
\end{Claim}
\pf Clear.
\endpf

Note that $\B_S,\C_S,i_{\D*}\O_{\D}$ are
$\{s\}$-exact for any one elements subset
$s\subset S$ (here $i_\D:\D\subset X^S$ is
the smallest diagonal. This immediately
implies that the functors $\i,\I,\delta$ are
exact.

\subsubsection{Filtration}
Let $\D_e\subset \D_f\subset \D_g$. Let
$p:S/g\to S/f$ and $q:S/f\to S/e$. For $x\in
S/e$ let $g_x$ be the equivalence relation
on $(qp)^{-1}x$ induced by $g$. We have a
projection
$$
p_x:(qp)^{-1}x\to (qp)^{-1}x/g_x\cong
q^{-1}(x)
$$
induced by $p$. We have a map
$$
J_{gfe}:\delta_{\D_g\D_f}\I_{\D_f\D_e}\to
\I_{\D_g\D_e}
$$
defined as follows:

\begin{eqnarray*}
\delta_{\D_g\D_f}\I_{\D_f\D_e}(M)\cong
i_{\D_f\D_g*}(i_{\D_e\D_f}^\ha(M)\otimes\boxtimes_{i\in
S/e}
C_{q^{-1}(x)})\\
\cong i_{\D_e\D_g}^\ha(M)\otimes
(i_{\D_f\D_g*}(\boxtimes_{i\in S/e}
C_{q^{-1}(x)}))\\
\cong
i_{\D_e\D_g}^\ha(M)\otimes\boxtimes_{x\in
S_e}
i_{\D_{g_x}X^{(qp)^{-1}x}*}C_{\D_{g_x}}\\
\to
i_{\D_e\D_g}^\ha(M)\otimes\boxtimes_{x\in
S_e}C_{(qp)^{-1}x}.
\end{eqnarray*}

This map is injective for all $M$. Indeed,
this needs to be checked only for the last
arrow, which follows from:

(1) injectivity of the arrow
$$i_{\D_{g_x}X^{(qp)^{-1}x}*}C_{\D_{g_x}}\\
\to \boxtimes_{x\in S_e}C_{(qp)^{-1}x};
$$

(2) both terms of this arrow are $\O_T$
exact where $T\subset S/e$ is such that
$T\to S/e\to S/g$ is a bijection, as well as
the cokernel of this arrow.

The above results imply that these
inclusions, for all $f$ such that $g\geq
f\geq e$, define a perfect filtration on
$\I_{\D_g\D_e}$. We denote by
$(\I_{\D_g\D_e})_f$ the correponding term of
this filtration.
\subsubsection{This filtration is perfect}
Let $p:S/g\to S/e$ be the projection. For
$i\in S/e$, let $(S/g)_i=p^{-1}i$. The
ordered  set of equivelence relations $f$
such that $g\geq f\geq e$ is isomorphic to
the product $\prod_{x\in  S/e}\Dg((S/g)_x)$.
Denote this ordered set by $[e,g]$. The
filtration on $\I_{\D_g\D_e}$ is induced by
the perfect $\prod_{x\in
S/e}\Dg((S/g)_x)=[e,g]$- filtrations on $
\boxtimes_{x\in S/e}\C_{(S/g)_x}. $ Denote
by
$$
F:[e,g]\to\DDD_{X^{S/g}}
$$
the functor determined by these filtrations:

$$
F(\{f_x\}_{x\in S/e})=\boxtimes_{x\in S/e}
 \C_{(S/g)_x,f_x}.
$$

The statement we are proving follows
immediately from the following one:

Let $T\subset S/g$ be a subset such that the
through map $T\to S/g\to S/e$ is a bijection
and
 $$
A:([e,g])^{\text{op}}\to \vect.
$$
Then $F\otimes_ {[e,g]}A$ is $T$-exact.

Let us prove this statement. Indeed, we have
seen that $F$ has a filtration $F_n$ such
that each $F_n/F_{n+1}$ is perfect.
 Furthermore,
as it follows from \ref{filtr},
$$
F_n/F_{n+1}\otimes_{[e,g]} A=
\bigoplus_{t\in (\prod_{x\in
S/e}\Eq((S/g)_x)))_n} \G_t\otimes A(t),
$$
where each $\G_t$ is $T$-free.

Therefore, since each $F_n/F_{n+1}$ is
perfect, the filtration on $F$ induces a
filtration on $F\otimes_{[e,g]} A$, its
associated graded quotient being isomorphic
to $F_n/F_{n+1}\otimes A$, which are, as we
have seen $T$-free. Therefore, $F\otimes A$
is also $T$-free.
\endpf
\subsubsection{}
Thus, $\I_{\D_g\D_e}$ is a perfect functor
from the category $\D_{\D_e}$-modules to the
category of $[e,g]$ -perfectly filtered
 $\D_{\D_g}$-modules.
We have a canonical isomorphism
$$
\Gr_f(\I_{\D_g\D_e}(M))\cong
\delta_{\D_g\D_f}\i_{\D_f\D_e}(M)$$

\subsubsection{Asymptotic decompositions}
The asymptotic maps $\as_{S,e}$ from
(\ref{asympt}) define maps:
$$
\as_{fge}:\I_{\D_g\D_e}\to
\I_{\D_g\D_f}\i_{\D_f\D_e}
$$
in the obvious way.

The compatibility of $\as_{S,e}$ with the
filtration implies that the map $\as_{fge}$
is compatible with the filtrations in the
following sence:
$$
\as_{fge}(\I_{\D_g\D_e})_{f'}=0
$$
if $f'\notin [f,g]$. Otherwise
\begin{equation}\label{b}
\as_{fge}(\I_{\D_g\D_e})_{f'}\subset
(\I_{\D_g\D_f})_{f'}\i_{\D_f\D_e}.
\end{equation}

Furthermore, we have
$$
(\I_{\D_g\D_e})_{f'}\cong
\delta_{\D_g\D_{f'}}\I_{\D_{f'}\D_e}
$$
and the above inclusion (\ref{b}) is given
by the map
\begin{eqnarray*}
\delta_{\D_g\D_{f'}}\I_{\D_{f'}\D_e}\to
\delta_{\D_g\D_{f'}}\I_{\D_{f'}\D_f}\i_{\D_{f}\D_e}\\
\cong (\I_{\D_g\D_f})_{f'}\i_{\D_{f}\D_e}.
\end{eqnarray*}

compute the associated graded map
\begin{equation}\label{asgr}
\Gr_{f'} \I_{\D_g\D_e}\to (\Gr_{f'}
\I_{\D_g\D_{f}})\i_{\D_f\D_e}\end{equation}

We have $$ \Gr_{f'} \I_{\D_g\D_e}\cong
\delta_{\De_g\De_{f'}} \i_{\De_{f'}\De_{e}};
$$

therefore the map (\ref{asgr}) is given by
the map
$$
\delta_{\De_{g}\De_{f'}}\i_{\De_{f'}\De_{e}}
\to
\delta_{\De_{g}\De_{f'}}\i_{\De_{f'}\De_{f}}
\i_{\De_{f}\De_{e}}
$$
induced by the asymptotic map
$$
\i_{f'e}\to \i_{f'f}\i_{fe}.
$$
\subsubsection{} Let $S_\alpha,\alpha\in A$ be a finite
family of finite sets. Let
 $e_\alpha,f_\alpha\in\Dg(S_\alpha)$;
 $e_\alpha\leq f_\alpha$; let
$M_\alpha$ be  some
$\DD_{\D_{e_\alpha}}$-modules. Let
$S=\disjoint_{\alpha\in A}S_\alpha$;
$e:=\disjoint_{\alpha\in A} e_\alpha$;
$f=\disjoint_{\alpha\in A} f_\alpha$.

We then  have a natural map
$$
\boxtimes_{\alpha \in A}
\I_{\D_{f_\alpha}\D_{e_\alpha}}(M_\alpha)
\to \I_{\D_{f}\D_e}(\boxtimes_{\alpha\in A}
M_\alpha).
$$

\section{Resolution}\label{endres}
We will focus our study on  the functors
$\i$ and $\I$. We will need the following
properties.
 Let $\D_1\subset \D_2\subset \D_3\subset\D_4$ be
a flag of diagonals.

1. We have natural transformations:
$$
\I_{\D_i\D_j}\to \i_{\D_i\D_j},\quad i>j;
$$
$$
\i_{\D_3\D_1}\to \i_{\D_3\D_2}\i_{\D_2\D_1};
$$
$$
\I_{\D_3\D_1}\to \I_{\D_3\D_2}\i_{\D_2\D_1};
$$

$$
\boxtimes_{\alpha \in A}
\I_{\D_{f_\alpha}\D_{e_\alpha}}(M_\alpha)
\to \I_{\D_{f}\D_e}(\boxtimes_{\alpha\in A}
M_\alpha).
$$

$$
\boxtimes_{\alpha \in A}
\i_{\D_{f_\alpha}\D_{e_\alpha}}(M_\alpha)
\to \i_{\D_{f}\D_e}(\boxtimes_{\alpha\in A}
M_\alpha).
$$
2. The properties are:

a) The functors $\i$ satisfy the axioms of
system (see (\ref{systems}).

b) The following diagrams commute:

\begin{equation}
\xymatrix{\I_{\De_4\De_1}\ar[r]\ar[d]&
 \I_{\De_4\De_2}\i_{\De_2\De_1}\ar[r]&
 \I_{\De_4\De_3}\i_{\De_3\De_2}\i_{\De_2\De_1}\\
\I_{\De_4\De_3}\i_{\De_3\De_1}\ar[rru]&&};
\end{equation}

\begin{equation}
\xymatrix{\boxtimes_{\alpha\in A}
\I_{\De_{3\alpha}\De_{1\alpha}}(M_\alpha)\ar[r]\ar[d]&
\boxtimes_{\alpha\in A}
\I_{\De_{3\alpha}\De_{2\alpha}}
\i_{\De_{2\alpha}
\De_{1\alpha}}(M_\alpha)\ar[r]&
\I_{\De_3\De_2}\boxtimes_{\alpha\in
A}\i_{\De_{2\alpha}
\De_{1\alpha}}(M_\alpha)\ar[r]&
\I_{\De_3\De_2}\i_{\De_2\De_1}M\\
\I_{\De_3\De_1}(M)\ar[rrru]&&&}
\end{equation}
\subsubsection{}Let $f>e$ be equivalence relations
on $S$ Let $\N(f,e)$ be the ordered set
defined as follows. Elements of $\N(S)$ are
sequences
$s:=(e_1i_{12}e_2i_{23}e_3i_{34}\cdots
i_{n-1n}e_n)$, where
$f=e_1>e_2>\cdots>e_n=e$ is a flag of
equivalence relations and each $i_{pp+1}$ is
one of the symbols $\i$ or $\I$. Let
$s'=(e'_1i'_{12}\cdots e'_{n'})$ be another
element of $\N(f,e)$. We write $e\geq e'$
if:

1) for all $k=1,2,\ldots,n'$, there exists
$n_k$ such that $e'_k=e_{n_k}$ (in
particular, $n_1=1; n_{n'}=n$;

2) if $i'_{kk+1}=\i$, then $i_{pp+1}=\i$ for
all $p=n_k,n_k+1,\ldots, n_{k+1}-1$;

3)if $i'_{kk+1}=\I$, then $i_{pp+1}=\i$ for
all $p=n_k+1,\ldots,n_{k+1}-1$ (it is
possible that $i_{n_kn_k+1}=\I$).

Let $\j^i_{\D_1\D_2}=\i_{\D_1\D_2}$ if
$i=\i$ and $\j^i_{\D_1\D_2}=\I_{\D_1\D_2}$
if $i=\I$. For $s\in\N(f,e)$ Write
$$
\j(s):=\j^{i_{12}}_{\D_1\D_2}\j^{i_{23}}_{\D_2\D_3}\cdots
\j^{i_{n-1n}}_{\D_{n-1}\D_n}.
$$

The above properties imply that $\j$ is a
functor from the category determined by the
ordered set $\N(f,e)$ to the category of
functors
$\DD_{\D_e}\text{-mod}\to\DD_{\D_f}\text{-mod}$;
our agreement is that whenever  $x'\leq x$,
$x',x\in \N(f,e)$, we have an arrow from
$\j(x')\to \j(x)$.
\subsection{Filtration on  the functor $\j:\N\to
\Funct(\DD_{\D_e}\text{-mod},\DD_{\D_f}\text{-mod})$}\label{Filt}
 To define such a filtration  we need
some combinatorics.
\subsubsection{}
Define the  ordered set $\segments(f,e)$. To
this end, we need a notion of {\em segment}
in an arbitrary ordered set $X$, which is
just an arbitrary pair of elements $x,y\in
X$ such that $x>y$. We denote such a segment
by $[x,y]$. Given two segments $[a,b]$ and
$[c,d]$, we say that $[a,b]>[c,d]$ iff
$b\geq c$ (in which case $a> b\geq c>d$.
Define the set $\segments(X)$ whose elements
are arbitrary flags of segments
$$
[a_0,b_0]>[a_1,b_1]>\cdots >[a_n,b_n],
$$
 Of
course, this simply means that $$
a_0>b_0\geq a_1>b_1\geq
a_2>b_2\geq\cdots\geq a_n>b_n.
$$

Introduce an order on the set $\segments(X)$
according to the following rule.

Let $$ u=([a_0,b_0]>[a_1,b_1]>\cdots
>[a_n,b_n]),
$$
and
$$
v=([a'_0,b'_0]>[a'_1,b'_1]>\cdots
>[a'_m,b'_m]),
$$
be elements in $\segments(X)$. We say that
$u\leq v$ iff for every segment
$[a'_i,b'_i]$  there exists a segment
$[a_j,b_j]$ such that $a_j=a'_i>b'_i\geq
b_j$.

Let $f\geq e $ be equivelence relations on
$S$. Let $\Dg(f,e)$ be the set of all
equivalence relations $g$ such that $f\geq
g\geq e$.

Set
$$
\segments(e,f):=\segments(\Dg(f,e)).
$$

For $s\in \zebra(f,e)$, where
$s=e_1i_{12}\cdots e_n$, we will  define an
element $\nu(s)\in \segments(f,e)$ by
setting
$$
\nu(s)=([e_{k_1},e_{k_1+1}]>[e_{k_2},e_{k_2+1}]>\cdots
>[e_{k_r},e_{k_{r}+1}]),
$$
where $k_1<k_2<\cdots<k_{r}$ is a sequence
of all numbers such that $i_{k,k+1}=\I$.

\subsubsection{} Let $s\in \zebra(f,e)$,
$$s=e_1i_{12}\cdots e_n,
$$
and let $t\in\segments(f,e)$ be an element
such that  $t\geq \nu(s)$. Let
$$
t=([a_1,b_1]>[a_2,b_2]>\cdots >[a_k,b_k]).
$$

Assume that $i_{p,p+1}=\I$. Then there are
two possibilities:

1) either there exists $p'$ such that
$e_p=a_{p'}$, $e_p=a_{p'}<b_{p'}\leq
e_{p+1}$. In this case write
$$
j'_p=\delta_{\D_{a_{p'}b_{p'}}}
\I_{b_{p'}e_{p+1}};
$$

2) there are no segments $[a_{p'},b_{p'}]$
as in 1).

We then set
$$
j'_p=\j^{i_{pp+1}}_{\D_{e_p}\D_{e_{p+1}}}.
$$

 Define:
$$
F^t\j(s)=j'_1j'_2\cdots j'_p.
$$

if it is not true that $t\geq \nu(s)$, we
then  set $F^t\j(s)=0$.
\begin{Claim}
For every $s$, $F$ is a perfect filtration
on $\j(s)$.
\end{Claim}
\pf
 Let $\segments(f,e)_s=\{t\in \segments[f,e]| t\geq s\}.$
 Let $u\in \segments(f,e)$
 We see that $\j(s)_u=0$ whenever $u\notin N_s$.
Therefore,
$$
\j(s)\otimes_{\segments(f,e)}A\cong
 \j(s)\otimes_{\segments(f,e)_s} A
$$
for every $A:\segments(f,e)\to\vect$. Thus,
it suffices to show that the
$\segments(f,e)_s$-filtration on $\j(s)$ is
perfect. Let
$$
s=([a_1,b_1]>[a_2,b_2]>\cdots>[a_nb_n]).
$$
be an element in $\segments(f,e)$. Then we
have an isomorphism
$$\segments(f,e)_s\cong \prod_{i}[a_i,b_i]_{\Dg(S)}.
$$
Indeed, let $a_i\geq u_i\geq b_i$. Let
$i_1>i_2>\cdots >i_r$ be the subsequence of
all numbers such that $a_i>u_{i_k}$. Then
the corresponding flag of segments is given
by the formula
$$
[a_{i_1},u_{i_1}]>[a_{i_2},u_{i_2}]>\cdots>
[a_{i_r},u_{i_r}].
$$
Consider two cases.

Case 1 $a_1=f$. Define an element $s'\in
\segments(b_1,e)$ by the formula
$$
s'=([a_2,b_2]>\cdots>[a_nb_n]).
$$
We then have $$
\j(s)=\I_{\D_{a_1}\D_{b_1}}\j(s_1)
$$
We have $$\segments(f,e)_s\cong
[a_1,b_1]_{\Dg(S)} \times
\segments(b_1,e)_{s'},
$$
where a pair $(u,r)$, where $a_1\geq u\geq
b_1$ and $r\in \segments(b_1,e)_{s'}$,
$$
r=[a'_2,b'_2]>[a'_3,b'_3]>\cdots>[a'_l,b'_l]
$$
determine the flag of segments $f$, where
$$
f=([a_1,u]>[a'_2,b'_2]>\cdots>[a'_n,b'_n])
$$
if $a_1>u$

and

$$
f=([a'_2,b'_2]>\cdots>[a'_n,b'_n])
$$
if $u=b_1$.

the filtration on $\j(s)$ is induced by the
corresponding filtrations on
$\I_{\D_f\D_{b_1}}$ and $\j(s')$.

We are going to use induction, so we can
assume that we have already proven that
  the filtration on
$\j(s')$ is  perfect. We denote by the same
letter  the functors determined by the
corresponding filtrations on $\j(s),\j(s')$
and $\I_{\D_f\D_{b_1}}$.

Denote $j_n:=\I_{\D_f\D_{b_1}n}\j(s')$.
Since the quotient
$$
\I_{\D_f\D_{b_1}n}/\I_{\D_f\D_{b_1}(n+1)}
$$
is $T$-free for every finite set $T\subset
S/f$ such that $T\to S/f\to S/b_1$ is
bijection, we have:
$$
j_n/j_{n+1}\cong
(\I_{\D_f\D_{b_1}n}/\I_{\D_f\D_{b_1}(n+1)})\j(s').
$$

Using (\ref{filtr}), we obtain
$$
j_n/j_{n+1}\otimes_{\segments(f,e)_s} A\cong
\oplus_{t\in [b_1f]_n} (G_t j'(s'))
\otimes_{\segments(b_1,f)_{s'}}A(t,s')$$
therefore, $j_n/j_{n+1}$ is perfect, hence
$j(s)$ is also perfect.

Case 2 $f>a_1$ is similar.
\endpf
\subsection{Desctription of $\Gr^t\j$.}
\subsubsection{}
Let $t\in \segments(f,e)$; let
$$
\N(f,e)_t=\{g\in\N(f,e)| \nu(g)=t\}.
$$
We consider $\N(f,e)_t$ as an ordered subset
of $\N(f,e)$. Let $i:\N_t\to\N(f,e)$ be the
inclusion.

 Let $\Funct(\N(f,e)_t,\C)$ be the category of functors from
$\N_t$ to an arbitrary abelian category
$\C$. We have the restriction functor
$$
i^{-1}:\Funct(\N(f,e),\C)\to\Funct(\N(f,e)_t,\C).
$$
Let $i_*$ be the right adjoint functor. It
can be constructed as follows.

Let $F:\N(f,e)_t\to \C$ and $s\in \N(f,e)$.
There are two cases:

1) It is false that $\nu(s)\leq t$, then
$i_*F(s)=0$;

2) $\nu(s)\leq t$. Then there exists the
least element $s_t\in \N(f,e)_t$ among the
elements in $\N(f,e)_t$ which are $\geq s$
(we will show it in the next paragraph). Set
$i_*F(s)=F(s_t)$. It is clear that if
$\nu(s_1),\nu(s_2)\leq t$ and $s_1\leq s_2$,
then $(s_1)_t\leq (s_2)_t$. This determines
the functor structure on $i_*F$.

We will now construct the  element $s_t$.
Let $$ s=(e=e_1i_{12}e_2i_{23}\cdots
i_{n-1n}e_n=f).
$$

Let $$
t=([a_1,b_1]>[a_2,b_2]>\cdots>[a_m,b_m]).
$$

The condition $\nu(s)\leq t$ means that for
every $\mu=1,2,\ldots, m$ there exists a
number $k_i$ such that
$e_{k_\mu}=a_\mu>b_\mu\geq e_{k_{\mu+1}}$
and $i_{k_(k_\mu+1)}=\I$.

Let $e=u_1>u_2>\cdots>u_N=f$ be a flag of
equivalence relations determined by the
condition
$$
\{u_1,u_2,\ldots,u_N\}=\{e_1,e_2,\ldots,e_n\}\cup
\{a_1,b_1,a_2,b_2,\ldots,a_m,b_m\}.
$$

Define the symbols $i'_{k,k+1}$, where
$k=1,2,\ldots, N-1$, according to the rule:

if $u_k=e_m$ and $u_{k+1}=e_{m+1}$, then
$i'_{k,k+1}=\i$;

if $u_k=e_m=a_r$ and $u_{k+1}=b_r$, then
$i'_{k,k+1}=\I$;

if $u_k=b_r$ and $u_{k+1}=e_l$, then
$i'_{k,k+1}=\i.$
 As we exhausted all the possibilities,  we can now define
 $$
 s_t=(u_1i'_{12}u_2i'_{23}\cdots u_N).
 $$

\subsubsection{}
We have $\Gr^t\j\cong i_*i^{-1}\Gr^t(\j)$
and it remains to describe
$G:=i^{-1}\Gr^t(\j)$. Let $s$ be such that
$\nu(s)=t$; $s=(e_1i_{12}\ldots e_n)$.

Set $c^{u}=i$ if $u=\i$; $c^u=\delta$ if
$u=\I$. Then
$$
G(s)=c^{i_{12}}_{\D_{e_1}\D_{e_2}}\cdots
c^{i_{n-1n}}_ {\D_{e_{n-1}}\D_{e_n}}.
$$
\subsubsection{The functors $\P_{\De_f\De_e}$.}
\label{Filt1} We will study the functor
$$
\P_{\De_f\De_e}:=\liminv_{s\in
\N(f,e)}\j(s).
$$

Set $$ \P_{\De_f,\De_e,t }:=\liminv_{s\in
\N(f,e)}\j(s)_t,
$$
where $t\in \segments(f,e)$. Our goal is to
show that the functor $t\mapsto
\P_{\De_f,\De_e,t}$ is

1) a filtration on $\P_{\De_f,\De_t}$;

2) a perfect functor on the category
$\segments(f,e)$.

Since these properties are the case for the
functor $t\mapsto \j_t$; it suffices to show
that the derived functors
$R^i\liminv_{\N(f,e)} $, $i\geq 1$, vanish
on $\Gr_t\j$. This is what we are going to
do.

\subsubsection{}
Let $I$ be a small category and
 $H:I\to \C$ be a functor,where $\C$ is an arbitrary
$k$-linear category. Let $I^-$ be the
abelian category of functors $I\to \vect$
Let $h_H(X):=\hom_{I}(X,H)\in \C$, where
$X\in I^-$.

$H$ is called {\em flabby} if the functor
$h_H$ is exact. It is clear that flabby
functors are adjusted
 to the functor $\varprojlim_{I}$ and that there are
enough flabby objects in the abelian
category of functors $I\to \C$. The functor
$i_*$ is exact and maps flabby functors to
flabby (this follows from the existence of
an exact left adjoint functor $i^{-1}$,
therefore
$$
h_{i_*H}(X)=\hom_{\zebra(f,e)}(X,i_*H)\cong
\hom_{ \zebra(f,e)_t}(i^{-1}X,H),
$$
which implies that $i_*H$ is flabby).

 Therefore,
$$R\liminv_{\N(f,e)}
\Gr^t(j)\cong R\liminv_{\N(f,e)_t}G.
$$
The category $\N(f,e)_t$ has an initial
object $t_i$, which is
$$
t_i=(f\i a_1\I b_1\i a_2\I\b_2\i\cdots b_n\i
f),
$$
where we assume that in the case $e=a_1$, or
$b_i=a_{i+1}$, or $b_n=f$, the fragment $f\i
a_1$ (resp. $b_i\i a_{i+1}$, resp. $b_n\i
f)$ is repaced with $f$ (resp. $b_i$, resp.
$f$).

 Therefore,
$$
R\liminv_{\N(f,e)}\Gr^t(j)=G(t_i).
$$
\subsubsection{Conclusion}\label{filtrP}
 As was mentioned above,
these facts imply that we have a
filtration on $\P_{\D_f\D_e}$ by subfunctors
$\P_{\D_f\D_e,t}$ and that this filtration
is perfect. We will also denote
$F^t\P_{\D_f\D_e}:=\P_{\D_f\D_e}$.

\subsubsection{Lemma}\label{lemmaforbp} We will prove a Lemma which will only
be used in the next section. We have an
element $e\I f\in \zebra(f,e)$. Let
$$\zebra^0(f,e):=\zebra(f,e)\backslash\{e\I
f\}.$$ Let
$\P^0_{fe}:=\liminv_{\zebra^0(f,e)}\j$. We
have natural maps
 \begin{equation}\label{sek}
0\to \delta_{fe}\to \P_{fe}\to \P^0_{fe}\to
0.
\end{equation}

\begin{Lemma} The sequence (\ref{sek}) is exact.
\end{Lemma}
\pf It is easy to check that the composition
of the arrows is zero. Let us now prove the
exactness.  Let $t\in \segments(f,e)$. Set
$$
\P^0(f,e)_t:=\liminv_{s\in
\zebra^0(f,e)}\j(s)_t.
$$

The same argument as above shows that:

1) $t\mapsto \P^0(f,e)_t$ is a filtration on
$\P^0(f,e)_t$;

2)  the lowest element of the filtration is
zero: $\P^0(f,e)_{[f,e]}=0$.

3) the induced map $$\Gr^t\P(f,e)\to
\Gr^t\P^0(f,e)$$ is an isomorphism for all
$t\neq [f,e]$. If $t=0$, then the induced
map is a surjection (onto zero).

The lemma then follows easily.
\endpf
\subsection{Formalism $\delta,\P$}
We are going to describe a structure
possessed by  the functors $\delta,\P$.  Let
us first
 introduce the elements of this structure
 and then describe
 their properties.
 \subsubsection{Decompositions}
Define a map
$$
\alpha: \P_{\D_1\D_3}\to
\P_{\D_1\D_2}\P_{\D_2\D_3}
$$
as follows. Let $s_1\in \N(\D_1,\D_2)$ and
$s_2\in \N(\D_2,\D_3)$. Let $(s_1s_2)\in
\N(\D_1,\D_2)$ be the obvious concatenation.
 Set
$$
(p_{s_1}\times p_{s_2})\alpha=p_{(s_1s-2)}.
$$
It is immediate that this defintion is
correct.

 \subsubsection{Concatenations}
Let $s\in \N(f,e)$; $s=(e_1i_{12}\ldots
e_n)$.
 Define a map
$$
c:\P_{\D_1\D_2}\delta_{\D_2\D_3}\P_{\D_3\D_4}\to\P_{\D_1\D_4}
$$
by setting
$$
p_sc=0
$$
if the following is wrong:

There exists an $m$ such that $i_{mm+1}=I$
and
 $$
\D_{e_m}=\D_2\overset\D_3\overset\D_{e_{m+1}}.
$$

Let $\NN_{\D_2\D_3}\subset\N(e,f)$ be the
set of $s$ for which this condition is true.
 For $s\in \NN_{\D_2\D_3}$ we define the elements
$s_1\in \N(\D_1\D_2)$ and $s_2\in
\N(\D_3,\D_4)$ (we do not distinguish
between a diagonal in $X^S$ and an
equivalence relation  on $S$ by which it is
determined) according to the rule:

$$
s_1=e_1i_{12}\ldots e_m
$$
and
$$
s_2=\D_3Ie_{m+1}\ldots e_n.
$$

We then have a composition:
$$
c_s:\P_{\D_1\D_2}\delta_{\D_2\D_3}\P_{\D_3\D_4}\to
j_{s_1}\delta_{\D_2\D_3}\j_{s_2}\to \j_s.
$$
Define $c$ by the condition $p_sc=c_s$. Show
that this definition is correct.

Let $t>s$. Let $a_{ts}:\j_{s}\to\j_t$ be the
induced map. We need to check that
$c_t=a_{ts}c_s.$ There are several cases.

Case 1. $s\notin \NN_{\D_2\D_3}$. Since
$t>s$, $t\notin \NN_{\D_2\D_3}$;the
correctness is obvious;

Case 2. $t\notin \NN_{\D_2\D_3}$; $s\in
\NN_{\D_2\D_3}$. This means that $t$
contains an element $\rho$ such that
$$
\D_2=e_m>\rho>e_{m+1},
$$
but it is not true that
$$
\D_3\geq  \rho.
$$
In this case, the composition
$$
\delta_{\D_2\D_3}\I_{\D_3e_{m+1}}\to
\I_{\D_2e_{m+1}}\to \I_{\D_2\rho}\i_{\rho
e_{m+1}}
$$
is zero, therefore $a_{ts}c_s=0$, and the
correctness condition is satisfied.

Case 3: $s,t\in \NN_{\D_2\D_3}$ ---
straightforward.
\subsubsection{Factorization maps}
Let $S_a$,$a\in A$ be a finite family of
finite sets. Let $f_a\geq e_a$ be
equivalence relations on ${S_a}$. Let
$S=\disjoint_a S_a$; $f=\disjoint_a f_a$;
$e=\disjoint_a e_a$; $f\geq e$. Let $M_a\in
\DD_{\D_{e_a}}$ Define a natural
transformation
$$
\mu:\boxtimes_{a}\P_{f_ae_a}(M_a)\to
\P_{fe}(\boxtimes_a M_a).
$$

Let $g$ be  such that $f\geq G\geq e$. Any
such an equivalence relation can be
represented as $g=\disjoint_a G_a$, where
$f_a\geq g_a\geq e_a$.

Let
 $$
\Phi=(g_1i_{12}g_2\ldots i_{n-1n}g_n);
$$
$\Phi\in \N(f,e)$. Let $g_r=\disjoint_a
g_{ra}$. Let
 $$
\Phi'_a=(g_{1a}i_{12}g_{2a}\ldots
i_{n-1n}g_{na}).
$$
After deletion of  repeating terms we get an
element $\Phi_a\in \N(f_a,e_a)$. Define
\begin{eqnarray*}
p_\Phi\mu:\boxtimes_{a}\P_{f_ae_a}(M_a))\to
\boxtimes_{a}p_{\Phi_a}\P_{f_ae_a}(M_a)\\
\cong \boxtimes_a(
\i^{i_{12}}_{g_{1a}g_{2a}}\i^{i_{23}}_{g_{2a}g_{3a}}
\ldots\i^{i_{n-1n}}_{g_{n-1a}g_{na}}(M_a))\\
\cong
(\boxtimes_a\i^{i_{12}}_{g_{1a}g_{2a}})
(\boxtimes_a\i^{i_{23}}_{g_{2a}g_{3a}})\ldots
(\boxtimes_a\i^{i_{n-1n}}_{g_{n-1a}g_{na}})
(\boxtimes_a M_a)\\
\to
\i^{i_{12}}_{g_{1}g_{2}}\i^{i_{23}}_{g_{2}g_{3}}
\ldots\i^{i_{n-1n}}_{g_{n-1}g_{n}})(\boxtimes_a
M_a).
\end{eqnarray*}
This defines the map $\mu$. This completes
the description of the elements. Now let us
pass to the properties.
\subsubsection{Concatenation+factorization}
\def\a{{\alpha}}

It follows that the map
$$
\boxtimes_a
(\P_{f_ae_a}(M_a))\to\P_{fe}(\boxtimes_a
M_a) \to \P_{eg}\P_{gf}(\boxtimes_a M_a)
$$
is equal to the map
\begin{eqnarray*}
\boxtimes_a (\P_{f_ae_a}(M_a))\to
\boxtimes_a (\P_{f_ag_a}\P_{g_ae_a}(M_a))\\
\to \P_{fg}\boxtimes_a \P_{g_ae_a}(M_a)\to
\P_{fg}\P_{ge}(\boxtimes_a M_a).
\end{eqnarray*}
\subsubsection{}
\def \F{{\cal F}}
The map
\begin{eqnarray}\label{sochet}
\F:\P_{fg'}\delta_{g'g''}\P_{g''e}(M)\boxtimes
(\P_{f_1e_1}(M_1))\\ \nonumber \to
\P_{fe}(M)\boxtimes (\P_{f_1e_1}(M_1))\\
\nonumber \to \P_{f\disjoint f_1,e\disjoint
e_1}(M\boxtimes M_1)
\end{eqnarray}
is equal to the sum of the maps $\F_{g_1}$,
where $f_1\geq g_1\geq e_1$:
\begin{eqnarray*}
f_{g}:
\P_{fg'}\delta_{g'g''}\P_{g''e}(M)\boxtimes
{\big(}\P_{f_1e_1}(M_1){\big)}\\
\to
\P_{fg'}\delta_{g'g''}\P_{g''e}(M)\boxtimes
{\big(}\P_{f_1g_1}
\P_{g_1e_1}(M_1){\big)}\\
\to \P_{f\disjoint f_1,g'\disjoint g_1}
{\big(}\delta_{g'g''} \P_{g''e}(M) \boxtimes
\P_{g_1e_1}(M_1){\big)}
\\
\to \P_{f\disjoint f_1,g'\disjoint
g_1}\delta_{g'\disjoint
 g_1,g''\disjoint g_1}
\P_{g''\disjoint g_1,e\disjoint e_1}(M
\boxtimes M_1)\\
\to \P_{f\disjoint f_1,e\disjoint
e_1}(M\boxtimes M_1).
\end{eqnarray*}

Let us prove  this statement. We need to
show that for every $s\in \zebra(f\disjoint
f_1,e\disjoint e_1)$,
$$
p_s \F=p_s\sum_{g_1} \F_{g_1}.
$$

Let $$ s=\Big( (k_0,h^0)>(k_1,h_1)>\cdots
 (k_n,h_n)\Big),
 $$
where $$ f=k_0\geq k_1\geq\cdots\geq k_n=e
$$
and
$$
f_1=h_0\geq h_1\geq\cdots\geq h_n=e_1.
$$

Let
$$
s'=(f=k'_0\geq k'_1\geq\cdots\geq k'_{n'}=e)
$$
and
$$
s''=(f_1=h'_0\geq h'_1\geq\cdots\geq
h'_{m'}=e_1.)
$$
be obtained from
$$
f=k_0\geq k_1\geq\cdots\geq k_n=e
$$
and
$$
f_1=h_0\geq h_1\geq\cdots\geq h_n=e_1.
$$
by deleting  repeating terms.

Let us compute $p_s\F$. To this end we first
compute the composition
\begin{equation}\label{str}
\P_{fg'}\delta_{g'g''}\P_{g''e}\to
\P_{fe}\to
\I_{k'_0k'_1}\I_{k'_1k'_2}\cdots\I_{k'_{n'-1}k'_n},
\end{equation}
which does not vanish only if there exists
an index $\a'$ such that $k'_{\a'}=g'\geq
g''\geq k'_{\a'+1}$. This  is equivalent to
existence of an index $\a$ such that
$$
k_\a=g'\geq g''\geq k_{\a+1}.
$$
The composition \ref{str} is then equal to
\begin{eqnarray*}
\P_{fg'}\delta_{g'g''}\P_{g''e}\to
\P_{fe}\to
\I_{k_0k_1}\cdots\I_{k_{\a-1}k_{\a}}\delta_{g'g''}
\I_{g''k'_{\a+1}}\I_{k_{\a+1}k_{\a+2}}\\
\cdots
\I_{k_{n-1}k_{n'}}\\
\to \I_{k_0k_1}\I_{k_1k_2}\cdots\I_{k_{n-1}k_n}\\
\cong
\I_{k'_0k'_1}\I_{k'_1k'_2}\cdots\I_{k'_{n-1}k'_n}
\end{eqnarray*}

The projection $p_s\F$ is then equal to:
\begin{eqnarray*}
\P_{fg'}\delta_{g'g''}\P_{g''e}(M)\boxtimes
\P_{f_1e_1}(M_1)\\
\to \Big(\I_{k_0k_1}\cdots
\I_{k_{\a-1}k_a}(\delta_{g'g''}
\I_{g''k_{a+1}}\cdots
\I_{k_{n-1}k_n}(M)\Big)\boxtimes
\Big(\I_{h_0h_1}\I_{h_1h_2}\cdots
\I_{h_{n-1}h_n}(M_1)
\Big)\\
\to \I_{(k_0,h_0)(k_1,h_1)}\cdots
\I_{(k_{\a-1},h_{\a-1}),(k_\a,h_\a)} \Big\{
      \big(
          (\delta_{g'g''}
             \I_{g''k_{a+1}})\cdots \I_{k_{n-1}k_n}(M)
       \big)
      \boxtimes
      \big(
            \I_{h_\a h_{\a+1}}\cdots \I_{h_{n-1}h_n}(M_1)
      \big)
\Big\}\\
\to \I_{(k_1,h_1)}\cdots
\I_{(k_n,h_n)}(M\boxtimes M_1).
\end{eqnarray*}

Where the last map is induced by the map
\begin{eqnarray*}
       \big(
          (\delta_{g'g''}
             \I_{g''k_{\a+1}})\cdots \I_{k_{n-1}k_n}(M)
       \big)
      \boxtimes
      \big(
            \I_{h_\a h_{\a+1}}\cdots \I_{h_{n-1}h_n}(M_1)
      \big)\\
      \to\big(
      \I_{g'k_{\a+1}}\cdots \I_{k_{n-1}k_n}(M)
       \big)
      \boxtimes
      \big(
            \I_{h_\a h_{\a+1}}\cdots \I_{h_{n-1}h_n}(M_1)
      \big)\\
\to
\I_{(k_\a,h_\a)(k_{\a+1},h_{\a+1})}\cdots
\I_{(k_{n-1},h_{n-1})(k_{n},h_{n})}(M\boxtimes
M_1).
\end{eqnarray*}
This map
 is  equal to the  map:
\begin{eqnarray*}
       \big(
          (\delta_{g'g''}
             \I_{g''k_{\a+1}})\cdots \I_{k_{n-1}k_n}(M)
       \big)
      \boxtimes
      \big(
            \I_{h_\a h_{\a+1}}\cdots \I_{h_{n-1}h_n}(M_1)
      \big)\\
      \to\big(
      \delta_{g'g''}\I_{g''k_{\a+1}}\cdots \I_{k_{n-1}k_n}(M)
       \big)
      \boxtimes
      \big(
            \I_{\h_a\h_a}\I_{h_\a h_{\a+1}}\cdots \I_{h_{n-1}h_n}(M_1)
      \big)\\
\to
\delta_{(g',h_\a)(g'',h_{a}}\I_{(g'',h_\a)(k_{\a+1},h_{\a+1})}\cdots
\I_{(k_{n-1},h_{n-1})(k_{n},h_{n})}(M\boxtimes M_1)\\
\to \I_{(g',h_\a)(k_{\a+1},h_{\a+1})}\cdots
\I_{(k_{n-1},h_{n-1})(k_{n},h_{n})}(M\boxtimes
M_1)
\end{eqnarray*}
The map $p_s\F$ can be then rewritten as
follows:
\begin{eqnarray}\label{mappa}
\P_{fg'}\delta_{g'g''}\P_{g''e}(M)\boxtimes
\P_{f_1e_1}(M_1)\\\nonumber \to
\Big(\I_{k_0k_1}\cdots\I_{k_{\a-1}k_{\a}}\delta_{g'g''}
\I_{g''k_{\a+1}}\cdots\I_{k_{n-1}k_n}(M)\Big)
\boxtimes
\Big(\I_{h_0h_1}\cdots\I_{h_{\a-1}h_\a}\I_{h_\a
h_\a} \I_{h_\a
h_{\a+1}}\cdots\I_{h_{n-1}h_n}(M_1)\Big)\\\nonumber
\to
\I_{(k_0,h_0)(k_1h_1}\cdots\I_{(k_{\a-1}h_{\a-1})
(k_\a h_\a)}\delta_{(g',h_\a)(g'',h_\a)}
\I_{(g'',h_\a)(k_{\a+1},h_{\a+1})}\cdots
\I_{(k_{n-1}h_{n-1})(k_nh_n)}(M\boxtimes
M_1).
\end{eqnarray}

Let us now compute $p_s\F_{g_1}$. It follows
that such a composition is not zero only if
there exists  an $\a$  such that
$$
(k_\a,h_\a)=(g',g_1)\geq (g'',g_1)\geq
(k_{\a+1},h_{\a+1}).
$$

Since $g'>g''$, $k_\a>k_{\a+1}$. There
exists at most one $\a$ such that
$k_\a>k_{\a+1}$ and $k_a=g'$. Then $g_1$ is
uniquely determined and equals $g_1$.

In other words, there exists at most one
$g_1$ such that $p_s\F_{g_1}\neq 0$. If such
a $g_1$ does not exist, then there is no
$\a$ such that $k_\a=g'\geq g''\geq
k{\a+1}$, therefore, $p_s\F=0$. Thus, in
this case $p_s\F=p_s\sum F_{g_1}$.

If there exists an $\a$ such that
$k_\a=g'\geq g''\geq k_{\a+1}$, then
$p_s\sum \F_{g_1}=p_s\F_{h_a}$. It is not
hard to see that $p_s\F_{h_a}$ coincides
with the map (\ref{mappa}) which is the same
as $p_s\F$, whence the statement.

\subsubsection{Concatenation+concatenation}
Let $\D_1\overset \D_2\overset
\cdots\overset \D_6$. Consider the following
maps
$$
a_1:\P_{\D_1\D_2}\delta_{\D_2\D_3}
\P_{\D_3\D_4}\delta_{\D_4\D_5}\P_{\D_5\D_6}\to
\P_{\D_1\D_4}\delta_{\D_4\D_5}\P_{\D_5\D_6}\to
\P_{\D_1\D_6}
$$
and

$$
a_2:\P_{\D_1\D_2}\delta_{\D_2\D_3}
\P_{\D_3\D_4}\delta_{\D_4\D_5}\P_{\D_5\D_6}\to
\P_{\D_1\D_2}\delta_{\D_2\D_3}\P_{\D_3\D_6}\to
\P_{\D_1\D_6}.
$$

In the case when $\D_3=\D_4$ we also have a
map
$$
a_3:\P_{\D_1\D_2}\delta_{\D_2\D_3}
\P_{\D_3\D_4}\delta_{\D_4\D_5}\P_{\D_5\D_6}\to
\P_{\D_1\D_2}\delta_{\D_2\D_5}\P_{\D_5\D_6}\to
\P_{\D_1\D_6}.
$$
In the case $\D_3\neq  \D_4$ set $a_3=0$.

\begin{Proposition}
We have:
$$
a_2=a_1+a_3.
$$
\end{Proposition}
\pf Straightforward
\endpf

\subsubsection{concatenation+decomposition}
Let $\D_1\overset\cdots\D_4$ be diagonals
and let $E$ be another diagonal such that
$\D_1\overset E\overset \D_4$. Compute the
composition:
$$
\P_{\D_1\D_2}\delta_{\D_2\D_3}\P_{\D_3\D_4}\to
\P_{\D_1\D_4}\to \P_{\D_1E}\P_{E\D_4}.
$$
\begin{Proposition}
If $\D_1\overset E\overset \D_2$, then this
composition is equal to the composition:
$$
\P_{\D_1\D_2}\delta_{\D_2\D_3}\P_{\D_3\D_4}\to
\P_{\D_1E}\P_{E\D_2}\delta_{\D_2\D_3}\P_{\D_3\D_4}\to
\P_{\D_1E}\P_{E\D_4};
$$

if $\D_3\overset E\overset \D_4$, then this
composition is equal to the composition:
$$
\P_{\D_1\D_2}\delta_{\D_2\D_3}\P_{\D_3\D_4}\to
\P_{\D_1\D_2}\delta_{\D_2\D_3}\P_{\D_3E}\P_{E\D_4}\to
\P_{\D_1E}\P_{E\D_4};
$$

otherwise this composition is zero.
\end{Proposition}
\pf Straightforward.
\endpf

\subsection{Filtrations}
We will study the relationship of the above
intriduced structure with the filtration on
the functors
 $\P_{E\De}$ (see
\ref{filtrP}), whenever $E\supset \De$. We
will see how it interacts with the maps
introduced in the previous section.

Let $t_1\in \NN_{\De_1\De_2}$ and $t_2\in
\NN_{\De_3\De_4}$. Let
$$
t_1=(e_1a_{12}e_2a_{23}\ldots e_n)
$$
and
$$
t_2=(e'_1a'_{12}e'_2a'_{23}\ldots e'_{n'}).
$$

Define the concatenation
$$
t_1\delta t_2:=(e_1a_{12}\ldots
e_n\delta_{\D_2 \D_3}e'_1a'_{12}\ldots
e'_{n'}.
$$

We say that $t_2$ starts with $\delta$ if
$a'_{12}=\delta$. In this case we define one
more concatenation
$$
t_1 \delta\circ t_2 := (e_1a_{12}\ldots
e_{n}\delta_{\D_{e_{n}}\D_{e'_2}}e'_{2}
a'_{23}e'_3\ldots e'_{n'}).
$$
\begin{Proposition}
If $t_1$ does not terminate in $\delta$,
then
$$
c(F^{t_1}\P_{\D_1\D_2}\delta_{\D_2\D_3}F^{t_2}\P_{\D_3\D_4})
\subset F^{t_1\delta t_2}\P_{\D_1\D_4};
$$
if $t_1$ terminates in $\delta$, then
$$
c(F^{t_1}\P_{\D_1\D_2}\delta_{\D_2\D_3}F^{t_2}\P_{\D_3\D_4})
\subset F^{t_1\delta t_2}\P_{\D_1\D_4}+
F^{t_1\delta\circ t_2}\P_{\D_1\D_4}.
$$
\end{Proposition}
\pf Straightforward.
\endpf

We have the inuced maps

$$
\Gr^{t_1}\P_{\D_1\D_2}\delta_{\D_2\D_3}\Gr^{t_2}\P_{\D_3\D_4}
\to \Gr^{t_1\delta t_2}\P_{\D_1\D_4},
$$
if $t_2$ does not start with $\delta$; and

$$
c(\Gr^{t_1}\P_{\D_1\D_2}\delta_{\D_2\D_3}\Gr^{t_2}\P_{\D_3\D_4})
\to \Gr^{t_1\delta t_2}\P_{\D_1\D_4}\oplus
\Gr^{t_1\delta\circ t_2}\P_{\D_1\D_4},
$$
if $t_2$ starts with $\delta$.

We see that
$$
 \Gr^{t_1}\P_{\D_1\D_2}\delta_{\D_2\D_3}\Gr^{t_2}\P_{\D_3\D_4}\cong \Gr^{t_1\delta t_2}\P_{\D_1\D_4}\cong
 \Gr^{t_1\delta\circ t_2}\P_{\D_1\D_4},
$$
 whenever $t_2$ starts with
$\delta$; otherwise we have only the first
isomorphism in this chain.

The above map (on the graded components) is
induced by this isomorphism in the case when
$t_2$ does not start with $\delta$;
otherwise the above map is induced by the
direct sum of our isomorphisms.

\subsection{Resolution}\label{resdif}
Fix two equivalence relations  $f\geq  e$ on
$S$. We are going to construct a resolution
of $\i_{\D_f\D_e}$.

Denote by $\Flags(f,e)$ the set of all
'non-strict' flags of the form
$$
f=a_0\geq b_0\geq a_1\geq b_1\cdots \geq
b_n=e,
$$
where $n\geq 0$ and $b_i\neq a_{i+1}$ for
all $i$. For $i=0,1,\ldots, n-1$ set
$$
A_{(k)}=a_0b_0\cdots
a_{k-1}b_{k-1}a_kb_{k+1}
a_{k+2}b_{k+2}\cdots a_{n}b_{n},
$$
 (we delete $b_k$ and $a_{k+1}$).
In the case $a_{k+1}=b_{k+1}$ set
$$
A_{[k]}:=a_0b_0\cdots
a_kb_ka_{k+2}b_{k+2}\cdots a_nb_n,
$$
where we  delete $a_{k+1}$ and $b_{k+1}$.

Denote $|A|:=n$ and set
$$
\PP(A):=\P_{a_1b_1}\delta_{b_1a_2}\P_{a_2b_2}
\cdots \P_{a_nb_n}.
$$
Let $\PP_n=\bigoplus_{|A|=n}\PP(A)$.

Denote
$$
x_k:\P_{a_kb_k}\delta_{b_ka_{k+1}}
\P_{a_{k+1}b_{k+1}}\to \P_{a_kb_{k+1}}
$$
Let $X_k:\PP(A)\to \PP(A_{(k)})$ be the map
induced by $x_k$.

We also need maps $Y_k$ defined as follows.
In the case when $a_{k+1}=b_{k+1}$ we have
an isomorphism $Y_k:\PP(A)\to \PP(A_{[k]})$.

In the case $a_{k+1}\neq b_{k+1}$ set
$Y_k=0$.

For example: Let $|A|=2$, then the above
theorem implies that $X_1(X_1+Y_1-X_2)=0$ as
a map $\PP(A)\to \P_{fe}$.

Define the map $d_n:\PP_n\to \PP_{n-1}$ by
the formula
$$
d_n=X_1+Y_1-X_2-Y_2+X_3+Y_3+\cdots
(-1)^{n}Y_{n-1}+ (-1)^{n+1}X_n.
$$
The above identity implies that
$d_{n-1}d_n=0$; thus, $(\PP_\bullet,d)$ is a
complex.

We have a  natural map
$\ve:\PP_0\to\i_{fe}$; we have $\ve d_0=0$.

\begin{Theorem}
1) The homology $H_i(\PP_\bullet,d)=0$ for
all $i>0$

2) The map $\ve$ identifies
$H_0(\PP_\bullet)$ with $\i_{fe}$.
\end{Theorem}
\pf We are going to consider the associated
graded complex with respect to a certain
filtration which we are going to define.

Define the set $\segments(fe)^0$ whose
elements are flags of segments
$$
[a_1,b_1]>[a_2,b_2]>\cdots>[a_n,b_n]
$$
such that
$$
f\geq a_1>b_1>a_2>b_2>a_3>\cdots>b_n\geq e.
$$
For each $t\in \segments(f,e)$,
$$
t=([a'_1,b'_1]>[a'_2,b'_2]>\cdots>[a'_n,b'_n])
$$
 define an element from $\segments(f,e)^0$,
$$
\nu(t)=([a_1,b_1]>[a_2,b_2]>\cdots>[a_k,b_k])
$$
according to the rule:

the sequence $a_1b_1a_2b_2\cdots b_k$ is
obtained from the sequence $a'_1b'_1\cdots
a'_nb'_n$ by deleting all its repeating
terms.

For an $s\in \segments(FE)^0$ we set
$$\Phi^s\PP_\bullet=\bigoplus_{\nu(t)=s} F^t(\PP_\bullet).
$$
We see that $F$ is a filtration on
$\PP_\bullet$ and that the associated graded
complex can be computed by the formula
$$
\Gr_\Phi^s\PP_\bullet=\bigoplus_{\nu(t)=s}
\Gr^t(\PP_\bullet).
$$
Let $f'>e'$ be a pair of equivalence
relations on $S$.
 Let $\oo\in \segments(f'e')^0$ be the least element,
 which
  is simply
$[f'e']$. Denote
$\RR_{f'e'}:=\Gr_F^{\oo}\PP_\bullet$. Let
$s\in \segments(fe)^0$  be an arbitrary
element;
$$
s=[a_1,b_1]>[a_2,b_2]>\cdots>[a_nb_n].
$$
We then have
$$
\PP_\bullet^s\cong
\i_{fa_1}\RR_{a_1b_1}\i_{b_1a_2}\RR_{a_2b_2}\cdots
\ldots \i_{\b_ne}.
$$

This implies that our task is reduced to
proving that $\RR_{fe}$ is acyclic, which
will be done in the next subsection.
\subsubsection{}

We see that the complex $\RR_{fe}$ is
isomorphic to the complex
$R_{fe\bullet}\otimes\delta_{fe}$, where
$R_{fe\bullet}$ is a complex of vector
spaces; the vector space  $R_{fen}$ has a
basis labelled by the elements
$$
H=(f=e_1u_{12}e_2u_{23}e_3\ldots e_N=e),
$$
where $e_1>\cdots > e_N$ each $u_{kk+1}$ is
either $p$ or $\delta$ and the total number
of deltas is $n$. Denote $|H|:=N$ The
differential is given by the sum of several
terms which we are now going to describe.
Let $A_kH$ be zero if $u_{kk+1}=p$ and let
it change $u_{kk+1}$ from $\delta$ to $p$
otherwise.

Let $B_kH$ be non-zero only if
$u_{kk+1}=\delta,u_{k+1k+2}=\delta$, in
which case it replaces the fragment
$u_{kk+1}e_{k+1} u_{k+1k+2}$ with
$\delta_{kk+2}$.

Let $C_kH$ be non-zero only if
$u_{kk+1}=\delta$ and $u_{k+1k+2}=p$, in
which case the fragment
$$u_{kk+1}e_{k+1}
u_{k+1k+2}
$$
is going to be replaced with $p$.

Denote by $d_k$ the number of symbols
$\delta$ before $a_{kk+1}$. It follows that
the differential on $R$ is given by
$$
d=\sum (-1)^{d_k}(A_k+B_k+C_k).
$$
Set $$ F_NR:=\oplus_{|H|\leq N} {\Bbb C}H.
$$
It is clear that $F$ is a filtration on $R$
the assoiated graded complex has the basis
labelled by the same elements, the
differential is given by
$$
d'=\sum (-1)^{d_k}A_k.
$$

Let $$ \Phi=(f= e_1> e_2>\cdots> e_n=e)
$$
be a flag and let $R_\Phi\subset \Gr_F R$ be
the subcomplex spanned by the elements
$$
H=(e_1u_{12}\ldots e_n),
$$
with arbitrary $u_{ii+1}$ (it is clear that
it is a subcomplex.)

We have $\Gr_F R=\oplus_F R_F$. Furthermore,
let $V={\Bbb C}<\delta,p>$ be  a complex in
which $|\delta|=1$; $|p|=0$ and $d\delta=0$.
Then $R_F\cong T^NV$ and is therefore
acyclic.
\subsection{The structure of system on the collection
of functors $\PP_{fe}$}
\subsubsection{} Let $f\geq g\geq e$ be a sequence of
equivalence relations. Define the
decomposition map
$$\as_{fge}:\R_{fe}\to\R_{fg}\R_{ge}.$$
Let
$$
A=(f=a_1b_1a_2\ldots a_nb_n=e)
$$
be an elements of $\Flags(fe)$. If there
exists $k$ such that $a_k\overset  g\overset
b_k$, then set
$$
a:\P(A)\to \P(a_1b_1\ldots
b_{k-1}a_kg)\P(gb_k\cdots a_nb_n)
$$
is induced by the  decomposition map
$$
\P_{a_kb_k}\to\P_{a_kg}\P_{gb_k}.
$$

Otherwise we set

$a|_{P(A)}=0$.
\subsubsection{Factorization maps}
We will first study

\subsection{Factorization maps for $\R$ }
\def\ba{{\frak a}}\def\bb{{\frak b}}
\subsubsection{}
We keep the notations of the previous
subsection. Let $F_\a\in \Flags(f_\a,e_\a)$,
$\a\in A$
 and $F\in \Flags(f,e)$.
We are going to define the map
$$
\mu(\{F_\a\}_{\a\in
A};F):\boxtimes_\a\R(F_\a)\to \R(F).
$$

This map is zero for all $F_a,F$ except
those determined by the following
conditions.

Let
$$
F=(f=\ba_1\geq \bb_1>a_2\geq \bb_2\cdots
\ba_n\geq \bb_n=e).
$$
 We then require that

1)For every $i$: $\bb_{i\a}=\ba_{i+1\a}$ for
all $\a$ except exactly one (denote it by
$\a_i$);

2) Fix $\a$ and consider the sequence
$$
\ba_{1\a}\geq
 \bb_{1\a}\geq \ba_{2\a}\geq\cdots\geq \ba_{n\a}
 \geq \bb_{n\a}.
$$
Construct a subsequence
$$
F(\a)=(\ba_{M_1\a}\bb_{N_1\a}\ba_{M_2\a}\bb_{N_2\a}\cdots
\ba_{M_r\a}\bb_{N_r\a})
$$
according to the following rule: we delete
every pair $\bb_{i\a}\geq \ba_{i+1\a}$ in
which $\bb_{i\a}=\ba_{i+1\a}$. We have:
$M_1=1$; $N_r=n$; $\bb_{N_i\a}\neq
\ba_{M_{i+1}\a}$; $M_{r+1}=N_r+1$.
Therefore, $F(\a)\in \Flags(f_\a,e_\a)$. Our
second condition is then $F_\a=F(\a)$ for
all $\a$.
\subsubsection{}
We have a natural map
$$
r_\a:\R(F(\a))\to \P_{\ba_{1\a}\bb_{1\a}}
\delta_{\bb_{1\a}\ba_{2\a}}
\cdots\P_{\ba_{n\a}\bb_{n\a}},
$$
induced by the maps
\begin{eqnarray*}
\P_{\ba_{M_{i\a}}\bb_{N_{i\a}}}\to
\P_{\ba_{M_{i\a}}\bb_{M_{i\a}}}\P_{\ba_{M_{i\a}+1}
\bb_{M_{i_a}+1}}\cdots
\P_{\ba_{N_{i\a}}\bb_{N_{i\a}}}
\end{eqnarray*}
which induce maps
\begin{eqnarray*}
\P_{\ba_{M_{1\a}}\bb_{N_{1\a}}}\delta_{\bb_{N_{1\a}}\ba_{M_{2\a}}}
\cdots\P_{\ba_{M_{r\a}}\bb_{N_{r\a}}}\to\\
\P_{\ba_{M_{1\a}}\bb_{M_{1\a}}}
\P_{\ba_{M_{1\a}+1}\bb_{M_{1\a}+1}}\cdots
\P_{\ba_{N_{1\a}}\bb_{N_{1\a}}}\delta_{\bb_{N_{1\a}}
\ba_{M_{2\a}}}\\
\P_{\ba_{M_{2\a}}\bb_{M_{2\a}}}
\P_{\ba_{M_{2\a}+1}\bb_{M_{2\a}+1}}\cdots
\P_{\ba_{N_{2\a}}\bb_{N_{2\a}}}\delta_{\bb_{N_{2\a}}
a_{M_{3\a}}}\cdots\\
\P_{\ba_{M_{r\a}}\bb_{M_{r\a}}}
\P_{\ba_{M_{r\a}+1}\bb_{M_{r\a}+1}}\cdots
\P_{\ba_{N_{r\a}}\bb_{N_{r\a}}}\\
\cong
\P_{\ba_1\bb_1}\delta_{\bb_1\ba_2}\P_{\ba_2\bb_2}
\cdots\P_{\ba_n\bb_n}.
\end{eqnarray*}

We then define
\begin{eqnarray*}
\boxtimes_\a \R(F_\a)\stackrel{\prod_\a r_\a}\to\\
\boxtimes_\a\P_{\ba_{1\a}\bb_{1a}}\delta_{\bb_{1\a}
\ba_{2\a}}
\cdots\P_{\ba_{n\a}\bb_{n\a}}\\
\to\P_{\ba_1\bb_1}\delta_{\bb_1\ba_2}\cdots\P_{\ba_n
\bb_n}\\
\cong \R(F).
\end{eqnarray*}

\subsubsection{Signs} The function $i\mapsto s_i$
defines a partition of the  set
$\{1,2,\ldots,n\}$ Fix an orientation of $S$
and denote  by
 $s(F)$  the sign of this partition.
\subsubsection{Definition of the map}
Define
$$
\mu=\sum_F s(F)\mu(\{F(\a)\}_{\a\in A},F).
$$
\subsubsection{} We are going to check that
$\mu$ commutes with the differential.

This follows from the several statements we
are going to formulate.

We assume that $F$ satisfies the conditions
from the previous section. 1) Let
$$
F<i>=(\ba_1\bb_1\ldots
\ba_i\bb_{i+1}\ba_{i+2}\ldots \bb_n);
$$
Let $X_i^F:\R(F)\to \R(F<i>)$ be induced by
the map
$$\P_{\ba_i\bb_i}\delta_{\bb_i\ba_{i+1}}
\P_{\ba_{i+1}\bb_{i+1}}
\to\P_{\ba_i\bb_{i+1}}.
$$

Let $U(F,i)$  be the set of all $F'\in
\Flags(f,e)$ which are obtained from $F$ by
changing $\bb_i,\ba_{i+1}$ only in such a
way that $\a_i$, $\bb_{i\a_i}$, and
$\aa_{i+1\a_i}$ do not change.

This means that every $F'$ is of the form
$$
\ba_1\bb_1\ba_2\bb_2\cdots\ba_{i}\bb'_i\ba'_{i+1}
\bb_{i+1}\ba_{i+2}\cdots\ba_n\bb_n,
$$
where $\bb'_{i,\a_i}=\bb_{i,\a_i}$;
$\ba'_{i+1,\a_i}=\ba_{i+1,\a_i}$, and for
all $\a\neq \a_i$,
$\bb'_{i\a}=\ba'_{i+1\a}$.

Let $j$ be such that $N_{j\alpha_i}=i$ (such
a $j$ always exists and is unique because
$\bb_{i\a_i}\neq \ba_{i+1}$). We then have:
$$
\sum\limits{F'\in U(F,i)}
X_i^{F'}\mu(F'(\a)_{\a\in A},F')=
\mu\Big(\{F(\a)_{\a\neq \a_i};
F(\a_i)<j>\};F<i>\Big)X_j^{F_{\a_i}}.
$$

To check this identity it suffices to
consider the case $n=2$, in which case the
statement follows immediately
 from
(\ref{sochet}).

2.Let $$
F=(\ba_1\bb_1\ba_2\bb_2\cdots\ba_n\bb_n).
$$
Assume that $\ba_i=\bb_{i}$ and set
$$
F[i]=(\ba_1\bb_1\ldots
\bb_{i-1}\ba_{i+1}\bb_{i+1}
\ldots\ba_n\bb_n.
$$
 We then have a natural map
$$
Y_i:\R(F)\to \R(F[i]).
$$

There are two cases:

Case 1: $\a_{i-1}=\a_{i}$.  Let $j$ be such
that $N_{j\a_i}=i$. In this case we have:
$$
Y_i\mu({F(\a)_{\a\in A}};F)=
\mu\Big(\{F(\a)_{\a\neq
\a_i},F(\a)[j]\};F[i]\Big) Y_j^{F(\a_i)}.
$$

Case 2.

$\a_{i-1}\neq \a_{i}$ In this case define
$$
F'=(\ba_1\ldots \bb_{i-1}\ba'_i\bb'_{i}
\ba_{i+1}\bb_{i+1}\ldots \ba_n\bb_n)
$$

in  such a way that
$\a'_i:=\a^{F'}_{i}=\a_{i+1}$,
$\a'_{i+1}=\a_i$ and
$\ba'_{i\a'_i}=\ba_{i\a_{i+1}}$.

We then have
$$
Y_i^F\mu(\{F(\a)_{\a\in A}\},F)=Y_i^{F'}
\mu(\{F'(\a)_{\a\in A}\},F').
$$

These facts imply that the factorization map
commutes with  the differential.
\subsection{}
The factorization  commutes with the
asymptotic
decomposition. We omit the proof as it is straightforward.

\subsection{The system $\m$ and a map $<\R>\to
<\m>$} Was discussed in detail above ...
\def\bt{\boxtimes}
\def\hl{\m}
\def\cc{C}\def\bt{\boxtimes}\def\Dmod{\text{D-mod}}
\def\Diff{{\cal D}}
\def\ii{\i}
\def\jj{\R}

\section{Bogolyubov-Parasyuk theorem}
Let $<\jj>$ be the resolution of the system
$<\ii>$ constructed in the previous section
and let $M$ be a cofibrant dg-$\D_X$-sheaf
endowed with an OPE-product over $<\ii>$.
\begin{Theorem} There exists an OPE
structure on $M$ over $<\jj>$ which lifts
that over $<\ii>$.
\end{Theorem}

The proof will occupy the rest of the
section.

\subsection{Unfolding the definition
of an OPE-algebra over $<\jj>$}

Let $p:S\to T$ be a surjection of finite
sets and $N$ a $\D_{X^T}$-module. We have
$\jj_p^0(N)\cong \P_p(N)\to \I_p(N)$. This
produces a natural transformation

(whose differential is not zero):
$$
\pi_p:\jj_p\to \I_p.
$$

Thus, we have an induced map $M^{\boxtimes
S}\to \I_p(M^{\boxtimes T})$.

We also have a map of systems
$$
\jj\to \hl$$
 which induces a strong homotopy *-Lie
 algebra structure on $M$. It turns out that
 the maps $\pi_p$ and the *-SHLA  structure
 on $M$ completely determine the
 OPE-structure on $M$. The precise
 formulation will be given below.
\subsubsection{} Suppose that for every surjective map
$p_S:S\to \pt$,  we are given a map
$$\aa_S:M^{\boxtimes S}\to \I_{p_S}(M),
$$
such that:

 - for $\#S=1$ we have: $\aa_S=\Id$;

 - $\aa_S$ is
equivariant with respect to bijections of
finite sets.

Assume, in addition, that we are given some
maps
$$
\cc_S:M^{\boxtimes S}\to \delta_{p_S}(M)
$$
of degree 1, where $\#S>1$, $\cc_S$ are
equivariant with respect to bijections of
finite sets.

We shall impose certain conditions on these
maps which will allow us to construct an
$OPE$-structure on $M$ using these maps.

\def \ap{{\frak ap}}
\def \bu{{\frak j}}
\def\bfp{{\bf p}}
\subsubsection{Condition 1}Let $q:S\to T$ be
a surjection of finite sets. As usual, the
product of maps $\aa_S$ gives rise to maps
$$\aa_q:M^{\boxtimes S}\to \I_q(M^{\boxtimes
T}).$$

Our first condition is as follows.

\begin{Condition}{\bf C1}
Let
$$
S\stackrel r\to R\stackrel s\to T
$$
be a sequence of surjections and $q=sr$.
Then the following diagram should commute:
$$
\xymatrix{ M^{\bt S}\ar[r]\ar[d]&
\I_r(M^{\bt R})\ar[r]&
\I_r\I_s(M^{\bt T})\ar[d]\\
\I_q(M^{\bt T})\ar[rr]&&\I_r\ii_s(M^{\bt
T})}
$$
\end{Condition}

One sees that it suffices to check this
condition for all $p:S\to \pt$.

This condition implies the following fact.
Let
$$\aa'_p:M^{\bt S}\to \I_pM^{\bt T}\to
\ii_pM^{\bt T}.
$$

Let $\bfp:=(p_i:S_i\to S_{i+1})$,
$i=0,\ldots, n-1$, be a sequence of
surjections, where $S_0=S$, $S_n=T$, and
$$
p_{n-1}p_{n-2}\cdots p_0=p.
$$

Let $\aa^\I:=\aa$ and let $\aa^\ii:=\aa'$.
Let $\bu:=(\bu_1,\bu_2,\ldots, \bu_{n-1})$
be an arbitrary sequence of elements from
the set $\{\I,\ii\}$.

Define the map
$$
\bu_\bfp:M^{\bt S}\to
(\bu_1)_{p_1}(\bu_{2})_{p_2}\cdots
(\bu_{n-1})_{p_{n-1}}(M^{\bt T})
$$
by the formula:
$$
\xymatrix{ M^{\bt
S}\ar[r]^{\aa^{\bu_1}_{p_1}}&
(\bu_1)_{p_1}(M^{S_1})\ar[r]^{\aa^{\bu_2}_{p_2}}&
(\bu_1)_{p_1}(\bu_2)_{p_2}(M^{\bt
S_2})\ar[r]&\cdots\ar[r]&(
\bu_1)_{p_1}(\bu_{2})_{p_2}\cdots
(\bu_{n-1})_{p_{n-1}}(M^{\bt T})}
$$

Condition C1 implies that the collection of
maps $\bu_{\bfp}$ for all $\bu$ and $\bfp$
determines a map $\ope_p:M^{\bt S}\to
\P_p(M^{\bt T})$.

\subsubsection{Condition 2} Let us now formulate the
condition on the collection of maps $\cc_S$
which is equivalent to the fact that this
collection endows $M[-1]$ with a structure
of *-SHLA.

We will formulate this condition in a
slightly unusual way. Let $p:S\to T$ be a
surjection. Define the map
$$
\cc_p:M^{\bt S}\to \delta_p(M^{\bt T})
$$
according to the following rule.

1. The map $\cc_p$ is not equal to zero only
if there exists a unique  $t_p\in T$ such
that $\#(p^{-1}t_p)>1$ (in which case
$\#(p^{-1}(t))=1$ for all $t\neq t_p$.)

2. If the above condition holds,  then
$\cc_p$ is defined as follows. Let
$S_p:=p^{-1}t_p$ and $S':=S\backslash S_p$.
Then $\cc_p$ is defined as the composition:
\begin{eqnarray*}
 M^{\bt S}\cong  M^{\bt S_p}\bt
M^{\bt S'}\\
\stackrel{\cc_{S_p}\bt \Id}{\longrightarrow}
\delta_{S_p}(M)\bt M^{\bt S'}\cong
\delta_p(M^{\bt T})
\end{eqnarray*}
where the last arrow is constructed via the
natural identification $T\cong S'\disjoint
\pt$.

Let now $p:S\to T$ be a surjection and let
$\Sigma_p$ be the set of all isomorphism
classes of splittings
$$\xymatrix{
S\ar[r]^q & R\ar[r]^r & T},
$$
where $r,q$ are surjections and $p=rq$. Then
the *-SHLA axiom can be formulated as
follows:
\begin{Condition}{\bf C2}
For any surjection $p:S\to T$ we have:

$$d\cc_p+\sum\limits_{(r,q)\in
\Sigma_p}\cc_r\cc_q=0,
$$
where we pick one representative for each
element in $\Sigma_p$.
\end{Condition}

It is clear that if this condition is
satisfied for all $p:S\to \pt$, then it is
satisfied for all $p$.
\subsubsection{Condition 3}
This condition describes the differential of
the maps $\aa_p$. Let $p,q,r$ be the same as
in the previous subsection. We have the
natural transformation
$$
f_{qr}:\delta_q\I_r\to\I_{qr}.
$$
Using this transformation, define a map
$$\xymatrix{
\phi_{qr}: M^{\bt
S}\ar[r]^{\cc_q}&\delta_q(M^{\bt
R})\ar[r]^{\aa_r}& \delta_q\I_r(M^{\bt
T})\ar[r]^{f_{qr}}& \I_{p}(M^{\bt T})}$$

\begin{Condition}{\bf C3} For every
surjection $p:S\to T$ we have:

$$
d\aa_p+\sum\limits_{(q,r)\in
\Sigma_p}\phi_{q,r}=0.
$$
\end{Condition}

As in the previous subsection, if   this
condition  holds for  all $p:S\to \pt$, then
it holds for all $p$.

 \subsubsection{}
We will show how, having the maps
$\aa_S,\cc_S$ satisfying conditions C1-C3,
one can construct an OPE structure on $M$
over $<\R>$.

The definition of $<\R>$ implies that to
define an OPE-structure over $<\R>$, we have
to prescribe maps
\begin{equation}\label{maptop}
M^{\bt S}\to
\P_{a_1}\delta_{b_1}\P_{a_2}\delta_{b_2}\P_{a_3}\cdots
\P_{a_{n}}\delta_{b_n}\P_{a_{n+1}}(M^{T}),
\end{equation}
where $a_i:S_{2(i-1)}\to S_{2i-1}$;
$b_{i}:S_{2i-1}\to S_{2i}$ are surjections;
$S_0=S$, $S_{2n+1}=T$, and $b_i$ are not
bijections. We define the map (\ref{maptop})
as  the composition
$$
\xymatrix{ M^{\bt S}\ar[r]^{\ope_{a_1}}&
\P_{a_1}M^{\bt S_1}\ar[r]^{\cc_{b_1}}&
\P_{a_1}\delta_{b_1}M^{\bt
S_2}\ar[r]&\cdots\ar[r]&
\P_{a_1}\delta_{b_1}\cdots
\P_{a_{n}}\delta_{b_n}\P_{a_{n+1}}(M^{T})}
$$

One checks straightforwardly that all the
conditions are satisfied.

\subsection{Proof of the Bogolyubov-Parasyuk
theorem}

We are going to use induction. To this end
introduce a notion  of $N$-OPE-{\em
structure on} $M$ (over $<\jj>)$, where
$N\geq 2$ is an integer. This means that the
maps $\aa_S,\cc_S$ are only defined when
$\#S\leq N$ and the conditions C1-3 are
satisfied for all surjections $p$ such that
$\forall i \#(p^{-1}(i))\leq N$.

The theorem follows from two statements:

1) (base of induction). There exists a 2-OPE
structure on $M$ such that the composition

$$\xymatrix{
M\boxtimes
M\ar[rr]^{\aa_{\{1,2\}}}&&\I_{\{1,2\}}(M)\ar[r]&
\ii_{\{1,2\}}(M)}
$$
equals to $\ope_{\{1,2\}}$.

2) (transition) Assume there exists an
$N$-OPE structure on $M$ such that for every
finite set $S$ with $\#S\leq N$ the
composition

\begin{equation}\label{pn}
\xymatrix{ M^{\bt S}\ar[r]^{\aa_S}&
\I_S(M)\ar[r]& \ii_S(M)}
\end{equation}
coincides with $\ope_S$. Then there exists
an $(N+1)$-OPE-structure on $M$ such that
for all $S$ with $\# S\leq N$ the maps
$\aa_S,\cc_S$ coincide with the existing
ones and the composition (\ref{pn})
coincides with $\ope_S$ for all $S$ with
$\#S\leq N+1$.

The statement 1 follows from  surjectivity
of the map $\I_{\{1,2\}}(M)\to
\ii_{\{1,2\}}(M)$ (because $M$ is
cofibrant). Therefore, the induced map
$$
r:\hom(M\boxtimes
M,\I_{\{1,2\}}(M))^{S_{\{1,2\}}}\to
\hom(M\boxtimes
M,\ii_{\{1,2\}}(M))^{S_{\{1,2\}}}
$$
is also surjective. Let $\aa_{\{1,2\}}$ be
any lifting of $\ope_{\{1,2\}}$. Then
$r(d\aa_{\{1,2\}})=0$, therefore, the image
of $d\aa_{\{1,2\}}$ is
$\delta_{\{1,2\}}(M)$. Set
$\cc_{\{1,2\}}=-da_{\{1,2\}}$. It is clear
that $(\aa_{\{1,2\}},\cc_{\{1,2\}})$
determine a 2-OPE-structure.
\bigskip

Statement 2. Let  $\j$  be the functor from
the category $\zebra(p_S)$ to the category
of functors $\Dmod_X\to \Dmod_{X^S}$ as in
(\ref{lemmaforbp})

 and let
$$
\P^0:=\P^0_{S}=\liminv_{s\in
\zebra^0(p_S)}\j(s).
$$
Let
$$
\P:=\P_{S}=\liminv_{s\in \zebra(p_S)}\j(s).
$$

 The existing $N$-OPE product defines an
equivariant  map
$$
\aa^0_S:M^{\bt S}\to \P^0(M).
$$

According to the Lemma from
(\ref{lemmaforbp}), the map
$$
\P(M)\to \P^0(M)
$$
is surjective. Therefore, there exists an
equivariant lifting
$$
\aa^1: M^{\bt S}\to \P(M)
$$
of $\aa^0$. Define $ \aa_S $ as the
composition
$$
M^{\bt S}\to \P(M)\to \j(M).
$$
 The condition C1 is then automatically
 satisfied.  The map $\cc_S$ can be uniquely
 found from the condition $C3$. Indeed, let
 $p_S:S\to \pt$. Let
 $$\Sigma_S^0:=
 \Sigma_{p_S}\backslash \{(p_S,\Id_S)\}.
 $$
 Then C3 reads as:
$$
\phi_{p_S,\Id_S}=-d\aa_S-\sum_{(q,r)\in
\Sigma_{S}^0}\phi_{q,r}.
$$
The right hand side is uniquely determined
by the  existing $N$-OPE structure and by
the chosen map $\aa_S$. It is only the left
hand side that depends on $\cc_S$. One can
find a unique $\cc_S$ satisfying C3 iff the
right hand side is a map whose image is
contained in $\delta_S(M)\subset \I_S(M)$.
Let us show that this is indeed the case.
Denote the map specified by the right hand
side by $u:M^{\boxtimes S}\to \I_S(M)$. The
image of $u$ lies in $\delta_S$ iff for
every $(q,r)\in\Sigma_S^0$, the through map
$$\xymatrix{
M^{\boxtimes S}\ar[r]^{u} &\I_S(M)\ar[r]&
\I_q\ii_r(M)}
$$
is zero. This can be checked directly.

 With such a choice  of $\cc_S$ the
 condition C3 is satisfied.

 The condition $C2$ is satisfied as well, as
 follows from the direct computation.

 Bogolyubov-Parasyuk theorem is proven.

\def\disjoint{{\sqcup}}
\def\bs{\backslash}
\def\A{{\frak A}}
\def\Aut{\text{Aut}}
\def\intr{\int}
\def\iI{{\frak iI}}
\def\ins{{\frak ins}}
\def\inr{{\frak int}}
\def\la{\lambda}
\def\An{{\frak An}}
\def\P{{\cal P}}
\def\limdir{\text{limdir}}
\def\ev{\text{ev}}
\def\notsubset{\text{notsubset}}
\def\identical{=}
\def\ins{{\frak ins}}
\def\dirlim{\limdir}
\def\usual{{\text{usual}}}
\def\A{{\cal A}}
\def\B{{\cal B}}
\def\F{{\frak F}}
\def\a{{\frak a}}
\def\b{{\frak b}}
\def\e{{\frak e}}
\def\jk{{\frak jk}}
\def \q{{\frak q}}
\def \qq{\q}
\def\sym{{\text{sym}}}
\def\res{{\text{res}}}
\def\f{{\fun}^{\sim}}
\def \Im{\text{Im}}
\def\jj{{\frak jj}}
\def\fp{({\fun'})^{\sim}}
\def\r{{\frak r}}
\def \into{\to}
\def \comm{\text{\bf comm}}
\def\lis{{\bf lis}}
\def\cat{{\bf cat}}
\def\red{\text{red}}
\def\lism{{\bf lism}}
\def\complexes{\text{complexes}}
\def\op{\text{op}}
\def \ores{{\O^{\text{res}}}}
\def\asm{{\frak ope'}}
\def\ens{{\bf ens}}
\def\ske{{\bf ske}}
\def\ve{\varepsilon}
\def\vect{{\bf vect}}
\def \id{\text{Id}}
\def\Id{\id}
\def \as{{\frak as}}
\def\ope{{\frak ope}}
\def\fun{{\frak fun}}
\def\funp{({\frak fun}^{\psys})}
\def \Fun{{\frak Fun}}
\def\funct{{\frak Funct}}
\def\system{{{\cal SYS}}}
\def\psys{{\cal SYS}[p]}

\def\lie{{\frak lie}}
\def\li{\lie^1}
\def \ll{{{\li}_M}}
\def\psystem{{{\cal SYS}[p,\ll]}}
\def\psystemg{{{\cal SYS}[p,\lie_M]}}
\def\OPE{{\cal OPE}}

\def \OPER{{\cal OPE}^{\text{res}}}
\def\psystemb{\OPE[p,\lie_M]}
\def\psystemg{\psystemb}
\def\psis{\psystemb^{\text{res}}}

\def\univ{\text{univ}}
\def\psim{\psystemb^{\text{sym}}}
\def\psimr{\psystem^{\text{sym},\text{res}}}

\def\lisp{{[p,\lie_M,\S]}}
\def\lispq{{\lisp^{\text{res}}}}
\def\psml{\system[p,\lie_M,\S]}
\def\psmq{\psml^{\text{res}}}
\def\psmlh{\psml[[\h]]}
\def\psmqh{\psmq[[\h]]}
\def\coker{\text{Coker}}
\def\M{{\cal M}}
\def\gen{{\frak gen}}
\def\u{{\frak u}}
\def\Mg{{\frak M}}
\def \hh{{\frak h}}
\def \kk{{\frak kk}}
\def \II{{\frak II}}
\def \JJ{{\frak JJ}}
\def \GG{{\frak GG}}
\def \FF{{\frak FF}}
\def \m{{\frak m}}
\def \z{{\frak z}}
\def \N{{\cal N}}
\def \B{{\cal B}}
\def \G{{\cal G}}
\def \M{{\cal M}}
\def \O{{\cal O}}
\def \R{{\Bbb R}}
\def \Z{{\frak Z}}
\def \Co{{\Bbb C}}
\def \C{{\cal C}}
\def \i{{\frak i}}
\def \j{{\frak j}}
\def \k{{\frak k}}
\def \ii{{\frak ii}}
\def \l{{\frak l}}

\def \L{{\frak L}}
\def \DDD{{\cal D}}
\def \D{{\cal D}}
\def \DD{\Delta}
\def \g{{\frak g}}
\def \ha{^\wedge}

\def \liminv{\text{liminv}}
\def \p{{\frak p}}
\def \pp{{\pi_p}}
\def \h{h}
\def \I{{\cal I}}
\def \s{{\frak s}}
\def \S{{\frak S}}
\def \U{{\cal U}}
\def \vs{\sigma}
\def\pt{\text{pt}}
\def\Hom{\text{Hom}}
\def \holi{{\frak holi}}
\def \La{{\holi_M}}
\def \Lab{{ {\frak holie}_M}}
\def \psl{\OPER[p,\Lab]}
\def \ps {\psl}

\section{The maps $\R_{f\disjoint g}\to
\R_{f\disjoint \Id}\delta_{\Id\disjoint g}$}
\subsection{Notations}\label{Notat}
\subsubsection{} Let $\phi:S\to T$, $g:A\to
B$ be surjections. Define a functor
$\iI_{\phi\disjoint g}$ from the category of
$\D_{X^{T\disjoint B}}$-modules to the
category of $\D_{X^{S\disjoint A}}$-modules
 by:
$$
\iI_{\phi\disjoint g}(M)=i_{\phi\disjoint
g}^{\ha}(M)\otimes_{\O_{X^{S\disjoint A}}}
(\B_\phi\boxtimes \C_g).
$$

One can also define $\iI_{\phi\disjoint g}$
as a quotient of $\I_{\phi\disjoint g}$ with
by the sum of images of all maps
$$
\delta_{\phi_1\disjoint
\Id}\I_{\phi_2\disjoint g}\to
\I_{\phi\disjoint g},
$$
where $\phi=\phi_2\phi_1$, $\phi_1,\phi_2$
are surjections,and $\phi_1$ is not
bijective. \subsubsection{} We then have
natural maps
\begin{equation}
\I_{\psi\phi\disjoint
g}\to\I_{\phi}\iI_{\psi\disjoint g },
\end{equation}
which shall be denoted by $a_{\psi\phi\times
g}$.
\subsection{Map $
\xi(\phi,g):\P_{\phi\disjoint  g}\to
\I_{\phi\disjoint \Id}\delta_{\Id\disjoint
g}$}Let
$$\phi:S\to T;\quad g:A\to B$$
be
surjections. We shall define a map
$$ \xi(\phi,g):\P_{\phi\disjoint  g}\to
\I_{\phi\disjoint \Id}\delta_{\Id\disjoint
g}$$ recursively. The parameter of the
recursion will be $|g|=\#A -\# B$. Since $g$
is surjective, $|g|\geq 0$. To describe the
recursive procedure we need to introduce
some notation.

 Suppose we are given an (arbitrary)
collection of maps
$$\xi(\phi,g)$$ for all $\phi$ and all $g$
with $|g|<N$. Fix a $g$ with $|g|=N$. We
then construct a map
$$
X(\phi,g):\P_{\phi\disjoint g}\to
\I_{\phi\disjoint \Id} \I_{\Id\disjoint g}
$$
by means of the formulas:
$$X(\phi,g)=U(\phi,g)-\sum\limits_{g=g_1\circ g_2,\ g_1\neq g}F(\phi,g_1,g_2),
$$
where
$$
U(\phi,g):\P_{\phi\disjoint g}\to
\I_{\phi\disjoint g}\to \I_{\phi\disjoint
\Id} \I_{\Id\disjoint g};
$$
\begin{equation}\label{Fg1g2}
F(\phi,g_1,g_2):\P_{\phi\disjoint g}\to
\P_{\phi\disjoint g_1}\P_{\Id\disjoint
g_2}\to \I_{\phi\disjoint \Id
}\delta_{\Id\disjoint g_1}\I_{\Id\disjoint
g_2}\to \I_{\phi\disjoint
\Id}\I_{\Id\disjoint g}.
\end{equation}

The recursive procedure will be now
described by means of:
\begin{Definition-Proposition}
There exists a unique collection of maps
$\xi(\phi,g)$ for all surjections $\phi,g$
such that

1) If $|g|=0$, i.e  $g$ is a bijection, then
$\xi(\phi,g)$ is the natural isomorphism
induced by $g$.

2) The composition $$\P_{\phi\disjoint g}\to
\I_{\phi\disjoint \Id}\delta_{\Id\disjoint
g}\to \I_{\phi\disjoint \Id}\I_{\Id\disjoint
g}
$$
equals $X(\phi,g)$.
\end{Definition-Proposition}
\pf We shall prove by  induction in $|g|$
that given a natural $N$, the required maps
$\xi(\phi,g)$ can be constructed for all $g$
with $|g|\leq N$.

The base of induction, $N=0$, is evident.
Let us now pass to the transition. Pick a
$g$ with $|g|=N$ and assume that our
statement is the case for all $g'$ with
$|g'|<N$.

 We will then show
that for every  decomposition $g=lk$, where
$k,l$ are proper surjections (i.e.
surjections but not bijections) the through
map
\begin{equation}\label{delt}
X(\phi,g):\P_{\phi\times g}\to
\I_\phi\I_g\to \I_\phi\I_k\i_l\end{equation}
is zero. Indeed, we have the following
commutative diagrams:

I.$$ \xymatrix{ P_{\phi\disjoint
g}\ar[r]^{U(\phi\disjoint g)}\ar[d]&
\I_{\phi\disjoint \Id}\I_{\Id\disjoint g}\ar[r]&
\I_{\phi\disjoint \Id}\I_{\Id\disjoint k}\i_{\Id\disjoint l}\\
\P_{\phi\disjoint k}\I_{\Id \disjoint l}\ar[urr]^{U(\phi,k)}\\
}
$$

II. The composition
\begin{equation}\label{diagII}
 \xymatrix{ \P_{\phi\disjoint
g}\ar[rr]^{F(\phi,g_1,g_2)} &&
\I_{\phi\disjoint \Id}\I_{\Id\disjoint
g}\ar[r]&
\I_{\phi\disjoint\Id}\I_{\Id\disjoint
k}\i_{\Id\disjoint l}}
\end{equation}
does not vanish only if one can decompose
$g=lug_1$ in such a way that $g_2=lu$ and
$k=ug_1$. In this case the map
(\ref{diagII}) is equal to the composition:
$$
\xymatrix{ \P_{\phi\disjoint g}\ar[r]&
\P_{\phi\disjoint g_1}\I_{\Id\disjoint
lu}\ar[r]& \P_{\phi\disjoint \Id}
\delta_{\Id\disjoint g_1}\I_{\Id\disjoint
u}\i_{\Id\disjoint l}\ar[r]&
\I_{\phi\disjoint \Id}\I_{\Id\disjoint
g_1u}\i_{\Id\disjoint l}}
$$
Therefore, the composition (\ref{delt}) is
equal to
$$\xymatrix{
\P_{\phi\disjoint g}\ar[r]&
\P_{\phi\disjoint k}\P_{\Id\disjoint
l}\ar[r]^{V}& \I_{\phi\disjoint \Id}\I_{\Id
\disjoint k}\i_{\Id\disjoint l}}
$$
where the arrow $V$ is induced by the map
$$
W:\P_{\phi\disjoint k}\to
\I_{\phi\disjoint\Id}\I_{\Id\disjoint l}
$$
given by the formula
$$
W=U_{\phi,k}-\sum_{k=g_2g_1,\ g_1\neq
g}F(\phi,g_1,g_2)-\xi_{\phi,k}.
$$
The induction assumption implies $W=0$,
therefore the map (\ref{delt}) vanishes as
well.

Thus, the map $X(\phi,g)$ actually passes
through
$\I_{\phi\disjoint\Id}\delta_{\Id\disjoint
g}$ thus defining a map
$$
\xi(\phi,g):\P_{\phi\disjoint
g}\to\I_{\phi\disjoint
\Id}\delta_{\Id\disjoint g}.
$$
 This
accomplishes the definition  of
$\xi(\phi,g)$.

\subsubsection{Claim}
Define maps $$
F(\phi,\psi,g_1,g_2):\P_{\phi\psi\times
g}\to \P_{\phi\times g_1}\P_{\psi\times g_2}
\to \I_{\phi}\delta_{g_1}\iI_{\psi\times
g_2}\to \I_\phi\iI_{\psi\times g};
$$
Let
$$
a(\phi,\psi,g):\P_{\phi\psi\times
g}\to\I_\phi\iI_{\psi\times g}
$$
be
 as in Sec \ref{Notat}.

 \begin{Claim}
 $$a(\phi,\psi,g)=\sum_{g=g_2g_1}F(\phi,\psi,g_1,g_2).
 $$
 \end{Claim}

\pf 1) If $\psi$ is bijective, then the
statement follows directly from the
Definition-Proposition.

2) For an arbitrary $\psi$, let
$$
D(\phi,\psi,g):=a(\phi,\psi,g)-\sum_{g=g_1g_2}F(\phi,\psi,g_1,g_2).
 $$
be the difference. It then suffices to show
that the composition
$$
\xymatrix{ \P_{\phi\psi\times g}\ar[r]^D &
\I_{\phi\disjoint\Id}\iI_{\psi\disjoint
g}\ar[r]&
\I_{\phi\disjoint\Id}\I_{\Id\disjoint
g}\i_{\psi\disjoint \Id}}
$$
vanishes, in virtue of injectivity of the
map
$$
\iI_{\psi\disjoint g}\to \I_{\Id\disjoint g}
\i_{\psi\disjoint \Id}.$$

 We have the following facts.

I: The diagram $$\xymatrix{
\P_{\psi\phi\disjoint
g}\ar[r]^{a(\phi,\psi,g)}\ar[d] &
\I_{\phi\disjoint \Id}\iI_{\psi\disjoint
g}\ar[r]& \I_{\phi\disjoint
\Id}\I_{\Id\disjoint g}
\i_{\Id\disjoint\psi}\\
\P_{\phi\disjoint g}\P_{\psi\disjoint \Id
}\ar[r]&\I_{\phi\disjoint
g}\i_{\psi\disjoint\Id}\ar[ur]&&}
$$
 commutes.

 II. The following diagram is commutative:

$$\xymatrix{
 \P_{\psi\phi\disjoint
g}\ar[rr]^{F(\phi,\psi,g_1,g_2)}\ar[d]&&
\I_{\phi\disjoint \Id}\iI_{\psi\disjoint
g}\ar[r]& \I_{\phi\disjoint \Id}\I_{\Id\disjoint g}
\i_{\psi\disjoint \Id}\\
\P_{\phi\disjoint g_1}
\P_{g_2}\P_{\psi\disjoint \Id}\ar[r]&
\I_\phi\delta_{\Id\disjoint
g_1}\I_{\Id\disjoint g_2}\i_{\psi\disjoint
\Id}\ar[rru]&&}
$$

Using $I,II$ we see that the statement
follows from the case when $\psi=\Id$.

%






\subsection{Claim}
Introduce a terminology. Let $g:A\to B$ be a
surjection.  Let $e$ be an equivalence
relation on $A$ determined by $g$. {\em A
decomposition $g=g_k\cdots g_2g_1$} is by
definition a diagram
$$
S\stackrel{g_1}\to S/e_1\stackrel{g_2}\to
S/e_2\cdots \stackrel{g_{k-1}}\to
S/e_{k-1}\stackrel{g_k}\to T,
$$
where $$e_1\ge e_2\ge\cdots\ge e_{k-1}\ge  e
$$
are equivalence relations on $S$ and $q_i$
are natural maps.

 Let $g=g_2g_1$ be a decomposition.
Define a map
$$
Y(\phi,\psi,g_1,g_2):\P_{\psi\phi\disjoint
g}\to \P_{\phi\disjoint
g_1}\P_{\psi\disjoint g_2}\to
\I_{\phi\disjoint \Id}\delta_{\Id\disjoint
g_1} \I_{\psi\disjoint
\Id}\delta_{\Id\disjoint g_2}\to
\I_{\phi\disjoint \Id}\I_{\psi\disjoint \Id}
\delta_{\Id\disjoint g}.
$$
Set
$$
Z(\phi,\psi,g)=\sum\limits_{g=g_2g_1}Y(\phi,\psi,g_1,g_2)
$$
\begin{Claim} The map $Z(\phi,\psi,g)$
coincides  with the composition
$$
\xymatrix{ \P_{\psi\phi\disjoint g}\ar[r] &
\I_{\psi\phi\disjoint\Id}\delta_{\Id\disjoint
g}\ar[r]& \I_{\phi\disjoint
\Id}\i_{\psi\disjoint
\Id}\delta_{\Id\disjoint g} }
$$
\end{Claim}
\pf Denote this composition by
$W(\phi,\psi,g)$. We shall use induction in
$|g|$.
 Let $g=g_2g_1$. Define a map
$$
Z(\phi,\psi,g_1,g_2):\P_{\phi\psi\times
g}\to \I_{\phi}\i_\psi\I_g
$$
as follows:
$$
\xymatrix{
Z(\phi,\psi,g_1,g_2):\P_{\psi\phi\disjoint
g}\ar[r]& \P_{\psi\phi\disjoint
g_1}\P_{\Id\disjoint
g_2}\ar[rr]^{Z(\phi,\psi,g_1)}&&
\I_{\phi\disjoint \Id
}\i_{\psi\disjoint\Id}\delta_{\Id\disjoint
g_1}\I_{\Id\disjoint g_2} \ar[r]&
\I_{\phi\disjoint \Id} \i_{\psi\disjoint
\Id}\I_{\Id \disjoint g}.}
$$
Define maps $W(\phi,\psi,g_1,g_2)$ in the
similar way (using $W(\phi,\psi,g)$ instead
of $Z(\phi,\psi,g)$. By the induction
assumption,
$$
Z(\phi,\psi,g_1,g_2)=W(\phi,\psi,g_1,g_2)$$
whenever $g_2g_1=g$ and $g_1\neq g$.
Therefore, it suffices to show that
$$
\sum\limits_{g_2g_1=g}Z(\phi,\psi,g_1,g_2)=
\sum\limits_{g_2g_1=g}W(\phi,\psi,g_1,g_2).
$$
Let $L$ be the sum on the LHS and $R$ be the
sum on the RHS. It follows that $L$ equals
the sum, over all decompositions
$g=g_3g_2g_1$, of the following maps

\begin{eqnarray*}
\P_{\psi\phi\disjoint g}\to
\I_{\phi\disjoint \Id}\delta_{\Id\disjoint
g_1}\P_{\psi\disjoint g_3g_2}\stackrel\alpha\to\\
\I_{\phi\disjoint \Id}\delta_{\Id\disjoint
g_1} \P_{\psi\disjoint g_2}\I_{\Id\disjoint
g_3}\to\I_{\phi\disjoint
\Id}\delta_{\Id\disjoint g_1}
\I_{\psi\disjoint \Id}\delta_{\Id\disjoint
g_2} \I_{\Id\disjoint g_3}\stackrel\omega\to\\
\I_{\phi\disjoint \Id}\delta_{\Id\disjoint
g_1} \I_{\psi\disjoint \Id} \I_{\Id\disjoint
g_3g_2}\to
 \I_{\phi\disjoint \Id}
\I_{\psi\disjoint \Id} \delta_{\Id\disjoint
g_1}\I_{\Id\disjoint g_3g_2}\to\\
\I_{\phi\disjoint \Id} \i_{\psi\disjoint
\Id}\I_{\Id\disjoint g}
\end{eqnarray*}

Fix $g_1$ and set $g^2=g_3g_2$. The previous
claim implies that the sum  of the
compositions of arrows from $\alpha$ to
$\omega$, over all decompositions
$g^2=g_3g_2$, equals the following
composition:
$$
\I_{\phi\disjoint \Id}\delta_{\Id\disjoint
g_1}\P_{\psi\disjoint g^2}\to
\I_{\phi\disjoint \Id}\delta_{\Id\disjoint
g_1}\I_{\psi\disjoint \Id}\I_{\Id\disjoint
g^2}.
$$
Therefore, $L$ equals the sum over all
decompositions $g=g_2g_1$ of the following
maps:
$$
\P_{\psi\phi\disjoint g}\to
\P_{\phi\disjoint g_1}\P_{\psi\disjoint g_2}
\to
\I_{\phi\disjoint\Id}\delta_{\Id\disjoint
g_1}\I_{\psi\disjoint g_2}\to
\I_{\phi\disjoint \Id}\I_{\psi\disjoint
g}\to \I_{\phi\disjoint
\Id}\I_{\psi\disjoint \Id}\I_{\Id\disjoint
g}\to \I_{\phi\disjoint
\Id}\i_{\psi\disjoint \Id}\I_{\Id\disjoint
g}.
$$
This can be rewritten as follows:
$$
\P_{\psi\phi\disjoint g}\to
\P_{\phi\disjoint g_1}\P_{\psi\disjoint g_2}
\to \I_{\phi\disjoint \Id}
\delta_{\Id\disjoint g_1}\I_{\psi\disjoint
g_2}\to \I_{\phi\disjoint
\Id}\I_{\psi\disjoint g}\to
\I_{\phi\disjoint \Id}\iI_{\psi\disjoint
g}\to \I_{\phi\disjoint
\Id}\i_{\psi\disjoint \Id}\I_{\Id\disjoint
g}.
$$
According to the previous statement, the sum
of these maps equals the following
composition:
$$
\P_{\psi\phi\disjoint g}\to
\I_{\psi\phi\disjoint g} \to
\I_{\psi\disjoint\Id}\iI_{\phi\disjoint
g}\to \I_{\phi\disjoint
\Id}\i_{\psi\disjoint \Id}\I_{\Id\disjoint
g}.
$$
This composition, in turn, is equal to:
$$
\P_{\psi\phi\disjoint g}\to
\I_{\psi\phi\disjoint g}
\to\I_{\phi\disjoint \Id}\I_{\psi\disjoint
\Id }\I_{\Id\disjoint g}\to
\I_{\phi\disjoint \Id}\i_{\psi\disjoint \Id}
\I_{\Id\disjoint g}.
$$
It easily follows that this sum equals $R$.
This completes the proof.
\subsubsection{Compatibility with $\boxtimes$}
\begin{Claim} The following diagram is
commutative:
$$
\xymatrix{ \P_{f_1\disjoint
g_1}(M_1)\boxtimes \P_{f_2\disjoint
g_2}(M_2)\ar[r]\ar[d]&\I_{f_1\disjoint
\Id}\delta_{Id\disjoint g_1}(M_1)\boxtimes
\I_{f_2\disjoint \Id}\delta_{\Id\disjoint
g_2}(M_2)\ar[r]&\I_{f_1\disjoint
f_2\disjoint\Id}\delta_{\Id\disjoint
g_1\disjoint g_2}(M_1\boxtimes M_2)\\
\P_{f_1\disjoint f_2\disjoint g_1\disjoint
g_2}(M_1\boxtimes M_2)\ar[rru]&&\\ }
$$
\end{Claim}\pf  We shall use induction.
 The composition

\begin{eqnarray*}
\P_{f_1\disjoint g_1}(M_1)\boxtimes
\P_{f_2\disjoint g_2}(M_2)\to
\P_{f_1\disjoint f_2\disjoint g_1\disjoint
g_2}(M_1\boxtimes M_2)\to\\
\I_{f_1\disjoint
f_2\disjoint\Id}\delta_{\Id\disjoint
g_1\disjoint g_2}(M_1\boxtimes M_2)\to
\I_{f_1\disjoint
f_2\disjoint\Id}\I_{\Id\disjoint
g_1\disjoint g_2}(M_1\boxtimes M_2)
\end{eqnarray*}
equals the negative of the sum over all
decompositions $g_1=h_2h_1$, $g_2=h_4h_3$,
$(h_2\disjoint h_4\neq \Id)$ of the
following maps:
\begin{eqnarray*}
\P_{f_1\disjoint g_1}(M_1)\boxtimes
\P_{f_2\disjoint g_2}(M_2)\to
\P_{f_1\disjoint f_2\disjoint g_1\disjoint
g_2}(M_1\boxtimes M_2)\to\\
\P_{f_1\disjoint f_2\disjoint h_1\disjoint
h_2 }\I_{\Id\disjoint h_3\disjoint
h_4}(M_1\boxtimes M_2)\to \I_{f_1\disjoint
f_2\disjoint\Id}\delta_{\Id\disjoint
h_1\disjoint h_2 }\I_{\Id\disjoint
h_3\disjoint h_4}(M_1\boxtimes M_2)\to
\I_{f_1\disjoint
f_2\disjoint\Id}\I_{\Id\disjoint
g_1\disjoint g_2}(M_1\boxtimes M_2),
\end{eqnarray*}

which is  (due to the induction assumption)
the same as:

\begin{eqnarray*}
\P_{f_1\disjoint g_1}(M_1)\boxtimes
\P_{f_2\disjoint g_2}(M_2)\to
\P_{f_1\disjoint h_1}\I_{h_3}(M_1)\boxtimes
\P_{f_2\disjoint h_2}\I_{h_4}(M_2)\to\\
\P_{f_1\disjoint f_2\disjoint h_1\disjoint
h_2}\I_{\Id\disjoint h_3\disjoint
h_4}(M_1\boxtimes M_2)\to \I_{f_1\disjoint
f_2\disjoint\Id}\delta_{\Id\disjoint
h_1\disjoint h_2 }\I_{\Id\disjoint
h_3\disjoint h_4}(M_1\boxtimes M_2)\to
\I_{f_1\disjoint
f_2\disjoint\Id}\I_{\Id\disjoint
g_1\disjoint g_2}(M_1\boxtimes M_2),
\end{eqnarray*}
which, in turn, equals:
\begin{eqnarray*}
\P_{f_1\disjoint g_1}(M_1)\boxtimes
\P_{f_2\disjoint g_2}(M_2)\to
\P_{f_1\disjoint h_1}\I_{\Id\disjoint
h_3}(M_1)\boxtimes \P_{f_2\disjoint
h_2}\I_{\Id \disjoint h_4}(M_2)\to\\
\I_{f_1\disjoint \Id}\delta_{\Id\disjoint
h_1}\I_{\Id\disjoint h_3}(M_1)\boxtimes
\I_{f_2\disjoint\Id}\delta_{\Id\disjoint
h_2} \I_{\Id\disjoint h_4}(M_2)\to
\I_{f_1\disjoint \Id} \I_{\Id\disjoint
g_1}(M_1) \boxtimes \I_{f_2\disjoint
\Id}\I_{\Id\disjoint g_2}(M_2)\to\\
\P_{f_1\disjoint f_2\disjoint \Id}\I_{
\Id\disjoint g_1\disjoint g_2}(M_1\boxtimes
M_2).
\end{eqnarray*}

The sum of all such maps over all
decompositions $g_1=h_3h_1,g_2=h_4h_2$, is
zero. Therefore, the negative of the sum
over all decompositions with $h_3\disjoint
h_4\neq \Id$ is equal to the map in which
$h_1=g_1$, $h_2=g_2$, $h_3=\Id$, $h_4=\Id$,
which immediately implies the commutativity
of the diagram in question.
\endpf

\subsection{Maps
$c(\phi,g):\P_{\phi\disjoint g}\to
\P_{\phi\disjoint \Id}\delta_{\Id\disjoint
g}$} First, we define maps
$$
\xi(\phi_1,\phi_2,\ldots,\phi_n,g):
\P_{\phi\disjoint  g}\to \I_{\phi_1\disjoint
\Id }\I_{\phi_2\disjoint \Id}\cdots
\I_{\phi_n\disjoint \Id}\delta_{\Id\disjoint
g},
$$
where $\phi=\phi_n\phi_{n-1}\cdots \phi_1$,
as the sum over all decompositions
$g=g_ng_{n-1}\cdots g_1$ of the maps:
\begin{eqnarray*}
\P_{\phi\disjoint g}\to \P_{\phi_1\disjoint
g_1}\P_{\phi_2\disjoint
g_2}\cdots\P_{\phi_n\disjoint g_n}\\
\to
\I_{\phi_1\disjoint\Id} \delta_{\Id\disjoint
g_1}\I_{\phi_2\disjoint\Id}
\delta_{\Id\disjoint g_2}\cdots
\I_{\phi_n\disjoint\Id} \delta_{\Id\disjoint
g_n}\to \I_{\phi_1\disjoint \Id
}\I_{\phi_2\disjoint \Id
}\cdots\I_{\phi_n\disjoint
\Id}\delta_{\Id\disjoint g}.
\end{eqnarray*}

The previous claim implies that the
collection of maps
$\xi(\phi_1,\phi_2,\ldots,\phi_n,g)$ for all
decompositions
$\phi=\phi_n\phi_{n-1}\cdots\phi_1$ gives
rise to a map
$$
c(\phi, g):\P_{\phi\disjoint g}\to
\P_{\phi\disjoint\Id}\delta_{\Id\disjoint
g}.
$$
\subsubsection{Claim}
\begin{Claim}\label{decp}
The composition
$$\xymatrix{
\P_{\psi\phi\disjoint g}\ar[r] &
\P_{\psi\phi\disjoint\Id}\delta_{\Id\disjoint
g}\ar[r]&
 \P_{\phi\disjoint \Id}
 \P_{\psi\disjoint\Id}\delta_{\Id\disjoint g}}$$
 is equal to the sum, over all decompositions
 $g=g_2g_1$, of the maps
$$\xymatrix{
\P_{\phi\psi\disjoint
g}\ar[r]&\P_{\phi\disjoint
g_1}\P_{\psi\disjoint g_2}\ar[r]&
\P_{\phi\disjoint\Id}\delta_{\Id\disjoint
g_1} \P_{
\psi\disjoint\Id}\delta_{\Id\disjoint g_2}
\ar[r]& \P_{\phi\disjoint
\Id}\P_{\psi\disjoint\Id}\delta_{\Id\disjoint
g}.}
$$
\end{Claim}
\pf Clear.
\subsection{Composition}\label{cmp}
\begin{Claim}
The following diagram  is commutative:
$$\xymatrix{
\P_{\phi\disjoint g\disjoint h}\ar[r]\ar[d]&
\P_{\phi\disjoint\Id\disjoint\Id}\delta_{\Id\disjoint
g\disjoint h}\\
\P_{\phi\disjoint
g\disjoint\Id}\delta_{\Id\disjoint\Id\disjoint
h}\ar[ur]&}
$$
\end{Claim}
\pf First, let us prove that the diagram
\begin{equation}\label{what}
\xymatrix{ \P_{\phi\disjoint g\disjoint
h}\ar[r]\ar[d]&
\I_{\phi\disjoint\Id\disjoint\Id}\delta_{\Id\disjoint
g\disjoint h}\\
\P_{\phi\disjoint g\disjoint \Id
}\delta_{\Id\disjoint\Id\disjoint
h}\ar[ur]\ar[ur]&}
\end{equation}
is commutative. Denote the composition
$$
\xymatrix{ \P_{\phi\disjoint g\disjoint
h}\ar[r]&\P_{\phi\disjoint
g\disjoint\Id}\delta_{\Id\disjoint
\Id\disjoint h}\ar[r]&
\I_{\phi\disjoint\Id\disjoint\Id}\delta_{\Id\disjoint
g\disjoint h}\ar[r]&
\I_{\phi\disjoint\Id\disjoint\Id}\I_{\Id\disjoint
g \disjoint
\Id}\delta_{\Id\disjoint\Id\disjoint h} }
$$
by $U(\phi,g,h)$.

Let $g=g_2g_1$ be a decomposition. Define a
map
$$
U(\phi,g_1,g_2,h):\P_{\phi\disjoint
g\disjoint
h}\to\I_{\phi\disjoint\Id\disjoint\Id}\I_{\Id\disjoint
g \disjoint
\Id}\delta_{\Id\disjoint\Id\disjoint h}
$$
as the following composition:
\begin{eqnarray*}
U(\phi,g_1,g_2,h):\P_{\phi\disjoint
g\disjoint h}\to \P_{\phi\disjoint
g\disjoint\Id}\delta_{\Id\disjoint\Id\disjoint
h}\to\\
 \P_{\phi\disjoint
g_1\disjoint\Id}\P_{\Id\disjoint
g_2\disjoint\Id}\delta_{\Id\disjoint\Id\disjoint
h}\to
\I_{\phi\disjoint\Id\disjoint\Id}\delta_{\Id\disjoint
g_1\disjoint\Id}\I_{\Id\disjoint
g_2\disjoint\Id}\delta_{\Id\disjoint\Id\disjoint
h}\to\\
\I_{\phi\disjoint\Id\disjoint\Id}\I_{\Id\disjoint
g\disjoint\Id}\delta_{\Id\disjoint\Id\disjoint
h.}
\end{eqnarray*}

 Let also
$$
A(\phi,g,h):\P_{\phi\disjoint g\disjoint
h}\to \I_{\phi\disjoint
g\disjoint\Id}\delta_{\Id\disjoint\Id\disjoint
h}\to
\I_{\phi\disjoint\Id\disjoint\Id}\I_{\Id\disjoint
g\disjoint \Id}\delta_{\Id\disjoint
\Id\disjoint h}.
$$

Then, by definition,
\begin{equation}\label{AU}
U(\phi,g,h)=A(\phi,g,h)-
\sum\limits_{g=g_2g_1,g_1\neq g
}U(\phi,g_1,g_2,h).
\end{equation}

The map $U(\phi,g_1,g_2,h)$ equals, in turn,
the sum over all decompositions $h=h_2h_1$
of the maps:
\begin{eqnarray*}
U(\phi,g_1,g_2,h_1,h_2):\P_{\phi\disjoint
g\disjoint h}\to\\
 \P_{\phi\disjoint
g_1\disjoint h_1}\P_{\Id\disjoint
g_2\disjoint
h_2}\stackrel{U(\phi,g_1,h_1)}\longrightarrow
\I_{\phi\disjoint\Id\disjoint\Id}
\delta_{\Id\disjoint
g_1\disjoint\Id}\delta_{\Id\disjoint\Id\disjoint
h_1}\I_{\Id\disjoint g_2\disjoint \Id
}\delta_{\Id\disjoint \Id\disjoint h_2}\to\\
\I_{\phi\disjoint\Id\disjoint \Id
}\I_{\Id\disjoint g\disjoint
\Id}\delta_{\Id\disjoint \Id\disjoint h}
\end{eqnarray*}

Then
\begin{equation}
U(\phi,g,h)=A(\phi,g,h)-
\sum\limits_{g=g_2g_1,h=h_2h_1,g_1\neq g
}U(\phi,g_1,g_2,h_1,h_2).
\end{equation}

Set
$$
U(\phi,g,h_1,h_2):=U(\phi,g,\Id,h_1,h_2):\P_{\phi\disjoint
g\disjoint h}\to
\I_{\phi\disjoint\Id\disjoint
\Id}\I_{\Id\disjoint g\disjoint
\Id}\I_{\Id\disjoint\Id\disjoint  h}
$$
to be \begin{eqnarray*}
 \P_{\phi\disjoint g\disjoint h}\to
\P_{\phi\disjoint g\disjoint h_1
}\P_{\Id\disjoint\Id\disjoint
h_2}\stackrel{U(\phi,g,h_1)}\longrightarrow\\
\I_{\phi\disjoint\Id\disjoint \Id
}\I_{\Id\disjoint g\disjoint \Id}
\delta_{\Id\disjoint\Id\disjoint
h_1}\I_{\Id\disjoint\Id\disjoint h_2}\to
\I_{\phi\disjoint\Id\disjoint
\Id}\I_{\Id\disjoint g\disjoint
\Id}\I_{\Id\disjoint\Id \disjoint h}
\end{eqnarray*}
Similarly, let
\begin{eqnarray*}
A(\phi,g,h_1,h_2):\P_{\phi\disjoint
g\disjoint h}\to \P_{\phi\disjoint
g\disjoint h_1 }\P_{\Id\disjoint\Id\disjoint
h_2}\stackrel{A(\phi,g,h_1)}\longrightarrow\\
\I_{\phi\disjoint\Id\disjoint \Id
}\I_{\Id\disjoint g\disjoint \Id}
\delta_{\Id\disjoint\Id\disjoint
h_1}\I_{\Id\disjoint \Id\disjoint h_2}\to
\I_{\phi\disjoint\Id\disjoint\Id}\I_{\Id\disjoint
g\disjoint \Id}\I_{\Id\disjoint \Id
\disjoint h}
\end{eqnarray*}
and
\begin{eqnarray*}
U(\phi,g_1,g_2,h_1,h_2,h_3):\P_{\phi\disjoint
g\disjoint h}\to \P_{\phi\disjoint
g\disjoint h_2h_1}
\P_{\Id\disjoint\Id\disjoint
\h_3}\stackrel{U(\phi,g_1,g_2,h_1,h_2)}\longrightarrow\\
\I_{\phi\disjoint\Id\disjoint
\Id}\I_{\Id\disjoint g\disjoint \Id}
\delta_{\Id\disjoint\Id\disjoint
h_2h_1}\I_{\Id\disjoint \Id\disjoint h_3}\to
\I_{\phi\disjoint \Id\disjoint \Id} \I_{\Id
\disjoint g\disjoint \Id}\I_{\Id\disjoint
\Id\disjoint h}
\end{eqnarray*}

The equation (\ref{AU}) implies that
\begin{equation}\label{AUU}
\sum\limits_{h=h_2h_1}U(\phi,g,h_1,h_2)=
\sum\limits_{h=h_2h_1}A(\phi,g,h_1,h_2)-
\sum\limits_{g=g_2g_1,g_1\neq g;h=h_3h_2h_1}
U(\phi,g_1,g_2,h_3,h_2,h_1).
\end{equation}

The map
$$
\sum_{h=h_2h_1}A(\phi,g,h_1,h_2)
$$ is equal to the following one:
$$
X(\phi,g,h):\P_{\phi\disjoint g\disjoint
h}\to \I_{\phi\disjoint g\disjoint h }\to
\I_{\phi\disjoint\Id\disjoint\Id}\I_{\Id\disjoint
g \disjoint \Id}\I_{\Id\disjoint\Id\disjoint
h}.
$$

The map
$$
Y(\phi,g_1,g_2,h_1,h_2):=\sum_{h^2=h_3h_2}U(\phi,g_1,g_2,h_1,h_2,h_3)
$$
equals: \begin{eqnarray*} \P_{\phi\disjoint
g\disjoint h}\to \P_{\phi\disjoint
g_1\disjoint
h_1}\P_{\Id\disjoint\g_2\disjoint
h^2}\stackrel{U(\phi,g_1,h_1)}\to \\
\I_{\phi\disjoint\Id\disjoint\Id}\delta_{\Id\disjoint
g_1\disjoint\Id}\delta_{\Id\disjoint\Id\disjoint
h_1}\I_{\Id\disjoint g_2\times h_2 }\to
\I_{\phi\disjoint\Id\disjoint\Id}\I_{\Id\disjoint
g\disjoint \Id}\I_{\Id\disjoint\Id\disjoint
h}
\end{eqnarray*}
Therefore,
$$Y(\phi,g_1,g_2,h_1,h_2)=U(\phi,g_1,g_2,h_1,h_2).$$
The equation (\ref{AUU}) can be now
rewritten as:
\begin{equation}\label{AUUU}
\sum_{h=h_1h_2}U(\phi,g,h_1,h_2)=X(\phi,g,h)-
\sum\limits_{g=g_1g_2,h=h_1h_2,g_1\neq g
}U(\phi,g_1,g_2,h_1,h_2).
\end{equation}
Note that
$$
U(\phi,g,h_1,h_2)=U(\phi,g,\Id,h_1,h_2).
$$
Therefore (\ref{AUUU}) implies that
$$
U(\phi,g,h,\Id)=X(\phi,g,h)-\sum_{g=g_2g_1,h=h_2h_1,g_1\disjoint
h_1\neq g_2\disjoint h_2
}U(\phi,g_1,g_2,h_1,h_2).
$$
 The induction assumption implies that
$U(\phi,g_1,h_1)=c(\phi,g_1\disjoint h_1)$
if $g_1\disjoint h_1\neq g\disjoint h$. This
implies that the right hand side equals the
composition:
$$\xymatrix{
\P_{\phi\disjoint g\disjoint
h}\ar[r]^{c(\phi,g\disjoint
h)}&\I_{\phi\disjoint\Id\disjoint\Id}\delta_{\Id\disjoint
g\disjoint
h}\ar[r]&\I_{\phi\disjoint\Id\disjoint\Id}
\I_{\Id\disjoint g\disjoint
\Id}\I_{\Id\disjoint \Id\disjoint h}.}
$$

By definition, the left hand side  equals
the composition:
\begin{eqnarray*}
\P_{\phi\disjoint g\disjoint
h}\to\P_{\phi\disjoint
g\disjoint\Id}\delta_{\Id\disjoint\Id\disjoint
h}\to\\
\I_{\phi\disjoint\Id\disjoint\Id}\delta_{\Id\disjoint
g\disjoint \Id}\delta_{\Id\disjoint
\Id\disjoint h}\to
\I_{\phi\disjoint\Id\disjoint\Id}\I_{\Id\disjoint
g\disjoint \Id}\I_{\Id\disjoint \Id\disjoint
h}
\end{eqnarray*}

Therefore, the diagram (\ref{what}) is
commutative. The  original statement  can be
now proven straightforwardly using
\ref{decp}).
\endpf
\subsection{Compositions $\P\delta\P\to
\P\to \I\delta$}
\begin{Claim}

 The composition
 \begin{equation}\label{xx}
 \P_{\phi_1\disjoint
 g_1}\delta_{\phi_2\disjoint
 g_2}\P_{\phi_3\disjoint g_3}\to
 \P_{\phi\disjoint
 g}\to \I_{\phi\disjoint\Id}\delta_{\Id\disjoint
 g}
\end{equation}
vanishes if $\phi_1\disjoint g_1\neq  \Id$
and: $\phi_1\neq \phi$ or $g_3\neq \Id$. In
the cases when it does not vanish we have
the following rules:

2. In the case $\phi_1\disjoint g_1=\Id$,
this composition equals:
$$
\delta_{\phi_2\disjoint
 g_2}\P_{\phi_3\disjoint g_3}\to
 \delta_{\phi_2\disjoint
 g_2}\I_{\phi_3\disjoint \Id}\delta_{\Id\disjoint g_3}\to
 \I_{\phi\disjoint\Id}\delta_{\Id\disjoint g}.
 $$

 3. In the case $\phi_1=\phi$, $g_3=\Id$, this
 composition equals $-A$, where
 $$
A:\P_{\phi\disjoint\g_1}\delta_{\Id\disjoint
g_2}\to \I_{\phi\disjoint \Id}
\delta_{\Id\disjoint
g_1}\delta_{\Id\disjoint g_2}\to
\I_{\phi\disjoint \Id}\delta_{\Id \disjoint
g}.
$$
\end{Claim}
\pf

We shall use induction in $g$. Compute the
composition
\begin{equation}\label{dF}\xymatrix{
 \P_{\phi_1\disjoint
 g_1}\delta_{\phi_2\disjoint
 g_2}\P_{\phi_3\disjoint
 g_3}\ar[r]&\P_{\phi\disjoint
 g}\ar[r]^{F(\phi,g^1,g^2)}
 &\I_{\phi\disjoint \Id}\I_{\Id\disjoint g}}
\end{equation}
Where $F(\phi,g^1,g^2)$ is as in
(\ref{Fg1g2}). This composition is equal to:
$$
\xymatrix{
 \P_{\phi_1\disjoint
 g_1}\delta_{\phi_2\disjoint
 g_2}\P_{\phi_3\disjoint
 g_3}\ar[r]&\P_{\phi\disjoint
 g}\ar[r] &\P_{\phi\disjoint g^1}\P_{\Id\disjoint
 g^2}\ar[r]
 &\I_{\phi\disjoint \Id}\delta_{\Id\disjoint g^1}
 \I_{\Id\disjoint g^2}\ar[r]&
 \I_{\phi\disjoint \Id}\I_{\Id\disjoint g}}
$$

This composition does not vanish only if

A:$g_3=g^2u$, $g^1=ug_2g_1$;

B:$\phi_1=\phi$; $g^1=g_1$, $g^2=g_3g_2$.

Consider several cases.

1. $\phi_1\neq \phi$ and $\phi_1\times
g_1\neq \Id$.  The induction assumption
implies that the composition (\ref{dF})
vanishes whenever $g_2\neq \Id$. Therefore,
the composition

$$
\xymatrix{
 \P_{\phi_1\disjoint
 g_1}\delta_{\phi_2\disjoint
 g_2}\P_{\phi_3\disjoint
 g_3}\ar[r]&\P_{\phi\disjoint
 g}\ar[r] &\I_{\phi\disjoint\Id}\delta_{\Id\disjoint g}
 \ar[r]&\I_{\phi\disjoint\Id}\I_{\Id \disjoint g}}
$$
equals
$$
\xymatrix{
 \P_{\phi_1\disjoint
 g_1}\delta_{\phi_2\disjoint
 g_2}\P_{\phi_3\disjoint
 g_3}\ar[r]&\I_{\phi\disjoint
 g}\ar[r]
 &\I_{\phi\disjoint \Id}\I_{\Id\disjoint g}}
$$

This composition   vanishes  because
$\phi_1\disjoint g_1 \neq \Id$.

Case 2. $\phi_1\disjoint g_1=\Id$, $g\neq
\Id$. The case $B$ is again excluded.  By
the inductive assumption, the composition
(\ref{dF}) equals
$$\xymatrix{
A(u,g^2):\delta_{\phi_2\disjoint
g_2}\P_{\phi_3\disjoint
g_3}\ar[r]^{F(\phi_3,u,g^2)}&
\delta_{\phi_2\disjoint
g_2}\I_{\phi_3}I_{\Id\disjoint g^2u}\ar[r]&
\I_{\phi\disjoint \Id}\I_{\Id\disjoint g}.}
$$

The composition
$$\xymatrix{
\delta_{\phi_2\disjoint
g_2}\P_{\phi_3\disjoint g_3}\ar[r]&
\I_{\phi\disjoint  g}\ar[r]&
\I_{\phi\disjoint\Id}\I_{\Id\disjoint g}}
$$
equals:
$$
B:\xymatrix{ \delta_{\phi_2\disjoint
g_2}\P_{\phi_3\disjoint g_3}\ar[r]&
\delta_{\phi_2\disjoint g_3}
\I_{\phi_3\disjoint\Id}\I_{\Id\disjoint
g_3}\ar[r]& \I_{\phi\disjoint\Id}\I_{\Id
\disjoint g}}
$$
 Therefore, the composition
$$\xymatrix{
\delta_{\phi_2\disjoint
g_2}\P_{\phi_3\disjoint g_3}\ar[r]&
\P_{\phi\disjoint g}\ar[r]&
\I_{\phi\disjoint\Id}\delta_{\Id\disjoint
g}\ar[r]&\I_{\phi\disjoint
\Id}\I_{\Id\disjoint g}}
$$
equals
$$B-\sum_{g_3=g^2u,g^2\neq \Id}A(u,g^2).$$
It follows that this composition equals:
$$
\xymatrix{ \delta_{\phi_2\disjoint
g_2}\P_{\phi_3\disjoint
g_3}\ar[r]&\delta_{\phi_2\disjoint g_2}
\I_{\phi_3\disjoint \Id}\delta_{\Id\disjoint
g_3}\ar[r]&
\I_{\phi\disjoint\Id}\I_{g\disjoint \Id}}
$$
which is what is predicted by 2.

We have the last remaining case
$\phi_1=\phi$,  $g\neq \Id$, where we have
to add contributions from A and B.

Then the contribution from $A$ is equal to
zero if $g_3=\Id$. Otherwise, according to
the inductive assumption, it equals to $-C$,
where

\begin{eqnarray*} C:\P_{\phi\disjoint
g_1}\delta_{\Id\disjoint
g_2}\P_{\Id\disjoint g_3}\to
\I_{\phi\disjoint\Id}\delta_{\Id\disjoint
g_1}\delta_{\Id\disjoint
g_2}\I_{\Id\disjoint g_3}\to
\I_{\phi\disjoint \Id}\I_{\Id\disjoint g}.
\end{eqnarray*}

 We see that the contribution from $B$
 equals $C$. Therefore, the composition
 (\ref{dF}) is zero if $g_3\neq \Id$ and $C$
 otherwise. Note that the composition
 $$
\xymatrix{ C:\P_{\phi\disjoint
g_1}\delta_{\Id\disjoint
g_2}\P_{\Id\disjoint
g_3}\ar[r]&\I_{\phi\disjoint  g }}
$$
is always zero. Therefore,  the map
(\ref{xx}) is zero if $g_3\neq \Id$ and $-C$
otherwise. This completes the proof. \endpf

\subsection{ Compositions $\P\delta\P\to \P\to \P\delta$}
\label{propc}
\begin{Claim}
Consider the composition
\begin{equation}\label{f2d}
\P_{\phi_1\disjoint
g_1}\delta_{\phi_2\disjoint
g_2}\P_{\phi_3\disjoint g_3}\to
\P_{\phi\disjoint g}\to
\P_{\phi\disjoint\Id}\delta_{\Id\disjoint
g}.
\end{equation}

If $\phi_2\neq \Id$, this composition is
equal to the following composition:
$$
\P_{\phi_1\disjoint
g_1}\delta_{\phi_2\disjoint
g_2}\P_{\phi_3\disjoint g_3}\to
\P_{\phi_1\disjoint\Id}\delta_{\Id\disjoint
g_1}\delta_{\phi_2\disjoint
g_2}\P_{\phi_3\disjoint\Id}\delta_{\Id\disjoint
g_3}\to
\P_{\phi_1\disjoint\Id}\delta_{\phi_2\disjoint\Id}
\P_{\phi_3\disjoint\Id}\delta_{\Id\disjoint
g}\to \P_{\phi\disjoint
\Id}\delta_{\Id\disjoint g}.
$$

Otherwise, this composition  is equal to
zero except the following cases:

a) $\phi_1\times g_1=\Id$, in which case our
composition equals
\begin{equation}\label{f3d}
\delta_{\Id\disjoint g_2}\P_{\phi\disjoint
g_3}\to \delta_{\Id\disjoint
g_2}\P_{\phi\disjoint\Id}\delta_{\Id\disjoint
g_3}\to \P_{\phi\disjoint
\Id}\delta_{\Id\disjoint g};
\end{equation}

b) $\phi_3\disjoint g_3=\Id$,  in which case
the composition is equal to
\begin{equation}\label{f4d}
-C,
\end{equation}
 where
$$
C:\P_{\phi\disjoint g_1}\delta_{\Id\disjoint
g_2}\to
\P_{\phi\disjoint\Id}\delta_{\Id\disjoint
g_1}\delta_{\Id\disjoint g_2}\to
\P_{\phi\disjoint \Id}\delta_{\Id\disjoint
g}.
$$

\end{Claim}
\pf

We will prove the statement by induction in
$\phi$.

It suffices to check that

 The composition (\ref{f2d}) coincides with
the maps (\ref{f3d}), (\ref{f4d}) after
composing each of them

1) with the map
$$
P_{\phi\disjoint \Id}\delta_{\Id\disjoint
g}\to
\I_{\phi\disjoint\Id}\delta_{\Id\disjoint g}
$$

2) with the maps
$$
P_{\phi\disjoint \Id}\delta_{\Id\disjoint
g}\to
\P_{\phi_1\disjoint\Id}\P_{\phi_2\disjoint\Id}\delta_{\Id\disjoint
g},
$$
where $\phi_1,\phi_2\neq \Id$.

1) can be checked straightforwardly:

1a) $\phi_2\neq \Id$.

If $\phi_1\neq \Id$, then both compositions
are immediately zero.

If $\phi_1=\Id$, $g_1\neq \Id$, then again
both compositions are zero (the composition
(\ref{f3d}) is zero because the
corresponding map $\P_{\phi_1\disjoint
g_1}\to
\P_{\phi_1\disjoint\Id}\delta_{\Id\disjoint
g_1}$ is zero.)

If $\phi_1\disjoint g_1=\Id$, then the two
compositions coincide.

1b) $\phi_2=\Id$. If none of $\phi_1,\phi_3$
is identity, then both compositions are
clearly  equal to zero.

If $\phi_1=\Id $ and $g_1\neq \Id$, then
both compositions are zero.

If $\phi_1\disjoint g_1=\Id$, then both
compositions do clearly coincide.

If $\phi_3=\Id$ and $g_3\neq \Id$, then both
compositions are zero.

If $\phi_3=\Id$ and $g_3=\Id$, then both
compositions coincide.

2a)  $\phi_2\neq \Id$. Compute the
composition
$$
\P_{\phi_1\disjoint
g_1}\delta_{\phi_2\disjoint g_2}
\P_{\phi_3\disjoint g_3}\to
\P_{\phi\disjoint\Id}\delta_{\Id\disjoint
g}\to
\P_{\phi^1\disjoint\Id}\P_{\phi^2\disjoint\Id}
\delta_{\Id\disjoint  g}
$$

This composition vanishes except the
following two cases:

i) $\phi_1=u\phi^1$.

ii) $\phi_3=\phi^2u$.

In both cases the coincidence is obvious.

2b) $\phi_2=\Id$,

$\phi_1,\phi_3\neq \Id$. If $\phi^1\neq
\phi_1$, both compositions are obviously
zero.

Assume $\phi_1=\phi^1$, $\phi_3=\phi^2$.
Then, according to (\ref{decp}) the
composition
$$
\P_{\phi_1\disjoint g_1}\delta_{\Id\disjoint
g_2} \P_{\phi_3\disjoint g_3}\to
\P_{\phi\disjoint\Id}\delta_{\Id\disjoint
g}\to
\P_{\phi_1\disjoint\Id}\P_{\phi_3\disjoint\Id}
\delta_{\Id\disjoint g}
$$
is equal to the sum of two maps which
annihilate each other. The second
composition is also zero.

If $\phi_1=\Id$ or $\phi_3=\Id$, then the
two compositions do clearly coincide.

This completes the proof.
\endpf
\subsubsection{Compatibility with $\boxtimes$}

\begin{Claim} The following diagram is
commutative:
$$
\xymatrix{ \P_{f_1\disjoint
g_1}(M_1)\boxtimes \P_{f_2\disjoint
g_2}(M_2)\ar[r]\ar[d]&\P_{f_1\disjoint
\Id}\delta_{Id\disjoint g_1}(M_1)\boxtimes
\P_{f_2\disjoint \Id}\delta_{Id\disjoint
g_2}(M_2)\ar[d]\\
\P_{f_1\disjoint f_2\disjoint g_1\disjoint
g_2}(M_1\boxtimes
M_2)\ar[r]&\P_{f_1\disjoint
f_2\disjoint\Id}\delta_{\Id\disjoint
g_1\disjoint g_2}(M_1\boxtimes M_2)}
$$
\end{Claim}

Direct application of formulas yields the
commutativity of the diagram

$$
\xymatrix{ \P_{f_1\disjoint
g_1}(M_1)\boxtimes \P_{f_2\disjoint
g_2}(M_2)\ar[r]\ar[d]&\P_{f_1\disjoint
\Id}\delta_{Id\disjoint g_1}(M_1)\boxtimes
\P_{f_2\disjoint \Id}\delta_{Id\disjoint
g_2}(M_2)\ar[d]\\
\P_{f_1\disjoint f_2\disjoint g_1\disjoint
g_2}(M_1\boxtimes M_2)\ar[d]&
\P_{f_1\disjoint
f_2\disjoint\Id}\delta_{\Id\disjoint
g_1\disjoint g_2}(M_1\boxtimes M_2)\ar[d]\\
\I_{u_1\disjoint v_1\disjoint
\Id}\I_{u_2\disjoint v_2\disjoint \Id
}\cdots \I_{u_n\disjoint v_n\disjoint \Id
}\delta_{\Id \disjoint g_1\disjoint
g_2}\ar[r] & \I_{u_1\disjoint v_1\disjoint
\Id}\I_{u_2\disjoint v_2\disjoint \Id}\cdots
\I_{u_n\disjoint v_n\disjoint
\Id}\delta_{\Id\disjoint g_1\disjoint g_2} }
$$

which proves the statement.

\subsection{The maps
$s(\phi,g):\R_{\phi\disjoint g}\to
\R_{\phi\disjoint \Id}\delta_{\Id\disjoint
g}$}
\subsubsection{Definition}
Define a map $$ S(\phi_1\disjoint
g_1,\phi^1\disjoint
g^1,\ldots,\phi_n\disjoint  g_n):
\P_{\phi_1\disjoint
g_1}\delta_{\phi^1\disjoint g^1}\cdots
\P_{\phi_n\disjoint g_n}\to
\P_{\phi_1\disjoint\Id}\delta_{\phi^1\disjoint\Id}\cdots
\delta_{\phi^{n-1}\disjoint\Id}
\P_{\phi_n\disjoint\Id}\delta_{\Id\disjoint
g},$$ where $g=g_ng_{n-1}\cdots g_1$, as
follows:
$$
\P_{\phi_1\disjoint
g_1}\delta_{\phi^1\disjoint g^1}\cdots
\P_{\phi_n\disjoint g_n}\to
\P_{\phi_1\disjoint\Id}\delta_{\Id\disjoint
g_1}\delta_{\phi^1\disjoint\Id}\cdots
\delta_{\phi^{n-1}\disjoint\Id}\P_{\phi_n\disjoint\Id}
\delta_{\Id\disjoint \phi^n}\to
\P_{\phi_1\disjoint\Id}\delta_{\Id\disjoint
g^1}\cdots \delta_{\Id\disjoint\phi^{n-1}}
\P_{\phi_n\disjoint\Id}\delta_{\Id\disjoint
g}.$$ Let $S'(\phi_1\disjoint
g_1,\phi^1\disjoint
g^1,\ldots,\phi_n\disjoint  g_n)=0$ if at
least one of $\phi^i$ is identity. Otherwise
set
$$
S'(\phi_1\disjoint g_1,\phi^1\disjoint
g^1,\ldots,\phi_n\disjoint
g_n)=S(\phi_1\disjoint g_1,\phi^1\disjoint
g^1,\ldots,\phi_n\disjoint g_n).
$$

The sum of all possible $S'(\phi_1\disjoint
g_1,\phi^1\disjoint
g^1,\ldots,\phi_n\disjoint g_n)$ produces a
map
$$s(\phi,g):\R_{\phi\disjoint g}\to
\R_{\phi\disjoint\Id}\delta_{\Id\disjoint
g}.
$$
 Let us study its properties.

 \subsubsection{}
Denote
$$s(g_1,\phi,g_2):\R_{\phi\disjoint
g}\to \R_{\Id\disjoint g_1}\R_{\phi\disjoint
g_2}\to \delta_{\Id\disjoint
g_1}\R_{\phi\disjoint\Id}\delta_{\Id\disjoint
g_2}\to
\R_{\phi\disjoint\Id}\delta_{\Id\disjoint
g};
$$
$$
s(\phi,g_1,g_2):\R_{\phi\disjoint g}\to
\R_{\phi\disjoint g_1}\R_{\Id\disjoint
g_2}\to
\R_{\phi\disjoint\Id}\delta_{Id\disjoint
\g_1}\delta_{\Id\disjoint g_2}\to
\R_{\phi\disjoint\Id}\delta_{\Id\disjoint
g}.
$$
\begin{Claim}\label{diffs}
$$
ds(\phi,g)=\sum\limits_{g=g_2g_1}(s(g_1,\phi,g_2)-s(\phi,g_1,g_2))
$$
\end{Claim}
\pf Follows directly from \ref{propc}.
\endpf
\subsubsection{}
\begin{Claim}
The following diagram is commutative:
$$
\xymatrix{ \R_{\phi\disjoint g\disjoint
h}\ar[r]\ar[rd] & \R_{\phi\disjoint
g\disjoint\Id}\delta_{\Id\disjoint\Id\disjoint
h}
 \ar[d]\\
&
\R_{\phi\disjoint\Id\disjoint\Id}\delta_{\Id\disjoint
g\disjoint h}}
$$
\end{Claim}
\pf Follows directly from \ref{cmp}.
\endpf
\subsubsection{}
\begin{Claim} Assume that $\phi$ is not
bijective.  Then the  composition
$$
\P_{\phi\disjoint g}\to
\P_{\phi\disjoint\Id}\delta_{\Id\disjoint
g}\to \delta_{\phi\disjoint
\Id}\delta_{\Id\disjoint g}$$ equals
$$
\P_{\phi\disjoint g}\to
\delta_{\phi\disjoint  g}.
$$
\end{Claim}\subsubsection{}
Introduce a map
$$K(\phi_1,\phi_2,g_1,g_2):\R_{\phi_2\phi_1\disjoint
g_2g_1 }\to \R_{\phi_1\disjoint
g_1}\R_{\phi_2\disjoint g_2}\to
\R_{\phi_1\disjoint\Id}\delta_{\Id\disjoint
g_1}\R_{\phi_2\disjoint\Id}\delta_{\Id\disjoint
g_2}\to
\R_{\phi_1\disjoint\Id}\R_{\phi_2\disjoint\Id}
\delta_{\Id\disjoint g_2g_1}.
$$
\begin{Claim} The map
 $$\R_{\phi_2\phi_1\disjoint g}\to
\R_{\phi_2\phi_1\disjoint\Id}\delta_{\Id\disjoint
g}\to
\R_{\phi_1\disjoint\Id}\R_{\phi_2\disjoint\Id}
\delta_{\Id\disjoint g}
$$
is equal to
$$
\sum_{g_2g_1=g} K(\phi_1,\phi_2,g_1,g_2).
$$
\end{Claim}
\pf

It suffices to check that the two maps
coincide when compose with the maps

1)
$$
\R_{\phi_1\disjoint
\Id}\R_{\phi_2\disjoint\Id}\delta_{\Id\disjoint
g}\to
\delta_{\phi_1\disjoint\Id}\delta_{\phi_2\disjoint\Id}
\delta_{ \Id\disjoint g}
$$

2) $$
\R_{\phi_1\disjoint\Id}\R_{\phi_2\disjoint\Id}
\delta_{\Id\disjoint g}\to
\R_{\psi_1\disjoint\Id}\R_{\psi_2\disjoint\Id}
\R_{\phi_2\disjoint\Id}\delta_{\Id\disjoint
g},
$$
where $\psi_2\psi_1=\phi_1$ and
$\psi_1,\psi_2\neq \phi_1$;

3) $$
\R_{\phi_1\disjoint\Id}\R_{\phi_2\disjoint\Id}
\delta_{\Id\disjoint g}\to
\R_{\phi_1\disjoint\Id}\R_{\chi_1\disjoint\Id}
\R_{\chi_2\disjoint\Id}\delta_{\Id\disjoint
g},
$$
where $\chi_2\chi_1=\phi_2$ and
$\chi_1,\chi_2\neq \phi_2$.

Let us check 1). The composition
$$
\P_{\psi_1\disjoint
g_1}\delta_{\psi^1\disjoint g^1}\cdots
\P_{\psi_n\disjoint g_n} \to
\R_{\psi\disjoint g}\to
\R_{\phi\disjoint\Id}delta_{\Id\disjoint
g}\to
\R_{\phi_1\disjoint\Id}\R_{\phi_2\disjoint\Id}
\delta_{\Id\disjoint g}\to
\delta_{\phi_1\disjoint \Id
}\delta_{\phi_2\disjoint
\Id}\delta_{\Id\disjoint g}
$$
does not vanish iff the leftmost term is
$$
\P_{ \Id}\delta_{\phi_1\disjoint
g_1}\P_{\Id}\delta_{\phi_2\disjoint g_2}
\P_{\Id},
$$
in which case it is

a) zero if $\phi_1=\Id$ or $\phi_2=\Id$;

b)  identity  otherwise.

Let us now examine the composition:
$$
\R_{\phi\disjoint g}\to \R_{\phi_1\disjoint
g_1}\R_{\phi_2\disjoint g_2}\to
\R_{\phi_1\disjoint\Id}\delta_{\Id\disjoint
g_1}\R_{\phi_2\disjoint\Id}\delta_{\Id\disjoint
g_2}\to
\delta_{\phi_1\disjoint\Id}\delta_{\phi_2\disjoint\Id}
\delta_{\Id\disjoint
g_1}\delta_{\Id\disjoint g_2}
$$

According to the previous statement, this
composition vanishes if $\phi_1=\Id$ or
$\phi_2=\Id$. Otherwise, this composition
equals:

$$
\R_{\phi\disjoint g}\to \R_{\phi_1\disjoint
g_1}\R_{\phi_2\disjoint g_2} \to
\delta_{\phi_1\disjoint g_1}
\delta_{\phi_2\disjoint g_2}.
$$

We see that the two maps coincide.

2), 3) are immediate by induction.
\subsubsection{Compatibility with
$\boxtimes$}

\endpf

\begin{Claim} The following diagram is
commutative:
$$
\xymatrix{ \R_{f_1\disjoint
g_1}(M_1)\boxtimes \R_{f_2\disjoint
g_2}(M_2)\ar[r]\ar[d]&\R_{f_1\disjoint
\Id}\delta_{Id\disjoint g_1}(M_1)\boxtimes
\R_{f_2\disjoint \Id}\delta_{Id\disjoint
g_2}(M_2)\ar[r]&\R_{f_1\disjoint
f_2\disjoint\Id}\delta_{\Id\disjoint
g_1\disjoint g_2}(M_1\boxtimes M_2)\\
\R_{f_1\disjoint f_2\disjoint g_1\disjoint
g_2}(M_1\boxtimes M_2)\ar[rru]&&\\ }
$$
\end{Claim}
\pf Similar to the previous one
\endpf

\subsection{Direct images with respect to
projections} The reformulation of the
properties that were proven in the previous
subsections  in terms of direct image
functors with respect to projections is
given in \ref{defpresymm}. We are now
passing to giving an appropriate formalism
for description of structures that we have
encountered.

\textwidth 7in \textheight
10in \topmargin -1cm \voffset -1cm \hoffset
-3cm
\def\disjoint{{\sqcup}}
\def\bs{\backslash}
\def\A{{\frak A}}
\def\Aut{\text{Aut}}
\def\intr{\int}
\def\iI{{\frak iI}}
\def\ins{{\frak ins}}
\def\inr{{\frak int}}
\def\la{\lambda}
\def\An{{\frak An}}
\def\P{{\cal P}}
\def\limdir{\text{limdir}}
\def\ev{\text{ev}}
\def\notsubset{\text{notsubset}}
\def\identical{=}
\def\ins{{\frak ins}}
\def\dirlim{\limdir}
\def\usual{{\text{usual}}}
\def\A{{\cal A}}
\def\B{{\cal B}}
\def\F{{\frak F}}
\def\a{{\frak a}}
\def\b{{\frak b}}
\def\e{{\frak e}}
\def\jk{{\frak jk}}
\def \q{{\frak q}}
\def \qq{\q}
\def\sym{{\text{sym}}}
\def\res{{\text{res}}}
\def\f{{\fun}^{\sim}}
\def \Im{\text{Im}}
\def\jj{{\frak jj}}
\def\fp{({\fun'})^{\sim}}
\def\r{{\frak r}}
\def \into{\to}
\def \comm{\text{\bf comm}}
\def\lis{{\bf lis}}
\def\cat{{\bf cat}}
\def\red{\text{red}}
\def\lism{{\bf lism}}
\def\complexes{\text{complexes}}
\def\op{{\text{op}}}
\def\Ob{{\text{Ob}}}
\def \ores{{\O^{\text{res}}}}
\def\asm{{\frak ope'}}
\def\ens{{\bf ens}}
\def\ve{\varepsilon}
\def\vect{{\bf vect}}
\def \id{\text{Id}}
\def\Id{\id}
\def \as{{\frak as}}
\def\ope{{\frak ope}}
\def\fun{{\frak fun}}
\def\funp{({\frak fun}^{\psys})}
\def \Fun{{\frak Fun}}
\def\funct{{\frak Funct}}
\def\system{{{\cal SYS}}}
\def\psys{{\cal SYS}[p]}

\def\lie{{\frak lie}}
\def\li{\lie^1}
\def \ll{{{\li}_M}}
\def\psystem{{{\cal SYS}[p,\ll]}}
\def\psystemg{{{\cal SYS}[p,\lie_M]}}
\def\OPE{{\cal OPE}}

\def \OPER{{\cal OPE}^{\text{res}}}
\def\psystemb{\OPE[p,\lie_M]}
\def\psystemg{\psystemb}
\def\psis{\psystemb^{\text{res}}}

\def\univ{\text{univ}}
\def\psim{\psystemb^{\text{sym}}}
\def\psimr{\psystem^{\text{sym},\text{res}}}

\def\lisp{{[p,\lie_M,\S]}}
\def\lispq{{\lisp^{\text{res}}}}
\def\psml{\system[p,\lie_M,\S]}
\def\psmq{\psml^{\text{res}}}
\def\psmlh{\psml[[\h]]}
\def\psmqh{\psmq[[\h]]}
\def\coker{\text{Coker}}
\def\M{{\cal M}}
\def\gen{{\frak gen}}
\def\u{{\frak u}}
\def\Mg{{\frak M}}
\def \hh{{\frak h}}
\def \kk{{\frak kk}}
\def \II{{\frak II}}
\def \JJ{{\frak JJ}}
\def \GG{{\frak GG}}
\def \FF{{\frak FF}}
\def \m{{\frak m}}
\def \z{{\frak z}}
\def \B{{\cal B}}
\def \G{{\cal G}}
\def \M{{\cal M}}
\def \O{{\cal O}}
\def \R{{\Bbb R}}
\def \Z{{\frak Z}}
\def \Co{{\Bbb C}}
\def \C{{\cal C}}
\def \i{{\frak i}}
\def \j{{\frak j}}
\def \k{{\frak k}}
\def \ii{{\frak ii}}
\def \l{{\frak l}}

\def \L{{\frak L}}
\def \DDD{{\cal D}}
\def \D{{\cal D}}
\def \DD{\Delta}
\def \g{{\frak g}}
\def \ha{^\wedge}

\def \liminv{\text{liminv}}
\def \p{{\frak p}}
\def \pp{{\pi_p}}
\def \h{h}
\def \I{{\cal I}}
\def \s{{\frak s}}
\def \S{{\frak S}}
\def \U{{\cal U}}
\def \vs{\sigma}
\def\pt{\text{pt}}
\def\Hom{\text{Hom}}
\def \holi{{\frak holi}}
\def \La{{\holi_M}}
\def \Lab{{ {\frak holie}_M}}
\def \psl{\OPER[p,\Lab]}
\def \ps {\psl}

\def\Sk{{\frak Ske}}
\def\body{{\bf body}}
\def\bodyas{{\frak AS}}
\def\bodya{{\cal B}_{\text{\bf presymm}}}
\def\bodyb{{\cal B}_{\text{\bf symm}}}
\def\supersur{{\bf super-sur}}
\def\setf{{\bf setf}}
\def\S{{\cal S}}
\def \Sk{{\bf Ske}}
\def\B{{\cal B}}
\def\bodyas{{{\cal B}_{as}}}
\def\bodya{{\cal B}_{\text{\bf presymm}}}
\def\bodyb{{{\cal B}_{\text{\bf symm }}}}
\def\fact{{\bf fact}}
\def\bfu{{\bf u}}
\def\bfv{{\bf v}}
\def\l{{\cal L}}
\def\lB{{\l[B]}}
\def\lbodya{{\l[\bodya]}}
\def\lbodyb{{\l[\bodyb]}}
\def\pseudo{{{\bf body}_\otimes}}
\def\N{\Rs}
\def\Asm{{\text{\bf presymm}}}
\def\Bsm{{\text{\bf symm}}}
\section{Formalism for description of
different structures on  a collection of
functors}
\subsection{Definition of skeleton}
\subsubsection{} Let $\C$ be a category (for example,
the category of finite sets).
 We consider it as a 2-category with
trivial 2-morphisms.

 A skeleton  over $\C$ is a 2-category $\S$
with the following features:

  objects  of $\S$ are the same as in $\C$;

  all
categories $\S(S,T)$ are groupoids;

we have a 2-functor $P:\S\to \C^{\op}$.

Let us decode this definition. Note that $P$
induces maps of groupoids
$$
P(T,S):\S(T,S)\to \C^{\op}(T,S).
$$
For $F:S\to T$ being an arrow in $\C$, let
$\S(F):=P(T,S)^{-1}(F)$. Since
$\setf^\op(S,T)$ is a trivial groupoid (with
only identity morphisms), we have an
isomorphism of groupoids:
$$
\S(T,S)=\disjoint_{F:S\to T} \S(F).
$$

The rest of the structure can be
reformulated as follows:

For every pair of $\C$-arrows $F:S\to R$ and
$G:R\to T$ there should be given composition
functors
$$\circ(F,G):\S(F)\times \S(G)\to \S(GF);
$$

for every triple of $\C$-arrows $F:S\to R$,
$G:R\to P$, $H:P\to T$, there should be
given isomorphism $i(F,G,H)$ of functors
$$
\S(F)\times \S(G)\times \S(H)\to\S(GF)\times
\S(H)\to \S(HGF)
$$
and
$$
\S(F)\times \S(G)\times \S(H)\to \S(F)\times
\S(HG)\to \S(HGF).
$$
These isomorphisms should satisfy the
pentagon axiom.

Namely, let $$
\xymatrix{S\ar[r]^F&P\ar[r]^G&Q\ar[r]^H& R
\ar[r]^K& S}$$ be a sequence of maps of
finite sets. Every bracketing of the product
$KHGF$ specifies a functor $$\B(F)\otimes
\B(G)\otimes \B(H)\otimes \B(K)\to
\B(KHGF):$$ for example, the bracketing
$(KH)(GF)$ corresponds to the functor
$$\xymatrix{[(KH)(GF)]:&&&&&\\
\B(F)\otimes \B(G)\otimes \B(H)\otimes
\B(K)\ar[rrr]^{\circ(G,F)\otimes\circ(K,H)}&&&
\B(GF)\otimes
\B(KH)\ar[rr]^{\circ(KH,GF)}&&\B(KHGF)}
$$
The other bracketings produce the
corresponding functors in a similar way.
Total there are 5 such bracketings. The
associativity maps induce isomorphisms
between these functors as shown on the
following diagram:
$$
\xymatrix{ &K((HG)F)\ar[rr]    && K(H(GF ))  &\\
           (K(HG))F\ar[ur]&    &&
                               &(KH)(GF)\ar[ul]\\
            &      & ((KH)G)F\ar[ull]\ar[urr]&
            &}
$$
The pentagon axioms requires that this
pentagon be commutative.

\subsection{Body} A body  $\B$ built on a skeleton
$\S$ is an arbitrary dg- 2-category with the
following features:

 Objects of $\B$ are the same as in $\C$;

 $\Ob\B(T,S)=\Ob\S(T,S)$;

There exists a 2-functor
$$
s:\S\to \B$$  identical on objects and on
$\Ob\S(T,S)$ for all  $T,S$;

There exists a 2-functor
$P_B:\B\to\setf^\op$ such that $P_B\s=P$.

This definition is equivalent to the
following one.

A body $\B$ is a collection of dg-categories
$\B(F)$ for all $\C$-arrows $F:S\to T$ with
the following features:

1) $\Ob \B(F)=\Ob \S(F)$;

2) There are given functors
$\s:=\s(F):k[\S(F)]\to \B(F)$ identical on
objects;

3) There are given functors
$\circ_\B(G,F):\B(F)\times \B(G)\to \B(GF) $
which coincide on the level of  objects with
$\circ(G,F)$ and such that
$$\circ_\B(G,F)(\s(a)\times s(b))=\s(\circ(F,G)(a\times
b));$$, where $a$ is an arrow in $\S(G)$ and
$b$ is an arrow in $\S(F)$.

4) There are given associativity constraints
$c_\B(H,G,F)$ which satisfy the pentagon
axiom and are compatible with $c(F,G,H)$ in
the obvious way, that is: given arrows
$a,b,c$ in resp. $\S(H),\S(G),\S(F)$, one
has:
$$
\s (c_\S(a,b,c))=c_\B(\s(a),\s(b),\s(c)).
$$

\def\Com{{\bf FULL}}

\subsubsection{} To define a body one has to
prescribe complexes $\B(X,Y)$ for all
$X,Y\in \S(F)$ and certain poly-linear maps
between these complexes. Assume that the
isomorphism classes of $\C$ and the
isomorphism classes of $\S(F)$ form a
(countable) set for any $\C$-arrow $F$. Then
it is clear that the structure of a body
with skeleton $\S$ is equivalent to a
structure of an algebra over a certain
colored operad with a (countable) set of
colors. Denote this colored operad by
$\body(\S)$ The countability hypothesis will
be always the case in our constructions.

Thus, given a fixed skeleton, we have
notions of a free body, a body generated by
generators and relations etc.
\subsubsection{Example} In this example hte
objects of $\S(F)$ won't form a set.

For $F:S\to T$ we set $B(F)$ to be the
category of all functors from the category
of $\D_{X^T}$-modules to the category of
$\D_{X^S}$-modules. Let $\S(F)$ be the
groupoid of isomorphisms of $B(F)$. The rest
of the structure is defined in an obvious
way. Denote such a body by $\Com$.
\subsubsection{} A map of bodies is
naturally defined; a map $B\to \Com$ is
referred to as {\em a representation}.

\subsection{Construction of a skeleton}
We will mainly use a skeleton $\Sk$, which
will be now described. We set $\C:=\setf$ to
be the category  of finite sets. Let $F:S\to
T$ Objects of $\S(F)$ are sequences
$$
\xymatrix{ S\ar@{^{(}->}[r]^i&
U_0\ar@{>>}[r]^{p_1}&U_1\ar@{>>}[r]^{p_2}&
U_2\ar@{>>}[r]^{p_3}&\cdots\ar@{>>}[r]^{p_{n}}&
U_{n}=T}
$$
where $p_n\cdots p_1i=F$ and each $p_k$ is a
proper surjection (i.e. is not a bijection)
Such objects will be also denoted by
$$
\p_i\R_{p_1}\R_{p_2}\cdots \R_{p_{n}}.
$$
We shall also use a notation
$$\R_{p_1}\R_{p_2}\ldots \R_{p_n}$$
instead of
$$
\p_{\Id_{U_0}}\R_{p_1}\R_{p_2}\ldots
\R_{p_n}.$$

We do not exclude the case $n=0$, in which
case the corresponding object will be
written simply as $\p_i$.

\subsubsection{}\label{isog}
 Define isomorphisms in this groupoid. Let
$$Y=\p_j\R_{q_1}\R_{q_2}\cdots\R_{q_m}
$$
be another object in $\Sk(F)$, where
$q_j:V_{j-1}\to V_j$ and $V_m=T$.

The set $\Sk(F)(X,Y)$ is non-empty  only if
$n=m$

An isomorphism $p:X\to Y$ is a collection of
bijections  $b_{k}:U_k\to V_k$ for all $k$
satisfying the following natural
compatibility properties:

 1) $b_{n}=\Id_T$;

 2) For $k<k'$  set
 $$p_{k'k}:=p_{k'}p_{k'-1}\cdots p_{k+1};$$
 $$q_{k'k}:=q_{k'}q_{k'-1}\cdots q_{k+1};$$

 Then the diagram
 $$
 \xymatrix{U_k'\ar[r]^{b_{k'}}\ar[d]^{p_{k'k}}& V_{k'}\ar[d]^{q_{k'k}}\\
            U_{k}\ar[r]^{b_k}&V_{k}}
            $$
            commutes;

3) The diagrams
$$
\xymatrix{S\ar[r]^{i}\ar[dr]^{j}&
                         U_k\ar[d]^{b_k}\\
                    &V_k}
$$
commute.

The composition law is obvious.

\subsubsection{}
Let
$$\xymatrix{S\ar[r]^F&T\ar[r]^G&R}.$$
 The composition morphisms
$$
\circ_\ske(G,F)
$$
are defined as follows.

Let
$$
X=\p_i\R_{p_0}\R_{p_1}\cdots \R_{p_n},$$

where $i:S\to U_0$, $p_k:U_k\to U_{k+1}$,
$U_{n+1}=T$, $F=p_np_{n-1}\cdots p_0i$. Let
$$
Y=\p_j\R_{q_0}\R_{q_1}\cdots \R_{q_m},$$
where $j:T\to V_0$, $q_k:V_k\to V_{k+1}$,
$V_{m+1}=R$, $G=q_mq_{m-1}\cdots q_0j$.

Let $Z:=V_0\backslash j(T)$. Set
$U_k':=U_k\disjoint Z$, $p'_k:=p_k\disjoint
\Id_Z$, $j':U'_{n+1}\to V_0$,
$j':=j\disjoint i_Z$, where $i_Z:Z\to V_0$
is the natural embedding, $j'$ is then
bijective. Set $i':S\to U_0\to U'_0$ to be
the natural map.

Set $\circ_\ske(G,F)(X,Y)$ to be
$$\p_{i'}\R_{p'_0}\R_{p'_1}\cdots
\R_{jp'_n}\R_{q_0}\R_{q_1}\cdots \R_{q_m}.
$$
We shall write $XY$ instead of
$\circ_\ske(G,F)(X,Y)$.
 \subsection{Bodyes
$\bodyas$, $\bodya$, $\bodyb$}\label{babb}
We are going to define the bodies which
axiomatize the situations we are working
with: those of a system ($\bodyas$); of a
pre-symmetric system ($\bodya$) and of a
symmetric system ($\bodyb$).
 All these
bodies are constructed on the skeleton
$\Sk$.
\subsubsection{Body
$\bodyas$} Is generated by the maps
$\as(q,p):\R_{pq}\to \R_q\R_p$ of degree
zero with zero differential, where $q:S\to
R$ and $p:R\to T$ and  the relation:

The compositions
$$\R_{pqr}\stackrel{\as(r,qp)}\to \R_r\R_{pq}\stackrel{\as(q,p)}\to
\R_r\R_q\R_p
$$
and
$$
\R_{pqr}\to \R_{qr}\R_{p}\to \R_r\R_q\R_p
$$
coincide.
\subsection{Explicit description of
the complexes $\hom_{\bodyas(F)}(X,Y)$} Let
$$
X=\p_i\R_{p_0}\R_{p_1}\cdots \R_{p_n},$$

where $i:S\to U_0$, $p_k:U_k\to U_{k+1}$,
$U_{n+1}=T$, $F=p_np_{n-1}\cdots p_0i$. Let
$$
Y=\p_j\R_{q_0}\R_{q_1}\cdots \R_{q_m},$$
where $j:S\to V_0$, $q_k:V_k\to V_{k+1}$,
$V_{m+1}=T$, $F=q_mq_{m-1}\cdots q_0j$.

The space $\hom_{\bodyas(F)}(X,Y)$ is
non-empty only if for every $U_k$ there
exists a $V_l$ such that $\#U_k=\#V_l$.
Define the set $S(X,Y)$ whose each element
$f$ is a collection of bijections
$f_{kl}:U_k\to V_l$, whenever $\#U_k=\#V_l$
satisfying  all the properties from
\ref{isog}. Set
$$\hom_{\bodyas(F)}(X,Y):=k[S(X,Y)].$$

The composition law in $\bodyas(F)$ and the
inclusion functor $\Sk(F)\to \bodyas(F)$ are
immediate.

\subsubsection{The body $\bodya$}
It is generated over $\bodyas$ by the
elements of two types:

Type 1. Consider a commutative triangle
$$\xymatrix{S\ar@{>>}[r]^p& T\\
            R\ar@{^{(}->}[u]^i\ar@{^{(}->}[ur]^j&}
            $$
in which $i,j$ are injections and $p$ is a
proper surjection. We then have a degree +1
map
$$
L(i,p):\p_i\R_p\to  \p_j.
$$

Type 2. Consider a commutative square
\begin{equation}
\xymatrix{ R\ar@{->>}[r]^p & T\\
           S\ar@{->>}[r]^q\ar@{^{(}->}[u]^i& P\ar@{^{(}->}[u]^j}
\end{equation}
in which $i,j$ are injections and $p,q$ are
proper surjections. Call such a square {\em
suitable} if the following is satisfied: Let
$T_1=T\backslash T_2$ be the subset of all
$t\in T$ such that $p^{-1}t\cap i(S)$
consists of $\geq 2$ elements. Then
$p^{-1}(T_1)\subset i(S)$, i.e.:
$$
\#(p^{-1}t\cap i(S))\geq 2\follows
p^{-1}(t)\subset i(S).
$$

 We then have a degree zero
map
$$
A(i,p,j,q):\p_i\R_p\to \R_q\p_j,
$$
where $\R_q\p_j:=\circ_\Sk(\R_q,\p_j)$.

\subsubsection{Relations}

1. Let $$\xymatrix{R\ar@{>>}[r]^p&  T\\
        S\ar@{^{(}->}[u]^i\ar@{>>}[r]^q
                                 &P\ar@{^{(}->}[u]^j }
                                 $$
        be a suitable square and $q=q_2q_1$,
where $q_1,q_2$ are surjections.

Define the  set $X(q_1,q_2)$ of isomorphism
classes of  commutative diagrams
$$\xymatrix{R\ar@{>>}[r]^{p_1}&
U\ar@{>>}[r]^{p_2}& T\\
        S\ar@{^{(}->}[u]^i\ar@{>>}[r]^{q_1}
        \ar@{^{(}->}[u]^j&
        V\ar@{^{(}->}[u]^{j'}\ar@{>>}[r]^{q_2}&
         P\ar@{^{(}->}[u]^j }$$
         We will refer to such a diagram as
         $(p_1,p_2,j')$.
Both squares in every such  a diagram are
automatically suitable. Therefore, every
element $x:=(p_1,p_2,j')\in X(q_1,q_2)$
determines a map
$$
m_x:\p_i\R_p\to \p_i\R_{p_1}\R_{p_2}\to
\R_{q_1}\p_{j'}\R_{p_2}\to
\R_{q_1}\R_{q_2}\p_j.
$$

         Then the relation says that the
        composition
        $$
        \p_i\R_p\to \R_q\p_j\to
        \R_{q_1}\R_{q_2}\p_j$$
 equals
 $$
 \sum_{x\in X(q_1,q_2)} m_x
 $$

 2. Consider the following commutative
 diagram
 $$\xymatrix{ R\ar[r]^p              &T\\
S_1\ar[u]^{i_2}\ar[r]^q   &P_1\ar[u]^{j_2} \\
S\ar[u]^{i_1}\ar[r]^r     & P\ar[u]^{j_1}}
$$
in which both small squares are suitable.
Then the large square is also suitable and
the following maps coincide:
$$
\p_{i_2i_1}\R_p\to \R_r\p_{j_2j_1}
$$
and
$$
\p_{i_2i_1}\R_p\to \p_{i_1}\p_{i_2}\R_p\to
\p_{i_1}\R_q\p_{j_2}\to
\R_{r}\p_{j_1}\p_{j_2} \to \R_{r}\p_{j_2j_1}
$$

3. Consider the following commutative
diagram:
$$
\xymatrix{R\ar@{>>}[r]^p        &                 T\\
            S\ar@{>>}[r]^q\ar@{^{(}->}[u]^i &
                                        P\ar@{^{(}->}[u]^j\\
            Q\ar@{^{(}->}[u]^k\ar@{^{(}->}[ru]^l&}
            $$
where the upper square is suitable. Then the
following maps coincide:

$$
\p_{ik}\R_p\to \p_k\p_i\R_p\to
\p_k\R_q\p_j\to\p_{qk}\p_j=\p_l
$$
and $$\p_{ik}\R_p\to \p_{pik}=\p_l.$$

 4. Consider the following commutative
diagram
$$ \xymatrix{S\disjoint U\ar[r]&
T\disjoint U\\
S\ar[r]\ar[u]&  T\ar[u]}
$$
This diagram is suitable and we require that
the corresponding map $A(i,p,j,q)$ be equal
to the corresponding isomorphism in $\Sk$.

5. Let $$
\xymatrix{ R\ar@{->>}[r]^p & T\\
           S\ar@{->>}[r]^q\ar@{^{(}->}[u]^i& X\ar@{^{(}->}[u]^j}
$$
and

$$
\xymatrix{ R_1\ar@{->>}[r]^{p_1} & T_1\\
           S_1\ar@{->>}[r]^{q_1}\ar@{^{(}->}[u]^{i_1}& X_1\ar@{^{(}->}[u]^{j_1}}
$$

be suitable squares and let $s:S\to S_1$,
$r:R\to R_1$, $t:T\to T_1$, $x:X\to X_1$ be
bijections fitting the two squares into a
commutative cube. Then the map $A(i,p,j,q)$
coincides with the map
$$
\p_i\R_p\cong
\p_s\p_{i_1}\p_{r^{-1}}\p_r\R_{p_1}\p_{t_1^{-1}}\cong
\p_s\p_{i_1}\R_{p_1}\p_{t_1^{-1}}\to
\p_s\R_{q_1}\p_{j_1}\p_{t_1^{-1}}\cong
\p_s\R_{q_1}\p_{x^{-1}}\p_x\p_{j_1}\p_{t_1^{-1}}
\R_q\p_j.
$$

 \subsubsection{Differentials} The
differential of the map $L(i,p)$ is computed
as follows. Consider the set of all
equivalence classes of decompositions
$p=p_2p_1$, where $p_1,p_2$ are surjections
and $p_1i$ is injection. We then have a map
$$l(p_1,p_2):\p_i\R_p\to
\p_i\R_{p_1}\R_{p_2}\to \p_{p_1i}R_{p_2}\to
\p_{p_2p_1i}=\p_{pi}.
$$
We then have
$$
dL(i,p)+\sum_{(p_1,p_2)}l(p_1,p_2)=0.
$$

2. Let
$$
\xymatrix{Q:R\ar@{>>}[r]^p&T\\
 S\ar@{>>}[r]^q\ar@{^{(}->}[u]^i&P\ar@{^{(}->}[u]^j}
 $$
 be a suitable square. Define two sets
 $L(Q)$ and $R(Q)$ as follows. The set
 $L(Q)$ is the set of all isomorphism
 classes of diagrams:
$$
\xymatrix{R\ar@{>>}[r]^p_1& R_1\ar@{>>}[r]^{p_2}&T\\
 S\ar@{>>}[rr]^q\ar@{^{(}->}[u]^i\ar@{^{(}->}[ur]^{i_1}&&
 P\ar@{^{(}->}[u]^j}
 $$
such that $p=p_1p_2$. It is clear that the
internal commutative square in this diagram
is also suitable.

Define the set $R(Q)$ as the set of
isomorphisms classes of diagrams

$$
\xymatrix{R\ar@{>>}[r]^{p_1}& R_1\ar@{>>}[r]^{p_2}&T\\
 S\ar@{>>}[rr]^q\ar@{^{(}->}[u]^i&&P\ar@{^{(}->}[u]^j
 \ar@{^{(}->}[ul]^{j_1}}
 $$
where $p=p_1p_2$. The internal square in
such a diagram is always suitable as well.

Every element $l:=(p_1,p_2,i_1)\in L(Q)$
determines a map
$$
f_l:\p_i\R_p\to \p_i\R_{p_1}\R_{p_2}\to
\R_{q}\p_{i_1}\R_{p_2}\to
\R_{q}\p_{p_2i_1}=\R_{q}\p_j.
$$
Every element $r=(p_1,p_2,j_1)\in R(Q)$
determines a map
$$
g_r:\p_i\R_p\to \p_i\R_{p_1}\R_{p_2}\to
\p_{p_1i}\R_{p_2}=\p_{j_1}\R_{p_2}\to
\R_q\p_j.
$$

We then have
$$
dA(i,p,j,q)=\sum_{l\in L(Q)}f_l-\sum_{r\in
R(Q)}g_r
$$
This completes the definition. We need to
check that $d^2=0$ and that $d$ preserves
the ideal generated by the relations, which
is left to the  reader.

\subsubsection{} The system $<\R>$ with its
properties provides for a representation of
$\bodya$.

\subsubsection{Explicit description of the
categories
 $\bodya(F)$}
Consider two objects $X,Y$ in $A(F)$:
$$\xymatrix{
S\ar@{^{(}->}[r]^i & U_0\ar@{->>}[r]^{p_0} &
U_1\ar@{->>}[r]^{p_1}&\cdots&
U_{k}\ar@{->>}[r]^{p_k}&U_{k+1}}
$$
and
$$\xymatrix{
S\ar@{^{(}->}[r]^j & V_0\ar@{->>}[r]^{q_0} &
V_1\ar@{->>}[r]^{q_1}&\cdots&
V_{l-1}\ar@{->>}[r]^{q_{l}}&U_{l+1}}
$$

 Define the set $M(X,Y)$ whose each
element is a collection of injections
$$
j_r:V_{m_r}\hookrightarrow U_r,
$$
where $r=0,1,\ldots,k+1$, $m_0=0$, $0\leq
m_{r+1}-m_r\leq 1$, $m_{k+1}=l+1$. The
following conditions should be satisfied:

1) if  $m_{r+1}=m_r$, then the diagram
$$\xymatrix{U_{r}\ar[r]^{p_r} &U_{r+1}\\
V_{m_r}\ar[u]^{j_r}\ar[ur]^{j_{r+1}}&}$$
must be commutative

2)if $m_{r+1}=m_r+1$, then the diagram
$$\xymatrix{U_{r}\ar[r]^{p_r} &U_{r+1}\\
V_{m_r}\ar[u]^{j_r}\ar[r]^{q_r}&V_{m_r+1}\ar[u]^{j_{r+1}}}
$$
must be commutative and suitable.

Every element $m=(j_1,j_2,\ldots,j_{k+1})$
in $M(X,Y)$ defines a map
$$
A(m):\p_i\R_{p_0}\R_{p_1}\cdots \R{p_k}\to
\p_j\R_{q_0}\R_{q_1}\cdots \R_{q_l},
$$
where we set $\R_{\Id}=\Id$, as follows.
Define
$$q'_{m_r}:V_{m_k}\to V_{m_{r+1}}$$
to be $\Id_{V_{m_r}}$ if $m_r=m_{r+1}$ and
$q_{m_r}$ if $m_{r+1}=m_r+1$. We then have
maps
$$F_k:\p_{j_r}\R_{p_r}\to
\R_{q'_{m_r}}\p_{j_{r+1}},
$$
where $ F_k=L(j_r,p_r,j_{r+1})$ if
$m_{r+1}=m_r$, and
$F_k=C(j_r,p_r,j_{r+1},q'_{m_r})$ if
$m_{r+1}=m_r+1$.

Set
\begin{eqnarray*}
A(m):\p_i\R_{p_0}\R_{p_1}\cdots \R_{p_l}\to
\p_j\p_{j_0}\R_{p_0}\R_{p_1}\cdots
\R_{p_l}\stackrel{A'(j_0,p_0,j_1,q'_{m_0})}\longrightarrow\\
\p_j\R_{q'_{m_0}}\p_{j_1}\R_{p_1}\cdots
\R_{p_l}\stackrel{A'(j_1,p_1,j_2,q'_{m_1})}\longrightarrow
\p_j\R_{q'_{m_0}}\R_{q'_{m_1}}\p_{j_2}\R_{p_2}\R_{p_3}
\cdots \R_{p_l}\to\\
\cdots\to
\p_j\R_{q'_{m_0}}\R_{q'_{m_1}}\cdots
\R_{q'_{m_{k+1}}}\cong
\p_j\R_{q_0}\R_{q_1}\cdots \R_{q_l},
\end{eqnarray*}
where
$$A'(j_u,p_u,j_{u+1},q'_{m_u})=A(j_u,p_u,j_{u+1},q'_{m_u})$$
if $q'_{m_u}\neq \Id$. Otherwise
$$
A'(j_u,p_u,j_{u+1},q'_{m_u})=l(j_u,p_u,j_{u+1}).
$$
 Let $N(X,Y):=k[M(X,Y)]$.  Let
$$H(X,Y)=\oplus_{Z} N(Z,Y),
$$
where the sum is taken over all refinements
$Z$ of $X$. We have an obvious map
$H(X,Y)\to \hom(X,Y)$. Set
$A_F(X,Y):=H(X,Y)$. The relations given in
the previous section  provide us with a
composition law $H(X,Y)\otimes H(Y,Z)\to
H(X,Z)$ and a differential.
\subsubsection{Body $\bodyb$} The definition
of the body $\bodyb$ is exactly the same as
the one of the body $\bodya$ except that the
maps $A(i,p,j,q)$ are defined for all
commutative squares
$$\xymatrix{R\ar@{>>}[r]^p& T\\
            S\ar@{>>}[r]^q\ar@{^{(}->}[u]^i&P\ar@{^{(}->}[u]^j}
$$
not necessarily suitable; the relations are
the same except that we lift everywhere the
restriction of suitability; the formulas for
the differential remain the same.

It is clear that we have a map of bodies
$$
i:\bodya\to \bodyb
$$

We are going to study this map.
\subsubsection{} Explicit expression for
$\hom_{\bodyb(F)}(X,Y)$, where $X,Y\in
\Sk(F)$ is exactly the same as for $\bodya$.

The further study of $\bodyb$ is facilitated
by the statement we are going to consider
\subsubsection{}  Let
$$
\xymatrix{S\ar@{^{(}->}[r]^i
&R\ar@{>>}[r]^p& T}
$$
be a diagram. Call it {\em super-surjective}
if for every $t\in T$:

  either $p^{-1}t\cap
i(S)$ has at least two elements

 or
$p^{-1}t$ is a one-element subset of $i(S)$.

\begin{Claim} Let
$$\xymatrix{R\ar@{>>}[r]^p& T\\
            S\ar@{>>}[r]^q\ar@{^{(}->}[u]^i&
                                 P\ar@{^{(}->}[u]^j}
            $$
be a commutative diagram.

Then there exists a decomposition:
$$\xymatrix{R\ar@{>>}[r]^p& T\\
U\ar@{^{(}->}[u]^{i_2}\ar@{>>}[r]^r&
                            P\ar@{^{(}->}[u]^j\\
            S\ar@{>>}[ru]^q\ar@{^{(}->}[u]^{i_1}&
                                 }
            $$
where the diagram is commutative,
$i=i_2i_1$, the square $(i_2,p,j,r)$ is
suitable and the pair $(i_1,q)$ is
super-surjective.

Such a decomposition is unique up-to an
isomorphism.
\end{Claim}
\pf {\em Existence}. Call an element $t\in
T$ {\em good} if  $p^{-1}t$ satisfies the
condition of the definition. Let $G_T\subset
T$ be the subset of all good elements.
  Let
$U:=i(S)\cup p^{-1}G_T$. Let $i_1,i_2$ be
the natural inclusions. By definition, for
every $t\in G_t$, the intersection
$p^{-1}t\cap i(S)$ is non-empty. Hence,
$G_T\subset pi(S)=jq(S)$ and, therefore,
$p(U)= jq(S)$.

 Thus, $p(U)=
\Im j$, which implies that the map
$p|_U:U\to T$ uniquely decomposes as $jr$,
where $r:U\twoheadrightarrow P$. It is clear
that all the conditions are satisfied.

{\em Uniqueness} is also clear.
 \endpf

 \subsubsection{Corollary} Let
$X\in \Sk(F)$ be an object of the form
$$
\p_j\p_{i_1}\R_{p_1}\p_{i_2}\R_{p_2}\cdots
\p_{i_n}\R_{p_n},
$$
where every pair $(i_k,p_k)$, $i_k:U_k\to
A_k$, $p_k:A_k\to U_{k+1}$, is
super-surjective. Let
$$Y=\p_j\R_{p_1i_1}\R_{p_2i_2}\cdots
\R_{p_ni_n}.
$$
The maps $\p_{i_k}\R_{p_k}\to \R_{p_ki_k}$
induce a map $f_X: X\to Y$. Call $X$ {\em a
super-surjective decomposition} of $Y$. Let
$\supersur(Y)$ be the groupoid  of all
super-symmetric  decompositions of $Y$ and
their isomorphisms (i.e.  collections of
isomorphisms $U_k\to U_k'$ fitting into the
commutative diagrams... . It is clear that
if $a:X_1\to X_2$ is such an isomorphism,
then $f_{X_1}=f_{X_2}a$.

Let $Z,Y\in \Sk(F)$. Define a functor
$h_Z:\supersur(Y)\to \complexes$ by the
formula $h_Z(X)=\hom_{\bodya(F)}(Z,X)$.

The collection of maps $f_X$ induces a
functor

$$\limdir _{\supersur(Y)} h_Z\to
\hom_{\bodyb(F)}(Z,Y).
$$
\begin{Claim} This map is an isomorphism.
\end{Claim}
\subsubsection{One more lemma} Let $Y\in
\Sk(F)$. Let $F=F_2F_1$ be a decomposition
and assume that we have an isomorphism
$t:Y\to Y_2Y_1$, where $Y_i\Sk(F_i)$. We
then have a natural functor:
$$
\supersur(Y_2)\times \supersur(Y_1)\to
\supersur(Y).
$$
\begin{Lemma}\label{supsur} This functor is
an equivalence of groupoids.
\end{Lemma}
\pf Clear.
\endpf
\subsection{Pseudo-tensor bodies}
Let $\B$ be a body. A pseudo-tensor
structure on $\B$ is a collection of several
pieces of data, the first one being functors
$$\Psi(\{F_i\}_{i\in I}): \otimes_{i\in I}
\B(F_i) \otimes \B(F)^\op\to \complexes,
$$
where  $F=\disjoint_{i\in I F_i}$, for all
$n>0$ and all collections $F_i:S_i\to T_i$
 of maps of finite sets indexed by an
 arbitrary finite non-empty set $I$.
 Let $X_i\in \B(F_i)$ and $X\in \B(F)$.
We then denote
$$
\hom(\{X_i\}_{i\in I};X):=\Psi(\{F_i\}_{i\in
I})((\otimes_{i\in I}X_i)\otimes X).
$$
Let $\pi:I\to J$ be  a surjection of finite
sets. Let $F_i:S_i\to T_i$, $i\in I$. be
maps of finite sets. For a $j\in J$ set
$$F_j:=\disjoint_{i\in \pi^{-1}j}F_i.$$
Let $X_i\in \B(F_i)$, $Y_j\in \B(F_j)$,
where  $i\in I$, $j\in J$. Set
$$
\hom_\pi(\{X_i\}_{i\in I};\{Y_j\}_{j\in
J}):=\otimes_{j\in J}\hom(\{X_i\}_{i\in
I},\{Y_j\}_{j\in J}).
$$
An element $f$ in this complex will be also
written as
$$f:\{X_i\}_{i\in I}\to
\{Y_j\}_{j\in J}.$$

Let $X,Y\in \B(G)$, where $G:S\to T$ be a
map of finite sets. It is assumed that
$\hom(X,Y)=\hom_{\B(G)}(X,Y)$.
\subsubsection{Composition of the first type}
Let $F=\disjoint_{i\in I} F_i$.  The second
feature of a  pseudo-tensor structure is a
collection of {\em composition maps of the
first kind} :
$$C_1(\{X_i\}_{i\in I};\{Y_j\}_{j\in J};
Z):\hom_\pi(\{X_i\}_{i\in I};\{Y_j\}_{j\in
J})\otimes \hom(\{Y_j\}_{j\in J};Z)\to
\hom(\{X_i\}_{i\in I};Z).
$$

 Let $\sigma:J\to K$ be another surjection
 and set $F_k:=\disjoint_{j\in \pi^{-1} k}$.
 Pick objects $Z_k\in \B(F_k)$, $k\in K$.
 Define a map
 \begin{eqnarray*}
 C_1(\{X_i\}_{i\in I};\{Y_j\}_{j\in J};
\{Z_k\}_{k\in K}):\hom_\pi(\{X_i\}_{i\in
I};\{Y_j\}_{j\in J})\otimes
\hom_\sigma(\{Y_j\}_{j\in J};\{Z_k\}_{k\in
K})\\
\to \hom_{\sigma\pi}(\{X_i\}_{i\in
I};\{Z_k\}_{k\in K})
\end{eqnarray*}
as the tensor product $$ \otimes_{k\in
K}C_1(\{X_i\}_{i\in
(\sigma\pi)^{-1}k};\{Y_j\}_{j\in
\sigma^{-1}k}; Z_k).
$$
\subsubsection{Compositions of the second
kind} Let $F_i:S_i\to T_i$, $G_i:T_i\to
R_i$, $i\in I$ be a family of maps of finite
sets. Let $F=\disjoint_{i\in I} F_i$,
$G=\disjoint_{j\in J}G_j$.  Let $X_i\in\B(
F_i)$, $Y_i\in \B(G_i)$,  $Z\in \B(F)$,
$W\in \B(G)$. We then have objects
$Y_iX_i\in \B(G_iF_i)$, $WZ\in \B(GF)$. The
last feature of a pseudo-tensor structure is
a prescription of {\em composition maps of
the second kind}:
$$
C_2(\{X_i\}_{i\in I},Z;\{Y_i\}_{i\in I},W)
:\hom(\{X_i\}_{i\in
I};Z)\otimes\hom(\{Y_i\}_{i\in
I};W)\to\hom(\{Y_iX_i\}_{i\in I};WZ).
$$
Let $\pi:I\to J$ be a surjection.  Let
$F_j=\disjoint_{i\in \pi^{-1}j}F_i$;
$G_j=\disjoint_{i\in \pi^{-1}j}G_i$. Let
 $X_i\in \B(F_i)$,
$Y_in \B(G_i)$, $i\in I$; $Z_j\in \B(F_j)$,
$W_j\in \B(G_j)$, $j\in J$. Define a map
\begin{eqnarray*}
C_2(\{X_i\}_{i\in I};\{Y_i\}_{i\in
I};\{Z_j\}_{j\in J};\{W_j\}_{j\in
J}):\hom_\pi(\{X_i\}_{i\in I} ,\{Z_j\}_{j\in
J} )\otimes\hom_\pi(\{Y_i\}_{i\in
I},\{W_j\}_{j\in J})\\ \to
\hom_{\pi}(\{Y_iX_i\}_{i\in I}
;\{W_jZ_j\}_{j\in J}).
\end{eqnarray*}
as the tensor product
$$
\otimes_{j\in J}C_2(\{X_i\}_{i\in
\pi^{-1}j},Z_j;\{Y_i\}_{i\in
\pi^{-1}j},W_j). $$
\subsubsection{Axiom}
The only axiom is as follows. Let $I$ be a
finite set and consider an $I$- family of
chains of maps
$$\xymatrix{
S_i^0\ar[r]^{F_i^1}&
S_i^1\ar[r]^{F_i^2}&\cdots\ar[r]^{F_i^N}&S_i^N}$$
where $N$ is a fixed number. For $0\leq p <
q\leq N$, set $F_i^{qp}:S_i^p\to S_i^q$ to
be the composition $$ F^q_iF^{q-1}_p\cdots
F^{p+1}_i.
$$
We also set $F^{qq}_i:=\Id_{S_i^q}$.

 We
shall also need a chain of surjections
$$\xymatrix{
I=I_0\ar@{>>}[r]^{\pi_1}&
I_1\ar@{>>}[r]^{\pi_2}&\cdots\ar@{>>}[r]^{\pi_M}&
I_M}
$$

where $M$ is a fixed natural number. For
$0\leq u<v\leq M$, denote by
$\pi_{vu}:I_u\to I_v$ the composition
$$
p_{vu}=\pi_v\pi_{v-1}\cdots \pi_{u+1};
$$
set $\pi_{uu}:=\Id_{I_u}$. For a $j\in I_k$
define a subset
$$\ov{j}:=\pi_{k0}^{-1}(j)$$ of $I_0$. Take
the following disjoint unions
$$
S_j^r=\disjoint_{i\in \ov{j}} S_i^r;\quad
F_j^r=\disjoint_{i\in \ov{j}} F_i^r;
$$
$$
F_j^{qp}=\disjoint_{i\in \ov{j}} F_j^{qp}.
$$

Pick elements  $X_j^{k}\in \B(F_j^{k})$, for
all $j\in I_u$, $0\leq u\leq M$ and all
$k=1,\ldots, N$. For $0\leq p< q\leq N,$
 set
$$
X_j^{qp}=X_j^{q}X_j^{q-1}\cdots X_j^{p+1}
$$
so that $X_j^{qp}\in \B(F_j^{qp})$.

Iterating various compositions of the two
kinds in various ways, one can construct, a
priori, several maps
\begin{eqnarray*}
\bigotimes\limits_{u=1}^{M}
\bigotimes\limits_{k=1}^{N}
\hom_{\pi_{u}}(\{X_j^{k}\}_{j\in
I_{u-1}};\{X_l^{k}\}_{l\in I_{u}})\\
\to \hom_{\pi_{M,1}}(\{X_s^{N,0}\}_{s\in
I_1},\{X_v^{N,0}\}_{v\in I_{M}})
\end{eqnarray*}
 The
axiom says that all these maps should
coincide. Denote thus obtained unique map by
$\comp\{X_s^k\}$.

 \subsubsection{}  Given a fixed skeleton $\S$,
 a structure of a peudo-tensor body on this
 skeleton is equivalent to the one of
 algebra over a certain colored operad $\pseudo(\S)$.
 Therefore, pseudo-tensor bodies can be
 specified by means of generators and
 relations.
 \subsubsection{Example} Introduce a
 pseudo-tensor structure on $\Com$ as
 follows. Let $X_a\in \Com(F_a)$ and $Y\in
 \Com (F)$, where $F_a:S_a\to T_a$,
 $S=\disjoint_a S_a$, $F=\disjoint_a F_a$,
 etc.

 Let $$\boxtimes_a:\prod_a \Dmod_{X^{S_a}}\to
 \Dmod_{X^S}$$,
 $$\boxtimes_a:\prod_a \Dmod_{X^{T_a}}\to
 \Dmod_{X^T}$$

be the functors of the exterior tensor
 product. Set
 $$\hom_{\Com} (\{X_a\};Y):=
 \hom(\boxtimes_a X_a; Y\circ \boxtimes_a),
 $$
 where $\hom$ is taken in the category of
 functors:
 $$
 \Prod_a \Dmod_{X^{T^a}}\to \Dmod_{X^S}.
 $$

 \subsection{Maps of  pseudo-tensor bodies}
 Let $\B_1,\B_2$ be pseudo-tensor bodies
 over  skeletons  resp. $\S_1$ and $\S_2$.
 Our goal is to define a notion of a map $R:\B_1\to \B_2$.
We shall give two equivalent definitions.
The first definition is based on a notion of
\subsubsection{Induced skeleton}
Let $X:k[\S_1]\to \B_2$ be a 2-functor which
maps $k[\S_1](F) \to \B_2(F)$ for all $F$.
In the sequel we shall write $\S_1$ instead
of $k[\S_1]$.

 This structure is equivalent to the
 following one:

1)we have functors $X(F):\S_1(F)\to \B_2(F)$
for all $F$

2) for all composable pairs $F,G$, the
natural transformation $I_X(G,F)$, shown on
the diagram:
$$\xymatrix{ \S_1(F)\times \S_1(G)\ar[rr]^{\circ(G,F)}
\ar[dd]^{X(F)\times X(G)}&&
\S_1(GF)\ar[dd]^{X(GF)}\\
&&\ar@2{->}[dl]_{I_X(G,F)}\\ \B_2(F)\otimes
\B_2(G)\ar[rr]^{\circ(G,F)}&&\B_2(GF)}
$$

3) The transformations $I_X$ should be
compatible with the associativity
transformations of $\S_1$ and $\B_2$ in a
natural way.

Using such an $X$ we shall construct a body
$X^{-1}\B_2$ on the skeleton $\S_2$. First
of all, we set
$$
\hom_{X^{-1}\B_2}(\{Y_a\}_{a\in
A};Z):=\hom_{\B_2}(\{X(Y_a)\}_{a\in
A},X(Z)).
$$
The compositions of the first and the second
kinds on $X^{-1}\B_2$ are naturally  induced
by those on $\B_2$. Thus constructed
pseudo-tensor body is called {\em induced}.
\subsubsection{Definition of a map $f:\B_1\to
\B_2$} By definition, such a  map is given
by a 2-functor $X_f:\S_1\to \B_2$ as above
and by a map $f':\B_1\to X_f^{-1}\B_2$,
where the meaning of $f'$ is as follows:
since $\B_1$ and $X_f^{-1}\B_2$ are
pseudo-tensor bodies over the same skeleton
$\S_1$ they can be both interpreted as
algebras over the operad $\pseudo(\S_1)$;
$f'$ is by definition a map of such
algebras.

This definition will  be now decoded.

\subsection{More straightforward approach}
To define a map $\B_1\to \B_2$ one  has to
prescribe the following
 data:

 1) a collection of functors $R_F:\B_1(F)\to \B_2(F)$;

 2) for every sequence of maps of finite sets
 $F_i:S_{i-1}\to S_{i}$, $i=1,2,\ldots, N$,
 such that
 $F_NF_{N-1}\cdots F_1=F$,
 consider a diagram of functors:
 $$\xymatrix{
 \otimes_{i} \B_1(F_i)\ar[r]^{\otimes_i R_{F_i}}\ar[d]^{\circ(F_N,F_{N-1}\ldots,F_1)}&
  \otimes_i\B_2{F_i}\ar[d]^{\circ(F_N,F_{N-1},\ldots,F_1)}\\
  \B_1(F)\ar[r]^{R_F}&\B_2(F)}
  $$
  There should be specified an isomorphism
  $I(F_1,F_2,\ldots F_N)$ between the composition of
  the top arrow followed by the right arrow
  and the composition of the left arrow followed by the bottom arrow.
  As it is common in the theory of 2-categories,
  $I(F_1,F_2,\ldots,F_N)$ will be denoted by a double diagonal arrow:
$$\xymatrix{
 \otimes_{i} \B_1(F_i)\ar[rr]^{\otimes_i
 R_{F_i}}\ar[dd]^{\circ(F_1,\ldots,F_N)}&&
  \otimes_i\B_2{F_i}\ar[dd]^{\circ(F_1,F_2,\ldots,F_N)}\\
  &&\ar@2{->}[dl]_I\\
  \B_1(F)\ar[rr]^{R_F}&&\B_2(F)}
  $$

  3) For every $\{X_a\}_{a\in A}$, $X_a\in \B_1(F_a)$,
  and every $Y\in \B_2(F)$, where $F=\disjoint_{a\in A}F_a$,
  there should be given a map of complexes:
  $$
  \R_{\{X_a\}_{a\in A};Y}:\hom_{\B_1}(\{X_a\}_{a\in A};Y)\to
  \hom_{\B_2}(\{R_{F_a}(X_a)\}_{a\in A};R_F(Y)).
  $$

  The axioms are as follows:

  1) Associativity axiom for $I(F_N,\ldots,
  F_1)$.

  Pick a sequence $1=i_1\leq
  i_2\leq\cdots i_k=N$. Set
$$G^{r}:=F_{i_r}F_{i_r-1}\cdots
F_{i_{r-1}+1}.
$$
Let
$$\circ_r:=\circ(F_{i_r},F_{i_r-1},\cdots,
F_{i_{r-1}+1});
$$
let $I_r:=I(F_{i_r},F_{i_r-1},\cdots,
F_{i_{r-1}+1})$. Let
$I:=I(G_k,G_{k-1},\ldots G_1)$.
 We then have the following diagram:
$$\xymatrix{
\otimes_i \B_1(F_i)\ar[rrr]^{\otimes_i
R_{F_i}}
        \ar[ddd]^{\otimes_r \circ_r}
                               &&&\otimes_i \B_2(F_i)
 \ar[ddd]^{\otimes_r \circ_r}\\
 &&&\\
                    &&& \ar@2{->}[dl]_{\otimes_r I_r}                 \\
\otimes_r \B_1(G_r)\ar[rrr]^{\otimes_r R_r}
\ar[ddd]^{\circ(G_k,G_{k-1},\ldots, G_1)}
 &&& \otimes_r\B_2(G_r)
 \ar[ddd]^{\circ(G_k,G_{k-1},\ldots, G_1)}\\
 &&&\\
                      &&&  \ar@2{->}[dl]_I \\
\B_1(F)\ar[rrr]^{R_F}               &&&
\B_2(F)}
$$

We then see that the two squares of this
diagram are composable and the axiom
requires that the composition be equal to
$I(F_N,F_{N-1},\cdots,F_1)$.

  2) Compatibility  of  $R(\{X_a\}_{a\in A};Y)$
with compositions of the first type.

Let $p:A\to B$ be a surjection of finite
sets. Let $F_a:S_a\to T_a$ be an $A$-family
of maps of finite sets. Let
$$
S_b=\disjoint_{a\in p^{-1}b} S_a,\quad
T_b=\disjoint_{a\in p^{-1}b} T_a;
$$
$$F_b=\disjoint_{a\in p^{-1}b} F_a,
$$
so that $F_b:S_b\to T_b$. Let  $X_a\in
\B_1(F_a)$, $Y_b\in \B_1(F_b)$. Let
$$
R_p(\{X_a\}_{a\in A};\{Y_b\}_{b\in B}):
\hom_{B_1,p}(\{X_a\}_{a\in A};\{Y_b\}_{b\in
B})\to \hom_{B_2,p}(\{R_{F_a}(X_a)\}_{a\in
A};\{R_{F_b}(Y_b)\}_{b\in B})
$$
be the tensor product
$$
\otimes_{b\in B} R(\{X_a\}_{a\in
p^{-1}b};Y_b).
$$

Let, finally, $q:B\to C$ be another
surjection. Let $S_c=\disjoint_{a\in
(qp)^{-1}c}S_a$; let $T_c,F_c$ be similar
disjoint unions. Let $Z_c\in \B_1(F_c)$. We
then have the following  diagram:
$$\xymatrix{
\hom_{\B_1,p}(\{X_a\};\{Y_b\})\otimes\hom_{B_2,q}(\{Y_b\};\{Z_c\})\ar[r]^{C_1}\ar[d]^{R(\{X_a\};\{Y_b\})}
&
\hom_{\B_1,qp}(\{X_a\};\{Z_c\})\ar[d]^{R(\{X_a\};\{Z_c\})}\\
\hom_{\B_2,p}(\{R(X_a)\} ;\{R(Y_b)\})
\otimes\hom_{B_2,q}(\{R(Y_b)\};\{R(Z_c)\})
\ar[r]^>{C_1}&
\hom_{\B_2,qp}(\{R(X_a)\};\{R(Z_c)\}) }
$$

The axiom says that this diagram should be
commutative.

  3) Compatibility of  $R({\{X_a\}_{a\in A};Y})$
with compositions of the second type.

Let $F_a:S_a\to T_a$, $G_a:T_a\to R_a$ be
$A$-families of maps of finite sets. Let
$F=\disjoint_a F_a$ and $G=\disjoint_a G_a$.
Let $X_a\in \B_1(F_a)$, $Y_a\in \B_2(G_a)$.
Let
$$I_a:=I(G_a,F_a)(Y_a,X_a):R(Y_a)R(X_a)\to
R(Y_aX_a); $$
$$
I:R(Y)R(X)\to R(YX).
$$

We then have the  following diagram:
$$\xymatrix{
\hom_{\B_1}(\{X_a\};X)\otimes\hom_{\B_1}(\{Y_a\};Y)
\ar[r]\ar[d]& \hom_{\B_1}(\{Y_aX_a\};YX)\ar[d]\\
\hom_{\B_2}(\{R(X_a)\};R(X))\otimes\hom_{\B_2}(\{R(Y_a)\};R(Y))
\ar[d]& \hom_{\B_2}(\{R(Y_aX_a)\};R(YX))\ar[d]\\
\hom_{\B_2}(\{R(Y_a)R(X_a)\};R(Y)R(X))\ar[r]&\hom_{\B_2}(\{R(Y_a)R(X_a)\};R(YX))
}$$

The axiom requires the commutativity of this
diagram.

 \subsection{Pseudo-tensor structure on
 $\bodyas,\bodya,\bodyb$}
\subsubsection{$\bodyas$}
The pseudo-tensor body $\bodyas$ is
generated over the usual body $\bodyas$ by
the following generators and relations.

Generators: Let $f_k:R_\k\to T_k$, $k\in K$
be a family of surjections and $i_k:S_k\to
T_k$ be a family of injections. Let
$f=\disjoint_{k\in K} f_k$ and
$i=\disjoint_{k\in K} i_k.$ We then have a
generator
$$\fact(\{i_k,f_k\}_{k\in
K}):\{\p_{i_k}\R_{f_k}\}_{k\in K}\to
\p_i\R_f.
$$

Let $\pi:K\to L$ be a surjection.  For $l\in
L$ set
$$
i_l=\disjoint_{k\in \pi^{-1}l}i_k;
$$
$$
f_l=\disjoint_{k\in \pi^{-1}l}f_k.
$$
Set
$$
\fact(\{i_k,f_k\}_{k\in K},\{i_l,f_l\}_{l\in
L}):\{\p_{i_k}\R_{f_k}\}_{k\in K}\to
\{\p_{i_l}R_{f_l}\}_{l\in L}
$$
to be $$\otimes_{l\in L}
\fact(\{i_k,f_k\}_{k\in \pi^{-1}l};i_l,f_l).
$$

 Relations:

 1) Let $\sigma:L\to M$ be a surjection.
 For $m\in M$ set $$f_m=\disjoint_{k\in
 (\sigma\pi)^{-1}m} f_k;$$
 $$
i_m=\disjoint_{k\in
 (\sigma\pi)^{-1}m} i_k.
 $$
Then the composition of the first kind
$$
\{\p_{i_k}\R_{f_k}\}_{k\in K}\to
\{\p_{i_l}\R_{f_l}\}_{l\in L}\to
\{\p_{i_m}\R_{f_m}\}_{m\in M}
$$
equals $C_1(\{\p_{i_k}R_{f_k}\}_{k\in
K},\{\p_{i_l}R_{f_l}\}_{l\in L})$.

2) Let $f_k:S_k\to R_k$, $k\in K$ be
surjections and $i_k:R_k\to T_k$ be
injections. Let $Z_k:=T_k\backslash
i_k(R_k)$. Let $S'_k=S_k\disjoint Z_k$,
$R'_k=R_k\disjoint Z_k$; let $i'_k:S_k\to
S'_k$ be the natural inclusion. Let
$f'_k:S'_k\to T_k$, $f'_k=f_k\disjoint
i_{Z_k}$, where $i_{Z_k}:Z_k\to T_k$. We
then have  isomorphisms in $\Sk(i_kf_k)$:
$$
\R_{f_k}\p_{i_k}\to \p_{i'_k}\R_{f'_k}
$$
Let $f=\disjoint_k f_k$, $i=\disjoint_k
i_k$, $f'=\disjoint_k f'_k$, $i'=\disjoint_k
i'_k$. We then have an isomorphism in
$\Sk(if)$:
$$
\R_f\p_i\to \p_{i'}\R_{f'}.
$$

The relation says that the composition
$$
\{\R_{f_k}\p_{i_k}\}_{k\in K}\to
\{\R_f\p_i\}\to \{\p_{i'}\R_{f'}\}
$$
equals the following composition:
$$
\{\R_{f_k}\p_{i_k}\}_{k\in K}\to
\{\p_{i'_k}\R_{f'_k}\}\to \{\p_i\R_f\}
$$

3)  Let $f_k:R_k\to T_k$, $k\in K$ be
surjections. Let $f_k=g_kh_k$, where
$g_k,h_k$ are surjections. Let
$f=\disjoint_{k\in K} f_k$,
$g=\disjoint_{k\in K} g_k$,
$h=\disjoint_{k\in K} h_k$.

Then the composition
$$
\{\R_{f_k}\}_{k\in K}\to \R_f\to \R_h\R_g
$$
equals the composition
$$
\{\R_{f_k}\}_{k\in K}\to
\{\R_{h_k}\R_{g_k}\}_{k\in K}\to \R_h\R_g.
$$

\subsubsection{$\bodya$} The
pseudo-tensor structure on $\bodya$ is
generated by the same generators as on
$\bodyas$, and the relations  include those
in $\bodyas$ with an addition of the
following relations:

a) let $i_k,p_k,j_k,q_k$, $k\in K$ be a
collection of suitable squares. Let
$i=\disjoint_{k\in K} i_k$,
$p=\disjoint_{k\in K} p_k$,
$j=\disjoint_{k\in K} j_k$,
$q=\disjoint_{k\in K} q_k$. Then the square
$i,p,j,q$ is also suitable and the following
compositions coincide:
$$
\{\p_{i_k}\R_{p_k}\}_{k\in K}\to
\{R_{q_k}\p_{j_k}\}_{k\in K}\to \R_q\p_j
$$
and
$$
\{\p_{i_k}\R_{p_k}\}\to \p_i\R_p\to
\R_q\p_j.
$$

b) Let $i_k:S_k\to R_k$, $k\in K$ be
injections and $p_k:R_k\to T_k$, $k\in K$ be
surjections such that $j_k:=p_ki_k$ are
injections. Assume that at least two of the
maps $p_k$ are proper surjections. Then the
composition
$$
\{\p_{i_k}\R_{p_k}\}_{k\in K}\to \p_i\R_p\to
\R_j
$$
vanishes.

If only one of the surjections $p_k$ is
proper, say $p_\kappa$, $\kappa\in K$, then
the above composition equals
$$
\{\p_{i_k}\R_{p_k}\}_{k\in K}=
\{\p_{i_\kappa}\R_{p_\kappa},\p_{i_k}\R_{\{p_k\}_{k\in
K}} \}\stackrel{L(i_\kappa,p_\kappa)}\to
\{\p_{j_\kappa},\{\p_{i_k}\R_{p_k}\}_{k\in
K}\}\to \{p_{j_k}\}_{k\in K}\to \p_j.
$$

\subsubsection{$\bodyb$} This pseudo-tensor body is
generated by the same generators and
relations as $\bodya$ except that we lift
the condition of suitability. We have a
natural map
\begin{equation}\label{bodyatobodyb}
\bodya\to \bodyb
\end{equation}
\subsubsection{} It is clear that the system
$<\R>$ determines a map of pseudo-tensor
bodies
$$
\bodya\to \Com,
$$
any such a functor will be also called {\em
representation}
\subsection{Explicit form of pseudo-tensor
maps}
\subsubsection{Category of special maps}
Consider a family of objects
$$
X_a=\p_{i^a}\R_{p_1^a}\R_{p_2^a}\cdots
\R_{p_{n_a}^a}
$$
indexed by a finite set $A$, where all
$p^a_k$ are proper surjections and
$$
p_{n_a}^a\cdots p_2^ap_1^ai^a=F_a
$$
so that $X_a\in \Sk(F_a)$. Let $N\geq j\geq
i\geq 1$. Set $p_{ji}^a=p_j^ap_{j-1}^a\cdots
p_{i+1}^a$ (if $i=j$, then we set
$p_{ji}^a=\Id$).
 Let $\bfu:=\{u^a_k\}$ ,
where $a\in A$, $k=0,1,2,\ldots,N$, be a
sequence of numbers satisfying: $u^a_0=0$,
$$0\leq u^a_{k+1}-u^a_k\leq 1,
$$
and $u^a_N=n_a$. Set
$$p_k(\bfu):=\disjoint_{a\in A}
p^a_{u^a_ku^a_{k-1}}
$$
Call $\bfu$ {\em proper} if such are all
$p_k(\bfu)$.

For proper $\bfu$ we set
$$
X(\bfu):=p_{\disjoint_a
i_a}\R_{p_1(\bfu)}\R_{p_2(\bfu)}\cdots
\R_{p_N(\bfu)}.
$$

We have natural maps
$$\fact(\bfu):\{X_a\}_{a\in A}\to X(\bfu).
$$

For an $Y\in \Sk(\disjoint_a F_a)$ consider
the groupoid $G_Y$ whose objects are
collections $$(\{X_a\in \Sk(F_a)\}_{a\in A},
\bfu),m:Y\to X(\bfu)$$ where the meaning of
the ingredients is the same as above and $m$
is an isomorphism; the isomorphisms in $G_Y$
are isomorphisms of such collections. It is
clear that $G_Y$ is a trivial groupoid. Let
$Z_a\in \Sk(f_a)$. We have a natural map
$$
\liminv_{\{X_a\}_{a\in A}\in
G_Y}\otimes_{a\in A}\hom_?(Z_a,X_a)\to
\hom_?(\{Z_a\}_{a\in A},Y).
$$
We claim that this map is  an isomorphism,
where $?=\bodyas,\bodya,\bodyb$.
\subsubsection{} We shall also need another
form of decomposition of the pseudo-tensor
maps in $\bodyb$.

Let $\{X_a\}_{a\in A}$ be a family of
objects $X_a\in \Sk(F_a)$. Let
$F=\disjoint_{a\in A} F_a$. We then have the
following natural functors
$$h_{\{X_a\}_{a\in A}}^\Asm:\bodya(F)\to\complexes
$$
and
$$h_{\{X_a\}_{a\in A}}^\Bsm:\bodyb(F)\to\complexes
$$
defined by the  formulas:
$$
h_{\{X_a\}_{a\in
A}}^\Asm(Y)=\hom_{\bodya}(\{X_a\}_{a\in
A};Y);
$$
$$
h_{\{X_a\}_{a\in
A}}^\Bsm(Y)=\hom_{\bodyb}(\{X_a\}_{a\in
A};Y)
$$

We shall also need a functor $$
G^\Bsm:\bodya^\op(F)\otimes\bodyb(F)\to
\complexes,
$$
where $G^\Bsm(Z,U)=\hom_{\bodyb}(Z,U)$. We
then have a natural map:
$$
h_{\{X_a\}_{a\in A}}^\Asm\otimes_{\bodya(F)}
G^\Bsm\to h_{\{X_a\}_{a\in A}}^\Bsm
$$
\begin{Lemma} \label{ismph}
\label{if2if1} This map is an isomorphism of
functors.
\end{Lemma}
\pf Straightforward
\endpf

\subsection{Linear span of a body}

Let $\B$ be a pseudo-tensor body. We shall
construct a body $\lB$, over another
skeleton, as follows. Set $\lB(F)$ to be the
category of functors $\B(F)^\op\to
\complexes$.  We shall start with the
composition maps
$$
\circ:=\circ(G,F):\lB(F)\otimes \lB(G)\to
\lB(GF).
$$

Introduce an auxiliary functor
$$\D:=\D^{G,F}:\B(F)\otimes \B(G)\otimes
\B(GF)^\op\to complexes,
$$
where $$\D(X,Y,Z)=\hom_{\B(GF)}(Z,YX).
$$

Let now $U\in \lB(F)$, $V\in \lB(G)$. Define
$$
V\circ U:=\D\otimes_{\B(G)\otimes
\B(F)}V\boxtimes U.
$$

Let us  construct the associativity map.
Define
$$\D_3:=\D^{H,G,F}:\B(H)\otimes \B(G)\otimes \B(F)\otimes
\B(HGF)^\op\to complexes$$ by
$$
\D_3(Z,Y,X,U):=\hom_{\B(HGF)}(U,ZYX)
$$.

Ioneda's lemma combained with the
associativity maps implies  isomorphisms
$$
\D^{G,F}\otimes_{\B(GF)}\D^{H,GF}\stackrel\sim\to
\D_3;
$$
$$
\D^{H,G}\otimes_{\B(HG)}\D^{HG,F}\stackrel\sim\to
\D_3.
$$
Let $U\in \B(F)$, $V\in \B(G)$, and $W\in
\B(H)$. Set $$ (WVU):=W\boxtimes V\boxtimes
U\otimes_{\B(H)\otimes \B(G)\otimes \B(F)}
D_3.
$$

We then have isomorphisms
$$
(WV)U\stackrel\sim\to
(WVU)\stackrel\sim\rightarrow W(VU),
$$
which furnish the associativity isomorphism.

For $X_i\in \lB(F_i)$, $i\in I$ and $Y\in
\lB(\disjoint_{i\in I} F_i)$ set
$$
\hom_{\lB}(\{X_i\}_{i\in I};
Y):=\hom(\boxtimes_{i\in I}X_i,
\hom_{\B}(\{.\},.)\otimes_{\B_F}Y).
$$

Define the compositions of the first kind.
Let $\pi:I\to J$ be a surjection. Let
$F_i:S_i\to T_i$ be a family
 of maps of finite sets. Let
 $F_j=\disjoint_{i\in \pi^{-1}j} F_i$.
 Let $K_i\in \lB(F_i)$ and $L_j\in
 \lB(F_j)$.
 Let
 $$A(\pi):\otimes_{i\in
 I}\B(F_i)^\op\otimes_{j\in J}\B(F_j)\to
 complexes
 $$
 be given by:
 $$
 A(\pi)(\{X_i\};\{Y_j\})=\hom_{\B}(\{X_i\};\{Y_j\}).
 $$
 Let $\B_I:=\otimes_{i\in I}
\B(F_i)$; $\B_J:=\otimes_{j\in J} \B(F_j)$.
 Let $K:\B_I^\op\to \complexes$ be
 $\boxtimes_{i\in I} K_i$.
 Let $L:\B_J^\op\to \complexes$ be
 $\boxtimes_{j\in J} L_j$.
 We then have:
$$
\hom_{\lB,\pi}(\{K_i\}_{i\in
I};\{L_j\})\cong \hom_{\B_I^\op}(K,
L\otimes_{\B_J}A(\pi)).
$$

Let $\sigma:J\to K$ be the third surjection.
For $k\in K$, let $F_k:=\disjoint_{i\in
(\sigma\pi)^{-1}k} F_i$. Let $M_k\in
\B(F_k)$. Let $\B_K=\otimes_{k\in
K}\B(F_k)$; let $M:=\boxtimes_{k\in K} M_k$.
Then $$
\hom_{\lB,\sigma}(\{L_j\};\{M_k\})\cong
\hom_{\B_j^\op}(L,
M\otimes_{\B_K}A(\sigma)).
$$

To construct the composition of the first
kind we shall also need  an isomorphism
$$A(\sigma)\otimes_{\B_J} A(\pi)\to
A(\sigma\pi),
$$
where the isomorphism follows from the
Ioneda's lemma.

In view of the above isomorphisms, the
composition of the first kind reduces to:
\begin{eqnarray*}
\hom_{\B_I^\op}(K,L\otimes_{\B_J}
A(\pi))\otimes
\hom_{B_J^\op}(L,M\otimes_{\B_K}
A(\sigma))\to
\hom_{\B_I^\op}(K,M\otimes_{\B_K}
A(\sigma)\otimes_{\B_J}A(\pi))\\
\cong
\hom_{\B_I^\op}(K,M\otimes_{\B_k}
A(\sigma\pi)).
\end{eqnarray*}

Lastly, let us define the compositions of
the second kind. We shall keep the above
notation. Let $G_i:T_i\to R_i$ be another
family of maps of finite sets. Let $K'_i\in
\lB(G_i)$ and $L'_j\in \lB(G_j)$. Let $ a\in
\hom(K,L)$ and $a'\in \hom(K',L')$. Let
$A(\pi)$ be as above and let $A'(\pi)$
(resp. $A''(\pi)$) be the same as $A(\pi)$
but $F_i$ are all replaced with $G_i$ (resp.
$G_iF_i$). Let $\B'_I=\otimes_{i\in I}
\B(G_i)$; $\B''_I=\otimes_{i\in
I}\B(G_iF_i)$
 Construct the
composition $a'a\in \hom(K'K,L'L)$.

As was mentioned above, $a$ determines a map
$$
a:K\to L\otimes_{B_J} A(\pi),
$$
and $a'$ produces a map
$$
a':K'\to L'\otimes_{B'_J} A'(\pi).
$$

To construct the compositions $K'K,L'L$,
introduce functors

$$
O_I:\B_I\otimes\B'_I\otimes(\B''_I)^\op\to
\complexes;
$$

$$
O_J:\B_J\otimes\B'_J\otimes(\B''_J)^\op\to
\complexes;
$$
by setting
$$
O_I(\{X_i\};\{X'_i\};\{X''_i\}):=\hom(\{X''_i\};
\{X'_iX_i\});
$$
$$
O_J(\{X_j\};\{X'_j\};\{X''_j\}):=\hom(\{X''_j\};
\{X'_jX_j\}).
$$
Then $$K'K=K'\boxtimes K\otimes_{B_I\otimes
B'_I}O_I;
$$
$$L'L=L'\boxtimes L\otimes_{B_J\otimes
B'_J}O_J;
$$

Applicaton of  $a,a'$ yields a  map:
\begin{eqnarray*}
K'\boxtimes K\otimes_{\B_I\otimes
\B'_I}O_I\to L'\boxtimes
L\otimes_{\B_J\otimes \B'_J}
(A(\pi)\boxtimes
A'(\pi))\otimes_{\B_I\otimes\B'_I} O_I
\end{eqnarray*}

Next, by Ioneda's lemma, we have an
isomorphism.
$$
(A(\pi)\boxtimes
A'(\pi))\otimes_{\B_I\otimes\B'_I} O_I\to
O_J.
$$

If we apply this isomorphism to the previous
map, we will get the desired sekond kind
composition map:
$$
K'\boxtimes K\otimes_{\B_I\otimes
\B'_I}O_I\to L'\boxtimes
L\otimes_{\B_J\otimes \B'_J}O_J.
$$

This concludes the definition of the
structure. Checking the axioms is
straightforward.
\subsubsection{} A representation of a pseudo-tensor body
$B$ (i.e. a map $B\to \Com$) naturally
extends to a representation of $\lB$.

\subsection{Representation of a body in another body}
An arbitrary map of bodies $\B_1\to k[\B_2]$
will be called {\em a representation of
$\B_1$ in $\B_2$} We shall construct
\subsection{Representation of $\bodyb$ in
$\bodya$} By constructing such a
representation, we shall automatically
obtain a  map $\bodyb\to \Com$, i.e. a
symmetric system.

 First of all construct maps
$R_F:\bodyb(F)\to\lbodya(F)$ by assigning
$$
R_F(X)(Y):=\hom_{\bodyb(F)}(Y,X),
$$
On $R_F(X)$, we have a natural structure of
 functor from the category $\bodyb(F)^\op$ to the category
 $\complexes$
 given by the map $\bodya\to
\bodyb$ as in (\ref{bodyatobodyb}). Let
$X_i\in \bodya(F_i)$. We then have a natural
map
$$
I(F_2,F_1):R_{F_1}(X_1)R_{F_2}(X_2)\to
R_{F_2F_1}(X_1X_2).
$$
given by:
\begin{eqnarray*}
(R_{F_1}(X_1)R_{F_2}(X_2))(Z)=R_{F_1}\otimes
R_{F_2}\otimes_{\B(F_1)\otimes
\B(F_2)}\D^{F_2F_1}_{\bodya}(Z)\\
\to
R_{F_1}\otimes
R_{F_2}\otimes_{\B(F_1)\otimes
\B(F_2)}\D^{F_2F_1}_\bodyb(Z)\cong
R_{F_2F_1}(X_1X_2)(Z).
\end{eqnarray*}
Furthermore, as follows from the
decomposition (\ref{if2if1}), $I(F_2,F_1)$
is an isomorphism.

To defined the maps $$R_{\{X_a\}_{a\in
A};Y}: \hom_{\bodyb}(\{X_a\}_{a\in A};Y)\to
\hom_{\bodya}(\{R_{F_a}(X_a)\}_{a\in A}
;R_F(Y))
$$
we shall use Lemma \ref{ismph}. We have
$$
R_F(X)(Y)=G^\Bsm(Y,X),
$$ where $G^\Bsm$ is as in the statement of Lemma \ref{ismph}.
Let $$h^\Asm:\otimes_{a\in A}
\bodya(F_a)^{\op}\otimes\bodya(F)\to
\complexes.$$ be defined by the formula
$$
h^\Asm(\{X_a\}_{a\in
A};Y):=\hom_{\bodya}(\{X_a\}_{a\in A};Y).
$$
Then, by definition,
$$
\hom_{\bodya}(\{R_{F_a}(X_a)\}_{a\in A}
;R_F(Y))= \hom_{\otimes_{a\in
A}\bodya(F_a)^\op} (\boxtimes_{a\in
A}R_F(X_a); h^\Asm\otimes_{\bodya(F)}
G^\Bsm)$$

The latter term is, by Lemma \ref{ismph},
isomorphic to
$$
\hom_{\otimes_{a\in A}\bodya(F_a)^\op}
(\boxtimes_{a\in A}R_F(X_a); h^\Bsm),
$$
where
$$
h^\Bsm(\{X_a\}_{a\in
A};Y):=\hom_\bodyb(\{X_a\};Y).
$$
Lastly, we have a natural map
$$
\hom_{\otimes_{a\in A}\bodyb(F_a)^\op}
(\boxtimes_{a\in A}R_F(X_a); h^\Bsm)\to
\hom_ {\otimes_{a\in A}\bodya(F_a)^\op}
(\boxtimes_{a\in A}R_F(X_a); h^\Bsm),
$$
and the first space is, by Ioneda's lemma,
isomorphic to
$$
\hom_\bodyb(\{X_a\};Y).
$$
This completes the desired construction.
Checking the axioms is straightforward.
\subsubsection{}
As was mentioned above, the above
construction provides us with a symmetric
system. Denote it $<\N>$.  An explicit
construction of $<\N>$ is given in
\ref{explicitrsymm}. Checking that this
construction  produces the same system as in
the previous section is straightforward, and
we omit it.

\def\Re{{\Bbb R}}
\def\fact{{\text{fact}}}
\def\bodya{{\cal B}_{\text{\bf presymm}}}
\def\bodyb{{\cal B}_{\text{\bf symm}}}
\def\follows{\Rightarrow}
\def\disjoint{{\sqcup}}
\def\bs{\backslash}
\def\A{{\frak A}}
\def\intr{\int}
\def\ins{{\frak ins}}
\def\inr{{\frak int}}
\def\la{\lambda}
\def\An{{\frak An}}
\def\P{{\cal P}}
\def\limdir{\text{limdir}}
\def\ev{\text{ev}}
\def\notsubset{\text{notsubset}}
\def\identical{=}
\def\ins{{\frak ins}}
\def\dirlim{\limdir}
\def\usual{{\text{usual}}}
\def\A{{\cal A}}
\def\B{{\cal B}}
\def\F{{\frak F}}
\def\a{{\frak a}}
\def\b{{\frak b}}
\def\e{{\frak e}}
\def\jk{{\frak jk}}
\def \q{{\frak q}}
\def \qq{\q}
\def\sym{{\text{sym}}}
\def\res{{\text{res}}}
\def\f{{\fun}^{\sim}}
\def \Im{\text{Im}}
\def\jj{{\frak jj}}
\def\fp{({\fun'})^{\sim}}
\def\r{{\frak r}}
\def \into{\to}
\def \comm{\text{\bf comm}}
\def\lis{{\bf lis}}
\def\cat{{\bf cat}}
\def\red{\text{red}}
\def\lism{{\bf lism}}
\def\complexes{\text{complexes}}
\def\op{\text{op}}
\def \ores{{\O^{\text{res}}}}
\def\asm{{\frak ope'}}
\def\ens{{\bf ens}}
\def\ske{{\bf ske}}
\def\ve{\varepsilon}
\def\vect{{\bf vect}}
\def \id{\text{Id}}
\def\Id{\id}
\def \as{{\frak as}}
\def\ope{{\frak ope}}
\def\fun{{\frak fun}}
\def\funp{({\frak fun}^{\psys})}
\def \Fun{{\frak Fun}}
\def\funct{{\frak Funct}}
\def\system{{{\cal SYS}}}
\def\psys{{\cal SYS}[p]}

\def\lie{{\frak lie}}
\def\li{\lie^1}
\def \ll{{{\li}_M}}
\def\psystem{{{\cal SYS}[p,\ll]}}
\def\psystemg{{{\cal SYS}[p,\lie_M]}}
\def\OPE{{\cal OPE}}

\def \OPER{{\cal OPE}^{\text{res}}}
\def\psystemb{\OPE[p,\lie_M]}
\def\psystemg{\psystemb}
\def\psis{\psystemb^{\text{res}}}

\def\univ{\text{univ}}
\def\psim{\psystemb^{\text{sym}}}
\def\psimr{\psystem^{\text{sym},\text{res}}}

\def\lisp{{[p,\lie_M,\S]}}
\def\lispq{{\lisp^{\text{res}}}}
\def\psml{\system[p,\lie_M,\S]}
\def\psmq{\psml^{\text{res}}}
\def\psmlh{\psml[[\h]]}
\def\psmqh{\psmq[[\h]]}
\def\coker{\text{Coker}}
\def\M{{\cal M}}
\def\gen{{\frak gen}}
\def\u{{\frak u}}
\def\Mg{{\frak M}}
\def \hh{{\frak h}}
\def \kk{{\frak kk}}
\def \II{{\frak II}}
\def \JJ{{\frak JJ}}
\def \GG{{\frak GG}}
\def \FF{{\frak FF}}
\def \m{{\frak m}}
\def \z{{\frak z}}
\def \B{{\cal B}}
\def \G{{\cal G}}
\def \M{{\cal M}}
\def \O{{\cal O}}
\def \R{{\Bbb R}}
\def \Z{{\frak Z}}
\def \Co{{\Bbb C}}
\def \C{{\cal C}}
\def \i{{\frak i}}
\def \k{{\frak k}}
\def \ii{{\frak ii}}
\def \l{{\frak l}}

\def \L{{\frak L}}
\def \DDD{{\cal D}}
\def \D{{\cal D}}
\def \DD{\Delta}
\def \g{{\frak g}}
\def \ha{^\wedge}

\def \liminv{\text{liminv}}
\def \p{{\frak p}}
\def \pp{{\pi_p}}
\def \h{h}
\def \I{{\cal I}}
\def \s{{\frak s}}
\def \S{{\frak S}}
\def \U{{\cal U}}
\def \vs{\sigma}
\def\pt{\text{pt}}
\def\Hom{\text{Hom}}
\def \holi{{\frak holi}}
\def \La{{\holi_M}}
\def \Lab{{ {\frak holie}_M}}
\def \psl{\OPER[p,\Lab]}
\def \ps {\psl}
\def\ov{\overline}
\def\comp{{\bf comp}}

\def\N{\Rs}

\def\D{{\frak D}}
\def\DD{{\cal D}}
\section{realization of the system $<\N>$ in the spaces
of real-analytic functions}
\subsection{Conventions and notation}
We do not consider sheaves in this sections,
but only their global sections. By
$\DD_{X^S}$ we denote the algebra  of
polynomial differential operators on $X^S$.
By a $\DD_{X^S}$-module we mean a module
over the algebra $\DD_{X^S}$.

We denote by $\D_{Y^S}$ the space of
compactly supported infinitely-differntiable
top-forms on $Y^S$, and by $\D'_{Y^S}$ the
space of distributions (=generalized
functions) on $\D_{Y^S}$. $\D'_{Y^S}$ is a
left $\D_{X^S}$-module.

For simplicity, we fix a translation
invariant top form $\omega$ on $Y^S$, and
define $\omega_S$ to be a top form on $Y^S$
which is the exterior product of copies of
omega. Because $\dim Y$ is even, the order
in this product does not matter.

The space $\D_{Y^S}$ is then identified with
the space of compactly supported infinitely
differentiable functions on $Y^S$.

\subsection{Asymptotic decompositions of functions from
$\C_S$}
\subsubsection{The main theorem}
Let $S$ be a finite set. Let $T\subset S$ be
a subset. Let $R:=S\bs T$. Pick an element
$\tau\in T$. We shall refer to a point of
$Y^S$ as $(\{y_s\}_{s\in S})$, where $y\in
Y$.

For a positive real $\la$ set
$$
U_\la((\{y_s\}_{s\in S}) =
(\{y_\tau+\frac{y_t-y_\tau}\la\}_{t\in
T};\{y_r\}_{r\in R}).
$$
This determines an action of the Lie group
$\Re_{>0}^\times$ on $Y^S$.

Let $F\in \C_S$.
\begin{Claim}
For every $g\in \D_{Y^S}$ there exist:

constants $A(F)$, $B(F)$;

distributions $C_{n,k}^F\in \D'_{Y^S}$, for
every $n\geq A(F)$ and every $k$ such that
$0\leq k\leq B(F)$;

 such that for every $N$ and every $g\in \D_{Y^S}$,
 the following asymptotics takes place:

 $$
 <F, U_\la g>=\sum\limits_{ A(F)\leq n\leq N,0\leq k\leq B(F)}
 C_{n,k}^F(g)\la^n(\ln \la)^k+o(\la^N).
 $$

 \end{Claim}

\pf We shall use induction in $\#R$ to prove
even stronger statement:

There exist:

constants $A(F)$, $B(F), K(F)$;

distributions $C_{n,k}^F$ on the space of
compactly supported $K(F)$-times
differentiable functions, for every $n\geq
A(F)$ and every $k$ such that $0\leq k\leq
B(F)$;

 such that for every $N$ there exists a constant
 $L:=L(N,F)$ such that whenever
 $$\phi\in C^L(Y^S), g\in C^L_c(Y^S),$$

 the following asymptotics takes place:

 \begin{equation}\label{asympt1}
 <F, \phi U_\la g>=\sum\limits_{ A(F)\leq n\leq N,0\leq k\leq B(F)}
 C_{n,k}^F((U_{\la^{-1}}\phi)g)\la^n(\ln \la)^k+o(\la^N).
 \end{equation}

{\bf Remark} We have
$$
U_{\la^{-1}}\phi(\{x_s\}_{s\in S})=
\phi(\{x_\tau+\la(x_t-x_\tau)\}_{t\in
T},\{x_r\}_{r\in R}),
$$
Therefore, for $L$ large enough, we can
replace $U_{\la^{-1}}\phi$ in
(\ref{asympt1}) with  a finite sum
$$
\sum \la^k\phi_k.$$

Base: $\#R=0$.  We then have $<F,U_\la
g>=<U_\la^*F,g>$. Let us study the action
$U_\la^*$ on $\C_S$. It is clear that this
action preserves the filtration on $\C_S$
and that the associated graded action is
diagonalizable. It then follows that for
every $F\in C_S$,
$$
U_\la^*F=\sum_{n,k} \la^n(\ln\la)^k F_{nk},
$$
where the sum is finite and $F_{nk}\in C_S$.

The statement now follows immediately.

Let now $R$ be arbitrary, and assume that
the statement is the case whenever $R$ has a
smaller number of elements.

Let $R_1\subset R$ be an arbitrary non-empty
subset. Let $R_2=S\bs  R_1$. Assume that
\begin{equation}\label{splitgen}
F= F_1F_2,
\end{equation}
where $F_i\in C_{R_i}$.

We then claim that the required asymptotics
is the case. Indeed, we have $$ <F_1F_2,\phi
U_\la g>=<F_1,<F_2,\phi U_\la g>>
$$

and the statement follows from the
corresponding statement for $F_2$ (which
holds in virtue of the induction
assumption).

Let us generalize this result. Let
$I^{N,K}_{R',R''}\subset C^K(Y^S)$ be the
subspace consisting of functions
 which vanish
on each diagonal $x_{r_1}=x_{r_2}$, $r_i\in
R_i$ up-to the order $N$.

Let $$Q_{R_1R_2}:=\prod_{r_i\in
R_i}q(x_{r_1}-x_{r_2}).
$$
It is not hard to see that for every $M,L$
there exist $N,K$ such that  for every
$\chi\in I^{N,K}_{R',R''}$ we have:
$$
\chi=Q^M\psi,$$ where $\psi\in C^M(Y^S)$.
Thus, for $N,K$ sufficiently large, and
$\chi\in I^{N,K}$ we have
$$
<F,\chi\phi U_\la g>=<FQ^M, \psi\phi U_\la
g>.
$$

But for $M$ large enough, $FQ^M$ splits into
a sum of elements of the form
(\ref{splitgen}), whence the statement for
all $\chi\in I^{N,K}_{R',R''}$.

Let us now define the space $J^{N,K}\subset
C^K{Y^S}$, consisting of all functions which
vanish on the diagonal
$$
\forall t\in T:x_t=x_\tau
$$
up-to the order $N$. It is not hard to see
that for $N,K$ large enough,
$$
<F,\chi\phi U_\la g>=o(\la^n).
$$

Therefore, there are large enough $N,K$ such
that whenever
$$
\phi\in \sum_{R'} I^{N,K}_{R',R''} +J^{N,K},
$$
the required asymptotics holds.

Let $A^{N,K}\subset C^K(Y^S)$ be the
subspace of functions which vanish on the
main diagonal in $Y^S$ up-to the order $N$.

By the Nullstellensatz,  for some $N',K'$,
$$
A^{N',K'}\subset \sum_{R'} I^{N,K}_{R',R''}
+J^{N,K}.
$$

Therefore, the required asymptotics holds
whenever
$$
\chi\in A^{N',K'}.
$$

Let us now pass to the original statement.

Let
$$
\Xi_l= \sum_{t\in
T}(x_t-x_\tau)\frac{\partial}{\partial x_t}
$$
$$
\Xi_r=\sum_{r\in R}(x_r-x_\tau)
\frac{\partial}{\partial x_r}.
$$

The action of the vector field $\Xi_l+\Xi_r$
on the space $\C_S$ preserves the
filtration, and the induced action on the
associated graded quotients is
diagonalizeable, therefore we may assume
that $(\Xi_l+\Xi_r-n)^NF=0$ for some $n,N$.

Consider expressions
\begin{equation}\label{pln}
<F,P(\Xi_l,\Xi_r,\la\frac d{d\la})\phi U_\la
g)>,
\end{equation}
where $P$ is a polynomial.

Let $U_M(z)=z(z-1)(z-2)\cdots (z-M)$.

Consider the following ideals in the ring of
polynomials of three variables:

$$
A_M=(U_M(\Xi_l-\la\frac d{d\la})); $$
$$
B_M=(U_M(\Xi_r));
$$
$$
C=((\Xi_l+\Xi_r-n)^N).
$$

It is not hard to see that for $M$ large
enough, whenever $P$ is large enough, the
expression
 (\ref{pln}) has the required asymptotics.

 Indeed, consider for example the ideal
 $A_M$.  We have:
 $$ U_M(\Xi_l-\la\frac
d{d\la})\phi U_\la(g)=
(U_M(\Xi_l)\phi)U_\la(g),
$$

and it is easy to see that
$\chi:=U_M(\Xi_l)\phi$ has  at least the
$M$-th order of vanishing along the main
diagonal, whence the statement. The ideals
$B_M,C_M$ can be checked in a similar way.

 Next, we see, by the Nullstellensatz,
  that  for some $M,L$
 $$U_M(\la\frac d{d\la}-n)^L\in A_M+B_M+C.$$

 Therefore, we see that
 there exists the required asymptotics for
 $$
 U_M(\la\frac d{d\la}-n)^L<F, \phi U_\la g>.
 $$
 The theory of ordinary
 differential equations now implies the statement.
 \endpf

 \subsubsection{A claim about the
 distributions $C_{n,k}$}
 Let $G$ be a function on $Y^T$ which is invariant
 under translations by a vector from $Y$, with support
 compact  modulo the action of $Y$.

 Let $H$ be a function on $Y^{R\disjoint \{\tau\}}$
 with compact support.

 We then have $U_\la(GH)=HU_\la(G)$.
 \begin{Claim}
 We have:
 $$
 C_{n,k}(HG)=D_{n,k}(G)(H),
 $$
 where $D_{n,k}(G)\in \C_{R\disjoint \{\tau\}}$.

 Furthermore, for every $N$,  the distributions
 $D_{n,k}(G)$, where  $n\leq N$ and $k$ is arbitrary,
 span a finitely dimensional vector subspace.
 \end{Claim}

\pf Use induction. If $R$ is empty, there is
nothing to prove.

Otherwise, let us split $R=R_1\disjoint
R_2$, in a non-trivial way.

 We then see that  for $M$ large enough
$D_{n,k}(G Q_{R_1\disjoint T,R_2}^M)$
satisfy the statement in virtue of the
induction statement.

Also, for $M'$ large enough and any $G$
vanishing on the diagonal $\forall t\in T:
y_t=y_\tau$ up to the order $M'$,
 $D_{n,k}(G)=0$.
 The  Nullstellensatz
then implies that for $L$ large enough, one
can write
$$Q_{R_1,R_2}^L=
P_1Q_{T\disjoint R_1, R_2}^M+P_2
$$ where $P_1,P_2$ are polynomials and  $P_2$
vanishes  on the diagonal $\forall t\in T:
y_t=y_\tau$ up-to the order $M'$.

This implies that $D_{n,k}(Q_{R_1,R_2}^LG)=
D_{n,k}(G)Q_{R_1,R_2}^L$ satisfy the
statement.

Also, $D_{n,k}(G)$ are all translation
invariant, therefore $D_{n,k}(G)\in
\C_{R\disjoint \{\tau\}}$.

Let $\C'$ be the quotient of $\C_{R\disjoint
\{\tau\}}$ by the distributions supported on
the main diagonal.

It then follows that $D_{n,k}$ span a
finitely dimensional space in $\C'$.

Let $Q_R=\prod_{r,s\in R,r\neq
s}q(x_r-x_s)$. Then $D_{n,k}(G)Q_R$ also
span a finitely dimensional space. Hence,
$D_{n,k}(G)$ span a finitely-dimensional
space in  $\C_{R\disjoint \{\tau\}}$.
\def\ovt{\overline{\tau}}
\subsubsection{} Consider a
decomposition $S=S_1\disjoint S_2$ such that
$T\subset S_1$.  Consider an element $F\in
\C_S$ which decomposes as a product
$F=F_1F_2$, where $F_i\in \C_{S_i}$. We are
going to express $D_{n,k}^F$ in terms of
$D_{n,k}^{F_1}$.

Let $G$ be as above (i.e. an infinitely
differentiable function on $Y^T$ invariant
under shifts  by $Y$ and with compact
support modulo  these shifts.
\begin{Claim}\label{razl}
We then have
$$
D_{n,k}^F(G)(H)=<F_2,D_{n,k}^{F_1}(G)(H)>.
$$
\end{Claim}
\pf {Clear}
\endpf
\subsubsection{}
Let $S$  be a finite set with a marked point
$\vs\in S$.
 Let $\Xi_r $ be the dilation vector field
on $X^S$ given by:
$$
\sum\limits_{s\in S}
(x_s-x_\vs)\frac{\partial}{\partial x_s}.
$$

Denote $\D'_{Y^S,n}$  the generalized
eigenspace of $\Xi_r$ with eigenvalue $n$.
Let $\C_{S,n}:=\C_S\cap \D'_{Y^S,n}$.

We know that $$\C_S=\oplus_{n\in \Bbb Z}
\C_{S,n}.$$

Let us now come back to our situation in
which we have a finite set $S$, its subset
$T$ and a marked point $\tau\in T$.

Consider a subspace $\D'_{T,n}\subset
\D'_{Y^T,n}$ consisting of all elements
which are nilpotent under translations by
$Y$.

It is then not hard to see that
\begin{Lemma}
$$D_{n,k}\in \D'_{T,n}\otimes_{\O_X} \C_{R\disjoint \{\tau\},N-n}
$$

where the $\O_X$-action is on the $\tau$-th
components of both tensor
 factors.
\end{Lemma}
\subsubsection{} Let
$$
\D'_{T,\geq n}:=\oplus_{N\geq n}\D'_{T,N};
$$
let
$$\D'=\oplus_{N} \D'_{T}.$$

Let $$E(S,T):=\liminv_{n\to\infty}
 \D'/\D'_{T,\geq n}\otimes C_{R\disjoint \{\tau\},N-n}.
 $$

Given a function $G\in \D_{Y^S}$ and an
element
$$s_n\in  \D'_{T,\geq n}\otimes \C_{R\disjoint \{\tau\},N-n}
$$
one has: $<s_n,U_\la G> =o(\la^{n-1})$.
Therefore,  for every  $s\in E(S,T)$ and
$G\in \D_{Y^S}$  we have an asymptotic
series
$$
<s,U_\la G>.
$$

 \def\eps{\epsilon}
 \begin{Claim} There exists a map
 $$
 \eps:C_S\to E(S,T)$$
 uniquely determined by the condition that
 $<\eps(F),U_\la G> $ is an asymptotic series for
 $<F,U_\la G>$.
 \end{Claim}
 \subsubsection{} Let $\pi:S\to S/T$ be  the natural surjection.
 Define a functor $\A_{\pi}$ from the category
 of $\D_{X^{S/T}}$-modules to the category
 of $\D_{X^S}$-modules  by the formula

 $$
 \A_{\pi}(M)=
 \liminv_{n\to\infty}
 i_\pi\ha(M)\otimes_{\O_{X^T}} \D'/\D'_{T,\geq n}.
$$

Then the above result can be rewritten as a
map
$$\C_S\to \A_{\pi}(\C_{S/T}).
$$

\subsubsection{} Let $q:S/T\to P$ be an arbitrary
surjection. Let $\ovt\in S/T$ be the image
of $T$. Let $\chi=\pi(\ovt)$. For $p\in P$
set $S_p:=(q\pi)^{-1}p$. Let $\vs:S\to
S/S_chi$ be the natural projection.  We then
have induced maps $q_1:S/S_\chi\to P$;

\begin{Lemma}\label{lemmarazl} The  composition
$$
\C_S \to    \A_{\pi}\C_{S/T} \to
\A_{\pi}I_q(\B_{X^P})
$$

equals the following composition:
\begin{eqnarray*}
\C_S\to i_{q\pi}\ha(\B_{X^P})\otimes_{p\in
P}(\C_{S_p})\to
i_{q\pi}\ha(\B_{X^P})\otimes_{p\neq \chi
}(\C_{S_p})
\otimes ( \A_{\vs}\C_{S_\chi/T})\\
\to \A_{\pi}(i_{q}\ha(\B_{X^P})\otimes
\otimes_{p\in P}(
\C_{q^{-1}p}))=\A_\pi\I_q(\B_{X^P}).
\end{eqnarray*}
\end{Lemma}
\pf Pick an $F$ in $\C_S$ and show that its
images under the two maps coincide.

First of all we note the following thing.
Let $s_1,s_2\in S$ be such that
$q\pi(s_1)\neq q\pi(s_2)$. Then
$q(x_{s_1}-x_{s_2})$ is invertible in
$\A_\pi\I_q(\B_{X^P})$.  Let us multiply $F$
by a product of sufficiently large number of
such factors. We  then shall obtain an
element in $\otimes_{p\in P} \C_{S_p}$, and
it is sufficient to prove the statement for
only such elements, (because
$q(x_{s_1}-x_{s_2})$ are all invertible in
the target space).  In this case the
statement follows directly from Lemma
\ref{razl}.

\subsubsection{} Let $q:S\to P$ be an arbitrary surjection.
Define a functor $\A_p$ from the category of
$\D_{X^P}$-modules to the category of
$\D_{X^S}$-modules by the formula
$$
\liminv_{N\to\infty}
i_p\ha(M)\otimes_{\O_{X^S}} \otimes_{p\in P}
\D'_{q^{-1}p}/\D'_{X^{q^{-1}p},N}.
$$
\subsubsection{}\label{p1p2}

Let $p:S\to R$ be an arbitrary surjection.
For $r\in R$ let $S_r:=p^{-1}r$. Pick
non-empty subsets $T_r\subset S_r$. Let
$P:=\disjoint_r S_r/T_r$. We then have a
natural decomposition: Let $p=p_2p_1$, where
$p_1:S\to P$, $p_2:P\to R$. We also denote
$p_r:S_r\to S_r/T_r$.

The above constructions allow us to define a
map
$$
\I_p\to \A_{p_1}\I_{p_2}
$$
as follows.

\begin{eqnarray*}
\I_p(M)=i_{p}\ha(M)\otimes (\boxtimes_{r\in
R} \C_{S_r}) \to i_{p}\ha(M)\otimes
 (\boxtimes_{r\in R} \A_{p_r}\C_{S_r/T_r})\\
\to i_p\ha(M)\otimes
\A_{p_1}(\boxtimes_{r\in R} \C_{p_2^{-1}r})
\\
\to \A_{p_1}(i_{p_2}\ha(M)\otimes
\boxtimes_{r\in R}
\C_{p_2^{-1}r})=\A_{p_1}\I_{p_2}(M).
\end{eqnarray*}

\subsubsection{} \label{uslp1p2}
It follows that the map $$ \I_p\to
\A_{p_1}\I_{p_2}$$ is defined for all
decompositions $p=p_2p_1$ such that
$p_1,p_2$ are surjections and for every
$t\in \Im p$, $p^{-1}_2t$ contains at most
one element $u$ such that $p_1^{-1}u$
consists of more than one element.

\subsubsection{}
Let $p_2=q_2q_1$ be a decomposition, where
$q_2,q_1$ are surjections.
\begin{Claim}
The following diagram is commutative:
$$\xymatrix{
\I_{p}\ar[r]\ar[d]& \A_{p_1}I_{p_2}\ar[r]& \A_{p_1}\I_{q_1}\i_{q_2}\\
\I_{q_1p_1}\i_{q_2}\ar[urr]&&}
$$
\end{Claim}
\pf Follows from Lemma \ref{lemmarazl}.
\endpf

\subsubsection{Compositions
$\delta_{q_1}\I_{q_2}\to \I_p\to
\A_{p_1}\I_{p_2}$}\label{bpbq}

Let $q_2q_1=p$ be a decomposition of $p$ as
a product of two surjections. We will
investigate the composition
$$
\delta_{q_1}\I_{q_2}\to \I_p\to
\A_{p_1}\I_{p_2}.
$$

Let $a$ be a universal surjection among
those that $p_1$ and $q_1$ pass through $a$:
$p_1=p'_1a, q_1=q'_1a$. The surjection $a$
is uniquely determined by the condition
$a(x)=a(y)$ iff $p_1(x)=p_1(y)$ and
$q_1(x)=q_1(y)$. Let $b$ be the universal
surjection among those that
$b=b_pp_1=b_qq_1$.

Let us describe $b$ more concretely. For
$t\in T$, let $R_t=p_2^{-1}t$ and
$S_t=p^{-1}t$. Then there is at most one
$r_t\in R_t$ such that $\# p^{-1}_1r_t>1$.
If there is no such an element pick $r_t$
arbitrarily.

Let $P_t=p_1^{-1}r_t$. Then
$$S_t/P_t\stackrel\sim\to R_t.$$

Let $e$ be the equivalence relation on $S$
determined by $q_1$. The subsets $S_t$ are
not connected by this relation. Define the
equivalence relation $f$ which determines
$b$. Let $u,v\in S_t$ we say $u\sim_f v$ if
either $u\sim_e v$ or if there are $u',v'\in
R_t$ such that $u\sim_e u'$ and $v\sim_e
v'$.

We then have a commutative diagram:

$$\xymatrix{&&\ar[dr]_{b_p}\ar[drr]^{p_2}&&\\
S\ar[r]^a\ar[urr]^{p_1}\ar[drr]_{q_1}&\ar[ur]_{p_1'}\ar[dr]^{q_1'}&&\ar[r]^c&T\\
&&\ar[ur]^{b_q}\ar[urr]_{q_2}&&}
$$

We  see that there is a natural map
$$
\delta_{q_1'}\A_{b_q}\to\A_{p_1'}\delta_{b_p}.
$$

\begin{Claim}
The composition
$$
\delta_{q_1}\I_{q_2}\to \I_p\to
\A_{p_1}\I_{p_2}.
$$
coincides with the composition:
$$
\delta_{q_1}\I_{q_2}\to
\delta_{a}\delta_{q_1'}\I_{cb_q}\to
\delta_a\delta_{q_1'}\A_{b_q}\I_c\to
\delta_a\A_{p_1'}\delta_{b_p}\I_c\to
\A_{p_1}\I_{p_2}.
$$
\end{Claim}
\pf Clear.
\endpf

\subsection{Maps $\P_p\to \A_{p_1}\delta_{p_2}$.}

We always assume that $p,p_1,p_2$ are the
same as above.

We are going to define maps
$x(p_1,p_2):\P_p\to \A_{p_1}\delta_{p_2}$
using induction in $|p_2|:=\#R-\#P$.

The base is $|p_2|=0$, i.e a bijective
$p_2$. Without loss of generality we can
assume that $P=S/T$ and $p_2=\Id$. The map
$x(p_1,\Id)$ is then defined as a
composition
$$
\P_p\to \I_p\to \A_p.
$$

The transition is as follows. We begin with
construction of a map $\xi(p_1,p_2):\P_p\to
\A_{p_1}\I_{p_2}$. We then show that it
passes through a unique map
$x(p_1,p_2):\P_p\to \A_{p_1}\delta_{p_2}$.

1. Construction of $\xi(p_1,p_2)$. For every
decomposition $p_2=q_2q_1$, all the maps
being  properly surjective, we define a map
$$\xi(p_1,q_1,q_2):\P_p\to \A_{p_1}\I_{p_2}$$
as the composition:
$$
\P_p\to\P_{q_1p_1}\I_{q_2}\to
\A_{p_1}\delta_{q_1}\I_{q_2} \to
\A_{p_1}\I_{q}.
$$

We also set $K(p_1,p_2):\P_p\to
\I_{p_1}\I_{p_2}\to \A_{p_1}\I_{p_2}$;
$L(p_1,p_2):\P_p\to \I_p\to
\A_{p_1}\I_{p_2}$ to be the natural maps.

We finally  define a map
$$
\xi(p_1,p_2):\P_p\to \A_{p_1}\I_{p_2}$$ as:
$$
\xi(p_1,p_2)=L(p_1,p_2)-K(p_1,p_2)-
\sum_{q=q_2q_1}\xi(p_1,q_1,q_2).
$$

2. We will now show that all compositions
$$\xymatrix{
\P_p\ar[rr]^{\xi(p_1,p_2)}
&&\A_{p_1}\I_{p_2}\ar[r]^{a}
 &
\A_{p_1}\I_{q_1}\i_{q_2}}
$$
vanish, where $p_2=q_2q_1$ is an arbitrary
decomposition into a product of proper
surjection.  To show the vanishing,
introduce a notation. For a map
$L:\P_{q_1p_1}\to \A_{p_1}\I_{q_1}$ we set
$$
L_!:\P_p\to \P_{q_1p_1}\i_{p_2}\to
\A_{p_1}\I_{q_1}\i_{p_2}.
$$

We then have 1)  If $q_1=q_3q^1$ and
$q_3,q^1$ is proper,
$$a\xi(p_1,q^1,q_2q_3)=\xi(p_1,q^1,q_3)_!;$$

2) $$a\xi(p_1,q_1,q_2)=\xi(p_1,q_1)_!;$$

3) $a\xi(q^1,q^2)=0$ if $q_1$ does not pass
through $q^3$;

4) $$aK(p_1,p_2)=K(p_1,q_1)_!;$$

5) $$aL(p_1,p_2)=L(p_1,q_1)_!.$$

Therefore,
$$a\xi(p_1,p_2)=\xi(p_1,q_1)_!-L(p_1,q_1)_!+
K(p_1,q_1)_!+\sum_{q_1=q_3q^1}\xi(p_1,q^1,q_3)_!=0,
$$
in virtue of the induction assumption.

This implies that $\xi(p_1,p_2)$ passes
through $\A_{p_1}\delta_{p_2}$. This
completes the construction.

\subsection{Interaction with the maps $\P\delta\P\to \P$}
\label{inter} We are going to study the
compositions
\begin{equation}\label{pdpa}
\P_{p^1}\delta_{p^2}\P_{p_3}\to\P_p\to
\A_{p_1}\delta_{p_2},
\end{equation}
where $p,p_1,p_2$ are as above and
$p=p^3p^2p^1$ is an arbitrary decomposition
into a product of surjections and $p^2$ is
proper.

\subsubsection{} We first of all note that
the map $\P_p\to \A_{p_1}\delta_{p_2}$
passes through the direct sum of natural
maps
$$
\P_p\to\I_{p^1}\I_{p^2}\cdots \I_{p^k},
$$
where $p^kp^{k-1}\cdots p^1=p$ is a
decomposition into a product of proper
surjections, and $p^1=ap_1$ for a surjection
$a$.

This implies that the composition
(\ref{pdpa})
 vanishes except  the following  cases

 1) $p_1$ is bijective;
 2) $p_3$ is bijective;
 3) $p_1=ap^1$.

Consider these cases.

1) Investigate the composition

$$\delta_{q_1}\P_{q_2}\to \P_p\to \A_{p^1}\delta_{p^2}.
$$
 We shall use the notations from \ref{bpbq}.
We then claim that this composition equals:
$$
\delta_{q_1}\P_{q_2}\to
\delta_a\delta_{q_1'}\P_{cb_q}\to
\delta_a\delta_{q_1'}\A_{b_q}\delta_c\to
\delta_a\A_{p_1'}\delta_{b_p}\delta_c\to
\A_{p_1}\delta_{p_2}.
$$

 2) The composition $\P_{p_1}\delta_{p_2}\to \P_p\to
 \A_{p^1}\delta_{p^2}$ does not vanish only if
 $p_1=ap^1$ for some $a$, in which case this map equals:
 $$
 \P_{p_1}\delta_{p_2}\to \A_{p^1}\delta_a\delta_{p_2}\to
 \A_{p^1}\delta_{p^2}.
 $$

  3) In this case the  composition
 vanishes. We shall use induction in $|p^2|$.

 The base, i.e.  the case when $p^2$ is bijective is clear.

 Let us pass to the transition.
 We will show that the composition
 $$
 \P_{ap^1}\delta_{p_2}\P_{p_3}\stackrel{u}\to
 \P_p\stackrel{\xi(p^1,p^2)}\to
  \A_{p^1}\I_{p^2}
 $$
vanishes.

We first consider the case when $a$ is
proper.

We see that $L(p^1,p^2)u=K(p^1,p^2)u=0$. and
that $\xi(p^1,q_1,q_2)u=0$ unless $q_1,q_2)$
belong to the isomorphism class of $q_1=a$
or $q_1=ap_2$ in which cases these
compositions mutually annihilate each other.

In the case $a=\Id$, $K(p^1,p^2)=0$,
$\xi(p^1,q_1,q_2)=0$ every time except when
the isomorphism class of $(q_1,q_2)$ is
given by $q_1=p_1$. In this situation $L$
and $\xi(p^1,p_2,p_3)$ annihilate each
other.
\subsubsection{Composition $\P_{fg\disjoint h}\to
\P_{fg\disjoint \Id}\delta_{\Id\disjoint
h}\to \A_{g\disjoint \Id}\delta_{f\disjoint
h}$} \label{longcomp} We claim that this
composition coincides with the map
$$
\P_{fg\disjoint \Id}\to
\A_{g\disjoint\Id}\to \delta_{f\disjoint
\Id}.
$$

\subsection{Interaction with the maps with $\P\to \P\P$}
The collection of functors $\A_p$ does not
form a system because it may include very
bad singularities which do not admit the
required asymptotic decomposition.

One can, nevertheless, define a
"correspondence". That is, for every
decomposition $p=p_2p_1$ One can define a
functor $\Ga(p_1,p_2)$ such that
  $$
  \Ga(p_1,p_2)(M)\subset \A_p(M)\oplus
\A_{p_1}\A_{p_2}(M).
$$

This is what we are going to do.

\subsubsection{ A subspace $\Ga_p\subset \D'_S\oplus
\A_p\D'_T$} Let $p:S\to T$ be a surjection.
We shall construct a subspace  $\Ga_p\subset
\D'_S\oplus\A_p\D'_T$.

 Pick a splitting $i:T\to S$ so that
$pi=\Id_T$. For $\{x_s\}_{s\in S}\in Y^S$
and $\la>0$ we set
$$
V_\la(\{x_s\}_{s\in S})=
\{x_{ip(s)}+\frac{x_s-x_{ip(s)}}\la\}_{\la\in
S}.
$$

Pick an element $\tau\in T$; for  a point
$\{x_t\}_{t\in T}\in Y^T$  and $\mu>0$ we
set
$$
U_\mu(\{x_t\}_{t\in T})=
\{x_\tau+\frac{x_t-x_\tau}\mu\}_{t\in T}.
$$

Pick  $f\in \D_{Y^S}$, $g\in \D_{Y^T}$ and
$F\in \D'_{S}$. We then have a function
$A(\la,\mu):=<F,V_\la fU_\mu g>$ in two
variables $\la,\mu$. This function is smooth
for all $\la,\mu>0$.

Let now $F'\in \A_{p}\D'_T$. We then can
construct an element
$$
A':=<F',V_\la fU_\mu g>\in \Co[\ln
\mu,\mu^{-1},\mu]][\ln \la,\la^{-1},\la]]
$$
in the obvious way.

We say that $A'$ is an asymptotic series for
$A$ if for every $P,Q>0$ and every
sufficiently large partial sum $A''$ of $A'$

$A-A''=\la^{P+1}x(\la,\mu)+\la^P\mu^Qy(\la,\mu)$,
where $x(\la,\mu)$ is continuous for all
$\la\geq 0$, $\mu>0$, and $\la,\mu$ is
continuous for all $\la,\mu\geq 0$.

Define $\Gamma_p$ as the set of all pairs
$F,F'$ such that $A'$ is an asymptotic
series for $A$ for all $f,g$ and all
splittings $i$ (one can actually show that
if this is true for one splitting $i$, it is
also true for every such a splitting).

The map $\Gamma_p\to \D'_S$ is injective and
closed under the action of dilations
$U^S_\la$. We may, therefore, split
$\Gamma_p=\oplus_n \Gamma_{p,n}$ into the
direct sum of generalized eigenvalues of
$U^S_\la$.

Let $p_1:S\to R$, $p_2:R\to T$ be
surjections. For $t\in T$ let $S_t,R_t$ be
the preimages and let $p_{1t}:S_t\to R_t$ be
the induced maps.

Set
$$\Gamma(p_1,p_2)(M):=\liminv_{N}(\bigotimes_{t\in T}
\Gamma_{p_{1t}}/\Gamma_{p_{1t},N})\otimes
i_{p}\ha(M).
$$

The inclusions $$ \Gamma_{p_{1t}}\subset
\D'_p\oplus \A_{p_1}\D'_{R_t}
$$
induce the inclusions
$$
\Gamma(p_1,p_2)\subset \A_p\oplus
\A_{p_1}\A_{p_2}.$$

\subsubsection{}
Let $p,p_1,p_2$ be as above. We then have
maps
$$
a:C_S\to \A_p(C_T)$$ and
$$b:C_S\to \A_{p_1}C_R\to  \A_{p_1}\A_{p_2}C_T.
$$

\begin{Claim} The map $a\oplus b$ passes through
$\Gamma(p_1,p_2)C_T$.
\end{Claim}

\subsubsection{Asymptotic series modulo diagonals} We will need
a weaker version of the above definition. In
the setting of the previous section, we say
that $F'$ is {\em an asymptotic series for
$F$ modulo diagonals in $X^{S/T}$ (resp. in
$X^S$)} if for every $P,Q$ there exists an
$N$ such that whenever $g$  vanishes on all
generalized diagonals up-to the order $N$
(resp. $f$ and $g$ vanish on all generalized
diagonals up-to the order $N$ ) , we have

$A-A''=\la^{P+1}x(\la,\mu)+\la^P\mu^Qy(\la,\mu)$,
where $x(\la,\mu)$ is continuous for all
$\la\geq 0$, $\mu>0$,  $y(\la,\mu)$ is
continuous for all $\la,\mu\geq 0$, and
$A''$ is a partial sum of $A'$ with
sufficiently  many terms.

Define $\Gl(p_1,p_2)$ (resp $\Gll(p_1,p_2)$)
in the same way as $\Ga(p_1,p_2)$ but using
asymptotic series modulo diagonals in
$X^{S/T}$ (resp. $X^S$).

Let $a,b\in S$ be such that $p_1(a)\neq
p_1(b)$. $\Gl(p_1,p_2)$ "does not feel"
sections  supported on the diagonal $a=b$".
Formal meaning is as follows. Let $\pi:S\to
S/\{a,b\}$; let $p':S/\{a,b\}\to T$. Let
$H:i_{\pi*}\A_{p'}\to \A_p$ be the  natural
map. Then  the functor
$$
(H(i_{\pi*}\A_{p'}(M)),0)\subset
\Gl(p_1,p_2).$$

Similarly, let $p_2=p_4p_3$ be a
decomposition into a product of surjections,
where $p_3$ is proper. Let
$i_{p_3*}\A_{p_4}\to \A_{p_2}$ be the
natural map. Let
$$G:\A_{p_1}i_{p_3*}\A_{p_4}\to \A_{p_1}\A_{p_2}
$$
be the induced map. Then
$$
(0,G(\A_{p_1}i_{p_3*}\A_{p_4})(M))\in
\Gl(p_1,p_2).
$$

\subsection{Decomposition of the map
$\P_p\to \A_{p_2}\delta_{p_1}$}

Let $p=p_2p_1$ be as in \ref{uslp1p2}.
Choose a decomposition $p_2=q_2q_1$, where
$q_2,q_1$ are surjections.

We are going to construct a map
$$
\P_p\to \A_{q_1}\A_{q_2}\delta_{p_1}
$$ such that its direct sum with
the map
$$\P_p\to \A_{p_2}\delta_{p_1}$$
will pass through
$\Gl(q_1,q_2)\delta_{p_1}$.

\subsection{}For a surjection $u:A\to B$
iet $B_m(p)\subset B$ be  given by
$$
B_m(u)=\{x\in B|\#p_1^{-1}(x)>1\}.
$$
Let $A_m(u)=u^{-1}B_m(u)$.  Let
$B=B_m(u)\disjoint B_s(u)$,
 $A=A_m(u)\disjoint A_s(u)$ be
the decompositions. We then have
$u=u_m\disjoint u_s$, where $u_s$ is
bijective and $u_m$ is essentially
surjective, i.e. $\#u_m^{-1}x>1$ for all
$x\in B_m(u)$.

\subsubsection{}Let $p_1:S\to R$, $p_2:R\to T$.
Let $q_1:S\to U$, $q_2:U\to R$.

Let $S_m:=S_m(p_1)$, $S_s=S_s(p_1)$. Then
$q_1(S_s),p_1(S_s)$ are identified with
$S_s$. Using this identification, we may
assume that $U=U_m\disjoint S_s$ and that
$q_1=q_{1m}\disjoint \Id_{S_s}$;
$R=R_m\disjoint S_s$, $q_2=q_{2m}\disjoint
\Id_{S_s}$ (see the diagrams below)

We will work with isomorphism classes of
maps $v:U\to X$ which are

1) injective on $U_m$;

2) there exists $w:X\to  T$ such that
$wv=p_2q_2$.

We may therefore assume that
$$
X=U_m\disjoint Y
$$ and  that $v=\Id\disjoint v_s$.

We see that equivalently, one can define a
map $v$ by a prescription of a map
$v_s:S_s\to U_m\disjoint Y$ such that

1) $v_s(S_s)\supset Y$;

2) there exists a map $w_s:Y\to T$ (it is
then determined uniquely)
 such that  the diagram below commutes.

We then have $w=p_2q_{2m}\disjoint w_s$.

Let $Z=R_m\disjoint Y$. Let $w_1:X\to Z$,
$w_1=q_{2m}\disjoint \Id_Y$; let $w_2:Z\to
T$, $w_2=p_2|_{R_m}\disjoint w_s$.

Let $\sigma:R\to Z$ be given by
$\Id_{R_m}\disjoint v_s$.

%
$$\xymatrix{S\ar[rr]^{p_1}\ar[dr]^{q_1}&& R\ar[dr]^\sigma
\ar[rr]^{p_2}&& T\\
              &U\ar[ur]^{q_2}\ar[dr]^v&&Z\ar[ur]^{w_2}&\\
              && X\ar[ur]^{w_1}&&}
$$

$$\xymatrix{S_m\disjoint
S_s\ar[rr]^{p_{1m}\disjoint
\Id_{S_s}}\ar[dr]_{q_{1m}\disjoint
\Id_{S_s}} && R_m\disjoint
S_s\ar[dr]_{\Id_{R_m}\disjoint v_s}
\ar[rr]^{p_2}&& T\\
&U_m\disjoint
S_s\ar[ur]^{q_{2m}\disjoint\Id_{S_s}}
\ar[dr]_{\Id_{U_m}\disjoint v_s}
&&R_m\disjoint Y\ar[ur]_{p_2|_{R_m}\disjoint w_s}&\\
              && U_m\disjoint Y\ar[ur]^{q_{2m}\disjoint \Id_{Y}}&&}
$$

We also see that there is a natural
transformation:

$$\delta_v\A_{w_1}\to \A_{q_2}\delta_\sigma.$$

Therefore, one constructs a map

$$
\mu_v:\P_p\to \P_{vq_1}\P_{w_2w_1}\to
\A_{q_1}\delta_v\A_{w_1}\delta_{w_2}\to
\A_{q_1}\A_{q_2}\delta_\sigma\delta_{w_2}\cong
\A_{q_1}\A_{q_2}\delta_{p_2}.
$$

Define a map
$$
\mu(q_1,q_2):\P_p\to
\A_{q_1}\A_{q_2}\delta_{p_2}
$$ as a sum of $\mu_v$ over the set of all isomorphism classes
of maps $v$.

Let $\nu:\P_p\to\A_{p_1}\delta_{p_2}$.
\begin{Claim}
The map $\nu\oplus\mu$ passes through
$\Gl(q_1,q_2)\delta_{p_2}$.
\end{Claim}
\subsubsection{} We need a Lemma.

Let $$L:\P_p\to \I_p\to \A_{p_1}\I_{p_2}.
$$
 For $v$ as in the previous section, set
$$\lambda_v: \P_p \to \P_{vq_1}\P_w\to
\A_{q_1}\delta_v\I_w\to
\A_{q_1}\I_{p_2q_2}\to
\A_{q_1}\A_{q_2}\I_{p_2}.
$$

Let $\lambda=\sum_v \lambda_v$.
\begin{Lemma}
$L-\lambda$ passes through
$\Gl(q_1,q_2)\I_{p_2}$
\end{Lemma}
\pf

Let $ \Lambda$ be given by:
$$
\P_p\to \I_p\to \A_{q_1}\I_{p_2q_2} \to
\A_{q_1}\A_{q_2}\I_{p_2}.
$$
As we have seen above, $L-\Lambda$ passes
through
$$
\Gl(q_1,q_2).
$$

We can now focus on the difference
$\Lambda-\lambda$. It suffices to show that
it passes through
$\A_{q_1}\Al_{q_2}\I_{p_2}$.

Let $$ H:\P_p\to \I_p\to
\A_{q_1}\I_{p_2q_2};
$$
Let
$$
G_v:\P_p\to
\P_{vq_1}\P_w\to\A_{q_1}\delta_v\I_w\to
\A_{q_1}\I_{p_2q_2}.
$$
We see that $\Lambda-\lambda$ equals the
composition of $H-\sum_v G_v$ with the map
$$\A_{q_1}\I_{p_2q_2}\to
\A_{q_1}\A_{q_2}\I_{p_2}
$$

Let $p_2q_2=wv$ be a decomposition such that
$v$ is as above. Then it is not hard to see
that the compositions of $H-\sum_v G_v$ with
the map
$$
\I_{p_2q_2}\to \I_v\i_w$$ vanish. This
implies that the composition of $H-\sum_v
G_v$ with the map
$$\A_{q_1}\I_{p_2q_2}\to
\A_{q_1}\A_{q_2}\I_{p_2}
$$
passes through
$$
\A_{q_1}\Al_{q_2}\I_{p_2}.
$$
This implies the statement.

\subsubsection{Proof of the Claim}

We shall use induction with respect to
$|p|$. The base is clear. Let us pass to the
transition.

By definition the composition
$$
\xi(p_1,p_2):\P_p\to \A_{p_1}\delta_{p_2}\to
\A_{p_1}\I_{p_2}
$$
equals $-L+\sum \xi(p_1;r_1,r_2)$, where we
changed $q$ for $r$ to avoid a confusion,
and the sum is taken over all isomorphism
classes of decompositions $p_2=r_2r_1$ such
that $r_2$ is proper (so that $K$ is
included as the term corresponding to
$r_1=\Id$.)

Define the map $\eta(p_1;r_1,r_2)$ as the
composition:

$$
\P_p\to \P_{r_1p_1}\P_{r_2} \to
\A_{q_1}\A_{q_2}\delta_{r_1}\I_{r_2} \to
\A_{q_1}\A_{q_2}\I_{p_2}.
$$

According to the induction assumption, the
direct sum
$$
\xi(p_1,r_1,r_2)+\eta(p_1;r_1,r_2)
$$
passes through $\Gl(q_1,q_2)\I_{p_1}$.

Let $\lambda$ be as in the Lemma. We then
know that $L+\lambda$ passes through
$\Gl(q_1,q_2)$.

It now suffices to prove that
$-\sum_v\lambda_v+\sum\eta(p_1,r_1,r_2)=0 $

We, first of all see that
$$
\sum\eta(p_1,r_1,r_2)$$ equals the sum of
the maps of the form
\begin{eqnarray*}
\P_p\to
\P_{vq_1}\P_{r_2w'_2w_1}\to \P_{vq_1}\P_{w'_2w_1}\P_{r_2}\\
\stackrel{\xi(w_1,w'_2,r_2)}\to
\A_{vq_1}\A_{w_1}\I_{r_2w'_2}\\
\to \A_{q_1}\delta_v\A_{w_1}\I_{r_2w'_2}\to
\A_{q_1}\A_{q_2}\I_{p_2}.
\end{eqnarray*}
where $v,w_1,w_2:=r_2w'_2$ are as in the
revious section, and the decompositions
$w_2=r_2w'_2$ are arbitrary, not necessarily
proper.

The map $\lambda$ equals

$$
\P_p\to  \P_{vq_1}\P_{w}\to
\P_{vq_1}\I_{w}\stackrel{L}\to
\P_{vq_1}\A_{w_1}\I_{w_2}\to
\A_{q_1}\delta_{v}\A_{w_1}\I_{w_2}\to
\A_{q_1}\A_{q_2}\I_{p_2}.
$$

 The statement now follows immideately.

 \subsection{Maps $\p_i\P_p\to \A_q\p_j$}
\subsubsection{Definition}

Suppose we have a commutative square
$$
\xymatrix{ R\ar@{>>}[r]^p           & T\\
           S\ar@{^{(}->}[u]^i\ar@{>>}[r]^q   & Q\ar@{^{(}->}[u]^j}
           $$

Let  us define a map $\p_i\P_p\to \A_q\p_j$
in the following way.

Let $L:=R\bs i(S)$. We then have an
identification $R=S\disjoint L$. Let
$p_1:R\to Q\disjoint L$ be just
$$R\stackrel{cong}\to S\disjoint L
\stackrel{q\disjoint \Id_L}\longrightarrow
Q\disjoint L.
$$

Let $p_2:Q\disjoint L\to  T$ be given by
$j\disjoint p_{1L}$.  We then have
$p=p_2p_1$, where $p_2,p_1$ are surjections.
We see that they satisfy the conditions
which are
 necessary
to define  the map
$$\P_p\to \A_{p_1}\delta_{p_2}.
$$

Finally, let $i_Q:Q\to Q\disjoint L$ be the
inclusion. We then have a natural map
$$
\p_i\A_{p_1}\to \A_q\p_{i_Q}.
$$

The map $\p_i\P_p\to \P_q\p_j$ is now
defined as the composition:

$$
\p_i\P_p\to \p_i\A_{p_1}\delta_{p_2}\to
\A_q\p_{i_Q}\delta_{p_2}\to
\A_q\p_{p_2i_Q}=\A_p\p_j.
$$

\subsubsection{Properties} We shall translate the properties
of the maps $\P_{p_2p_1}\to
\A_{p_1}\delta_{p_2}$ into the language of
the maps
$$\p_i\P_p\to \A_q\p_j.$$

1. If the square $(i,p,j,q)$ is suitable,
then  the diagram
$$\xymatrix{
\p_i\P_p\ar[r]\ar[dr]& \A_q\p_j\\
               &  \I_q\p_j\ar[u]}$$
is  commutative.

2. Let $p=p_3p_2p_1$ be a decomposition into
a product of surjections, where $p_2$ is
proper.

3.1. The composition
$$
 \p_i\P_{p_1}\delta_{p_2}\P_{p_3}\to \p_i\P_p\to \P_q\delta_j
$$

vanishes unless $p_1$ or $p_3$ are
bijective.

3.2. Investigate the composition

$$
\p_i\delta_{p_2}\P_{p_3}\to \p_i\P_p\to
\P_q\p_j.
$$

We can uniquely decompose $ p_2i=j_2q_2,$
where $j_2$ is injective and $q_2$ is
bijective. Furthermore, we can decompose
$p_3j_2=jq'$ for a surjection $q'$

The above composition is then:
$$
\p_i\delta_{p_2}\to
\delta_{q_2}\p_{j_2}\P_{p_3}\to
\delta_{q_2}\A_{q'}\p_j\to\A_{q}\p_j.
$$

3.3. The composition
$$\p_i\P_{p_1}\delta_{p_2}\to \p_i\P_p\to
\P_q\p_j$$ does not vanish only if one can
decompose $p_1i=j_1q$, where $j_1$ is
injective, in which case it equals

$$
\p_i\P_{p_1}\delta_{p_2}\to
\P_q\p_{j_1}\delta_{p_2}\to \P_q\p_j.
$$

 4. Let $q=q^2q^1$ be a decomposition into a product of
 surjections.

 Consider the set of all isomorphism classes of the diagrams
 $$\xymatrix{
 R\ar[r]^{p^1} &R_1\ar[r]^{p^2} & T\\
 S\ar[u]^{i}\ar[r]^{q^1}& P_1\ar[u]^{j_1}\ar[r]^{q^2}& P\ar[u]^j}
 $$
 where $p^2p^1=p$. For every such a diagram $D$ we have
 a map

 $$
u_D: \p_i\P_p\to \p_i\P_{p_1}\P_{p_2}\to
 \A_{q_1}\p_{j_1}\P_{p_2}\to \A_{q_1}\A_{q_2}\p_j,
 $$

Let $u$ be the sum of $u_D$ taken over the
set of all diagrams $D$.

Let $v:\p_i\P_p\to \P_q\p_j$. Then the
direct sum $u\oplus v$ passes through
$\Gl(q_1,q_2)\p_j$.

\subsection{Maps $\int_q:\N_q\to \Al_q$}

We  define $\int_q=0$ on all terms of
cohomological degree $<0$.  The terms of
degree zero are all of the form
$p_{i*}\P_p$, where $pi=q$, $i$ is injective
and $P$ is surjective. We define
$\int_q|_{p_{i*}\P_p}$ as the composition:
$$
p_{i*}\P_p\to \A_q\to \Al_q.
$$

\begin{Claim} $d\int_q=0$.
\end{Claim}
\pf

We need to check that the composition
$$
(\N_p)^{-1}\stackrel d\to
(\N)^0_p\stackrel{\int_p} \to \A_p
$$
vanishes.

The functor  $(\N_p){-1}$ is a direct sum of
the terms
$\p_i\P_{p_1}\delta_{p_2}\P_{p_3}$, where
$p=p_1p_2p_3i$, where $i$ is injective and
$p_1,p_2,p_3$ are surjective and $p_2$ is
proper.

Consider several cases.

1) $p_1$ is bijective. We may think that
$p_1=\Id$. The restriction of the
differential onto this term equals the sum
$-D_1+D_2$, where
$$D_1:\p_i\delta_{p_2}\P_{p_3}\to \p_i\P_{p_3p_2};
$$
The map $D_2$ does not vanish only if
$i_1=p_2i$ is injective, in which case
$$
D_2:\p_i\delta_{p_2}\P_{p_3}\to
\p_{i_2}\P_{p_3}.
$$

The check now reduces to showing that the
diagram
$$\xymatrix{   &\p_i\P_{p_3p_2}\ar[dr]&\\
\p_i\delta_{p_2}\P_{p_3}\ar[ur]^{D_1}
\ar[dr]^{D_2}&&\Al_q\\
 &\p_{i_2}\P_{p_3}\ar[ur]&}
 $$
 is commutative which follows from the property 3.2.

2) $p_3$  is bijective. We may assume
$p_3=\Id$. In this case, the restriction of
the differential onto
$\p_i\P_{p_1}\delta_{p_2}$ equals
$-D_1+D_2$, where
$$D_1:\p_i\P_{p_1}\delta_{p_2}\to
\p_i\P_{p_2p_1}.
$$

The second term $D_2$ does not vanish only
if $p_1i=jq$, where $j$ is injective.  In
this case we can construct a commutative
diagram (uniquely up-to an isomorphism):
$$
\xymatrix{\ar[r]^{p_1}&\ar[r]^{p_2}&\\
          \ar[u]^{i_2}\ar[dr]^r&&\\
          \ar[u]^{i_1}\ar[r]^q&\ar[uu]^j\ar[uur]^{\Id}&}
          $$

in which the square $i_2,p,r,j$ is suitable.

The map $D_2$ is then:
\begin{eqnarray*}
\p_i\P_{p_1}\delta_{p_2}\stackrel\sim\to
\p_{i_1}\p_{i_2}\P_{p_1}\delta_{p_2}\\
\to \p_{i_1}\P_r\p_{j}\delta_{p_2}
\stackrel\sim\to \p_{i_1}\P_r
\end{eqnarray*}

The property 3.3, and \ref{longcomp} imply
that the diagram
$$\xymatrix{   &\p_i\P_{p_2p_1}\ar[dr]&\\
\p_i\P_{p_1}\delta_{p_2}\ar[ur]^{D_1}
\ar[dr]^{D_2}&&\Al_q\\
 &\p_{i_1}\P_{r}\ar[ur]&}
 $$

is commutative,  whence the statement

3) $p_1,p_3$ are proper. In this case the
restriction of the differential onto
$\P_{p_1}\delta_{p_2}\P_{p_3}$ simply
equals:
$$
\p_i\P_{p_1}\delta_{p_2}\P_{p_3}\to
\p_i\P_p.
$$

The composition
$$
\p_i\P_{p_1}\delta_{p_2}\P_{p_3}\to
\p_i\P_p\to \Al_q
$$
vanishes according to 3.1.

\endpf
\subsubsection{Interaction with the maps $\N\to \N\N$}
Let $p=p_2p_1$ be surjections. Let
$\int_{p_1p_2}:\N_p\to \N_{p_1}\N_{p_2}\to
\Al_{p_1}\Al_{p_2}$

\begin{Claim}  The map $\int_p\oplus \int_{p_1,p_2}$
passes through $\Gll(p_1,p_2)$.
\end{Claim}
\pf Compute the restriction of the map
$\N_p\to \N_{p_1}\N_{p_2}$ onto $\p_i\P_q$.

By definition, such a restriction equals the
sum of maps $m(q_1,q_2)$, where $q=q_2q_1$
and $q_1i=jp_2$, where $j$ is injective. In
this case one can construct a unique, up-to
an isomorphism, commutative diagram
$$\xymatrix{ \ar[r]^{q_1}&\ar[r]^{q_2}&\\
\ar[u]^{i_2}\ar[r]^{\pi}&\ar[u]^j\ar[ur]^{p_2}&\\
\ar[u]^{i_1}\ar[ur]^{p_1}&&}
$$
where the square $i_2,q_1,j,\pi$ is
suitable.

The map $m(q_1,q_2)$ is then given by:
$$
\p_i\P_q\to
\p_{i_1}\p_{i_2}\P_{q_1}\P_{q_2}\to
\p_{i_1}\P_{\pi}\p_j\P_{q_2}.$$

The composition
$$
\p_i\P_q\to
\p_{i_1}\p_{i_2}\P_{q_1}\P_{q_2}\to
\p_{i_1}\P_{\pi}\p_j\P_{q_2}\to
\A_{p_1}\A_{p_2}
$$
equals, in virtue of \ref{longcomp},
$$
u(q_1,q_2): \p_i\P_q\to
\p_i\P_{q_1}\P_{q_2}\to\A_{p_1}\p_j\A_{q_2}
\to \A_{p_1}\A_{p_2}.
$$

The sum of all $u(q_1,q_2)$ is  the map $u$
from 4. Therefore,  the direct sum of the
composition
$$
\p_i\P_q\to
\p_{i_1}\p_{i_2}\P_{q_1}\P_{q_2}\to
\p_{i_1}\P_{\pi}\p_j\P_{q_2}\to
\A_{p_1}\A_{p_2}
$$
with the  map

$$
\p_i\P_q\to \A_p$$

passes through $\Al(p_1,p_2)$, whence the
statement.
\endpf

\end{document}